\pdfoutput=1

\documentclass[oneside,reqno,11pt]{amsart}

\usepackage{preamble}

\myexternaldocument[II:]{corank1_ASW_II}
\myexternaldocument[III:]{corank1_ASW_III}
\myexternaldocument[IV:]{corank1_ASW_IV}

\title{Co-rank \texorpdfstring{$1$}{1} Arithmetic Siegel--Weil I: \\ Local non-Archimedean}
\author{Ryan C. Chen}
\date{May 2, 2024}
\address{Department of Mathematics, Massachusetts Institute of Technology, 182 Memorial Drive, Cambridge, MA 02139, USA}
\email{rcchen@mit.edu}

\begin{document}
    
    \begin{abstract}
This is the first in a sequence of four papers, where we prove the arithmetic Siegel--Weil formula in co-rank $1$ for Kudla--Rapoport special cycles on exotic smooth integral models of unitary Shimura varieties of arbitrarily large even arithmetic dimension. 
Our arithmetic Siegel--Weil formula implies that degrees of Kudla--Rapoport arithmetic special $1$-cycles are encoded in the first derivatives of unitary Eisenstein series Fourier coefficients.

The crucial input is a new local limiting method at all places. 
In this paper, we formulate and prove the key local theorems at all 
non-Archimedean places.
On the analytic side, the limit relates local Whittaker functions on different groups. On the geometric side at nonsplit non-Archimedean places, the limit relates degrees of $0$-cycles on Rapoport--Zink spaces and local contributions to heights of $1$-cycles in mixed characteristic.
\end{abstract}

    \maketitle
    
    \tableofcontents
    
    \clearpage


    \part*{Preliminary}

        \section{Introduction}
        \label{sec:intro}
            The landmark work of Gross and Zagier \cite{GZ86} showed that N\'eron--Tate heights of Heegner points on elliptic curves over $\mathbb{Q}$ are encoded in the first central derivatives of associated Rankin--Selberg $L$-functions. After the work of Gross and Keating \cite{GK93} on arithmetic intersection numbers for modular correspondences, Kudla proposed to recast such formulas in the language of special cycles on higher-dimensional Shimura varieties. This was originally formulated for integral models of orthogonal Shimura varieties in low dimensions by Kudla \cite{Kudla97a, Kudla97b, Kudla04} and the subsequent work of Kudla and Rapoport \cite{KR99,KR00}, where they pioneered the moduli definition of special cycles on integral models. Later, the attention was shifted to unitary Shimura varieties by Kudla--Rapoport in \cite{KR11,KR14}. Along with other closely related predictions about special cycles (e.g. modularity of generating series), these ideas are now called \emph{Kudla's program}. 
Kudla's program has played a role in a range of works, such as Gross--Zagier formulas on Shimura curves \cite{YZZ13}, the averaged Colmez conjecture \cite{AGHMP18, YZ18}, the arithmetic fundamental lemma \cite{Zhang21}, results on the Beilinson--Bloch conjecture \cite{LL22I}, and Picard rank jumps for K3 surfaces \cite{SSTT22}. We refer to Li's excellent surveys \cite{Li23,Li22IHES} for more.

Our work is about \emph{arithmetic Siegel--Weil formulas} in Kudla's program \cite[{Problem 6}]{Kudla04}, which (conjecturally) relate the first derivatives of Siegel Eisenstein series for unitary (resp. symplectic) groups with ``arithmetic theta series'' formed from special cycles on integral models of unitary (resp. orthogonal) Shimura varieties (see \cref{equation:intro:arith_Siegel_Weil:arithmetic_Siegel-Weil} below). These are closely parallel to the classical (resp. geometric) Siegel--Weil formulas, which state that special values of Eisenstein series encode representation numbers (resp. complex degrees) for lattices (resp. complex special cycles on Shimura varieties). For the interested reader, we refer to the introduction of our companion paper \cite{corank1_ASW_IV.pdf} where formulations of classical/geometric/arithmetic Siegel--Weil are compared in more detail.

Our main global results for this four-paper sequence are stated in \cref{ssec:intro:results}. The present paper focuses on the key ``local arithmetic Siegel--Weil'' theorems at all non-Archimedean places. Further comparisons with our companion papers \cite{corank1_ASW_II.pdf,corank1_ASW_III.pdf,corank1_ASW_IV.pdf} may be found in \cref{ssec:intro:outline}.

            \subsection{Eisenstein series}
            \label{ssec:intro:Eisenstein}
                In our four-paper sequence, we focus on the unitary/Hermitian case (literature review in \cref{ssec:intro:literature} also discusses the symplectic/quadratic case).
For the introduction, fix an imaginary quadratic field $F / \Q$ with ring of integers $\mc{O}_F$ and odd discriminant $\Delta$. Given $m \in \Z_{\geq 0}$ and an even integer $n \in \Z$, we consider the (normalized) \emph{Siegel Eisenstein series}
    \begin{align}\label{equation:intro:results:Eisenstein}
    E^*(z,s)^{\circ}_n & \coloneqq \Lambda_m(s)^{\circ}_n \sum_{\begin{psmallmatrix} a & b \\ c & d \end{psmallmatrix} \in P_1(\Z) \backslash SU(m,m)(\Z)} \frac{\det(y)^{s - s_0}}{\det(c z + d)^n |\det (c z + d)|^{2(s - s_0)}}
    \end{align}
for the group
    \begin{equation}\label{equation:intro:results:U(m,m)}
    U(m,m) \coloneqq \left \{ h \in \Res_{\mc{O}_F / \Z} \GL_{2m} : {}^t \overline{h} \begin{pmatrix} 0 & 1_m \\ -1_m & 0 \end{pmatrix} h = \begin{pmatrix} 0 & 1_m \\ -1_m & 0 \end{pmatrix} \right \}
    \end{equation}
where $\Lambda_m(s)^{\circ}_n$ is the normalizing factor
    \begin{align}
    \Lambda_m(s)^{\circ}_n &\coloneqq \frac{(2 \pi)^{m(m-1)/2}}{(-2 \pi i)^{nm}} \pi^{m(- s + s_0)} |\Delta|^{m(m-1)/4 + \lfloor m / 2 \rfloor (s + s_0)} 
    \\
    & \mathrel{\phantom{\coloneqq}} \cdot \left ( \prod_{j = 0}^{m - 1} \Gamma(s - s_0 + n - j) \cdot L(2s + m - j, \eta^{j + n})  \right ). \notag
    \end{align}
In \eqref{equation:intro:results:U(m,m)}, the notation $1_m$ stands for the $m \times m$ identity matrix, we wrote $SU(m,m) \subseteq U(m,m)$ for the determinant $1$ subgroup, and we set $P_1 \coloneqq P \cap SU(m,m)$ for the Siegel parabolic $P \subseteq U(m,m)$ (consisting of $m \times m$ block upper triangular matrices). The variable $s \in \C$ is a complex parameter, we set $s_0 = (n - m)/2$, and the element $z = x + i y$ lies in Hermitian upper-half space (i.e. $x \in \mrm{Herm}_m(\R)$ and $y \in \mrm{Herm}_m(\R)_{>0}$; the latter means that $y$ is positive definite).\footnote{Here, the notation $\mrm{Herm}_m$ denotes a scheme over $\Spec \Z$, e.g. $\mrm{Herm}_m(\R)$ denotes $m \times m$ complex Hermitian matrices, and $\mrm{Herm}_m(\Q)$ denotes $m \times m$ Hermitian matrices with entries in $F$. See conventions at the end of Section \ref{ssec:Hermitian_conventions:Hermitian_alternating_symmetric}.} The symbol $\eta$ denotes the quadratic character associated to $F / \Q$ (via class field theory).
The sum in \eqref{equation:intro:results:Eisenstein} is convergent for $\mrm{Re}(s) > m / 2$ and admits meromorphic continuation to all $s \in \C$. When $m = 1$, the expression in \eqref{equation:intro:results:Eisenstein} is a classical Eisenstein series on the usual upper-half plane. 

The normalized Eisenstein series has a symmetric functional equation 
    \begin{equation}
    E^*(z,s)^{\circ}_n = (-1)^{m (m - 1)(n - m - 1)/2} E^*(z,-s)^{\circ}_n,
    \end{equation}
see \crefext{IV:ssec:Eisenstein:global_normalized_Fourier:global_normalization}. Our definition of the normalizing factor $\Lambda_m(s)^{\circ}_n$ is motivated by symmetry of global and local functional equations, along with certain local special value formulas; see \crefext{IV:part:part_IV:Eisenstein} for further discussion. The function $\Lambda_m(s)^{\circ}_n$ should be closely related with the $L$-function of an Artin--Tate motive attached to the group $U(m,m)$, in the sense of Gross \cite{Gross97} (see \cite[{Remark 1.1.1}]{BH21}).

Given $T \in \mrm{Herm}_m(\Q)$, the Eisenstein series $E^*(z,s)^{\circ}_n$ has \emph{$T$-th Fourier coefficient}
    \begin{equation}\label{equation:intro:results:Fourier_coeff}
    E^*_T(y,s)^{\circ}_n \coloneqq 2^{m(m-1)/2} |\Delta|^{-m(m-1)/4} \int_{\mrm{Herm}_m(\Z) \backslash \mrm{Herm}_m(\R)} E^*(z, s)^{\circ}_n e^{-2 \pi i \mrm{tr}(Tz)} ~dx
    \end{equation}
for $z = x + i y$ in Hermitian upper-half space, where this integral is taken with respect to the Euclidean measure\footnote{The factor $2^{m(m-1)/2} |\Delta|^{-m(m-1)/4}$ disappears in the (usual) equivalent ad\`elic formulation, upon taking a certain self-dual Haar measure. The ad\`elic formulations of \eqref{equation:intro:results:Eisenstein} and \eqref{equation:intro:results:Fourier_coeff} are used in \crefext{IV:sec:setup}.} 
on $\mrm{Herm}_m(\R)$. The integral is convergent for $\mrm{Re}(s) > m / 2$, and admits meromorphic continuation to all $s \in \C$. When $\det T \neq 0$, there is a factorization into normalized \emph{local Whittaker functions}
    \begin{equation}\label{equation:intro:Eisenstein:Euler_product}
    E^*_T(y,s)^{\circ}_n = W^*_{T,\infty}(y,s)^{\circ}_n \prod_p W^*_{T,p}(s)^{\circ}_n
    \end{equation}
over all places, see \crefext{IV:ssec:Eisenstein:global_normalized_Fourier:global_normalization}. We suppress the definitions of local Whittaker functions in the introduction; see \cref{sec:part_I:local_Whittaker} for the non-Archimedean cases (and \crefext{IV:sec:local_Whittaker} more generally).
        
            \subsection{Arithmetic Siegel--Weil}
            \label{ssec:intro:arith_Siegel-Weil}
                \emph{Arithmetic Siegel--Weil formulas} predict that the derivative of $E^*_T(y,s)^{\circ}_n$ at $s = s_0$ should encode \emph{arithmetic degrees} of special cycles, where we replace the complex Shimura variety and its special cycles by integral models. For this, we restrict to the case of signature $(n - 1, 1)$, as is standard in the literature. 

Since the work of Kudla--Rapoport \cite{KR14} (also Rapoport--Smithling--Zhang \cite{RSZ21}), it has been customary to consider special cycles $\mc{Z}(T) \ra \mc{M}$ over (stacky) integral models $\mc{M} \ra \Spec \mc{O}_F$ for Shimura varieties associated to $G' \coloneqq \Res_{F / \Q} \G_m \times U(V)$, for signature $(n - 1, 1)$ non-degenerate $F / \Q$ Hermitian spaces $V$ with pairing $(-,-)$. In this paper, we assume $V$ contains a full-rank self-dual\footnote{We always mean \emph{self-dual} for the bilinear \emph{trace pairing} $\mrm{tr}_{F / \Q}(v,w)$ unless otherwise specified; see conventions in \cref{ssec:Hermitian_conventions:lattices}.} $\mc{O}_F$-lattice and we take $\mc{M} \ra \Spec \mc{O}_F$ to be the ``exotic smooth'' Rapoport--Smithling--Zhang (RSZ) integral model of relative dimension $n - 1$ \cite[{\S 6}]{RSZ21}.\footnote{Our assumption on $V$ forces $n \equiv 2 \pmod{4}$. We allow a slightly more general setup in our companion papers \cite{corank1_ASW_III.pdf,corank1_ASW_IV.pdf} for general $V$, at the cost of throwing out finitely many primes (particularly the ramified primes when $n$ is odd). This does not affect the essential ideas of our method, which is local in nature.} When $n = 2$, the stack $\mc{M}$ is essentially a disjoint union of (stacky) modular curves (\cref{example:integral_models:Serre_tensor_global}).

The stack $\mc{M}$ admits a moduli description: it parameterizes tuples $(A_0, \iota_0, \lambda_0, A, \iota, \lambda)$ where $A_0$ and $A$ are abelian schemes (dimensions $1$ and $n$ respectively) with $\mc{O}_F$-actions $\iota_0$ and $\iota$, and with compatible quasi-polarizations $\lambda_0$ and $\lambda$. The datum $(A_0, \iota_0, \lambda_0, A, \iota, \lambda)$ satisfies a few additional conditions, which we postpone to Section \ref{ssec:part_I:arith_intersections:integral_models}.

The moduli stack $\mc{M}$ carries a natural family of Hermitian $\mc{O}_F$-lattices 
    \begin{equation}
    \mc{L}at \ra \mc{M} \quad \quad \mc{L}at \coloneqq \underline{\Hom}_{\mc{O}_F}(A_0, A).
    \end{equation}
Given any $T \in \mrm{Herm}_m(\Q)$, the associated \emph{Kudla--Rapoport} special cycle $\mc{Z}(T) \ra \mc{M}$ is defined as the substack
    \begin{equation}\label{equation:intro:arith_Siegel-Weil:KR_cycle}
    \mc{Z}(T) \coloneqq \{\underline{x} \in \mc{L}at^m : (\underline{x}, \underline{x}) = T\} \subseteq \mc{L}at^m
    \end{equation}
consisting of $m$-tuples $\underline{x}$ with Gram matrix $T$.
More precisely, see \cref{ssec:part_I:arith_intersections:special_cycles}. This is in close analogy with classical Siegel--Weil: there one considers $\mc{O}_F$-lattices varying in a given \emph{genus} (and counts tuples with a fixed Gram matrix),
and here we are considering $\mc{O}_F$-lattices varying over the moduli stack $\mc{M}$. 
The morphism $\mc{Z}(T) \ra \Spec \mc{O}_F$ is smooth of relative dimension $n - 1 - \rank(T)$ in the generic fiber over $\Spec F$. If $T$ is not positive semi-definite, then $\mc{Z}(T)$ is empty.

An \emph{arithmetic Siegel--Weil formula} is an identity roughly of the shape
    \begin{equation}\label{equation:intro:arith_Siegel_Weil:arithmetic_Siegel-Weil}
    \frac{h_F}{w_F} \frac{d}{ds} \bigg|_{s = s_0} \frac{2 \Lambda_n(s-s_0)^{\circ}_n}{\kappa \Lambda_m(s)^{\circ}_n} E^*_T(y,s)^{\circ}_n \overset{?}{=} \widehat{\mrm{vol}}_{\widehat{\mc{E}}^{\vee}}([\widehat{\mc{Z}}(T)]).
    \end{equation}
Here we set $\kappa = 1$ (resp. $\kappa = 2$) if $m \neq n$ (resp. $m = n$). The right-hand side of \cref{equation:intro:arith_Siegel_Weil:arithmetic_Siegel-Weil} denotes an \emph{arithmetic volume}, which is a real number ``defined'' by an arithmetic intersection product
    \begin{equation}
    \widehat{\mrm{vol}}_{\widehat{\mc{E}}^{\vee}}([\widehat{\mc{Z}}(T)]) ``\coloneqq" \widehat{\deg}([\widehat{\mc{Z}}(T)] \cdot \widehat{c}_1(\widehat{\mc{E}}^{\vee})^{n - m})
    \end{equation}
in an arithmetic Chow ring $\arithCh^*(\mc{M})_{\Q}$ (roughly in the sense of Gillet--Soul\'e \cite{GS87})
for a certain metrized tautological bundle $\widehat{\mc{E}}^{\vee}$ on $\mc{M}$ (the bundle $\widehat{\mc{E}}^{\vee}$ is discussed in \cref{ssec:part_I:arith_intersections:tautological_bundle,ssec:part_I:arith_intersections:metrized_taut_bundle}). 
The notation $[\widehat{\mc{Z}}(T)]$ indicates a class in $\arithCh^m(\mc{M})_{\Q}$, which is expected to involve $\mc{Z}(T)$ and some additional Archimedean data (e.g. from a Green current on the complex Shimura variety), as appearing in arithmetic intersection theory.

An expected application of arithmetic Siegel--Weil formulas is in the theory of \emph{arithmetic theta lifting}. One expects to form automorphic \emph{arithmetic theta series} as generating series
    \begin{equation}\label{equation:intro:arith_Siegel-Weil:arith_theta}
    \widehat{\Theta} = \sum_T [\widehat{\mc{Z}}(T)] q^T
    \end{equation}
with ``Fourier coefficients'' $[\widehat{\mc{Z}}(T)]$ valued in the arithmetic Chow group $\arithCh^m(\mc{M})_{\Q}$. These should be analogous to (weighted averages of) classical theta series, as in the classical Siegel--Weil formula. In analogy with classical theta lifting, one expects to use $\widehat{\Theta}$ as an integral kernel to lift $U(m,m)$ automorphic forms to elements of $\arithCh^m(\mc{M})_{\Q}$. In analogy with the classical Rallis inner product formula, one expects to use the doubling method and arithmetic Siegel--Weil formulas to relate the derivative of an $L$-function with the arithmetic inner product of this arithmetic theta lift \cite[{Part III}]{Kudla04}. We refer to \cite{KRY06,BHKRY20II,LL22I,LL22II} for some cases where versions of this have been realized, with applications to Beilinson--Bloch. For modularity results on generating series of arithmetic divisors, see \cite{KRY06,BBK07,BHKRY20,Qiu22}.

We sketched the arithmetic Siegel--Weil formula as a rough expectation, because precise formulations remain open in the general case \cite[{Remark 4.4.2}]{Li22IHES}. In general, it is necessary to renormalize or modify the Eisenstein series in a way which is not completely understood. In fact, our normalization on the left-hand side of \cref{equation:intro:arith_Siegel_Weil:arithmetic_Siegel-Weil} is already nonstandard (new). 
We are not certain about this normalization for arithmetic Siegel--Weil in general, but our Theorem \ref{theorem:intro:results:main} (when $m = n$ for $T$ of co-rank $1$, and $m = n - 1$ for $T$ nonsingular; more discussion appears below) provides some evidence. The case of $m = n$ and $T$ nonsingular also holds, as can be extracted from known theorems in the literature (see discussion following \cref{equation:intro:results:conjecture} below). 

In general, posing a good (precise) definition of the arithmetic cycle class $[\widehat{\mc{Z}}(T)]$ is an open problem, especially for singular $T$ (due to arithmetic-intersection-theoretic difficulties), and particularly in the unitary case or over general totally real fields (due to a certain class number phenomenon), see \crefext{III:sec:arith_cycle_classes}. Previous works used $K$-theoretic methods to define special cycle classes (e.g. \cite{KR14} and \cite{HM22}), and the works by Feng--Yun--Zhang (moduli of shtukas) \cite{FYZ21,FYZ22SW} and Madapusi (Shimura varieties) \cite{Madapusi22} have employed derived algebro-geometric methods to define special cycle classes. As of now, these constructions do not incorporate the Archimedean place, which would be needed for arithmetic intersection theory (e.g. there seems to be no ``derived arithmetic intersection theory'' at the moment). Garcia and Sankaran have defined (Archimedean) Green currents associated to singular $T$ using the Mathai--Quillen theory of superconnections, but there has been no proposal to combine this with the non-Archimedean theory.

We first propose a method to construct the arithmetic special cycle classes $[\widehat{\mc{Z}}(T)]$ for arbitrary $T$. Our proposed definition mixes the work of Garcia--Sankaran with $K$-theoretic methods for positive characteristic contributions. Our construction may need adjustment on compactifications of integral models, but we expect it to apply in already-compact situations (e.g. the Rapoport--Smithling--Zhang \cite{RSZ21} setup for CM extensions of totally real fields $\neq \Q$). 

The first part of the main theorem in our four-paper sequence (\cref{theorem:intro:results:main}(1)) is a proof of (1) the arithmetic Siegel--Weil formula when $m = n$ and $T \in \mrm{Herm}_n(\Q)$ is singular of co-rank $1$. Most known (fully global) results concern special cycles $\mc{Z}(T) \ra \mc{M}$ which are empty in the generic fiber. These previous results include the non-Archimedean Kudla--Rapoport conjectures (for $T \in \mrm{Herm}_n(\Q)$ nonsingular) proved by Li--Zhang \cite{LZ22unitary} (and the ramified versions \cite{HLSY22,LL22II}), as well as the purely Archimedean results of Liu \cite{Liu11} and Garcia--Sankaran \cite{GS19}. Our theorem is the first (fully global) arithmetic Siegel--Weil result which involves mixed characteristic special cycles $\mc{Z}(T)$ on Shimura varieties of arbitrarily large dimension. We further discuss the comparison with previous literature in \cref{ssec:intro:literature}.

We also prove (2) the arithmetic Siegel--Weil formula when $m = n - 1$ and $T$ is nonsingular (\cref{theorem:intro:results:main}(2)). This is very closely related with our theorem for singular $T \in \mrm{Herm}_n(\Q)$ of co-rank $1$, as we explain further in \cref{ssec:intro:results}. This theorem implies that both the first derivative and the special value of a $U(n - 1, n - 1)$ Eisenstein series at the \emph{non-central} point $s = s_0 = 1/2$ have geometric meaning; see discussion following \cref{theorem:intro:results:main}.

As a byproduct of our methods, we prove (3) a version of the arithmetic Siegel--Weil formula (up to an volume constant which we did not calculate) for arbitrary $m$ when $T$ is nonsingular and not positive-definite (corresponding to a ``purely Archimedean'' arithmetic intersection number) (\cref{theorem:intro:results:Archimedean}). This purely Archimedean result is analogous to those in \cite{GS19} (there in a situation with compact Shimura varieties, which need not apply in the setup above), but our method of proof is completely different and is insensitive to compactness.

More importantly, we propose and apply a new uniform strategy to prove (1), (2), and (3). This is the key conceptual novelty in our work. Our strategy is a certain ``local limiting method'' at all places, Archimedean and non-Archimedean. We further sketch this strategy in Section \ref{ssec:intro:strategy_overview}, and at a finer level of detail in \cref{sec:part_I:sketch}.
        
            \subsection{Global results}
            \label{ssec:intro:results}
                We describe our global results in more detail, retaining the notation from \cref{ssec:intro:arith_Siegel-Weil}.

First, we propose a new candidate definition of arithmetic cycle classes
    \begin{equation}\label{equation:intro:results:arith_cycle_class}
    [\widehat{\mc{Z}}(T)] \coloneqq [\widehat{\mc{Z}}(T)_{\ms{H}}] + \sum_{p \text{ prime}} [{}^{\mathbb{L}}\mc{Z}(T)_{\ms{V},p}] \in \arithCh^m(\mc{M})_{\Q}
    \end{equation}
associated to arbitrary (possibly singular) $T$. Here, $[\widehat{\mc{Z}}(T)_{\ms{H}}]$ is intended to describe ``horizontal'' contributions and ${}^{\mathbb{L}}\mc{Z}(T)_{\ms{V},p}$ is intended to describe ``vertical'' contributions. 

In this paper, the vertical (positive characteristic) classes ${}^{\mathbb{L}}\mc{Z}(T)_{\ms{V},p}$ will be constructed in \cref{ssec:part_I:arith_intersections:vertical_classes}. For each prime $p$, we give a (new) definition of an element ${}^{\mathbb{L}}\mc{Z}(T)_{\ms{V},p} \in \mrm{gr}^m_{\mc{M}} K'_0(\mc{Z}(T)_{\F_p})_{\Q}$ (``vertical'') lying in the dimension $n - m$ graded piece of the Grothendieck group (tensor $\Q$) of coherent sheaves on $\mc{Z}(T)_{\F_p} \coloneqq \mc{Z}(T) \times_{\Spec \Z} \Spec \F_p$. Our construction is based on a certain ``$p$-local linear invariance'', and is explained in \cref{ssec:part_I:arith_intersections:vertical_classes}.

The construction of $[\widehat{\mc{Z}}(T)_{\ms{H}}]$ involves the flat part\footnote{Give an algebraic stack $\mc{X}$ over a Dedekind domain $R$, its \emph{flat part} or \emph{horizontal part} $\mc{X}_{\ms{H}}$ is the largest closed substack $\mc{X}_{\ms{H}} \subseteq \mc{X}$ which is flat over $\Spec R$. The stack $\mc{X}_{\ms{H}}$ is also the scheme-theoretic image of the generic fiber of $\mc{X}$. Given a formal algebraic stack $\mc{X}$ over $\Spf R$ for a complete discrete valuation ring $R$, its flat part of horizontal part $\mc{X}_{\ms{H}}$ is the largest closed substack $\mc{X}_{\ms{H}} \subseteq \mc{X}$ which is flat over $\Spf R$ (in the sense discussed in \crefext{III:footnote:flatness_formal_alg_stacks}).} $\mc{Z}(T)_{\ms{H}}$ of $\mc{Z}(T)$ and certain ``modified currents'' $g_{T,y}$ on the complex fiber $\mc{M}_{\C} \coloneqq \mc{M} \times_{\Spec \mc{O}_F} \Spec \C$ (choose either embedding $F \ra \C$), where $g_{T,y}$ is allowed to vary with a parameter $y \in \mrm{Herm}_m(\R)_{>0}$. We defer further discussion on the construction of $[\widehat{\mc{Z}}(T)_{\ms{H}}]$ 
to our companion papers \cite{corank1_ASW_II.pdf,corank1_ASW_III.pdf}.

As part of the expected automorphic behavior of $[\widehat{\mc{Z}}(T)]$, it is expected that these classes should satisfy a certain ``linear invariance'' property for the action\footnote{For any $\gamma \in \GL_m(\mc{O}_F)$ and any Hermitian matrix $T \in \mrm{Herm}_m(\Q)$, we say e.g. that $T$ and ${}^t \overline{\gamma} T \gamma$ are \emph{$\GL_m(\mc{O}_F)$-equivalent}, and that they lie in the same \emph{$\GL_m(\mc{O}_F)$-equivalence class}.} of $\GL_m(\mc{O}_F)$ on Hermitian matrices $T$. We verify this for the classes we define: for any $g_{T,y}$ satisfying 
    \begin{equation}\label{intro:results:Archimedean_linear_invariance}
    g_{T,y} = g_{{}^t \overline{\gamma} T \gamma, \gamma^{-1} y {}^t \overline{\gamma}^{-1}},
    \end{equation}
we show
    \begin{equation}
    [\widehat{\mc{Z}}(T)] = [\widehat{\mc{Z}}({}^t \overline{\gamma} T \gamma)]
    \end{equation}
where $[\widehat{\mc{Z}}(T)]$ is formed with respect to $y$ and $[\widehat{\mc{Z}}({}^t \overline{\gamma} T \gamma)]$ is formed with respect to $\gamma^{-1} y {}^t \overline{\gamma}^{-1}$.
In fact, we prove refined results: the vertical part at each prime $p$ is linearly invariant on the level of Grothendieck groups \cref{equation:arith_cycle_classes:vertical:linear_invariance}, and the horizontal part is linear invariant on its own \crefext{III:ssec:arith_cycle_classes:horizontal}.
The currents $g_{T,y}$ appearing in our main arithmetic Siegel--Weil results do satisfy the linear invariance property in \eqref{intro:results:Archimedean_linear_invariance}; see \crefext{III:ssec:Arch_uniformization:Archimedean}. Note that the Garcia--Sankaran Green currents in \cite[{(4.38)}]{GS19} also satisfy the same linear invariance property (but the modified currents in \cite[{Definition 4.7}]{GS19} do not).

Due to non-properness of $\mc{M} \ra \Spec \mc{O}_F$ in general, one should likely modify $[\widehat{\mc{Z}}(T)]$ on a suitable compactification of $\mc{M}$. 
If $\mc{Z}(T) \ra \Spec \mc{O}_F$ is proper, however, we consider certain ``arithmetic degrees without boundary contributions'' (a real number)
    \begin{align}\label{equation:intro:results:if_proper}
        \widehat{\deg}([\widehat{\mc{Z}}(T)] \cdot \widehat{c}_1(\widehat{\mc{E}}^{\vee})^{n-m}) & \coloneqq \left ( \int_{\mc{M}_{\C}} g_{T,y} \wedge c_1(\widehat{\mc{E}}^{\vee}_{\C})^{n-m} \right ) 
        \\
        & \hphantom{\coloneqq} + \widehat{\deg}((\widehat{\mc{E}}^{\vee})^{n - \rank(T)}|_{\mc{Z}(T)_{\ms{H}}}) \notag
        \\
        & \hphantom{\coloneqq} + \sum_{p \text{ prime}} \deg_{\F_p}({}^{\mathbb{L}}\mc{Z}(T)_{\ms{V},p} \cdot (\mc{E}^{\vee})^{n-m}) \log p \notag
    \end{align}
conditional on convergence of the integral, for a certain metrized tautological bundle $\widehat{\mc{E}}$ on $\mc{M}$ (\cref{ssec:part_I:arith_intersections:tautological_bundle,ssec:part_I:arith_intersections:metrized_taut_bundle}) (we do check convergence of the integral in the settings of our arithmetic Siegel--Weil results). The middle term is mixed characteristic in nature: for $\rank (T) = n - 1$, it is (essentially) a weighted sum of Faltings heights of abelian varieties (\crefext{IV:remark:arithmetic_Siegel-Weil:main_results:Faltings_height}).
For proper $\mc{Z}(T) \ra \mc{O}_F$,
the quantity in \eqref{equation:intro:results:if_proper} should coincide with the arithmetic degree (without boundary contributions) of a version of $[\widehat{\mc{Z}}(T)]$ on any reasonable compactification of $\mc{M}$. 

Our main theorems concern the $T$-th Fourier coefficients $E^*_T(y,s)^{\circ}_n$ of $E^*(z,s)^{\circ}_n$. As above, we write $h_F$ (resp. $w_F$) for the class number of (resp. number of roots of unity in) $\mc{O}_F$. The following theorem appears in the final paper of our four-paper sequence as \crefext{IV:theorem:arithmetic_Siegel-Weil:main_results:main}. 

\begin{theoremLetter}[Co-rank $1$ arithmetic Siegel--Weil]\label{theorem:intro:results:main}
Assume the prime $2$ splits in $\mc{O}_F$.
    \begin{enumerate}[(1)]
        \item For any $T \in \mrm{Herm}_n(\Q)$ with $\rank(T) = n - 1$ and any $y \in \mrm{Herm}_n(\R)_{>0}$, we have
            \begin{equation}\label{equation:intro:results:singular}
            \frac{h_F}{w_F} \frac{d}{d s} \bigg |_{s = 0} E^*_{T}(y,s)^{\circ}_n = \widehat{\deg}([\widehat{\mc{Z}}(T)]).
            \end{equation}
        \item For any $T^{\flat} \in \mrm{Herm}_{n - 1}(\Q)$ with $\det T^{\flat} \neq 0$ and any $y^{\flat} \in \mrm{Herm}_{n-1}(\R)_{>0}$, we have
            \begin{equation}\label{equation:intro:results:nonsingular}
            2 \frac{h_F}{w_F} \frac{d}{d s} \bigg|_{s = 0} \left ( \frac{\Lambda_n(s)_n^{\circ}}{\Lambda_{n - 1}(s + 1/2)^{\circ}_n} E^*_{T^{\flat}}(y^{\flat}, s + 1/2)^{\circ}_n \right ) = \widehat{\deg}([\widehat{\mc{Z}}(T^{\flat}) \cdot \widehat{c}_1(\widehat{\mc{E}}^{\vee})).
            \end{equation}
    \end{enumerate}
\end{theoremLetter}

Note that \cref{theorem:intro:results:main}(1) concerns the \emph{central} derivative of a $U(n,n)$ Eisenstein series, while part \cref{theorem:intro:results:main}(2) concerns a \emph{non-central} derivative of a $U(n - 1, n - 1)$ Eisenstein series.
For $n \equiv 0 \pmod {4}$, Theorem \ref{theorem:intro:results:main}(1) also holds in the sense that there is no self-dual $\mc{O}_F$-lattice of signature $(n - 1, 1)$ and the right-hand side is $0$ \crefext{IV:remark:arithmetic_Siegel-Weil:main_results:0_mod_4}.

\begin{remark}\label{remark:intro:results:main:value_and_derivative_both}
In the situation of Theorem \ref{theorem:intro:results:main}(2), there is also a ``geometric Siegel--Weil formula'' when we evaluate
    \begin{equation}\label{equation:part_I:intro:results:normalized_off-central_Eisenstein}
    2 \frac{h_F}{w_F} \frac{\Lambda_n(s)_n^{\circ}}{\Lambda_{n - 1}(s + 1/2)^{\circ}_n} E^*_{T^{\flat}}(y^{\flat}, s + 1/2)^{\circ}_n
    \end{equation}
at $s = 0$; the resulting expression is exactly $-\deg_{\C} \mc{Z}(T^{\flat})_{\C}$ (negative degree of complex fiber $\mc{Z}(T^{\flat})_{\C}$, which is a proper and quasi-finite Deligne--Mumford stack over $\Spec \C$). In other words, both the special value and the first derivative at $s = 1/2$ of the $U(n - 1, n - 1)$ Eisenstein series (normalized as in \cref{equation:part_I:intro:results:normalized_off-central_Eisenstein}) \emph{simultaneously} have arithmetic-geometric meaning.

The above ``geometric Siegel--Weil'' formula is also needed as an ingredient in our proof of Theorem \ref{theorem:intro:results:main}, and will be treated in our companion paper \cite{corank1_ASW_IV.pdf} via uniformization (Archimedean and non-Archimedean both work).
\end{remark}

We highlight the simplicity of the analytic side in Theorem \ref{theorem:intro:results:main}(1). It is expected that arithmetic Siegel--Weil for integral models with bad reduction should be corrected on the analytic side, e.g. by special values of other Eisenstein series. See for example \cite{HSY22a,HLSY22} for bad reduction in the nonsingular case $\det T \neq 0$ for the central derivative at $s = 0$ (i.e. $T$ is $n \times n$), or \cite{KRY06} for quaternionic Shimura curves. 
We do not know whether the analytic formulation \cite{HSY22a,HLSY22} is expected to hold for singular $T$. 

We argue that arithmetic Siegel--Weil formulas should be simplest to formulate on integral models with everywhere good reduction, as in our case. We thus propose a precise formulation of the analytic side of the central derivative arithmetic Siegel--Weil formula in our setup.

\begin{question*}[Arithmetic Siegel--Weil, central point]
Let $T \in \mrm{Herm}_n(\Q)$ be arbitrary. For a suitable current $g_{T,y}$, a suitable compactification of $\mc{M}$, and a possibly modified class $[\widehat{\mc{Z}}(T)]$ on the compactification, do we have
    \begin{equation}\label{equation:intro:results:conjecture}
    \frac{h_F}{w_F} \frac{d}{d s} \bigg|_{s = 0} E^*_T(y,s)^{\circ}_n \overset{?}{=} \widehat{\deg}([\widehat{\mc{Z}}(T)]).
    \end{equation}
\end{question*}

Our theorem verifies this proposed arithmetic Siegel--Weil formula for all singular $T \in \mrm{Herm}_n(\Q)$ of rank $n - 1$, in the sense of ``arithmetic degrees without boundary contributions''. The formula also holds (in the same sense) for all nonsingular $T \in \mrm{Herm}_n(\Q)$. This latter case (``central derivative nonsingular arithmetic Siegel--Weil'') is possibly considered known to experts up to a volume constant by collecting the local theorems in \cite{Liu11,LZ22unitary,LL22II}. This particular global statement does not appear in the literature, though other variants are available (e.g. for unramified CM fields $F / F_0$ with $F_0 \neq \Q$ \cite{LZ22unitary} or on integral models with bad reduction and correction terms by special values of other Eisenstein series \cite{HLSY22}). In one of our companion papers, we will compute the volume constant and explain how to extract the $\det T \neq 0$ case of \eqref{equation:intro:results:conjecture} from the literature \crefext{IV:remark:arithmetic_Siegel-Weil:main_results:nonsingular}.

A more optimistic version of \cref{equation:intro:results:conjecture} was given in \cref{equation:intro:arith_Siegel_Weil:arithmetic_Siegel-Weil} involving $T \in \mrm{Herm}_m(\Q)$ for arbitrary $m$, but we are less certain about the validity of that formulation in general.

In the general case of \cref{equation:intro:results:conjecture}, we expect the current $g_{T,y}$ to be essentially the currents of \cite[{Definition 4.7}]{GS19}, though $\mf{g}(T, \mbf{y}, \varphi_f)$ as defined in loc. cit. may need some modification (see the introductions of our companion papers \cite{corank1_ASW_II.pdf,corank1_ASW_III.pdf} for further discussion). Since our main theorems take a different approach to define $g_{T,y}$, we do not pursue this issue further.

Part (2) of Theorem \ref{theorem:intro:results:main} is the special case of part (1) when $T = \mrm{diag}(0,T^{\flat})$ and $y = \mrm{diag}(1, y^{\flat})$. The geometric sides will agree essentially by definition \cref{equation:intro:results:if_proper}. On the analytic side, the relation is provided by the formula
    \begin{align}
    E_T^*(y, s)^{\circ}_n & = \frac{\Lambda_n(s)^{\circ}_n}{\Lambda_{n - 1}(s + 1/2)^{\circ}_n} E^*_{T^{\flat}}(y^{\flat}, s + 1/2)^{\circ}_n -\frac{\Lambda_n(- s)^{\circ}_n}{\Lambda_{n - 1}(- s + 1/2)^{\circ}_n} E^*_{T^{\flat}}(y^{\flat}, s - 1/2)^{\circ}_n
    \end{align}
from \crefext{IV:corollary:Eisenstein:singular_Fourier:corank_1}, along with the functional equation $E^*_{T^{\flat}}(y^{\flat},s)^{\circ}_n = E^*_{T^{\flat}}(y^{\flat},-s)^{\circ}_n$. The general case of Theorem \ref{theorem:intro:results:main} is proved in a similar way as the special case $T = \mrm{diag}(0,T^{\flat})$, with an additional ``local diagonalizability argument'' (proof of \crefext{IV:theorem:arithmetic_Siegel-Weil:main_results:main}) where the identity is proved modulo $\sum_{\ell \neq p} \Q \cdot \log \ell$ for any given $p$ (and varying $p$ removes the ambiguity).

It is also possible to formulate and prove Theorem \ref{theorem:intro:results:main} in terms of Faltings heights (i.e. replacing the middle term in \eqref{equation:intro:results:if_proper} with the degree of the metrized Hodge bundle). The formulation in Theorem \ref{theorem:intro:results:main} seems more natural to us, but the version with Faltings heights is in \crefext{IV:remark:arithmetic_Siegel-Weil:main_results:Faltings_height}.

The simplest case of Theorem \ref{theorem:intro:results:main} is the case $n = 2$. When $\mc{O}_F^{\times} = \{ \pm 1 \}$, the Serre tensor construction gives an open and closed embedding $\ms{M}_0 \times_{\Spec \mc{O}_F} \ms{M}_{\text{ell}} \ra \mc{M}$, where $\ms{M}_0$ is the moduli stack of elliptic curves with signature $(1,0)$ action by $\mc{O}_F$ and $\ms{M}_{\text{ell}}$ is the moduli stack of all elliptic curves, base-changed to $\mc{O}_F$ \crefext{IV:ssec:arithmetic_Siegel-Weil:Serre_tensor}. In this case, the special cycle $\mc{Z}(j) \ra \mc{M}$ for $j \in \Z_{>0}$ pulls back to the $j$-th Hecke correspondence. Then the proof of \cref{theorem:intro:results:main} gives the following corollary (appearing later as \crefext{IV:corollary:arithmetic_Siegel-Weil:Serre_tensor}). One might think of this corollary as reformulating a result of Nakkajima--Taguchi \cite{NT91} (they compute Faltings heights of elliptic curves with CM by possibly non-maximal orders) by averaging over Hecke translates and expressing the result in terms of Eisenstein series Fourier coefficients. 

\begin{corollary}\label{corollary:intro:results:n=2}
Assume $2$ is split in $\mc{O}_F$. Fix any elliptic curve $E_0$ over $\C$ with $\mc{O}_F$-action. For any integer $j > 0$, we have
    \begin{equation}
    \sum_{E \xra{w} E_0} (h_{\mrm{Fal}}(E) - h_{\mrm{Fal}}(E_0)) = \frac{1}{2} \frac{d}{ds} \bigg|_{s = 1/2} \left ( j^{s + 1/2} \s_{-2s}(j) \right )
    \end{equation}
where the sum runs over degree $j$ isogenies $w \colon E_0 \ra E$ of elliptic curves.
\end{corollary}

The notation $h_{\mrm{Fal}}(E)$ denotes the (stable) Faltings height of the elliptic curve $E$ after descent to any number field, and similarly for $E_0$. The notation $\s_s(j) \coloneqq \sum_{d \mid j} d^s$ denotes the usual divisor function, and the quantity $j^{s+1/2} \s_{-2s}(j)$ is the product of the normalized non-Archimedean local Whittaker functions in the $j$-th Fourier coefficient $E^*_j(y,s)^{\circ}_2$ (with $m = 1$) as in \cref{equation:intro:Eisenstein:Euler_product}. The derivative of the Archimedean local Whittaker function $W^*_{j,\infty}(y,s)^{\circ}_2$ at $s = 1/2$ is also calculated explicitly and compared with its geometric counterpart (integral of Green function wedge Chern form on upper half-plane) in our companion paper \crefext{II:ssec:Archimedean_identity:case_n_is_2}. 

Our purely Archimedean result (for arbitrary $n$ and $m^{\flat} \geq 1$) is the following.

\begin{theoremLetter}[Archimedean arithmetic Siegel--Weil, nonsingular]\label{theorem:intro:results:Archimedean}
Consider any integer $m^{\flat}$ with $1 \leq m^{\flat} \leq n$, and consider any $T^{\flat} \in \mrm{Herm}_{m^{\flat}}(\Q)$ which is nonsingular and not positive definite.
    \begin{enumerate}[(1)]
        \addtocounter{enumi}{2}
        \item For any $y^{\flat} \in \mrm{Herm}_{m^{\flat}}(\R)_{>0}$, we have an equality of real numbers
            \begin{equation}\label{equation:intro:results:Archimedean}
            \widehat{\deg}([\widehat{\mc{Z}}(T^{\flat})] \cdot \widehat{c}_1(\widehat{\mc{E}}^{\vee})^{n - m^{\flat}}) \coloneqq \int_{\mc{M}_{\C}} g_{T^{\flat},y^{\flat}} \wedge c_1(\widehat{\mc{E}}^{\vee}_{\C})^{n - m^{\flat}} = (-1)^{n - m^{\flat}} C \cdot \frac{h_F}{w_F} \frac{d}{ds} \bigg|_{s = s_0^{\flat}} E^*_{T^{\flat}}(y^{\flat}, s)^{\circ}_n
            \end{equation}
        where $s_0^{\flat} \coloneqq (n - m^{\flat}) / 2$. Here $C \in \Q_{>0}$ is the volume constant from \crefext{IV:lemma:local_Siegel-Weil:uniformization_degree:main:1}, for the Hermitian space $V$ and $v_0 = \infty$ in the notation of loc. cit.. The constant $C$ may depend on $n$ and $m^{\flat}$ (and $F$), but does not otherwise depend on $T^{\flat}$.
    \end{enumerate}
\end{theoremLetter}

This appears (in stronger form) as \crefext{IV:theorem:arithmetic_Siegel-Weil:main_results:Archimedean}. That version applies for all $n$ (even or not) and arbitrary level, as it is a statement about the complex Shimura variety. We gave the weaker version here to avoid more notation. Due to non-properness of $\mc{M}_{\C} \ra \Spec \C$ for $n > 2$, the corresponding Archimedean Siegel--Weil result of \cite{GS19} does not apply here if $n > 2$.

When $m^{\flat} = n$, our Theorem \ref{theorem:intro:results:Archimedean} follows from Liu's result \cite[{Theorem 4.17}]{Liu11}. We do not have a new proof of this case. Instead, we deduce our general result from his by a certain limiting argument. This is also our method at non-Archimedean places (replacing Liu's Archimedean results with the non-Archimedean results of Li--Zhang \cite{LZ22unitary} and Li--Liu \cite{LL22II}). Our limiting method will be sketched further in \cref{ssec:intro:strategy_overview} below.

            \subsection{Previous work}
            \label{ssec:intro:literature}
                We summarize what was previously known on arithmetic Siegel--Weil formulas. These were originally formulated for $\mrm{GSpin}$ Shimura varieties (as opposed to the unitary Shimura varieties considered in \cref{ssec:intro:arith_Siegel-Weil}); we call these the \emph{orthogonal} and \emph{unitary} cases respectively. In both cases, we write $n$ for the arithmetic dimension of the Shimura varieties (i.e. complex dimension $n - 1$).

The problem was initially studied in low-dimensional situations.
For quaternionic Shimura curves, the full arithmetic Siegel--Weil formula has been proved in the influential work of Kudla--Rapoport--Yang \cite{KRY04,KRY06}. For modular curves, the formula has been proved in the papers \cite{Yang04,BF06,DY19,SSY22,Zhu23a,Zhu23b}. 

For Shimura varieties of complex dimension $> 1$, results on arithmetic Siegel--Weil formulas are currently incomplete. Most the available results concern the case $m = n$ and $\det T \neq 0$; we restrict to this case for the moment. Then $s_0 = 0$ is the central point and the special cycle $\mc{Z}(T) \ra \mc{M}$ is empty in the generic fiber. The arithmetic cycle class $[\widehat{\mc{Z}}(T)]$ is thus ``purely vertical'', i.e. either purely in positive characteristic (non-Archimedean), or with $\mc{Z}(T)$ being empty with possibly nontrivial Green current (Archimedean).

The purely Archimedean case (with $\det T \neq 0$ and $s_0 = 0$) was proved by \cite{Liu11,BY21} (unitary and orthogonal, respectively) using different methods. Garcia--Sankaran's Archimedean results apply here as well if the Shimura varieties are compact (more discussion below).

For unitary groups (with $\det T \neq 0$ and $s_0 = 0$), the purely non-Archimedean case for hyperspecial level was first proposed and studied by Kudla--Rapoport \cite{KR11,KR14} at an odd inert prime, where they proved the formula when $\mc{Z}(T)$ has dimension $0$ (reducing locally to the case $n = 2$). The case $n = 3$ at an odd inert prime was solved by Terstiege \cite{Terstiege13}. The case of arbitrary $n$ at odd inert primes was solved in the breakthrough work of Li and Zhang \cite{LZ22unitary} by an inductive ``uncertainty principle'' strategy. This strategy was later adapted to solve the analogous problem at odd ramified primes \cite{LL22II,HLSY22}. We mention that the problem formulation itself needed to be resolved at ramified primes in the presence of bad reduction, and this was done in \cite{HSY22a} for the Kr\"amer model. Split primes play a relatively trivial role when $\det T \neq 0$ and $s_0 = 0$. The timeline for non-Archimedean aspects of the $\mrm{GSpin}$ arithmetic Siegel--Weil formula is similar, i.e. results for $\mc{Z}(T)$ of dimension $0$ were obtained by Kudla--Rapport and Bruinier--Yang \cite{KR99,KR00,BY21}, the case $n = 3$ was resolved by Terstiege \cite{Terstiege11}, and the case of general $n$ at hyperspecial level was resolved by Li and Zhang using (a modified version of) their ``uncertainty principle'' strategy \cite{LZ22orthogonal}.

We now drop the restrictions $\det T \neq 0$ and $s_0 = 0$. For the purpose of arithmetic theta lifting \cref{equation:intro:arith_Siegel-Weil:arith_theta}, it is desirable to also understand the special cycle classes $[\widehat{\mc{Z}}(T)]$ when $\det T = 0$, to fill out the complete arithmetic theta series. Much less is known about this case, which presents new difficulties on both the analytic and geometric sides. It also presents new opportunities: our arithmetic Siegel--Weil result for singular $T$ relates Faltings heights and derivatives of Eisenstein series. Such formulas were observed by Kudla--Rapoport--Yang on Shimura curves \cite{KRY04}; our result applies on unitary Shimura varieties of arbitrarily high dimension. These mixed characteristic phenomena are not visible from arithmetic Siegel--Weil for nonsingular $T$ at the central point $s_0 = 0$ (which was ``purely vertical'').

We mention known partial results for singular $T$, besides the previously mentioned work on Shimura curves. There is concrete progress on the case $T = 0$, where the expected geometric side (``arithmetic volumes'') has been computed for certain levels in the work of H\"ormann and Bruinier--Howard \cite{Hormann14,BH21}, with some partial results on the comparison with Eisenstein series. In the general case, an important advance was made by Garcia and Sankaran \cite{GS19}, who defined Green currents via superconnections and proved a purely Archimedean version of the arithmetic Siegel--Weil formula on compact Shimura varieties (e.g. when $T$ is not positive semi-definite, giving an empty special cycle with possibly nontrivial Green current) via the classical Siegel--Weil formula.

Besides the partial results for $T = 0$, we are not aware of any previous arithmetic Siegel--Weil results which treat non-Archimedean (or combined Archimedean and non-Archimedean) aspects for singular $T$
on Shimura varieties of complex dimension $>1$. This is closely related to the following open problem: for Shimura varieties of complex dimension $> 1$, we are also unaware of any fully global arithmetic Siegel--Weil results (incorporating non-Archimedean places) at a non-central point $s_0 \neq 0$, besides the partial results in \cite[{Theorem C}]{BH21} (there for certain nonzero $1 \times 1$ matrices $T \in \Z$). As discussed at the end of \cref{ssec:intro:arith_Siegel-Weil}, our main theorems make new contributions in both of these directions.

            \subsection{Non-Archimedean local main theorems}
            \label{ssec:intro:non-Arch_local_theorems}
                Our proof of Theorems \ref{theorem:intro:results:main} and \ref{theorem:intro:results:Archimedean} is local in nature. In \cref{ssec:intro:strategy_overview}, we outline the key (new) strategy for proving our main local theorems via our new limit argument. In \cref{ssec:intro:local-to-global}, we outline some of the (new) ideas involved in decomposing our main global theorems into our main local theorems (whose statements and proofs are both new) at every place.

Our local theorems are stated in terms of local special cycles on the Hermitian symmetric domain in the Archimedean case (resp. Rapoport--Zink spaces in the non-Archimedean case). We defer our Archimedean ``local arithmetic Siegel--Weil'' main theorem to our companion paper \cite{corank1_ASW_II.pdf}; additional comparisons between the Archimedean and non-Archimedean cases may be found in \cref{ssec:intro:strategy_overview} and \cref{sec:part_I:sketch}.

The present paper focuses on the non-Archimedean case. 
In the rest of \cref{ssec:intro:non-Arch_local_theorems}, we restrict to the case of an odd prime $p$ inert in $\mc{O}_F$ (for illustration purposes).

Fix an embedding $F \ra \breve{\Q}_p$ into the completion of the maximal unramified extension of $\Q_p$. Set $F_p \coloneqq F \otimes_{\Q} \Q_p$. There is a space of ``local special quasi-homomorphisms'' $\mbf{V}_p$, which is a non-split $F_p / \Q_p$ Hermitian space of dimension $n$. Given any tuple $\underline{\mbf{x}} \in \mbf{V}_p^m$, there is an associated (Kudla--Rapoport) local special cycle $\mc{Z}(\underline{\mbf{x}}) \hookrightarrow \mc{N}$ on a Rapoport--Zink space $\mc{N}$.
These are locally Noetherian formal schemes over $\Spf \breve{\Z}_p$, and they represent certain moduli problems for $p$-divisible groups (\cref{sec:moduli_pDiv}). They appear in Rapoport--Zink uniformization of global special cycles (\crefext{III:sec:non-Arch_uniformization}) in a manner analogous to the complex uniformization of the unitary Shimura varieties and their special cycles. 
There is also an analogous tautological bundle $\mc{E}^{\vee}$ on $\mc{N}$.

There are associated \emph{derived local special cycle (classes)} ${}^{\mbb{L}} \mc{Z}(\underline{\mbf{x}}) \in \mrm{gr}_{\mc{N}}^m K'_0(\mc{Z}(\underline{\mbf{x}}))_{\Q} \coloneqq \mrm{gr}_{n - m} K'_0(\mc{Z}(\underline{\mbf{x}}))_{\Q}$ in codimension $m$. If all elements of the tuple $\underline{\mbf{x}} = [\mbf{x}_1, \ldots, \mbf{x}_m]$ are nonzero, then we have
    \begin{equation}
    {}^{\mbb{L}} \mc{Z}(\underline{\mbf{x}}) \coloneqq \mc{O}_{\mc{Z}(\mbf{x}_1)} \otimes^{\mbb{L}} \cdots \otimes^{\mbb{L}} \mc{O}_{\mc{Z}(\mbf{x}_m)}.
    \end{equation}
For any $x_i = 0$, is it usual to make a similar definition by replacing $\mc{O}_{Z(x_i)}$ with $[\mc{O}_{\mc{N}}] - [\mc{E}]$. For any proper closed subscheme $Z \subseteq \mc{N}$, there is a degree map $\deg_{\overline{\F}_p} \colon \mrm{gr}^n_{\mc{N}} K'_0(Z)_{\Q} \ra \Q$ given by the composition
    \begin{equation}
    \mrm{gr}_0 K'_0(Z)_{\Q} \xra{\sim} \mrm{gr}_0 K'_0(Z_{\overline{\F}_p})_{\Q} \ra \mrm{gr}_0 K'_0(\Spec \overline{\F}_p)_{\Q} = \Q
    \end{equation}
where the first arrow is induced by the d\'evissage pushforward isomorphism $K'_0(Z_{\overline{\F}_p}) \ra K'_0(Z)$ and the second arrow is pushforward along $Z_{\overline{\F}_p} \ra \Spec \overline{\F}_p$ (e.g. induced by taking Euler characteristics of coherent sheaves on $Z_{\overline{\F}_p}$).

Our main ``local arithmetic Siegel--Weil'' result (at an odd inert prime $p$) is the following.
    \begin{theorem*}[Non-Archimedean inert local version of \cref{theorem:intro:results:main}]
    For any nonsingular $T^{\flat} \in \mrm{Herm}_{n - 1}(\Q_p)$ and any $\underline{\mbf{x}}^{\flat} \in \mbf{V}_p^{n - 1}$ with Gram matrix $T^{\flat}$, we have
    \begin{equation}\label{equation:intro:strategy:inert_local_theorem}
    - \frac{d}{ds} \bigg|_{s = 1/2} W^*_{T^{\flat},p}(s)^{\circ}_n 
    = \left ( 2 \deg_{\overline{\F}_p}(\mc{E}^{\vee} \cdot {}^{\mbb{L}} \mc{Z}(\underline{\mbf{x}}^{\flat})_{\ms{V}}) 
    + 2 \sum_{\mc{Z} \hookrightarrow \mc{Z}(\underline{\mbf{x}}^{\flat})_{\ms{H}}} \deg(\mc{Z}) \cdot \delta_{\mrm{tau}}(\mc{Z}) \right) \cdot \log p.
    \end{equation}
    \end{theorem*}
The left-hand side is the derivative of a certain normalized (non-Archimedean) local Whittaker function $W^*_{T^{\flat},p}(s)^{\circ}_n$. On the right, the notation ${}^{\mbb{L}} \mc{Z}(\underline{\mbf{x}}^{\flat})_{\ms{V}} \in \mrm{gr}^{n - 1}_{\mc{N}} K'_0(\mc{Z}(\underline{\mbf{x}}^{\flat})_{\overline{\F}_p})_{\Q}$ denotes a ``derived vertical local special cycle class'' and $\mc{Z}(\underline{\mbf{x}}^{\flat})_{\ms{H}}$ denotes the flat part of $\mc{Z}(\underline{\mbf{x}}^{\flat})$; the former is finite flat over $\Spf \breve{\Z}_p$. The sum runs over all irreducible components $\mc{Z}$ of (the scheme associated to) $\mc{Z}(\underline{\mbf{x}}^{\flat})_{\ms{H}}$. The quantity $\delta_{\mrm{tau}}(\mc{Z}) \in \Q$ is a certain ``local change of tautological height'' which arises from the reduction process from (global) mixed characteristic heights to local quantities.\footnote{Each $\mc{Z}$ is associated with a quasi-canonical lifting of some level $s \in \Z_{\geq 0}$ in the sense of Gross \cite{Gross86}; see \cref{ssec:can_and_qcan:qcan_cycles}. Our notation $\delta_{\mrm{tau}}(\mc{Z})$ here is the $\delta_{\mrm{tau}}(s)$ in \cref{equation:can_and_qcan:qcan:local_change_taut}.} The definition of $\delta_{\mrm{tau}}(\mc{Z})$ is somewhat involved, but some additional discussion may be found in \cref{ssec:intro:local-to-global}.

For a more concrete simple case, see \cref{example:intro:strategy_overview:n=2_sketch} in the next section. For the full precise formulation of our non-Archimedean local main theorems, we refer to \cref{sec:non-Arch_identity} in the body of this paper (there stated in terms of local densities) where the inert/split/ramified cases are treated in parallel.
        
            \subsection{Strategy}
            \label{ssec:intro:strategy_overview}
                We now describe our key local strategy: ``take a limit'' (Figure \ref{figure:intro:strategy:local_limit_method}).

In Figure \ref{figure:intro:strategy:local_limit_method} below, for a given place $v$ of $\Q$, we consider $T^{\flat} \in \mrm{Herm}_{n - 1}(\Q_v)$ with $\det T^{\flat} \neq 0$, and $T = \mrm{diag}(t, T^{\flat})$ for suitable nonzero $t \in \Q_v$. On the left, the limit refers to $t \ra 0$ in the $v$-adic topology (meaning the real topology if $v = \infty$).
The upper horizontal arrow should be understood as a local version of Theorem \ref{theorem:intro:results:main}(2), and the lower horizontal arrow should be understood as the (known) local version of \eqref{equation:intro:results:conjecture} when $\det T \neq 0$.

\begin{figure}[htb]
\centering
\begin{tikzpicture}[>=Stealth, node distance=3cm and 10cm, on grid, auto]
    \node (rect1) [rectangle, draw, text width=5cm, align=center] {Normalized local Whittaker functions $W^*_{T^{\flat},v}(s)^{\circ}_n$ for $U(n - 1, n - 1)$ at $s = 1/2$ (near-center)};
    \node (rect2) [rectangle, draw, text width=5cm, align=center, right=of rect1] {Local $1$-cycles, ``heights''};
    \node (rect3) [rectangle, draw, text width=5cm, align=center, below=of rect2] {Local $0$-cycles, degrees};
    \node (rect4) [rectangle, draw, text width=5cm, align=center, left=of rect3] {Normalized local Whittaker functions $W^*_{T,v}(s)^{\circ}_n$ for $U(n,n)$ at $s = 0$ (center)};

    \draw[->, dotted] (rect1) -- node[right, above] {Our main} node[right, below] {local theorems} (rect2);
    \draw[->, decorate, decoration={snake, amplitude=0.5mm, segment length=5mm}] (rect3) -- node[above, rotate=90] {limit} (rect2);
    \draw[->] (rect4) -- node[right, above] {Known local theorems} node[right, below] {\cite{Liu11,LZ22unitary,LL22II}} (rect3);
    \draw[->, decorate, decoration={snake, amplitude=0.5mm, segment length=5mm}] (rect4) -- node[above, rotate=90] {limit} (rect1);
\end{tikzpicture}
\caption{A local limiting method} \label{figure:intro:strategy:local_limit_method}
\end{figure}
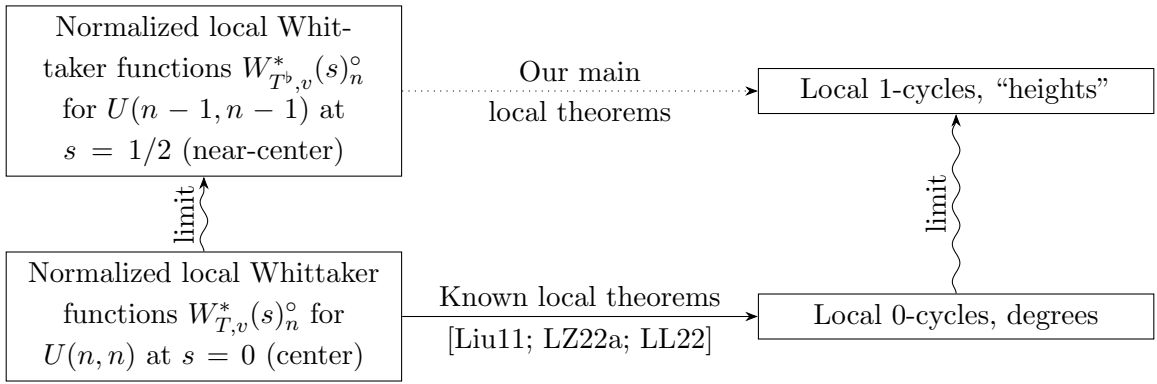
This limiting method is the main conceptual novelty in our work, and is the the key idea driving our main results. 
In Figure \ref{figure:intro:strategy:local_limit_method}, the left vertical arrow and upper horizontal arrow are new in this work. In the right vertical arrow, the relation between limits and Faltings heights is also new in this work.

It is striking that the limiting method plays a similar role at all places, Archimedean and non-Archimedean. In the purely Archimedean case, i.e. when $v = \infty$ with $T^{\flat}$ nonsingular and not positive definite, we are able to run our limiting argument for special cycles (currents) in arbitrary dimension. This is why our purely Archimedean result (Theorem \ref{theorem:intro:results:Archimedean}) applies in arbitrary codimension.

For non-Archimedean places, it is interesting to ask whether the limiting method in Figure \ref{figure:intro:strategy:local_limit_method} can be adapted to the case of higher dimensional special cycles (corresponding to $T^{\flat}$ of smaller rank). Some key difficulties are mentioned in Remark \ref{remark:intro:strategy_overview:higher_corank_difficulties}.

We also mention a slight difference if $v = p$ is a prime split in $\mc{O}_F$. The known local theorems \cite{Liu11,LZ22unitary,LL22II} apply in the Archimedean, inert, and ramified cases respectively. In the split case, the lower left corner of Figure \ref{figure:intro:strategy:local_limit_method} will involve the special value $W^*_{T,v}(0)^{\circ}_n$ (while the derivative at $s = 0$ appears in the Archimedean, inert, and ramified cases). In the split case, the ``known local theorem'' in Figure \ref{figure:intro:strategy:local_limit_method} refers to a certain vanishing statement for a certain contribution to $W^*_{T,v}(0)^{\circ}_n$ (``vertical part'' via an analogue of Cho--Yamauchi's formula; the vanishing is proved in \cref{lemma:non-Arch_identity:limits:split_vertical_den_vanishes}) and emptiness of local $0$-cycles.

In the next example, we sketch the strategy from \cref{figure:intro:strategy:local_limit_method} in a simple non-Archimedean case. The reader seeking a more detailed sketch of a more general setup (along with a comparison between the Archimedean and non-Archimedean strategies) may refer to \cref{sec:part_I:sketch}.

\begin{example}\label{example:intro:strategy_overview:n=2_sketch}
Suppose $p$ is a prime which is inert in $\mc{O}_F$. Take $T^{\flat}$ and $T$ as in \cref{figure:intro:strategy:local_limit_method}, and assume $t$ is such that $T$ defines a non-split Hermitian space. We prove a limiting formula (inert case)
    \begin{equation}\label{equation:intro:strategy:inert:Whittaker_limit}
    \frac{d}{ds} \bigg|_{s = -1/2} W^*_{T^{\flat},p}(s)^{\circ}_n = \lim_{t \ra 0} \left ( \frac{d}{ds} \bigg|_{s = 0} W^*_{T,p}(s)^{\circ}_n + (\log |t|_{p} - \log p) W^*_{T^{\flat},p}(-1/2)^{\circ}_n \right ),
    \end{equation}
which appears in the text as Proposition \ref{proposition:non-Arch_identity:limits} (there stated via local densities). This is the left vertical arrow in \cref{figure:intro:strategy:local_limit_method}. Here $|-|_p$ is the usual $p$-adic norm.

The right vertical arrow in \cref{figure:intro:strategy:local_limit_method} asserts that the analytic limit formula in \cref{equation:intro:strategy:inert:Whittaker_limit} has a geometric interpretation in terms of the local special cycles $\mc{Z}(\underline{\mbf{x}})$.

To illustrate a relatively simple case, consider the case $n = 2$ (and $p \neq 2$). We then have $\mc{N} \cong \Spf \psring{\breve{\Z}_p}{u}$ (non-canonically). Consider
    \begin{equation}
    T^{\flat} = \begin{pmatrix} p \end{pmatrix} \quad \quad T = \diag(t, T^{\flat}) \quad \quad t = p^e
    \end{equation}
for even integers $e \in \Z_{\geq 0}$, where $T^{\flat}$ is a $1 \times 1$ matrix. We have
    \begin{align}
    W^*_{T^{\flat},p}(s)^{\circ}_2 = p^{s+1/2} + p^{-s+1/2} \quad \quad W^*_{T,p}(s)^{\circ}_2 = p^{(e+1)s} - p^{-(e+1)s} + (1-p) p^{(e+1)s} \sum_{i = 1}^e (-q^{-2s})^i. \notag
    \end{align}
    
Set $\mc{O}_{F_p} \coloneqq \mc{O}_F \otimes_{\Z} \Z_p$. With notation as sketched in \cref{ssec:intro:non-Arch_local_theorems}, let $\underline{\mbf{x}} = [\mbf{x}, \mbf{x}^{\flat}] \in \mbf{V}_p^n$ be any tuple with Gram matrix $(\underline{\mbf{x}}, \underline{\mbf{x}}) = T$. Then $\mc{Z}(\mbf{x}^{\flat}) \cong \Spf \mc{O}_{\breve{E}}$ for a certain degree $p + 1$ extension $\breve{E} / \breve{\Q}_p$. Over $\mc{Z}(\mbf{x}^{\flat})$, the universal $p$-divisible group from $\mc{N}$ pulls back to $\mf{X}_1 \otimes_{\Z_p} \mc{O}_{F_p}$ (Serre tensor construction) where $\mf{X}_1 \ra \Spf \mc{O}_{\breve{E}}$ is a \emph{quasi-canonical lifting} of level $1$ in the sense of Gross \cite{Gross86}. If $\mf{X}_0$ denotes the canonical lifting, there is an isogeny $\psi_1 \colon \mf{X}_0 \ra \mf{X}_1$ of minimal degree (unique up to $\mc{O}_{F_p}^{\times}$), with $\deg \psi_1 = p$. We have
    \begin{align}
    \breve{\delta}_{\mrm{Fal}}(\psi_1) & \coloneqq \frac{1}{2} \log(\deg \psi_1) - \frac{1}{[\breve{E} : \breve{\Q}_p]} \mrm{length}_{\mc{O}_{\breve{E}}}(e^* \Omega^1_{\ker \psi_1 / \mc{O}_{\breve{E}}}) \log p 
    \\
    & = \left ( \frac{1}{2} - \frac{1}{p + 1} \right ) \log p
    \end{align}
where the second equality follows from a computation of Nakkajima--Taguchi \cite{NT91} (with $e^*$ denoting pullback along the identity section). We call the left-hand side a ``local change of Faltings height''; its relation with (global) Faltings height is sketched in \cref{ssec:intro:local-to-global}.

In this case, our geometric analogue of \cref{equation:intro:strategy:inert:Whittaker_limit} is the limit formula
    \begin{equation}\label{equation:intro:strategy:inert:geometric_limit}
    - \breve{\delta}_{\mrm{Fal}}(\psi_1) \cdot \deg \mc{Z}(\mbf{x}^{\flat}) = \lim_{\mbf{x} \ra 0} \left ( (\deg_{\overline{\F}_p} (\mc{O}_{\mc{Z}(\mbf{x})} \otimes^{\mbb{L}} \mc{O}_{\mc{Z}(\mbf{x}^{\flat})})) \cdot \log p - \frac{1}{2} (\log |t|_p - \log p) \cdot \deg \mc{Z}(\mbf{x}^{\flat}) \right )
    \end{equation}
where the limit $t \ra 0$ is $p$-adic.
We have
    \begin{equation}\label{equation:intro:strategy:inert:LZ22unitary}
    \frac{1}{2} \frac{d}{ds} \bigg|_{s = 0} W^*_{T,p}(s)^{\circ}_2 \overset{\text{\cite{LZ22unitary}}}{=} (\deg_{\overline{\F}_p} {}^{\mbb{L}} \mc{Z}(\underline{\mbf{x}})) \cdot \log p \quad \quad W^*_{T^{\flat},p}(-1/2)^{\circ}_2 \overset{\text{\cite{LZ22unitary}}}{=} \deg \mc{Z}(\mbf{x}^{\flat}).
    \end{equation}
These may be thought of as (nonsingular, central point) ``local arithmetic Siegel--Weil'' and (nonsingular, near-central) ``local geometric Siegel--Weil'' formulas, respectively.

Applying our limiting formulas in \cref{equation:intro:strategy:inert:Whittaker_limit,equation:intro:strategy:inert:geometric_limit} (along with the local functional equation $W^*_{T^{\flat},p}(s)^{\circ}_2 = W^*_{T^{\flat},p}(-s)^{\circ}_2$) produces the formula
    \begin{equation}
    \frac{d}{ds} \bigg|_{s = 1/2} W^*_{T^{\flat},p}(s)^{\circ}_2 = 2 \breve{\delta}_{\mrm{Fal}}(\psi_1) \cdot \deg \mc{Z}(\mbf{x}^{\flat}).
    \end{equation}
This is one form of our main local theorem in this simple special case (compare the global version, Theorem \ref{theorem:intro:results:main}(2), and the more general local formulation in \cref{ssec:intro:non-Arch_local_theorems}). 

In extremely impressionistic terms, the formula in \cref{equation:intro:strategy:inert:geometric_limit} states that the ``limit'' of the special divisor $\lim_{\mbf{x} \ra 0} \mc{Z}(\mbf{x})$ ``converges'' to some ``local part'' of a tautological bundle (essentially the Hodge bundle) when intersecting against $\mc{Z}(\mbf{x}^{\flat})$, after regularizing. 
We remark that the analogous ``numerical limit'' statement (without regularizing) is literally true if $\mc{Z}(\mbf{x}^{\flat})$ is replaced by any (proper) curve in the special fiber over $\mc{N}$. That case has a conceptual explanation: Grothendieck--Messing deformation theory. This discussion is continued in \cref{equation:intro_strategy:non-Arch_limit:vertical} and the text below in loc. cit..
\end{example}

            \subsection{Local-to-global}
            \label{ssec:intro:local-to-global}
                For the global-to-local reduction process, we use complex uniformization (Archimedean place) and Rapoport--Zink uniformization (non-Archimedean places). Unlike the previously known case $\det T \neq 0$ for $T \in \mrm{Herm}_n(\Q)$ (giving a purely vertical arithmetic special cycle class), we have a new mixed characteristic ``horizontal'' contribution in the Faltings height. While the Faltings height decomposes locally \emph{after picking a section} of the metrized Hodge bundle, it does not admit an obvious \emph{canonical} local decomposition. Such a canonical decomposition seems necessary for the comparison with local Whittaker functions (as appearing in Eisenstein series Fourier coefficients), which presumably does not retain information on which section was picked.

Instead, we do have a decomposition for the \emph{difference} between (stable) Faltings heights of any two abelian varieties $A_1, A_2$ (over a number field $E$) in a fixed isogeny class.
We may assume $A_1$ and $A_2$ have everywhere semi-abelian reduction after extending $E$. The difference of Faltings heights is then
    \begin{align}\label{equation:intro:local-to-global:change_Faltings_height}
    h_{\mrm{Fal}}(A_2) - h_{\mrm{Fal}}(A_1) & = \sum_{p \mid \deg \phi} a_p \log p = \frac{1}{[E : \Q]} \sum_p \sum_{\breve{w} \mid p} [\breve{E}_{\breve{w}} : \breve{\Q}_p] \breve{\delta}_{\mrm{Fal}}(\phi_{\breve{w}})
    \\
    \breve{\delta}_{\mrm{Fal}}(\phi_{\breve{w}}) & \coloneqq \frac{1}{2} \log(\deg \phi_{\breve{w}}) - \frac{1}{[\breve{E}_{\breve{w}} : \breve{\Q}_p]} \mrm{length}_{\mc{O}_{\breve{E}_{\breve{w}}}}(e^* \Omega^1_{\ker \phi/ \mc{O}_{\breve{E}_{\breve{w}}}}) \log p
    \label{equation:intro:local-to-global:change_Faltings_height_local}
    \end{align}
for some $a_p \in \Q$ and any choice of isogeny $\phi \colon A_1 \ra A_2$ (extend $\phi$ over $\mc{O}_E$). Here, the symbol $e$ means the identity section, the inner sum in \cref{equation:intro:local-to-global:change_Faltings_height} runs over all prime ideals $\breve{w}$ of $E \otimes_{\Q} \breve{\Q}_p$ (with associated residue field $\breve{E}_{\breve{w}}$), and $\phi_{\breve{w}}$ means the associated isogeny of $p$-divisible groups $A_1[p^{\infty}] \ra A_2[p^{\infty}]$ base-changed to $\Spf \mc{O}_{\breve{E}_{\breve{w}}}$.

The coefficients $a_p$ do not depend on the choice of $\phi$, by linear independence of $\log p$ for different $p$. This is also the reason why a difference of Faltings heights appears in Corollary \ref{corollary:intro:results:n=2}. 
We then argue that these numbers $a_p \in \Q$ (averaged over the special cycle) can be calculated in a purely local way, in terms of local special cycles on Rapoport--Zink spaces. The argument we give is somewhat delicate, as we wish to avoid writing down explicit (global) isogenies $\phi \colon A_1 \ra A_2$, so that we obtain a more local formulation. The quantities $\delta_{\mrm{Fal}}(\phi_{\breve{w}})$ are amenable to (local) calculation on Rapoport--Zink spaces, but may depend on the chosen isogeny of $p$-divisible groups $\phi_{\breve{w}}$. We show that certain ``minimal degree'' isogenies $\phi_{\breve{w}}$ lift to global isogenies of $p$-divisible groups (\cref{ssec:qcan_heights:minimal_isogenies}) and can be used to calculate local contributions to (differences of) Faltings height.

The global-to-local reduction process from $h_{\mrm{Fal}}(-)$ (Faltings heights) to $\delta_{\mrm{Fal}}(\phi_{\breve{w}})$ (local quantities calculate-able on Rapoport--Zink spaces) for ``minimal degree'' $\phi_{\breve{w}}$ is the content of \cref{sec:Faltings_and_taut,sec:qcan_heights}. The relation with global special cycles via uniformization is explained in \crefext{III:sec:non-Arch_uniformization,III:sec:Arch_uniformization}. There is also the issue that the tautological bundle $\widehat{\mc{E}}^{\vee}$ is not the same as the metrized Hodge bundle (but is known to behave similarly, as first observed by Gross \cite{Gross78} and studied further in \cite{BHKRY20II}), so the (more natural) version with ``tautological height'' needs additional argument. The ``tautological height'' and Faltings height are treated in parallel in \cref{sec:Faltings_and_taut,sec:qcan_heights}.

While previous work for special points on Shimura curves \cite{KRY04} also studied the change in Faltings heights along isogenies, our insistence on a purely local formulation is an important difference for our method. We only observe the limiting phenomena in Figure \ref{figure:intro:strategy:local_limit_method} (below) on a local level; this is what allows us to prove a theorem on Shimura varieties of arbitrarily large dimension.

\begin{remark}\label{remark:intro:strategy_overview:higher_corank_difficulties}
After our work, the arithmetic Siegel--Weil formula in \eqref{equation:intro:results:conjecture} remains open only for $T$ with $\rank(T) \leq n - 2$ (corresponding to special cycles $\mc{Z}(T)$ of dimension $\geq 1$ in the generic fiber, if nonempty). We mention some of the difficulties for these higher dimensional special cycles, from our perspective.

Our strategy in Figure \ref{figure:intro:strategy:local_limit_method} strongly emphasizes \emph{local} limits on both the analytic and geometric sides. On the geometric side at non-Archimedean places, this is possible in the rank $\geq n - 1$ case because e.g. the special cycles are contained inside supersingular loci at all nonsplit primes. This allows us to describe global special cycles in terms of local special cycles on a single Rapoport--Zink space (at each prime).

Higher dimensional special cycles may pass through several strata, so a single Rapoport--Zink space should not be enough to capture all information. We are not sure whether it is sensible to ``piece together'' the intersection-theoretic information coming from several Rapoport--Zink spaces. Moreover, the corresponding local special cycles on Rapoport--Zink spaces may be non-proper (so it is unclear how to extract local intersection numbers in the more general situation).

For our method, it is also important to locally decompose Faltings heights in a canonical way. In the rank $n - 1$ case, this was accomplished using the ``change of Faltings height along an isogeny'' formula. To generalize our limiting strategy to higher dimensional cycles, we may need a similar canonical local decomposition involving heights of higher dimensional cycles on the Shimura variety.
\end{remark}

            \subsection{Outline}
            \label{ssec:intro:outline}
                We briefly summarize the remaining content in this paper, and discuss the relation with our companion papers \cite{corank1_ASW_II.pdf,corank1_ASW_III.pdf,corank1_ASW_IV.pdf}.
Further explanations may be found at the beginning of some sections.

In \cref{sec:part_I:arith_intersections}, we define the (global) objects appearing in the ``arithmetic intersection numbers without boundary contributions'' considered in \cref{equation:intro:results:if_proper}. Some technical properties of global special cycles are only stated in this paper, with corresponding proofs appearing in our companion paper \cite{corank1_ASW_III.pdf}. In \cref{ssec:part_I:arith_intersections:vertical_classes}, we complete the ``vertical'' (positive characteristic) part of our proposed construction for arithmetic special cycle classes in \cref{equation:intro:results:arith_cycle_class}.

We focus on the local situation from \cref{sec:moduli_pDiv} onwards in this paper, and return to the global special cycles in our companion paper \cite{corank1_ASW_III.pdf}. Nevertheless, we have included these global definitions as they underlie the formulations of our main local theorems. 

The remaining sections are divided into \cref{part:local_cycles,part:local_main_theorems,part:local_change_heights} and an appendix. 
In \cref{part:local_cycles} ``Local special cycles'', we set up Kudla--Rapoport local special cycles on Rapoport--Zink spaces (inert/ramified/split). The case of split primes is less well-studied in the literature than the inert/ramified cases (we need uniformization in a non supersingular situation at split primes). Section \ref{sec:more_moduli_pDiv_split} contains some new results on decomposing local special cycles into quasi-canonical lifting cycles at split primes, which we need later. These are analogous to known results at inert and ramified primes (Section \ref{ssec:can_and_qcan:qcan_cycles}), though our method of proof is different.

In \cref{part:local_main_theorems} ``Local main theorems'', we formulate and prove our main local non-Archimedean theorems. 
\Cref{sec:non-Arch_identity} in particular contains the heart of our work. Here, we prove our ``local arithmetic Siegel--Weil'' theorems at inert/ramified/split non-Archimedean places via the local limiting method sketched in \cref{ssec:intro:strategy_overview}. The Archimedean analogue will appear in our companion paper \cite{corank1_ASW_II.pdf}; for comparison purposes, a brief sketch of the Archimedean strategy is included in \cref{ssec:part_I:sketch:Archimedean}.

In \cref{part:local_change_heights} ``Local change of heights'', we give a key argument in our reduction process from global heights in mixed characteristic to local quantities. We formulate these local quantities to be computable in terms of local special cycles and certain ``local change of heights'' along ``minimal isogenies'' (sketched in \cref{ssec:intro:strategy_overview}). The ``local change of height'' constants appearing in our non-Archimedean ``local arithmetic Sigel--Weil'' main theorems emerge from this global-to-local reduction process. \Cref{part:local_change_heights} is included mainly to motivate the formulation of our non-Archimedean local main theorems.

Appendix \ref{appendix:K0} explains the setup we use for $K_0$ groups of Deligne--Mumford stacks. This setup is roughly as in \cite[{Appendix A}]{YZ17}, but they work over a field (we need schemes over $\mc{O}_F$, or at least over discrete valuation rings).

In our companion paper \cite{corank1_ASW_III.pdf}, we combine the material from \cref{part:local_change_heights} ``Local change of heights'' (in the present paper) with Rapoport--Zink (non-Archimedean) uniformization and complex (Archimedean) uniformization to complete the reduction process from the global intersection numbers \cref{equation:intro:results:if_proper} to the geometric quantities appearing in our main local theorems from the present paper and our companion paper \cite{corank1_ASW_II.pdf}. In our companion paper \cite{corank1_ASW_IV.pdf}, we make an analogous local-to-global reduction for Eisenstein series and show how to patch together our local arithmetic Siegel--Weil theorems at all places to finish the proof of our main global arithmetic Siegel--Weil theorems (Theorems \ref{theorem:intro:results:main} and \ref{theorem:intro:results:Archimedean}); see \crefext{IV:sec:arithmetic_Siegel-Weil}. This final proof will use almost all preceding results from our four-paper sequence.

Our algebro-geometric conventions follow the Stacks project \cite{stacks-project} unless stated otherwise.

            \subsection{Acknowledgements}
            \label{acknowledgements}
                I thank my advisor Wei Zhang for suggesting this topic, for his dedicated support and constant enthusiasm, for insightful discussions throughout the entire course of this project, and for helpful comments on earlier drafts. I thank Tony Feng, Qiao He, Benjamin Howard, Ishan Levy, Chao Li, Keerthi Madapusi, Andreas Mihatsch, Siddarth Sankaran, Ananth Shankar, Yousheng Shi, Tonghai Yang, Shou-Wu Zhang, and Zhiyu Zhang for helpful comments or discussions.

This work was partly supported by the National Science Foundation Graduate Research Fellowship under Grant Nos. DGE-1745302 and DGE-2141064.
Parts of this work were completed at the Mathematical Sciences Research Institute (MSRI), now becoming the Simons Laufer Mathematical Sciences Institute (SLMath), and the Hausdorff Institute for Mathematics. I thank these institutes for their support and hospitality. The former is supported by the National Science Foundation (Grant No. DMS-1928930), and the latter is funded by the Deutsche Forschungsgemeinschaft (DFG, German Research Foundation) under Germany's Excellence Strategy – EXC-2047/1 – 390685813.
        
        \section{Conventions on Hermitian spaces and lattices}
        \label{sec:Hermitian_conventions}
            We will use the conventions from \cref{sec:Hermitian_conventions} throughout all of our four-paper sequence (i.e. this paper and its companions \cite{corank1_ASW_II.pdf,corank1_ASW_III.pdf,corank1_ASW_IV.pdf}).
            
            \subsection{Hermitian, alternating, symmetric}
            \label{ssec:Hermitian_conventions:Hermitian_alternating_symmetric}
                Consider a Dedekind domain $\mc{O}_{F_0}$ with fraction field $F_0$. Let $F$ be a finite \'etale $F_0$-algebra of degree $2$, i.e. $F$ is either a degree $2$ separable field extension of $F_0$, or $F = F_0 \times F_0$. Let $\mc{O}_F \subseteq F$ be the integral closure of $\mc{O}_{F_0}$ in $F$. Write $a \mapsto a^{\s}$ for the nontrivial involution of $F$ over $F_0$, and $\mrm{tr} \colon F \ra F_0$ for the trace map $a \mapsto a + a^{\s}$.

Assume that the different ideal $\mf{d}$ of $\mc{O}_F$ over $\mc{O}_{F_0}$ is principal, and choose a generator $u \in \mf{d}$ satisfying $u^{\s} = - u$. This is always possible if $\mc{O}_{F}$ is a free $\mc{O}_{F_0}$-module.

Let $L$ be a finite locally free $\mc{O}_F$-module of constant rank. If $F_0$ has characteristic $\neq 2$, the following data are equivalent.
    \begin{enumerate}[(1)]
        \item A \emph{Hermitian pairing} on $L$, i.e. a $\mc{O}_{F_0}$-bilinear map $(-,-) \colon L \times L \ra F$ satisfying
            \begin{equation}
            (x, a y) = a (x, y) \quad \quad (y, x) = (x, y)^{\s}
            \end{equation}
        for all $a \in \mc{O}_F$ and $x, y \in L$.
        
        \item An \emph{$\mc{O}_F$-compatible alternating pairing} on $L$, i.e. a $\mc{O}_{F_0}$-bilinear map $\langle -,- \rangle \colon L \times L \ra F_0$ satisfying
            \begin{equation}
            \langle a x, y \rangle = \langle x, a^{\s} y \rangle \quad \quad \langle y, x \rangle = - \langle x, y \rangle
            \end{equation}
        for all $a \in \mc{O}_F$ and $x, y \in L$.
        
        \item An \emph{$\mc{O}_F$-compatible symmetric pairing} on $L$, i.e. a $\mc{O}_{F_0}$-bilinear map $\lbrbrak -,- \rbrbrak \colon L \times L \ra F_0$ satisfying 
            \begin{equation}
            \lbrbrak a x, y \rbrbrak = \lbrbrak x, a^{\s} y \rbrbrak \quad \quad \lbrbrak y, x \rbrbrak = \lbrbrak x, y \rbrbrak
            \end{equation}
        for all $a \in \mc{O}_F$ and $x, y \in L$.
    \end{enumerate}
If $L$ is equipped with any of the equivalent data above, we say that $L$ is a \emph{Hermitian $\mc{O}_F$-lattice} (or \emph{Hermitian $\mc{O}_F$-module}). Note that our Hermitian pairings $(-,-)$ are conjugate linear in the first argument. We pass between these pairings using the formulas (depending on the choice of $u$)
    \begin{align*}
    & 2 (x,y) = \lbrbrak x, y \rbrbrak - u^{-1} \lbrbrak u x, y \rbrbrak
    \quad
    && \langle x, y \rangle = \lbrbrak u^{-1} x, y \rbrbrak
    \quad
    && \lbrbrak x, y \rbrbrak = \operatorname{tr}((x,y))    
    \\
    & 2 (x,y) = \langle u x, y \rangle - u \langle x, y \rangle
    \quad
    && \langle x, y \rangle = -\operatorname{tr}((x, y) u^{-1})
    \quad
    && \lbrbrak x, y \rbrbrak = \langle u x, y \rangle
    \end{align*}
and this will be freely used in the paper. The choice of $u$ plays a limited role for us, so we generally suppress it.

We say that $\lbrbrak -,- \rbrbrak$ is the associated \emph{trace pairing}, and otherwise avoid the notation $\lbrbrak -,- \rbrbrak$ outside of Section \ref{ssec:Hermitian_conventions:Hermitian_alternating_symmetric}.

Given any tuple $\underline{x} = [x_1, \ldots, x_m] \in L^m$, its \emph{Gram matrix}
is the matrix $(\underline{x}, \underline{x}) = T$ with $i,j$-th entry $T_{i,j} = (x_i, x_j)$.
We write $L_F \coloneqq L \otimes_{\mc{O}_F} F$ and say that a Hermitian $\mc{O}_F$-module $L$ is \emph{non-degenerate} if the Gram matrix for any $F$-basis of $L_F$ has nonzero determinant. Given non-degenerate Hermitian $F$-modules $V$ and $V'$ with Hermitian pairings $(-,-)$ and $(-,-)'$, there is a canonical $\s$-linear involution of $F$-modules
    \begin{equation}
    \begin{tikzcd}[row sep = tiny]
    \Hom_F(V, V') \arrow{r}{f \mapsto f^{\dagger}} & \Hom_F(V', V)    \\
    \end{tikzcd}
    \quad 
    \begin{array}{l}
    \text{such that } (f x, y')' = (x, f^{\dagger} y')
    \\
    \text{for all $x \in V$ and $y' \in V'$.}
    \end{array}
    \end{equation}
The notation $\Hom_F(V,V')$ and $\Hom_F(V', V)$ does not include any requirement on preserving Hermitian pairings.

Given a non-degenerate Hermitian $\mc{O}_F$-lattice $L$, we always form its \emph{dual lattice} $L^{*}$ with respect to the trace pairing $\lbrbrak -,- \rbrbrak$, i.e.
    \begin{equation}
    L^{*} \coloneqq \{ x \in L_F : \mrm{tr}(x,y) \in \mc{O}_{F_0} \text{ for all $y \in L$} \}.
    \end{equation}
The dual lattice $L^{\vee}$ with respect to $(-,-)$ is the same as the dual lattice for $\langle -,- \rangle$. We have $L^{\vee} = u L^{*}$ (as sublattices of $L_F$).
If the dual $L^{\vee}$ with respect to $(-,-)$ or $\langle -,- \rangle$ is intended, we will state this explicitly. We say that $L$ is \emph{self-dual} if $L = L^{*}$.

As a typical example of passing between $(-,-)$ and $\langle -,- \rangle$, suppose $\mc{O}_{F_0} = \Z$ and suppose $\mc{O}_F$ is the ring of integers in an imaginary quadratic field $F / \Q$. Let $(A, \iota, \lambda)$ be a \emph{Hermitian abelian variety} (Definition \ref{definition:ab_var:integral_models:Hermitian_abelian_scheme}) over an algebraically closed field $k$ of characteristic $\neq p$, i.e. $A$ is an abelian variety over $k$ with an action $\iota \colon \mc{O}_{F} \rightarrow \End(A)$, and $\lambda$ is an $\mc{O}_F$-compatible quasi-polarization on $A$. After picking a trivialization $\Z_p(1) \cong \Z_p$ of $p$-th power roots of unity over $k$, the polarization $\lambda$ induces an $(\mc{O}_F \otimes_{\Z} \Z_p)$-compatible alternating pairing on the Tate module $T_p(A)$, so we automatically view $T_p(A)$ as a Hermitian $(\mc{O}_F \otimes_{\Z} \Z_p)$-lattice without further mention. If $(A', \iota', \lambda')$ is another such Hermitian $p$-divisible group, note that the induced Hermitian pairing on $\Hom(T_p(A), T_p(A'))$ does not depend on the choice of trivialization $\Z_p(1) \cong \Z_p$ or the choice of $u$.

The notation $\mrm{Herm}_n(\mc{O}_{F_0})$ means the set of $n \times n$ Hermitian matrices with coefficients in $\mc{O}_F$ (i.e. $T \in M_{n,n}(\mc{O}_F)$ satisfying $T = {}^t \overline{T}$ where ${}^t \overline{T}$ means conjugate transpose). Here we are considering the subfunctor $\mrm{Herm}_n \subseteq \Res_{\mc{O}_F / \mc{O}_{F_0}} M_{n,n}$ of the Weil restriction (of $n \times n$ matrices $M_{n,n}$). We adhere strictly to this notation (when $\mc{O}_F$ is understood), e.g. $\mrm{Herm}_n(\R)$ will typically mean $n \times n$ complex Hermitian matrices when $\mc{O}_F / \mc{O}_{F_0} = \C / \R$ is understood.
        
            \subsection{Lattices for local fields}
            \label{ssec:Hermitian_conventions:lattices}
                Continuing in the setup of Section \ref{ssec:Hermitian_conventions:Hermitian_alternating_symmetric}, suppose $F_0$ is a local field. Let $\eta \colon F_0^{\times} \ra \{ \pm 1\}$ be the character associated to $F / F_0$ by local class field theory. Given a non-degenerate Hermitian $F$-module $V$ of rank $n$, define its \emph{local invariant}
    \begin{equation}
    \varepsilon(V) \coloneqq \eta((-1)^{n(n - 1)/2} \det T) \in \{ \pm 1 \}
    \end{equation}
where $T$ is the Gram matrix of any basis for $V$. This is normalized so that $\varepsilon(V) = 1$ for the the Hermitian $F$-module $V$ given by the antidiagonal unit Gram matrix. Rank $n$ non-degenerate Hermitian $F$-modules $V$ and $V'$ are isomorphic if and only if $\varepsilon(V) = \varepsilon(V')$. If $T \in \mrm{Herm}_n(F_0)$ is a Hermitian matrix (with entries in $F$) satisfying $\det T \neq 0$, we set $\varepsilon(T) \coloneqq \eta((-1)^{n(n - 1)/2} \det T)$.

Next, assume $F_0$ is non-Archimedean and that $\mc{O}_{F_0} \subseteq F_0$ is its ring of integers. Write $q$ for the residue cardinality of $\mc{O}_{F_0}$. If $q$ is even, we require $F / F_0$ to be unramified. Let $\varpi_0 \in \mc{O}_{F_0}$ and $\varpi \in \mc{O}_F$ be uniformizers (meaning $\varpi \in \varpi_0 \mc{O}_F^{\times}$ in the unramified cases) satisfying $\varpi^{\s} = - \varpi$. If a non-degenerate Hermitian $F$-module $V$ contains a full rank self-dual $\mc{O}_F$-lattice, then $\varepsilon(V) = 1$.

The ``norm'' $\norm{-}$ on a Hermitian $F$-module $V$ with pairing $(-,-)$ is given by
    \begin{equation}
    \norm{x} \coloneqq q^{-v_{\varpi_0}((x,x))/2}
    \end{equation}
where $v_{\varpi_0}$ is the $\varpi_0$-adic valuation, normalized so that $v_{\varpi_0}(\varpi_0) = 1$.

Given a non-degenerate Hermitian $\mc{O}_F$-lattice $L$ of rank $n$, we set $\varepsilon(L) \coloneqq \varepsilon(L_F)$. By a \emph{lattice} or \emph{sublattice} $L' \subseteq L$, we mean any $\mc{O}_F$-submodule which is finite free of constant rank (similarly for lattices or sublattices in $L_F$). If $L'$ has rank $n$, we say that $L'$ is \emph{full rank} in $L_F$. A sublattice $L' \subseteq L$ is \emph{saturated} if $a x \in L'$ with $a \in F^{\times}$ and $x \in L$ implies $x \in L'$ (equivalently, $L'$ is a direct summand of $L$).

We say that $L$ is \emph{integral} if $L \subseteq L^{*}$. If $F / F_0$ is nonsplit, we say that $L$ is \emph{almost self-dual} if $L \subseteq L^{*}$ and $\operatorname{length}_{\mc{O}_{F}}(L^{*} / L) = 1$. We say that a non-degenerate integral lattice $L$ is \emph{maximal integral} if any integral lattice $L' \subseteq L_F$ with $L \subseteq L'$ satisfies $L = L'$.

If $L$ is a non-degenerate Hermitian $\mc{O}_F$-lattice, we define the \emph{valuation} $\mrm{val}(L) \in \frac{1}{2} \Z$ such that
    \begin{equation}
    q^{-\mrm{val}(L)} = \mrm{vol}(L)
    \end{equation}
where $\mrm{vol}(L)$ is the volume of $L$ for the self-dual Haar measure on $L_F$ with respect to the pairing $x,y \mapsto \psi(\mrm{tr}{(x,y)})$ for any unramified (unitary) additive character $\psi \colon F_0 \ra \C^{\times}$.
If $L$ is integral, we have $q^{2 \mrm{val}(L)} = |L^* / L|$.
If $F / F_0$ is unramified, we have $\mrm{val}(L) \in \Z$. 
Given $x \in L$, we write $\langle x \rangle \subseteq L$ for the rank one $\mc{O}_F$-submodule generated by $x$. If $(x,x) \neq 0$, we set $\mrm{val}(x) \coloneqq \mrm{val}(\langle x \rangle)$ (and otherwise set $\mrm{val}(x) = \infty$).

Continuing to assume $L$ is non-degenerate and integral, we define its sequence of \emph{fundamental invariants} to be the unique sequence of integers $(a_1, \ldots, a_n)$ with $0 \leq a_1 \leq \cdots \leq a_n$ such that $L^* / L \cong \oplus_{i = 1}^n \mc{O}_F / \varpi^{a_i}$ (where $n$ is the rank of $L$). Two non-degenerate integral Hermitian $\mc{O}_F$-lattices of the same rank are isomorphic if and only if they have the same sequence of fundamental invariants (in the unramified case, this follows from diagonalizability of Hermitian lattices; in the ramified case, this follows from \cite[{Proposition 4.3, Proposition 8.1}]{Jacobowitz62} (see also \cite[{Lemma 2.12}]{LL22II})). We set
    \begin{equation}
    t(L) \coloneqq | \{a_i \in \{a_1, \ldots, a_n \} : a_i \neq 0\} | \in \Z \quad \quad a_{\mrm{max}}(L) \coloneq a_n
    \end{equation}
and refer to $t(L)$ as the \emph{type} of $L$. If $F / F_0$ is ramified, recall that $t(L)$, $2 \mrm{val}(L)$, and $n$ all have the same parity (follows from \cite[{Proposition 4.3, Proposition 8.1}]{Jacobowitz62}).

Given a finite length $\mc{O}_F$-module $M$, we define $\ell(M) \in \Z$ such that
    \begin{equation}
    q^{\ell(M)} = |M|.
    \end{equation}
where $|M|$ denotes the cardinality of $M$.

The above terminology is adapted from  e.g. \cite{LZ22unitary} (inert), \cite{FYZ21} (inert and split), \cite{LL22II} (ramified). We made slight modifications to give a uniform description (e.g. our $\mrm{val}(L)$ is half of the $\mrm{val}(L)$ appearing in \cite{LL22II}, and our $\ell(M)$ differs by a factor of $2$ from some of the references).

    \clearpage


        \section{Arithmetic intersection numbers}
        \label{sec:part_I:arith_intersections}
            Fix an imaginary quadratic field extension $F / \Q$ with ring of integers $\mc{O}_F$ and write $a \mapsto a^{\s}$ for the nontrivial automorphism $\s$ of $F$. We write $\Delta \in \Z_{<0}$ and $\sqrt{\Delta} \in \mc{O}_F$ (pick a square root) for (generators of the) discriminant and different, respectively.
    
            \subsection{Integral models}
            \label{ssec:part_I:arith_intersections:integral_models}
                \begin{definition}\label{definition:ab_var:integral_models:Hermitian_abelian_scheme}
Let $S$ be a scheme over $\Spec \mc{O}_F$. By a \emph{Hermitian abelian scheme} over $S$, we mean a tuple $(A, \iota, \lambda)$ where
    \begin{align*}
    & A & & \text{is an abelian scheme over $S$ of constant relative dimension $n$} 
    \\
    & \iota \colon \mc{O}_F \ra \End(A) & & \parbox[t]{0.75\textwidth}{is a ring homomorphism
    } 
    \\
    & \lambda \colon A \ra A^{\vee} & & \parbox[t]{0.75\textwidth}{is a quasi-polarization satisfying:
    \begin{itemize}[label={}]
    \item (Action compatibility) The Rosati involution $\dagger$ on $\End^{0}(A)$ satisfies $\iota(a)^{\dagger} = \iota(a^{\sigma})$ for all $a \in \mc{O}_F$.
    \end{itemize}
    }
    \end{align*}
\end{definition}
An \emph{isomorphism} of Hermitian abelian schemes is an isomorphism of abelian schemes which respects the $\mc{O}_F$-actions and polarizations (exactly). 

Consider any integer $r$ with $0 \leq r \leq n$. Given a pair $(A, \iota)$ where $A \ra S$ is an abelian scheme of relative dimension $n$ for a scheme $S$ over $\Spec \mc{O}_F$, and $\iota \colon \mc{O}_F \ra \End(A)$ is an action, we say that $(A, \iota)$ has \emph{signature $(n - r, r)$} if the following is satisfied:
    \begin{itemize}[label = {}]
    \item (Kottwitz $(n - r, r)$ signature condition) For all $a \in \mc{O}_F$, the characteristic polynomial of $\iota(a)$ acting on $\Lie A$ is $(x - a)^{n - r}(x - a^{\s})^r \in \mc{O}_S[x]$.
    \end{itemize}
Here we viewed $\mc{O}_S$ as an $\mc{O}_F$-algebra via the structure map $S \ra \Spec \mc{O}_F$. We say a Hermitian abelian scheme $(A, \iota, \lambda)$ has signature $(n - r, r)$ if $(A, \iota)$ has signature $(n - r, r)$.

Let
    \begin{equation}
    \ms{M}_0 \ra \Spec \mc{O}_F
    \end{equation}
be the moduli stack of signature $(1,0)$ CM elliptic curves. More precisely, if $S$ is a scheme over $\Spec \mc{O}_F$, then $\ms{M}_0(S)$ is the groupoid of signature $(1,0)$ Hermitian abelian schemes $(A_0, \iota_0, \lambda_0)$ where $A_0$ is relative $1$-dimensional and the quasi-polarization $\lambda_0$ is a principal polarization. Then $\ms{M}_0$ is a Deligne--Mumford stack and the structure morphism $\ms{M}_0 \ra \Spec \mc{O}_F$ is proper, quasi-finite,\footnote{Following the Stacks project \cite[\href{https://stacks.math.columbia.edu/tag/0CHU}{Definition 0CHU}]{stacks-project}, we require that finite morphisms of algebraic stacks are by definition (relatively) representable by schemes.
The morphism $\ms{M}_0 \ra \Spec \mc{O}_F$ is not finite in this sense, because $\ms{M}_0$ is not a scheme. Nevertheless, we continue to use terminology like ``representable by schemes and finite'' for morphisms of stacks which are not necessarily algebraic.} and \'etale by \cite[{Proposition 3.1.2}]{Howard12} or \cite[{Proposition 2.1.2}]{Howard15}.

We now require $2 \nmid \Delta$ and assume $n$ is even, for the rest of Section \ref{sec:part_I:arith_intersections}. We write $\ms{M}(n - r, r)^{\mrm{Kot},\circ} \ra \Spec \mc{O}_F$ for the moduli stack of signature $(n - r, r)$ Hermitian abelian schemes subject to the following condition:
    \begin{enumerate}[label = {}]
        \item (Polarization condition $\circ$) The quasi-polarization $|\Delta| \cdot \lambda$ is a polarization and we have $\ker(|\Delta| \cdot \lambda) = A[\sqrt{\Delta}]$.
    \end{enumerate}
Then $\ms{M}(n - r, r)^{\mrm{Kot},\circ}$ is a Deligne--Mumford stack which is separated and finite type over $\Spec \mc{O}_F$ (it admits a closed immersion to the moduli stack of signature $(n - 1, 1)$ Hermitian abelian schemes whose quasi-polarizations are polarizations of degree $|\Delta|^n$).
If $(A, \iota, \lambda)$ is a signature $(n - r, r)$ Hermitian abelian scheme over an algebraically closed field of characteristic $0$, then requiring the preceding polarization condition $\circ$ is equivalent to requiring that the Hermitian $\mc{O}_F \otimes_{\Z} \Z_p$-lattice $T_p(A)$ is self-dual (for the trace pairing) for all primes $p$. 

We now restrict to signature $(n - 1, 1)$ and $n \geq 2$. Let $\ms{M}(n - 1, 1)^{\circ}$ be the flat part of $\ms{M}(n - 1, 1)^{\mrm{Kot}, \circ}$, i.e. the scheme-theoretic image of the generic fiber. Equivalently, this is the largest closed substack which is flat over $\Spec \mc{O}_F$. 

This $\ms{M}(n - 1, 1)^{\circ}$ turns out to be smooth over $\Spec \mc{O}_F$ (also separated and finite type). This particular smoothness phenomenon, even over ramified primes, is typically called ``exotic smoothness'' in the literature. This moduli stack has been studied by Rapoport--Smithling--Zhang,\footnote{Strictly speaking, Rapoport--Smithling--Zhang normalize their polarization differently (i.e. our $\lambda$ is their $|\Delta|^{-1} \lambda$). Their convention is more common elsewhere in the literature, and is of course equivalent to our formulation. We prefer our normalization, which seems more natural for our main results on the comparison with Eisenstein series Fourier coefficients. A related remark is \cite[{Footnote 9}]{LL22II}.\label{footnote:pi_inverse_modular}} who gave a moduli description of $\ms{M}(n - 1, 1)^{\circ}$ as follows: given any scheme $S$ over $\Spec \mc{O}_F$, the groupoid $\ms{M}(n - 1, 1)^{\mrm{RSZ}}(S) \subseteq \ms{M}(n - 1, 1)^{\mrm{Kot}, \circ}(S)$ is the full subcategory consisting of tuples $(A, \iota, \lambda)$ such that the action $\iota \colon \mc{O}_F \ra \End(A)$ satisfies:
    \begin{enumerate}[(1)]
        \item (Pappas wedge condition) For all $a \in \mc{O}_F$, the action of $\iota(a)$ on $\Lie A$ satisfies
            \begin{equation*}
            \bigwedge^{2} (\iota(a) - a) = 0 \quad \text{ and } \quad \bigwedge^{n}(\iota(a) - a^{\s}) = 0.
            \end{equation*}
        \item (PRRSZ spin condition) For every geometric point $\overline{s}$ of $S$, the action of $(\iota(a) - a)$ on $\Lie A_{\overline{s}}$ is nonzero for some $a \in \mc{O}_F$.
    \end{enumerate}
The signature condition implies that the equation involving $\bigwedge^{n}$ in the wedge condition is automatic, and that the wedge condition is empty if $n = 2$.
The wedge and spin conditions are automatic (given the signature condition) over $\Spec \mc{O}_F[1/\Delta]$, i.e. $\ms{M}(n-1,1)^{\mrm{Kot}, \circ}[1/\Delta] = \ms{M}(n-1,1)^{\mrm{RSZ}}[1/\Delta]$. For closedness of the spin condition, we refer to the closedness assertion in \cite[{Theorem 5.4}]{RSZ21}. The acronym PRRSZ stands for Pappas, Rapoport, Richarz, Smithling, and Zhang (for the work of Rapoport--Smithling--Zhang and the earlier works of Pappas \cite{Pappas00} and Richarz \cite[{Proposition 4.16}]{Arzdorf09}).

\begin{example}\label{example:integral_models:Serre_tensor_global}
Suppose $\mc{O}_F^{\times} = \{ \pm 1 \}$ (i.e. further exclude $F = \Q[\sqrt{-3}]$). If $\ms{M}_{\text{ell}} \ra \Spec \mc{O}_F$ denotes the moduli stack of elliptic curves base-changed to $\Spec \mc{O}_F$, the \emph{Serre tensor construction} $E \mapsto E \otimes_{\Z} \mc{O}_F$ defines an open and closed immersion $i_{\mrm{Serre}} \colon \ms{M}_{\text{ell}} \ra \ms{M}(1,1)^{\circ}$ \crefext{IV:equation:arithmetic_Siegel-Weil:Serre_tensor}. If we replace $\mc{O}_F$ by (representatives of) fractional ideal classes for $\mc{O}_F$, we obtain an isomorphism $\coprod_{\mrm{Cl}(\mc{O}_F)} \ms{M}_{\text{ell}} \ra \ms{M}(1,1)^{\circ}$, where $\mrm{Cl}(\mc{O}_F)$ is the class group. This is \cite[{Proposition 14.4}]{KR14}. The local analogue (e.g. Lemma \ref{lemma:moduli_pDiv:Serre_tensor}) will play an important role in this work.
In \crefext{IV:ssec:arithmetic_Siegel-Weil:Serre_tensor}, we revisit this description of $\ms{M}(1,1)^{\circ}$ to restate our main theorem in the simplest case.
\end{example}

We now set
    \begin{equation}
    \mc{M} \coloneqq \ms{M}_0 \times_{\Spec \mc{O}_F} \ms{M}(n - 1, 1)^{\circ}.   
    \end{equation}
This moduli stack $\mc{M}$ is (essentially) a special case of the Rapoport--Smithling--Zhang Shimura variety integral models from \cite[{\S 6}]{RSZ21}, in an everywhere smooth situation. This stack $\mc{M}$ (and its special cycles) are the global geometric objects of main interest in our work.
    
            \subsection{Global special cycles}
            \label{ssec:part_I:arith_intersections:special_cycles}
                The stack $\mc{M}$ carries a family of special cycles as defined by Kudla and Rapoport \cite[{Definition 2.8}]{KR14} (there in a principally polarized situation).

    \begin{definition}[Kudla--Rapoport special cycles]\label{definition:global_special_cycles}
    Given an integer $m \geq 0$, let $T \in \mrm{Herm}_m(\Q)$ be a $m \times m$ Hermitian matrix (with coefficients in $F$). The \emph{Kudla--Rapoport (KR) special cycle} $\mc{Z}(T)$ is the stack in groupoids over $\Spec \mc{O}_F$ defined as follows: for schemes $S$ over $\Spec \mc{O}_F$, we take $\mc{Z}(T)(S)$ to be the groupoid
        \begin{equation}
        \mc{Z}(T)(S) \coloneqq 
        \left \{ (A_0, \iota_0, \lambda_0, A, \iota, \lambda, \underline{x}) : \begin{array}{c} (A_0, \iota_0, \lambda_0, A, \iota, \lambda) \in \mc{M}(S) \\ \underline{x} = [x_1,\ldots,x_m] \in \Hom_{\mc{O}_F}(A_0,A)^m \\ (\underline{x}, \underline{x}) = T \end{array} \right \}
        \end{equation}
    where $(\underline{x}, \underline{x})$ is the matrix with $i,j$-th entry given by $x_i^{\dagger} x_j \in \End^0_{\mc{O}_F}(A_0)$ (a quasi-endomorphism), with $\dagger$ denoting the Rosati involution. We sometimes refer to elements $x \in \Hom_{\mc{O}_F}(A_0, A)$ as \emph{special homomorphisms}.
    \end{definition}

    \begin{example}
    Suppose $2 \nmid \Delta$, and consider $L$ which is self-dual of signature $(1,1)$. Let $j \in \Z_{>0}$ be any positive integer. If $\mc{O}_F^{\times} = \{ \pm 1 \}$, consider the inclusion
        \begin{equation}
        \ms{M}_0 \times_{\Spec \mc{O}_F} \ms{M}_{\text{ell}} \xra{1 \times i_{\mrm{Serre}}} \ms{M}_0 \times_{\Spec \mc{O}_F} \ms{M}(1, 1)^{\circ} = \mc{M}
        \end{equation}
    with $i_{\mrm{Serre}}$ as in Example \ref{example:integral_models:Serre_tensor_global}. Then $\mc{Z}(j) \ra \mc{M}$ pulls back to the $j$-th Hecke correspondence over the left-hand side, parameterizing triples $(E_0, E, w)$ where $E_0$ and $E$ are elliptic curves, $E_0$ has $\mc{O}_F$ action of signature $(1,0)$, and $w \colon E \ra E_0$ is an isogeny of degree $j$.
    This is \cite[{Proposition 14.5}]{KR14}. 
    We revisit this example in \crefext{IV:ssec:arithmetic_Siegel-Weil:Serre_tensor}, where we restate our main (global) theorem in the simplest case via this description.
    \end{example}

    We state a few properties of the special cycles $\mc{Z}(T)$ which we need (proofs are available in our companion paper \cite{corank1_ASW_III.pdf}, there in greater generality). The forgetful map $\mc{Z}(T) \ra \mc{M}$ is representable by schemes, finite, and unramified \crefext{III:lemma:quasi-compactness_lemma}, which also implies that $\mc{Z}(T)$ is a Deligne--Mumford stack. Moreover, the generic fiber $\mc{Z}(T)_F \ra \Spec F$ is smooth of relative dimension $n - 1 - \rank(T)$ if nonempty \crefext{III:lemma:special_cycles_generically_smooth}. If $p$ is any prime, then the restriction $\mc{Z}(T)[1/p] \ra \mc{M}[1/p]$ becomes a disjoint union of closed immersions after a finite \'etale cover\footnote{We use the shorthand $\mc{M}[1/p]$ for $\mc{M} \times_{\Spec \mc{O}_F} \Spec \mc{O}_F[1/p]$, etc..} of $\mc{M}[1/p]$ \crefext{III:lemma:special_cycles_level_structure}, e.g. after adding enough level structure at $p$. Note that $\mc{Z}(T) \ra \mc{M}$ may not literally be a cycle (e.g. may not be a closed immersion), but we still call $\mc{Z}(T)$ a special cycle, as is standard abuse of terminology.

In the rest of Section \ref{ssec:part_I:arith_intersections:special_cycles}, we also collect some needed miscellaneous facts which follow more readily from the definition of $\mc{Z}(T)$. If $T_i \in \mrm{Herm}_{m_i}(\Q)$ for $i = 1, \ldots, j$ with $m \coloneqq m_1 + \cdots + m_j$ and all $m_i > 0$, then there is an identification
    \begin{equation}\label{equation:special_cycles_diagonal_entry_intersection}
    \mc{Z}(T_1) \times_{\mc{M}} \cdots \times_{\mc{M}} \mc{Z}(T_j) \cong \coprod_{\substack{T \in \mrm{Herm}_m(\Q) \\ \text{satisfying \eqref{equation:special_cycles_diagonal_entry_intersection_2}}}} \mc{Z}(T)
    \end{equation}
where the disjoint union runs over $T$ of the form
    \begin{equation}\label{equation:special_cycles_diagonal_entry_intersection_2}
    T = \begin{pmatrix} T_1 & & * \\ & \ddots & \\ * & & T_j \end{pmatrix}
    \end{equation}
(i.e. whose block diagonal entries are given by the $T_i$).

We write $\mc{Z}(T)_{\ms{H}} \subseteq \mc{Z}(T)$ for the largest closed substack flat over $\Spec \mc{O}_F$, and call $\mc{Z}(T)_{\ms{H}}$ a \emph{horizontal special cycle} or the \emph{flat part}. 
There is a decomposition of $\mc{Z}(T)$ as a scheme-theoretic union of closed substacks\footnote{By the \emph{scheme-theoretic union} of finitely many closed substacks $\mc{Z}_i$ of a Deligne--Mumford stack $\mc{Z}$, we mean the closed substack whose ideal sheaf is given by intersecting the ideal sheaves of $\mc{Z}_i$ on the small \'etale site of $\mc{Z}$.}
    \begin{equation}\label{equation:setup:ab_var:horizontal_vertical_decomp}
    \mc{Z}(T) = \mc{Z}(T)_{\ms{H}} \cup \bigcup_{p} \mc{Z}(T)_{\ms{V}, p}
    \end{equation}
where $\mc{Z}(T)_{\ms{V}, p} \coloneqq \mc{Z}(T) \times_{\Spec \Z} \Spec \Z/p^{e_p}\Z$ for a choice of $e_p \gg 0$ (notation $e_p$ is temporary). This follows from quasi-compactness of $\mc{Z}(T)$ (which also ensures that we may take $e_p = 0$ for all but finitely many $p$). We think of \eqref{equation:setup:ab_var:horizontal_vertical_decomp} as decomposition into a ``horizontal part'' and ``vertical parts''. A similar horizontal/vertical decomposition for local special cycles on Rapoport--Zink spaces is \cite[{\S 2.9}]{LZ22unitary} (inert case).

While the horizontal part $\mc{Z}(T)_{\ms{H}}$ is defined canonically, the vertical parts $\mc{Z}(T)_{\ms{V}, p}$ depend on $e_p$ as above. We are mostly interested in intersection numbers defined via ``derived vertical special cycle classes'' from Section \ref{ssec:part_I:arith_intersections:vertical_classes}, which do not depend on such a choice of $e_p$. 

Given any $T \in \mrm{Herm}_m(\Q)$ and $\gamma \in M_{n,n}(\mc{O}_F)$, there is a commutative diagram
    \begin{equation}\label{equation:ab_var:special_cycles:gamma_action}
    \begin{tikzcd}[row sep = tiny]
    (A_0, \iota_0, \lambda_0, A, \iota, \lambda, \underline{x}) \arrow[mapsto]{rr} & & (A_0, \iota_0, \lambda_0, A, \iota, \lambda, \underline{x} \cdot \gamma)
    \\
    \mc{Z}(T) \arrow{ddr} \arrow{rr} & & \mc{Z}({}^t \overline{\gamma} T \gamma) \arrow{ddl} 
    \\
    & &
    \\
    & \mc{M} &
    \end{tikzcd}
    \end{equation}
induced by $\gamma$. Below, we set $\mc{O}_{F,(p)} \coloneqq \mc{O}_F \otimes_{\Z} \Z_{(p)}$.

\begin{lemma}\label{lemma:ab_var:special_cycles:gamma_action}
Fix any $T \in \mrm{Herm}_m(\Q)$ and $\gamma \in M_{n,n}(\mc{O}_F)$. Consider the induced map
    \begin{equation}
    \mc{Z}(T) \ra \mc{Z}({}^t \overline{\gamma} T \gamma).
    \end{equation}
    \begin{enumerate}[(1)]
        \item This map is a finite morphism of algebraic stacks. If moreover $\gamma \in \GL_m(F)$ (resp. $\gamma \in \GL_m(\mc{O}_F)$) then the map is a closed immersion (resp. isomorphism).
        \item Assume $\gamma \in \GL_m(F)$, and let $N \in \Z$ be the product of primes $p$ such that $\gamma \not \in \GL_m(\mc{O}_{F,(p)})$. Then the restriction $\mc{Z}(T)[1/N] \ra \mc{Z}({}^t \overline{\gamma} T \gamma)[1/N]$ is an open immersion.
    \end{enumerate}
\end{lemma}
\begin{proof}
(1) The map $\mc{Z}(T) \ra \mc{Z}(T)({}^t \overline{\gamma} T \gamma)$ is finite because the projections to $\mc{M}$ are finite (finiteness for morphisms of algebraic stacks may be checked fppf locally on the target, so we reduce to the case of schemes). If $\gamma \in \GL_m(F)$, then $\mc{Z}(T) \ra \mc{Z}(T)({}^t \overline{\gamma} T \gamma)$ is a monomorphism of algebraic stacks (check via the moduli description), and any proper monomorphism of algebraic stacks is a closed immersion. If $\gamma \in \GL_m(\mc{O}_F)$, there is an inverse $\mc{Z}({}^t \overline{\gamma} T \gamma) \ra \mc{Z}(T)$ sending $\underline{x} \mapsto \underline{x} \cdot \gamma^{-1}$.

(2) Consider the substack $\mc{Z} \subseteq \mc{Z}({}^t \overline{\gamma} T \gamma)$ consisting of tuples $(A_0, \iota_0, \lambda_0, A, \iota, \lambda, \underline{w})$ such that $\underline{w} \cdot \gamma^{-1} \in \Hom_{\mc{O}_F}(A_0, A)^m$ (i.e. that the tuple of quasi-homomorphisms $\underline{w} \cdot \gamma^{-1}$ is a tuple of homomorphisms). The closed immersion $\mc{Z}(T) \ra \mc{Z}({}^t \overline{\gamma} T \gamma)$ maps isomorphically onto $\mc{Z}$. 

A quasi-homomorphism $f \colon B \ra B'$ of abelian schemes (over any base scheme $S$) is a homomorphism if and only the induced quasi-homomorphisms of $p$-divisible groups $f[p^{\infty}] \colon B[p^{\infty}] \ra B'[p^{\infty}]$ are homomorphisms for all primes $p$. This is a closed condition on $S$ for each prime $p$ (e.g. the quasi-homomorphism version of \cite[{Proposition 2.9}]{RZ96}). This is also an open condition for any prime $p$ which is invertible on $S$ (by \'etaleness of the $p$-divisible groups). If $p$ is a prime such that $\gamma \in \GL_m(\mc{O}_{F,(p)})$, then the tuple $\underline{w} \cdot \gamma^{-1}$ induces a tuple of quasi-homomorphisms $A_0[p^{\infty}] \ra A[p^{\infty}]$ consisting of homomorphisms.
\end{proof}

            \subsection{Hermitian vector bundles}
            \label{ssec:part_I:arith_intersections:Hermitian_bundles}
                Given a smooth algebraic stack $\mc{X}$ over $\Spec \C$, a \emph{Hermitian vector bundle} $\widehat{\mc{E}}$ on $\mc{X}$ is the following functorial assignment: for every morphism $S \ra \mc{X}$ with $S$ a scheme smooth over $\Spec \C$, the assignment gives a vector bundle on $S$ with a smooth Hermitian metric on the analytification. We sometimes write $\widehat{\mc{E}} = (\mc{E}, \norm{-})$ where $\mc{E}$ is the underlying line bundle on $\mc{X}$ and $\norm{-}$ is the norm associated with the smooth Hermitian metric (on the analytification), defined functorially.

Let $(R,\Sigma,c_{\infty})$ be an arithmetic ring in the sense of Gillet--Soul\'e \cite[{\S 3.1}]{GS90}, i.e. $R$ is an excellent regular Noetherian integral domain (e.g. Dedekind domains with fraction field of characteristic $0$ or fields), $\Sigma$ is a finite nonempty set of injective homomorphisms $\t \colon R \ra \C$, and $c_{\infty} \colon \C^{\Sigma} \ra \C^{\Sigma}$ is a conjugate-linear involution of $\C \otimes_{\Z} R$-algebras. Write $K$ for the fraction field of $R$. 

Suppose $\mc{X}$ is an algebraic stack which is flat and finite type over $\Spec R$. Write $\mc{X}_{\t} \coloneqq \mc{X} \times_{\Spec K, \t} \Spec \C$. Assume that the generic fiber $\mc{X}_K$ is smooth over $\Spec K$. A \emph{Hermitian vector bundle} on $\mc{X}$ is a vector bundle $\mc{E}$ on $\mc{X}$ equipped with a smooth Hermitian metric on $\mc{E}|_{\mc{X}_{\t}}$ for each $\t \in \Sigma$, such that this collection of metrics is conjugation invariant (meaning $c_{\infty}$-invariant). We write $\widehat{\Pic}(\mc{X})$ for the group of (isomorphism classes of) Hermitian line bundles, with group structure given by the tensor product. We use the notation $\widehat{\mc{L}} = (\mc{L}, \norm{-})$ for a Hermitian line bundle with underlying line bundle $\mc{L}$ and $\norm{-}$ its norm (meaning a collection of norms $\{\norm{-}_{\t}\}_{\t \in \Sigma}$). We also write $\norm{-}_{\infty} = \prod_{\t \in \Sigma} \norm{-}_{\t}$.

Next, we discuss stacky degrees of Hermitian line bundles. We fix the arithmetic ring $(R, \Sigma, c_{\infty})$ associated with $R = \mc{O}_K[1/N]$ for a number field $K$ and an integer $N \in \Z$, i.e. $\Sigma$ is the set of all embeddings $\t \colon K \ra \C$, and $c_{\infty}$ is induced by complex conjugation. For the rest of Section \ref{ssec:part_I:arith_intersections:Hermitian_bundles}, we assume that 
    \begin{equation}\label{equation:stacky_degrees_stack_hypotheses}
    \parbox[t]{.8 \textwidth}{$\mc{X}$ is a reduced $1$-dimensional Deligne--Mumford stack which is proper and flat over $\Spec R$.}
    \end{equation}
Here, dimension is used in an absolute sense rather than a relative one (e.g. we could have $\mc{X} = \Spec R$).

Let $\widehat{\mc{L}} = (\mc{L}, \norm{-})$ be a Hermitian line bundle on $\mc{X}$. For each complex place $v$ of $\mc{O}_K$, we set $\norm{-}_v \coloneqq \norm{-}_{\t_1} \norm{-}_{\t_2}$ where $\t_1, \t_2 \colon K \ra \C$ are the two embeddings corresponding to $v$.
The Arakelov \emph{arithmetic degree} $\widehat{\deg}(\widehat{\mc{L}})$ of $\widehat{\mc{L}}$ on $\mc{X}$ is valued in $\R_N \coloneqq \R/(\sum_{p \mid N} \Q \cdot \log p)$, and may be described as follows. If $\mc{X} = \Spec \mc{O}_E[1/N]$ for a number field $E$, we have the standard definition
    \begin{equation}
    \widehat{\deg}(\widehat{\mc{L}}) \coloneqq - \sum_{v \nmid N} \log \norm{s}_v \quad \text{with} \quad \norm{s}_v \coloneqq q_v^{-\mrm{ord}_v(s)} \text{ if $v < \infty$}
    \end{equation}
where the sum runs over all places $v$ of $E$ (Archimedean included) with $v \nmid N$, the quantity $q_v$ is the cardinality of the residue field at $v$, and $s$ is any rational section of $\mc{L}$.

If $\mc{X}$ is integral (equivalently, reduced and irreducible by quasi-separatedness), select any number field $E$ with a finite surjection $\Spec \mc{O}_E[1/N] \ra \mc{X}$ and set
    \begin{equation} \label{equation:stacky_degree_Hermitian_bundle}
    \widehat{\deg}(\widehat{\mc{L}}) \coloneqq \frac{1}{\deg(E/\mc{X}_K)} \widehat{\deg}(\widehat{\mc{L}}|_{\Spec \mc{O}_{E}[1/N]}) 
    \end{equation}
where $\deg(E/\mc{X}_K)$ denotes the degree of the finite \'etale morphism $\Spec E \ra \mc{X}_K$. This definition does not depend on the choice of $E$ or the morphism $\Spec \mc{O}_E[1/N] \ra \mc{X}$, since any two such choices may be covered by finite surjections from a third such choice $\Spec \mc{O}_{E'}[1/N]$ (and these finite surjections can be made compatible with the maps to $\mc{X}$).

\begin{remark}\label{remark:existence_finite_covers_stack}
Existence of $E$ as in the preceding paragraph follows from some general theory. Indeed, a general fact about Noetherian Deligne--Mumford stacks with separated diagonal \cite[{Th\'eor\`eme 16.6}]{LMB00} implies that there exists a scheme $Z$ with a morphism $Z \ra \mc{X}$ which is finite, surjective, and generically \'etale (in the sense that $Z_{\mc{U}} \ra \mc{U}$ is \'etale for a dense open substack $\mc{U} \subseteq \mc{X}$). Selecting an irreducible component of $Z$ which surjects onto $\mc{X}$, we may assume that $Z$ is also integral. Thus $Z$ is a $1$-dimensional integral scheme which is proper and flat over $\Spec R$. Such a map $Z \ra \Spec R$ must be quasi-finite, hence also finite. If $E$ denotes the fraction field of $Z$, the normalization of $Z$ must be $\tilde{Z} = \Spec \mc{O}_{E}[1/N]$.
\end{remark}

One can check that the definition of stacky arithmetic degree in \eqref{equation:stacky_degree_Hermitian_bundle} recovers the definition of \cite[{(4.4)}]{KRY04} and \cite[{\S 2.1}]{KRY06} which counts geometric points weighted by orders of automorphism groups.\footnote{In loc. cit., the functorial definition of Hermitian line bundles should also include the additional assumption of complex conjugation invariance as above.}

In the general case when $\mc{X}$ is not necessarily irreducible, we take
    \begin{align*}
    \widehat{\deg}(\widehat{\mc{L}}) \coloneqq \sum_{\xi} \widehat{\deg}(\widehat{\mc{L}}|_{\mc{X}_{\xi}})
    \end{align*}
where the sum runs over generic points $\xi$ of irreducible components $\mc{X}_{\xi}$ of $\mc{X}$.

The preceding discussion showed that $\mc{X}$ admits a finite surjection from a scheme which is finite over $\Spec R$, hence $\mc{X} \ra \Spec R$ is also quasi-finite (in the sense of \cite[\href{https://stacks.math.columbia.edu/tag/0G2M}{Definition 0G2M}]{stacks-project}).

Suppose $\mc{X}'$ and $\mc{X}$ are Deligne--Mumford stacks which both satisfy the hypotheses from \crefext{equation:stacky_degrees_stack_hypotheses}, and consider a morphism $f \colon \mc{X}' \ra \mc{X}$ over $\Spec R$. Let $\widehat{\mc{L}}$ be a Hermitian line bundle on $\mc{X}$. First consider the case where $\mc{X}$ is irreducible. We say that the morphism $\mc{X}'_K \ra \mc{X}_K$ has \emph{degree} $\deg(\mc{X}'_K/\mc{X}_K) \coloneqq \deg(\mc{X}'_K/K) / \deg(\mc{X}_K/K)$ (with stacky degrees of $0$-cycles over fields as in \crefext{equation:zero_cycle_degree_over_field}). We have $\widehat{\deg}(f^* \widehat{\mc{L}}) = \deg(\mc{X}'_K / \mc{X}_K) \widehat{\deg}(\widehat{\mc{L}})$.
Next, consider the general case where $\mc{X}$ is not necessarily irreducible. We say that $\mc{X}'_K \ra \mc{X}_K$ has \emph{constant degree $d$} if $\mc{X}'_K \times_{\mc{X}_K} \mc{X}_{\xi,K} \ra \mc{X}_{\xi,K}$ has degree $d$ for every irreducible component $\mc{X}_{\xi}$ of $\mc{X}$. In this case, we have
    \begin{equation} \label{equation:stacky_degrees_pullback_equation}
    \widehat{\deg}(f^* \widehat{\mc{L}}) = d \cdot  \widehat{\deg}(\widehat{\mc{L}}).
    \end{equation}
    
            \subsection{Tautological bundle}
            \label{ssec:part_I:arith_intersections:tautological_bundle}
                In order to extract arithmetic intersection numbers (``heights'') from special cycles of positive expected dimension, we will need a certain metrized tautological bundle $\widehat{\mc{E}}^{\vee}$ on $\mc{M}$.

Given $(A, \iota, \lambda) \in \ms{M}(n - 1, 1)^{\circ}(S)$ over an $\mc{O}_F$-scheme $S$, we set
    \begin{equation}\label{equation:ab_var:integral_models:Lie_decomp_ram}
    (\Lie A)^+ \coloneqq \bigcap_{a \in \mc{O}_F} \ker (\iota(a) - a)|_{\Lie A}.
    \end{equation}

\begin{notation}
Given a commutative ring $R$ with an automorphism $\s \colon R \ra R$ (e.g. $R = \mc{O}_F$), given a presheaf of modules $\mc{F}$ on a scheme $S$ over $\Spec R$, and given an action $\iota \colon R \ra \End(\mc{F})$ (with $\mc{F}$ viewed as a presheaf of abelian groups), we say the $R$ action via $\iota$ is \emph{$R$-linear} (resp. \emph{$\s$-linear}) if $\iota(a) = a$ (resp. $\iota(a) = a^{\s}$) for all $a \in R$. Here we view $\mc{O}_S$ as an $R$-algebra via the structure morphism $S \ra \Spec R$.
\end{notation}
    
\begin{lemma}\label{lemma:unique_Kramer_hyperplane}
For objects $(A, \iota, \lambda) \in \ms{M}(n - 1, 1)^{\circ}(S)$, the subsheaf $(\Lie A)^+ \subseteq \Lie A$ is a local direct summand of rank $n - 1$ whose formation commutes with arbitrary base change. The $\mc{O}_F$ action via $\iota$ on $(\Lie A)^+$ (resp. the line bundle $(\Lie A)/(\Lie A)^+$) is $\mc{O}_F$-linear (resp. $\s$-linear).
\end{lemma}
\begin{proof}
This lemma (and its proof) is a global analogue of \cite[{Lemma 2.36}]{LL22II} (the latter is an analogous statement on a Rapoport--Zink space).

Fix $a \in \mc{O}_F$ such that $\{1, a\}$ forms a $\Z$-basis of $\mc{O}_F$.
We have exact sequences
    \begin{equation}\label{equation:lemma:unique_Kramer_hyperplane:1}
    \begin{tikzcd}[row sep = small]
    0 \arrow{r} & (\Lie A)^+ \arrow{r} & \Lie A \arrow{r}{\iota(a) - a} \arrow{r} &  \operatorname{im}(\iota(a) - a) \arrow{r} & 0
    \\
    0 \arrow{r} & \operatorname{im}(\iota(a) - a) \arrow{r} & \Lie A \arrow{r} & \coker(\iota(a) - a) \arrow{r} & 0
    \end{tikzcd}
    \end{equation}
of quasi-coherent sheaves on $S$. The wedge and spin conditions imply that $(\Lie A)^+$ has rank $n - 1$ if $S = \Spec k$ for a field $k$. If $S$ is an arbitrary reduced scheme, the rank constancy of $\coker(\iota(a) - a)$ on geometric points implies that $\coker(\iota(a) - a)$ is finite locally free of rank $n - 1$ (e.g. by \cite[\href{https://stacks.math.columbia.edu/tag/0FWG}{Lemma 0FWG}]{stacks-project}). Hence, when $S$ is reduced, every sheaf appearing in \eqref{equation:lemma:unique_Kramer_hyperplane:1} is finite locally free, with $(\Lie A)^+$, $\operatorname{im}(\iota(a) - a)$, and $\coker(\iota(a) - a)$ having ranks $n-1$, $1$, and $n-1$ respectively. Thus the exact sequences of \eqref{equation:lemma:unique_Kramer_hyperplane:1} remain exact after pullback along any morphism of schemes $S' \ra S$ (where $S$ is reduced but $S'$ is not necessarily reduced). For arbitrary $S$ (not necessarily reduced), the morphism
$S \ra \ms{M}(n-1, 1)^{\circ}$ corresponding to $(A, \iota, \lambda)$ factors through a regular (hence reduced) locally Noetherian scheme fppf locally on $S$ (since the moduli stack $\ms{M}(n-1,1)^{\circ}$ is smooth over $\Spec \mc{O}_F$). These considerations show that $(\Lie A)^+ \subseteq \Lie A$ is a local direct summand of rank $n - 1$ whose formation commutes with arbitrary base change.

It is clear that the $\mc{O}_F$ action via $\iota$ on $(\Lie A)^+$ is $\mc{O}_F$-linear. To show that the action on $(\Lie A)/(\Lie A)^+$ is $\s$-linear, it is enough to check the case where $S$ is an integral scheme (argue fppf locally as above). When $S$ is integral, the $\s$-linearity follows from the $(n-1,1)$ signature condition on $\Lie A$.
\end{proof} 

We write $\Omega_0^{\vee}$ for the line bundle on $\ms{M}_0$ given by the functorial association $(A_0, \iota_0, \lambda_0) \mapsto \Lie A_0$ for objects $(A_0, \iota_0, \lambda_0) \in \ms{M}_0(S)$.

\begin{definition}\label{definition:ab_var:exotic_smooth:tautological_bundle}
By the \emph{tautological bundle} $\ms{E}$ on $\ms{M}(n - 1, 1)^{\circ}$, we mean the rank one locally free sheaf $\ms{E}$ (for the fppf topology) whose dual is $\ms{E}^{\vee} \coloneqq (\Lie A)/(\Lie A)^+$ for $(A, \iota, \lambda) \in \ms{M}(n - 1, 1)^{\circ}(S)$ for $\mc{O}_F$-schemes $S$.

By the \emph{tautological bundle} on $\mc{M}$, we mean the rank one locally free sheaf $\mc{E}$ whose dual is
    \begin{equation}\label{equation:ab_var:exotic_smooth:tautological_bundle}
    \mc{E}^{\vee} \coloneqq \Omega_0^{\vee} \otimes \ms{E}^{\vee}.
    \end{equation}
\end{definition}

On the right-hand side of Equation \eqref{equation:ab_var:exotic_smooth:tautological_bundle}, the bundles $\Omega_0^{\vee}$ and $\ms{E}^{\vee}$ denote pulled-back bundles along the projections $\mc{M} \ra \ms{M}_0$ and $\mc{M} \ra \ms{M}(n - 1, 1)^{\circ}$, respectively (by abuse of notation).

We write $\Omega^{\vee}$ for the rank $n$ vector bundle on $\ms{M}(n - 1, 1)^{\circ}$ given by the functorial association $(A, \iota, \lambda) \mapsto \Lie(A)$, and also write $\Omega^{\vee}$ for its pullback along the forgetful map $\mc{M} \ra \ms{M}(n - 1, 1)^{\circ}$.

The natural surjection $\Omega^{\vee} \ra \ms{E}^{\vee}$ admits a unique $\mc{O}_F$-linear splitting $\ms{E}^{\vee} \ra \Omega^{\vee}$ over $\ms{M}(n - 1, 1)^{\circ}[1/\Delta]$: indeed, if $S$ is an $\mc{O}_F[1/\Delta]$-scheme and $(A, \iota, \lambda) \in \ms{M}(n - 1, 1)^{\circ}(S)$, then there is a unique direct sum decomposition
    \begin{equation}
    \Lie A = (\Lie A)^+ \oplus (\Lie A)^-
    \end{equation}
where the $\iota$ action on $(\Lie A)^+$ (resp. $(\Lie A)^-$) is $\mc{O}_F$-linear (resp. $\s$-linear).
    
            \subsection{Metrized tautological bundle}
            \label{ssec:part_I:arith_intersections:metrized_taut_bundle}
                Given an abelian scheme $\pi \colon A \ra S$ of relative dimension $n$, we consider the \emph{Hodge bundles} $\Omega_A \coloneqq \pi_* \Omega^1_{A/S}$ and $\omega_A \coloneqq \pi_* \bigwedge^n \Omega_A$. Here $\Omega_A$ and $\omega_A$ are locally free of ranks $n$ and $1$ respectively. If $e \colon S \ra A$ denotes the identity section, there are canonical isomorphisms $\Omega_A \cong e^* \Omega_{A/S} \cong (\Lie A)^{\vee}$. These Hodge bundles are (contravariantly) functorial in $A$, and their formation commutes with arbitrary base change \cite[{Proposition 2.5.2}]{BBM82}.

When $S$ is a scheme which is smooth over $\Spec \C$, the analytification of $\omega_A$ may be equipped with a Hermitian metric (\emph{Faltings} or \emph{Hodge} metric), normalized as follows. The fiber of $\omega_A$ over any $s \in S(\C)$ is canonically identified with $H^0(A_s, \omega_{A_s})$. Viewing $H^0(A_s, \omega_{A_s})$ as the $1$-dimensional $\C$-vector space of holomorphic $n$-forms on $A_s(\C)$, we take the metric on $\omega_A$ at $s$ to be
    \begin{equation}\label{equation:arith_cycle_classes:Hodge_bundles:Faltings_metric}
    (\a, \b) = \left ( \frac{i}{2 \pi} \right ) ^n \int_{A_s(\C)} \b \wedge \bar{\a}
    \end{equation}
for $\a, \b \in H^0(A_s, \omega_{A_s})$. We call the resulting Hermitian line bundle $\widehat{\omega}_A \coloneqq (\omega_A, \norm{-})$ a \emph{metrized Hodge bundle}. This metric on $\omega_A$ is functorial for isomorphisms between abelian schemes $A$ over $S$.

Next, suppose $S = \Spec \C$ and suppose $\lambda \colon A \ra A^{\vee}$ is a quasi-polarization. There is an associated Hermitian metric on $\Omega_A^{\vee} \cong \Lie(A)$ which we normalize as follows: 
if $\lambda$ is a polarization and $\lambda(a) = (t_a^* \mc{L}) \otimes \mc{L}^{-1}$ for an ample line bundle $\mc{L}$ on $A$ (where $t_a$ is translation by $a$), then the Chern class $c_1(\mc{L}) \in H^2(A,\Z)$ defines a $\Z$-valued alternating form $\psi$ on $H_1(A,\Z)$ (following the usual conventions, \cite[{\S I}]{Mumford85}, except our Hermitian pairings are conjugate linear in the first variable) and a positive definite Hermitian pairing
    \begin{equation}\label{equation:arith_cycle_classes:Hodge_bundles:Hermitian_normalization}
    (x,y) \coloneqq \pi \sqrt{|\Delta|} \big ( \psi(ix, y) - i \psi(x,y) \big )
    \end{equation}
on $\Lie A$. 
If instead $\lambda$ is a quasi-polarization, the associated Hermitian pairing is $m^{-1}$ times the Hermitian pairing of $m \lambda$ for any $m \in \Z_{>0}$ such that $m \lambda$ is a polarization. If $\dim A = 1$ and if $\lambda$ is the unique principal polarization of $A$, the induced dual metric on $\Omega_A \cong (\Lie A)^{\vee}$ is $\sqrt{|\Delta|}^{-1}$ times the Faltings metric (cf. the proof of \cite[{Lemma 5.1.4}]{BHKRY20II}).

Over an arbitrary smooth scheme $S \ra \Spec \C$, any quasi-polarization $\lambda \colon A \ra A^{\vee}$ defines a smooth Hermitian metric on $\Lie A$, given fiberwise by the construction above. We also call the resulting Hermitian vector bundle $\widehat{\Omega}_A^{\vee}$ a \emph{metrized Hodge bundle}.

Next, let $(R, \Sigma, c_{\infty})$ be an arithmetic ring, and suppose $\mc{X}$ is an algebraic stack which is flat and finite type over $\Spec R$, and whose generic fiber is smooth. Suppose we are given a relative abelian scheme $\mc{A}$ over $\mc{X}$ (equivalently, a functorial assignment of abelian schemes $A \ra S$ to objects $x \in \mc{X}(S)$, for $R$-schemes $S$). Formation of the metrized Hodge bundle $\widehat{\omega}_{A}$ is functorial, hence defines a \emph{metrized Hodge bundle} $\widehat{\omega}$ on $\mc{X}$. If $\mc{A}$ is equipped with a quasi-polarization, then there is similarly a \emph{metrized Hodge bundle} $\widehat{\Omega}^{\vee}$ on $\mc{X}$.

With $\Omega^{\vee}$ the vector bundle on $\ms{M}(n - 1, 1)^{\circ}$ (described at the end of Section \ref{ssec:part_I:arith_intersections:tautological_bundle}), we equip $\Omega^{\vee}$ with the Hermitian metric which is 
    \begin{equation}\label{equation:arith_cycle_classes:Hodge_bundles:4_pi_gamma_normalization}
    (x,y) \coloneqq 4 \pi^2 e^{\gamma} \sqrt{|\Delta|} \big ( \psi(ix, y) - i \psi(x,y) \big )
    \end{equation}
in complex fibers (i.e. multiply the Hermitian metric from \eqref{equation:arith_cycle_classes:Hodge_bundles:Hermitian_normalization} by $4 \pi e^{\gamma}$) where $\gamma$ is the Euler--Mascheroni constant. We remark that the normalization constant $4 \pi e^{\gamma}$ has appeared previously in similar contexts, e.g. \cite[{(0.4)}]{KRY04} \cite[{\S 7}]{KRY06} \cite[{\S 7.2}]{BHKRY20}. We refer to loc. cit. for possible conceptual explanations of this constant.

The $\mc{O}_F$-linear splitting $\ms{E}^{\vee} \ra \Omega^{\vee}$ over $\ms{M}(n - 1, 1)^{\circ}$ then endows $\ms{E}^{\vee}$ with a Hermitian metric as well, given by restriction from $\Omega^{\vee}$. We write $\widehat{\ms{E}}^{\vee}$ for the resulting Hermitian line bundle. With $\Omega_0^{\vee}$ denoting the dual line bundle to $\Omega_0$ on $\ms{M}_0$, we equip $\Omega_0^{\vee}$ with the metric given fiberwise by \eqref{equation:arith_cycle_classes:Hodge_bundles:Hermitian_normalization}. We write $\widehat{\Omega}_0^{\vee}$ for the resulting metrized line bundle.

By the \emph{metrized dual tautological bundle} on $\mc{M}$, we mean
    \begin{equation}
    \widehat{\mc{E}}^{\vee} \coloneqq \widehat{\Omega}_0^{\vee} \otimes \widehat{\ms{E}}^{\vee}
    \end{equation}
where we have suppressed pullbacks from notation. This definition of $\widehat{\mc{E}}^{\vee}$ is similar to \cite[{\S 2.4, \S 7.2}]{BHKRY20}, though in a different setup (we are considering not-necessarily principal polarizations). Taking a dual gives the \emph{metrized tautological bundle} $\widehat{\mc{E}}$.

Suppose $A_0 \ra \Spec \mc{O}_E$ is any (relative) elliptic curve with $\mc{O}_F$-action, where $E$ is a number field. If $\widehat{\omega}_{A_0}$ denotes the associated metrized Hodge bundle (normalized as in \eqref{equation:arith_cycle_classes:Hodge_bundles:Faltings_metric}), we recall that the \emph{Faltings height} of $A_0$ is 
    \begin{equation}\label{equation:part_I:arith_cycle_classes:Hodge_bundles:CM_Faltings_height}
    h_{\mrm{Fal}}^{\mrm{CM}} \coloneqq \frac{1}{[E : \Q]} \widehat{\deg}(\widehat{\omega}_{A_0}) = \frac{1}{2} \frac{L'(1,\eta)}{L(1, \eta)} + \frac{1}{2} \frac{\Gamma'(1)}{\Gamma(1)} + \frac{1}{4} \log |\Delta| - \frac{1}{2} \log(2 \pi)
    \end{equation}
where $\eta$ is the quadratic character associated to $F / \Q$, and $\Gamma$ is the usual gamma function. This comes from the classical Chowla--Selberg formula (the statement above is as in \cite[{Proposition 10.10}]{KRY04}). It will be convenient to define the height constant
    \begin{equation}\label{equation:part_I:arith_cycle_classes:Hodge_bundles:taut_height_constants}
    h_{\mrm{tau}}^{\mrm{CM}} \coloneqq - h_{\mrm{Fal}}^{\mrm{CM}} + \frac{1}{4} \log|\Delta| - \frac{1}{2} \log(4 \pi e^{\gamma})
    \end{equation}
which will reappear in \cref{ssec:qcan_heights:minimal_isogenies}.

            \subsection{Derived vertical cycle classes}
            \label{ssec:part_I:arith_intersections:vertical_classes}
                We define vertical special cycle classes via $K_0$ groups. 

For our notation and definitions regarding $K_0$-groups for Deligne--Mumford stacks, we refer to \crefext{appendix:K0}. Note that the stacky $K_0$ groups we use are different from those used in \cite{HM22}. 

Fix any prime $p$ and set
    \begin{equation}
    \mc{M}_{(p)} \coloneqq \mc{M} \times_{\Spec \Z} \Spec \Z_{(p)} \quad \quad \mc{Z}(T)_{(p)} \coloneqq \mc{Z}(T) \times_{\Spec \Z} \Spec \Z_{(p)}
    \end{equation}
for any $T \in \mrm{Herm}_m(\Q)$. Since $\mc{M}_{(p)}$ admits a finite \'etale cover by a scheme (add away-from-$p$ level structure), we may consider filtrations for $K_0'$ groups as in \crefext{definition:stack_dimension_filtration}.

Let $T \in \mrm{Herm}_m(\Q)$ be any $m \times m$ Hermitian matrix (with entries in $F$). We first describe a ``$p$-local'' derived special cycle class ${}^{\mbb{L}} \mc{Z}(T)_{(p)} \in F^m_{\mc{M}_{(p)}} K'_0(\mc{Z}(T)_{(p)})_{\Q}$ (with $F^m_{\mc{M}_{(p)}} \coloneqq F_{n - m}$ denoting the $m$-th step of the codimension filtration) before extracting a ``vertical'' piece.

For any $t \in \Q$, we define ${}^{\mbb{L}} \mc{Z}(t)_{(p)} \in F^1_{\mc{M}_{(p)}} K'_0(\mc{Z}(t)_{(p)})_{\Q}$ to be the element
    \begin{align}
    {}^{\mbb{L}} \mc{Z}(t)_{(p)} \coloneqq
    \begin{cases}
    [\mc{O}_{\mc{Z}(t)_{(p)}}] & \text{if $t \neq 0$} \\
    [\mc{O}_{\mc{M}_{(p)}}] - [\mc{E}] & \text{if $t = 0$}. \notag
    \end{cases}
    \end{align}
Write $t_1,\ldots,t_m$ for the diagonal entries of $T$. Using the intersection pairing from \crefext{lemma:K0_relative}, we form the intersection ${}^{\mbb{L}} \mc{Z}(t_1)_{(p)} \cdots {}^{\mbb{L}} \mc{Z}(t_m)_{(p)}$ and define ${}^{\mbb{L}} \mc{Z}(T)_{(p)}$ by the restriction
    \begin{equation}\label{equation:setup:arith_cycle_classes:vertical:derived_definition} 
    \begin{tikzcd}[row sep = tiny]    
    F^m_{\mc{M}_{(p)}} K_0'(\mc{Z}(t_1)_{(p)} \times_{\mc{M}_{(p)}} \cdots \times_{\mc{M}_{(p)}} \mc{Z}(t_m)_{(p)})_{\Q} \arrow{r} & F^m_{\mc{M}_{(p)}} K'_0(\mc{Z}(T)_{(p)})_{\Q} \\
    {}^{\mbb{L}} \mc{Z}(t_1)_{(p)} \cdots {}^{\mbb{L}} \mc{Z}(t_m)_{(p)} \arrow[mapsto]{r} & {}^{\mbb{L}} \mc{Z}(T)_{(p)}.
    \end{tikzcd}
    \end{equation}
This displayed restriction map comes from the disjoint union decomposition in \eqref{equation:special_cycles_diagonal_entry_intersection}. We are using \cref{lemma:K0_relative_with_filtrations} for multiplicativity of the codimension filtration, since $\mc{Z}(T)_{(p)} \ra \mc{M}_{(p)}$ becomes a disjoint union of closed immersions after adding sufficient away-from-$p$ level structure \crefext{III:lemma:special_cycles_level_structure}.
We call ${}^{\mbb{L}} \mc{Z}(T)_{(p)}$ the \emph{$p$-local derived special cycle class}\footnote{This is the construction of \cite[{Definition 5.1.3}]{HM22} (there for orthogonal Shimura varieties). This construction also underlies the intersection numbers considered in \cite{KR14} for non-degenerate $T$. We differ slightly from those references by localizing at $p$, since we will only be interested in the ``vertical'' part of ${}^{\mbb{L}} \mc{Z}(T)_{(p)}$. The ``horizontal part'' will be accounted for by the class $[\widehat{\mc{Z}}(T)_{\ms{H}}]$ mentioned in \cref{ssec:intro:results} (described in detail in \crefext{III:ssec:arith_cycle_classes:horizontal}).
} associated with $T$. 

The following lemma is a ``$p$-local'' version of linear invariance, and is proved using a variant on ideas from \cite{Howard19,HM22}. The map in \eqref{equation:arith_cycle_classes:vertical:p-local_linear_invariance} was defined in \eqref{equation:ab_var:special_cycles:gamma_action}.

\begin{lemma}\label{lemma:arith_cycle_classes:vertical:p-local_linear_invariance}
Given any $T \in \mrm{Herm}_m(\Q)$ and any $\gamma \in \GL_m(\mc{O}_{F,(p)}) \cap M_{m,m}(\mc{O}_F)$, the pullback along
    \begin{equation}\label{equation:arith_cycle_classes:vertical:p-local_linear_invariance}
    \mc{Z}(T)_{(p)} \ra \mc{Z}({}^t \overline{\gamma} T \gamma)_{(p)}
    \end{equation}
sends ${}^{\mbb{L}} \mc{Z}({}^t \overline{\gamma} T \gamma)_{(p)}$ to ${}^{\mbb{L}} \mc{Z}(T)_{(p)}$.
\end{lemma}
\begin{proof}
By Lemma \ref{lemma:ab_var:special_cycles:gamma_action}, we know that \eqref{equation:arith_cycle_classes:vertical:p-local_linear_invariance} is an open and closed immersion.

The ring $\mc{O}_{F,(p)}$ is a Euclidean domain, with Euclidean function $\phi(a) \coloneqq \sum_{\mf{p}_i} v_{\mf{p}_i}(a) \cdot f_i$ for nonzero $a \in \mc{O}_{F,(p)}$ (summing over primes $\mf{p}_i$ in $\mc{O}_F$ lying over $p$, with residue cardinality $p^{f_i}$). Row reducing via the Euclidean algorithm shows that $\GL_m(\mc{O}_{F,(p)})$ is generated by elementary matrices.

Any $\gamma \in \GL_m(\mc{O}_{F,(p)})$ may thus be expressed as $\gamma = \gamma_1 \gamma_2^{-1}$ where each $\gamma_1$ and $\gamma_2$ are products of elementary matrices lying in $\GL_m(\mc{O}_{F,(p)}) \cap M_{m,m}(\mc{O}_F)$. If moreover $\gamma \in \GL_m(\mc{O}_{F,(p)}) \cap M_{m,m}(\mc{O}_F)$, the commutative diagram
    \begin{equation}
    \begin{tikzcd}
    & \mc{Z}({}^t \overline{\gamma}_1 T \gamma_1)_{(p)} & \\
    \mc{Z}(T)_{(p)} \arrow{rr} \arrow{ur} & & \arrow{ul} \mc{Z}({}^t \overline{\gamma} T \gamma)_{(p)}
    \end{tikzcd}
    \end{equation}
shows that it is enough to prove the lemma when $\gamma \in \GL_m(\mc{O}_{F,(p)}) \cap M_{m,m}(\mc{O}_F)$ is an elementary matrix.

If $\gamma$ is a permutation matrix, the lemma is clear. Next, consider $a \in \mc{O}_{F,(p)}^{\times} \cap \mc{O}_F$. For any $t \in \Q$, note that $\mc{Z}(t)_{(p)} \ra \mc{Z}(\overline{a} t a)_{(p)}$ is an open and closed immersion (by Lemma \ref{lemma:ab_var:special_cycles:gamma_action} again). This fact implies the present lemma for the case where $\gamma = \mrm{diag}(a, 1, \ldots, 1)$.

It remains to check the case where $\gamma$ is an elementary unipotent matrix. This case follows as in the analogous result \cite[{Proposition 5.4.1}]{HM22} (there for $\mrm{GSpin}$).\footnote{Strictly speaking, our setup for stacky $K'_0$ groups is slightly different from that of \cite{HM22}, see \crefext{appendix:K0}. This makes no difference in the proof of the cited result. Alternatively, one can replace $\mc{M}_{(p)}$ by a finite \'etale cover by a scheme to reduce to the case of schemes, where our setup agrees with \cite[{\S A.2}]{HM22}.} 
The latter is proved using methods from \cite{Howard19} (the analogous local linear invariance result on Rapoport--Zink spaces). We are also using global analogues of \cite[{Lemma 2.36, Lemma 2.37, Lemma 2.41}]{LL22II} (there about a tautological bundle on an exotic smooth Rapoport--Zink space) which may be proved similarly, e.g. our Lemma \ref{lemma:unique_Kramer_hyperplane} replaces \cite[{Lemma 2.36}]{LL22II} in the global setup.
Alternatively, linear invariance for $\gamma \in \GL_m(\mc{O}_F)$ should also follow from the derived algebro-geometric methods in \cite{Madapusi22}.
\end{proof}

Next, we define a \emph{derived vertical special cycle class}
    \begin{equation}
    {}^{\mbb{L}} \mc{Z}(T)_{\ms{V},p} \in \mrm{gr}^m_{\mc{M}} K'_0(\mc{Z}(T)_{\F_p})_{\Q}
    \end{equation}
at $p$, where $\mc{Z}(T)_{\F_p} \coloneqq \mc{Z}(T) \times_{\Spec \Z} \Spec \F_p$.

First consider the case where $\det T \neq 0$. Using \crefext{III:lemma:special_cycles_generically_smooth} and \cref{lemma:K'0_component_decomp} as well as \eqref{equation:setup:ab_var:horizontal_vertical_decomp}, we decompose
    \begin{equation}\label{equation:nonsingular_horizontal_vertical_decomp}
    \mrm{gr}^{m}_{\mc{M}_{(p)}} K'_0(\mc{Z}(T)_{(p)})_{\Q} = \mrm{gr}^{m}_{\mc{M}_{(p)}} K'_0(\mc{Z}(T)_{{(p)},\ms{H}})_{\Q} \oplus \mrm{gr}^{m}_{\mc{M}} K'_0(\mc{Z}(T)_{\F_p})_{\Q}
    \end{equation}
into a ``horizontal'' part and a ``vertical'' part. This uses nonsingularity of $T$ (via \crefext{III:lemma:special_cycles_generically_smooth}), so that $\mc{Z}(T)_{{(p)},\ms{H}} \cap \mc{Z}(T)_{\F_p}$ is of dimension $< n - m$.
We are also using the d\'evissage pushforward identification $K'_0(\mc{Z}(T)_{\F_{p}}) \xra{\sim} K'_0(\mc{Z}(T)_{\ms{V},p})$, with $\mc{Z}(T)_{\ms{V},p}$ as in \eqref{equation:setup:ab_var:horizontal_vertical_decomp}. The above decomposition of $\mrm{gr}^{m}_{\mc{M}_{(p)}} K'_0(\mc{Z}(T)_{(p)})_{\Q}$ is independent of the choice of $e_{p}$ in \eqref{equation:setup:ab_var:horizontal_vertical_decomp}.
We define ${}^{\mbb{L}} \mc{Z}(T)_{\ms{V},p}$ to be given by the projection
    \begin{equation}
    \begin{tikzcd}[row sep = tiny]
    \mrm{gr}^{m}_{\mc{M}_{(p)}} K'_0(\mc{Z}(T)_{(p)})_{\Q} \arrow{r} & \mrm{gr}^{m}_{\mc{M}} K'_0(\mc{Z}(T)_{\F_p})_{\Q} \\
    {}^{\mbb{L}} \mc{Z}(T)_{(p)} \arrow[mapsto]{r} & {}^{\mbb{L}} \mc{Z}(T)_{\ms{V},p}.
    \end{tikzcd}
    \end{equation}

If $T = \mrm{diag}(0,T^{\flat})$ (with $\det T^{\flat} \neq 0$), we set
    \begin{equation}
    {}^{\mbb{L}} \mc{Z}(T)_{\ms{V},p} \coloneqq (\mc{E}^{\vee})^{m - \rank(T)} \cdot {}^{\mbb{L}} \mc{Z}(T^{\flat})_{\ms{V},p} \in \mrm{gr}^m_{\mc{M}} K'_0(\mc{Z}(T)_{\F_p})_{\Q}
    \end{equation}
where $\mc{E}^{\vee}$ stands for the class $[\mc{O}_{\mc{M}_{(p)}}] - [\mc{E}] \in F^1_{\mc{M}_{(p)}}(\mc{M}_{(p)})_{\Q}$. Given arbitrary $T$ (not necessarily block diagonal), select any $\gamma \in \GL_m(\mc{O}_{F,(p)}) \cap M_{m,m}(\mc{O}_F)$ such that
    \begin{equation}
    {}^t \overline{\gamma} T \gamma = \mrm{diag}(0,T^{\prime\flat})
    \end{equation}
where $\det T^{\prime \flat} \neq 0$. Set $T' \coloneqq \mrm{diag}(0,T^{\prime\flat})$. We define ${}^{\mbb{L}} \mc{Z}(T)_{\ms{V},p}$ to be the pullback class
    \begin{equation}\label{equation:arith_cycle_classes:vertical}
    \begin{tikzcd}[row sep = tiny]
    \mrm{gr}^m_{\mc{M}}(\mc{Z}(T)_{\F_p})_{\Q} & \arrow{l} \mrm{gr}^m_{\mc{M}}(\mc{Z}(T')_{\F_p})_{\Q} \\
    {}^{\mbb{L}} \mc{Z}(T)_{\ms{V},p} & \arrow[mapsto]{l} {}^{\mbb{L}} \mc{Z}(T')_{\ms{V},p}
    \end{tikzcd}
    \end{equation}
along the open and closed immersion $\mc{Z}(T)_{\F_p} \ra \mc{Z}(T')_{\F_p}$ induced by $\gamma$.

By Lemma \ref{lemma:arith_cycle_classes:vertical:p-local_linear_invariance} (applied to $T^{\prime \flat}$, in the notation above), the preceding definition of ${}^{\mbb{L}} \mc{Z}(T)_{\ms{V},p}$ does not depend on the choice of $\gamma$. Moreover, the class $\mc{Z}(T)_{\ms{V},p}$ is \emph{linearly invariant} in the following sense: given any $\gamma \in \GL_m(\mc{O}_{F,(p)}) \cap M_{m,m}(\mc{O}_F)$ (no additional assumptions on ${}^t \overline{\gamma} T \gamma)$, the pullback along
    \begin{equation}\label{equation:arith_cycle_classes:vertical:linear_invariance}
    \mc{Z}(T)_{\F_p} \ra \mc{Z}({}^t \overline{\gamma} T \gamma)_{\F_p}
    \end{equation}
sends ${}^{\mbb{L}} \mc{Z}({}^t \overline{\gamma} T \gamma)_{\ms{V},p}$ to ${}^{\mbb{L}} \mc{Z}(T)_{\ms{V},p}$. This follows from the construction of ${}^{\mbb{L}} \mc{Z}(T)_{\ms{V},p}$.

We also see that ${}^{\mbb{L}}\mc{Z}(T)_{\ms{V},p} = 0$ for all but finitely many primes $p$. 
    
    \clearpage


    \part{Local special cycles}
    \label{part:local_cycles}

        \section{Moduli spaces of \texorpdfstring{$p$}{p}-divisible groups}
        \label{sec:moduli_pDiv}
            We review some unitary Rapoport--Zink spaces \cite{RZ96} and their special cycles, which appear in $p$-adic uniformization of the moduli stacks from Section \ref{sec:part_I:arith_intersections} (also \crefext{III:sec:ab_var}). Some (freely used) notation on $p$-divisible groups appears in \crefext{III:appendix:pDiv_prelim:terminology}.

Fix a prime $p$ and let $F / \Q_p$ be a degree $2$ \'etale algebra, i.e. $F / \Q_p$ is an inert quadratic extension, a ramified quadratic extension, or $F = \Q_p \times \Q_p$. If $F / \Q_p$ is ramified, we assume $p \neq 2$. Write $a \mapsto a^{\s}$ for the nontrivial automorphism $\s$ of $F$ over $\Q_p$. We write $\mc{O}_F$ for the integral closure of $\Z_p$ in $F$ (hence $\mc{O}_F = \Z_p \times \Z_p$ in the split case). If $F / \Q_p$ is ramified, we write $\varpi$ for a uniformizer of $\mc{O}_F$ satisfying $\varpi^{\s} = - \varpi$.

We use the usual notation $\breve{\Q}_p$ for the completion of the maximal unramified extension of $\Q_p$, with ring of integers $\breve{\Z}_p$. 

For $F / \Q_p$ nonsplit, let $\breve{F}$ be the completion of the maximal unramified extension of $F$. If $F / \Q_p$ is split, set $\breve{F} = \breve{\Q}_p$, and view $\breve{F}$ as an $F$-algebra by choosing one of the two morphisms of $\Q_p$-algebras $F \ra \breve{F}$. We also equip $\breve{F}$ with the structure of a $\breve{\Q}_p$-algebra (taking the identity map if $F / \Q_p$ is split).

In all cases, let $\mc{O}_{\breve{F}} \subseteq \breve{F}$ be the ring of integers and let $\overline{k}$ be the residue field of $\breve{F}$. There is a canonical map $\mc{O}_F \ra \mc{O}_{\breve{F}}$ (using the above choice of $F \ra \breve{F}$ when $F / \Q_p$ is split).

We write $\Delta \subseteq \Z_p$ (resp. $\mf{d} \subseteq \mc{O}_F$) for the discriminant ideal (resp. different ideal), which is $\Delta = \Z_p$ and $\mf{d} = \mc{O}_F$ in the split case.
We also abuse notation and write $\mf{d}$ for a chosen generator of the different ideal satisfying $\mf{d}^{\s} = -\mf{d}$, taking $\mf{d} = \varpi$ in the ramified case.

In the split case, let $e^+$ (resp. $e^-$) be the nontrivial idempotent in $\mc{O}_F$ which maps to $1 \in \mc{O}_{\breve{F}}$ (resp. $0 \in \mc{O}_{\breve{F}}$). Given an $\mc{O}_F$-module $M$, we write $M = M^+ \oplus M^-$ where $e^+$ projects to $M^+$ and $e^-$ projects to $M^-$. We use similar notation $f = f^+ \oplus f^-$ for morphisms $f \colon M \ra M'$ of $\mc{O}_F$-modules. We often use this for $p$-divisible groups $X$ with $\mc{O}_F$ action, e.g. $X = X^+ \times X^-$ (and similarly for $\mc{O}_F$-linear quasi-homomorphisms).

            \subsection{Rapoport--Zink spaces}
            \label{ssec:moduli_pDiv:RZ}
                \begin{definition}\label{definition:moduli_pDiv:RZ:Hermitian_p-divisible}
Let $S$ be a formal scheme and let $n \geq 1$ be an integer. By a \emph{Hermitian $p$-divisible group} over $S$, we mean a tuple $(X, \iota, \lambda)$ where
    \begin{align*}
    & X & & \parbox[t]{0.75\textwidth}{is a $p$-divisible group over $S$ of height $2n$ and dimension $n$
    }
    \\
    & \iota \colon \mc{O}_F \ra \End(X) & & \parbox[t]{0.75\textwidth}{is a ring homomorphism
    } 
    \\
    & \lambda \colon X \ra X^{\vee} & & \parbox[t]{0.75\textwidth}{is a quasi-polarization satisfying:
    \begin{enumerate}[(1)]
        \item (Action compatibility) The Rosati involution $\dagger$ on $\End^{0}(A)$ satisfies $\iota(a)^{\dagger} = \iota(a^{\sigma})$ for all $a \in \mc{O}_F$.
    \end{enumerate}
    }
    \end{align*}
\end{definition}

An \emph{isomorphism} of Hermitian $p$-divisible groups is an isomorphism of underlying $p$-divisible groups which respects the $\mc{O}_F$-actions and polarizations.

In Sections \ref{sec:moduli_pDiv}--\ref{sec:can_and_qcan}, we only consider Hermitian $p$-divisible groups over formal schemes $S$ equipped with a morphism $S \ra \Spf \mc{O}_{\breve{F}}$, and we assume that $X$ is supersingular (resp. ordinary) if $F / \Q_p$ is nonsplit (resp. split). 

We primarily discuss Hermitian $p$-divisible groups satisfying either of the following two conditions.
    \begin{enumerate}[(1)]
        \setcounter{enumi}{1}
        \item (Principal polarization) The quasi-polarization $\lambda$ is a principal polarization.
    \end{enumerate}
    \begin{enumerate}[(1$^{\circ}$)]
        \setcounter{enumi}{1}
        \item (Polarization condition $\circ$) Assume $n$ is even if $F / \Q_p$ is ramified. The quasi-polarization $\Delta \lambda$ is a polarization, and $\ker(\Delta \lambda) = X[\iota(\mf{d})]$.
    \end{enumerate}
In these cases, we say that $(X, \iota, \lambda)$ is \emph{principally polarized} or \emph{$\circ$-polarized} respectively.\footnote{The local analogue of Footnote \ref{footnote:pi_inverse_modular} applies as well.} If $F / \Q_p$ is ramified, we also use the alternative terminology \emph{$\varpi^{-1}$-modular}.

Given an integer $r$ with $0 \leq r \leq n$, we next consider
    \begin{enumerate}[(1)]
        \item (Kottwitz $(n - r, r)$ signature condition) For all $a \in \mc{O}_F$, the characteristic polynomial of $\iota(a)$ acting on $\Lie X$ is $(T - a)^{n-r}(T - a^{\s})^r \in \mc{O}_S[T]$.
    \end{enumerate}
for pairs $(X, \iota)$, i.e. $n$-dimensional $p$-divisible groups $X$ over a formal scheme $S$ with action $\iota \colon \mc{O}_F \ra \End(X)$.
Here we view $\mc{O}_S$ as an $\mc{O}_F$-algebra via $S \ra \Spf \mc{O}_{\breve{F}} \ra \Spec \mc{O}_F$. 

If $(X,\iota,\lambda)$ is a Hermitian $p$-divisible group of signature $(n-r,r)$, then $(X^{\s}, \iota^{\s}, \lambda^{\s})$ with
    \begin{equation}\label{equation:sigma_twist_triple}
    X^{\s} = X \quad \quad \iota^{\s} = \iota \circ \sigma \quad \quad \lambda^{\sigma} = \lambda
    \end{equation}
is a Hermitian $p$-divisible group of signature $(r, n-r)$. We use similar notation $(X, \iota) \leftrightarrow (X^{\s}, \iota^{\s})$ without the presence of a polarization. This allows us to switch between signatures $(n - r, r)$ and signature $(r, n - r)$ (e.g. for comparison with the literature).

From here on, we always implicitly restrict to signature $(n - 1, 1)$ (and even $n$) when discussing $\circ$-polarized Hermitian $p$-divisible groups for ramified $F / \Q_p$.
In this case, we also impose
    \begin{enumerate}[(1)]
        \setcounter{enumi}{1}
        \item (Pappas wedge condition) For all $a \in \mc{O}_F$, the action of $\iota(a)$ on $\Lie X$ satisfies
            \begin{equation*}
            \bigwedge^{2} (\iota(a) - a) = 0 \quad \text{ and } \quad \bigwedge^{n}(\iota(a) - a^{\s}) = 0.
            \end{equation*}
        \item (Pappas--Rapoport--Smithling--Zhang spin condition) For every geometric point $\overline{s}$ of $S$, the action of $(\iota(a) - a)$ on $\Lie X_{\overline{s}}$ is nonzero for some $a \in \mc{O}_F$.
    \end{enumerate}
The signature $(n - 1, 1)$ condition implies that the equation involving $\bigwedge^{n}$ in the wedge condition is automatic, and that the wedge condition is empty if $n = 2$.
    
We temporarily allow $p = 2$ even if $F / \Q_p$ is ramified. Recall that there exists a supersingular (resp. ordinary) $p$-divisible group $\mbf{X}_0$ of height $2$ and dimension $1$ over $\overline{k}$, and that $\mbf{X}_0$ is unique up to isomorphism (this also holds for any algebraically closed field $\kappa$ over $\overline{k}$).
In the supersingular case $\mbf{X}_0$ is given by a Lubin--Tate formal group law, and in the ordinary case we have $\mbf{X}_0 \cong \pmb{\mu}_{p^{\infty}} \times \underline{\Q_p / \Z_p}$. We have $\End(\mbf{X}_0) \cong \mc{O}_D$ (resp. $\End(\mbf{X}_0) \cong \Z_p \times \Z_p$) in the supersingular case (resp. ordinary case) where $\mc{O}_D$ is the unique maximal order in the quaternion division algebra $D$ over $\Q_p$ (e.g. \cite{Gross86} or \cite[{\S 1}]{Wewers07}). 

Quasi-polarizations on $\mbf{X}_0$ exist and are unique up to $\Q_p^{\times}$ scalar, and there exists a principal polarization $\lambda_{\mbf{X}_0}$ on $\mbf{X}_0$ (unique up to $\Z_p^{\times}$ scalar). 
The induced Rosati involution on $\End(\mbf{X}_0)$ is the standard involution (this can be verified on on the Dieudonn\'e module, see e.g. \cite[{Page 2205}]{RSZ17}), hence induces the nontrivial Galois involution on $F$ for any embedding $F \hookrightarrow \End^0(\mbf{X}_0)$ (if such an embedding exists).

From now on, we assume $\mbf{X}_0$ is supersingular (resp. ordinary) if $F / \Q_p$ is nonsplit (resp. split). Then there exists an embedding $j \colon \mc{O}_F \hookrightarrow \End(\mbf{X}_0)$. Given such a $j$, form $(\mbf{X}_0, j)$ and $(\mbf{X}_0^{\s}, j^{\s})$ as above. There is an $\mc{O}_F$-linear isogeny of degree $|\Delta|_p^{-1}$
    \begin{equation}\label{equation:moduli_pDiv:RZ:isogeny_from_tensor}
    \begin{tikzcd}[row sep = tiny]
    \mbf{X}_0 \otimes_{\Z_p} \mc{O}_F \arrow{r} & \mbf{X}_0 \times \mbf{X}_0^{\s}
    \\
    x \otimes a \arrow[mapsto]{r} & (j(a) x, j(a^{\s}) x).
    \end{tikzcd}
    \end{equation}
where $X \otimes_{\Z_p} \mc{O}_F$ is the Serre tensor $p$-divisible group (as in \crefext{III:equation:appendix:Serre_tensor}), with its Serre tensor $\mc{O}_F$-action. See also \cite[{Lemma 6.2}]{KR11} (there in the inert case for $p \neq 2$, but the version in \eqref{equation:moduli_pDiv:RZ:isogeny_from_tensor} allows for $p = 2$).

Suppose $\lambda_{\mbf{X}_0}$ is a principal polarization of $\mbf{X}_0$. Under the map in \eqref{equation:moduli_pDiv:RZ:isogeny_from_tensor}, the ($\mc{O}_F$-action compatible) product polarization $\lambda_{\mbf{X}_0} \times \lambda_{\mbf{X}_0}^{\s}$ on $\mbf{X}_0 \times \mbf{X}_0^{\s}$ pulls back to the polarization 
    \begin{equation}\label{equation:moduli_pDiv:RZ:isogeny_from_tensor:polarizations}
    \lambda_{\mbf{X}_0} \otimes \lambda_{\mrm{tr}} \colon \mbf{X}_0 \otimes_{\Z_p} \mc{O}_F \ra (\mbf{X}_0 \otimes_{\Z_p} \mc{O}_F^*) \cong \mbf{X}_0^{\vee} \otimes_{\Z_p} \mc{O}_F^*
    \end{equation}
where $\mc{O}_F^* \coloneqq \Hom_{\Z_p}(\mc{O}_F, \Z_p)$ and $\lambda_{\mrm{tr}} \colon \mc{O}_F \ra \mc{O}_F^*$ is induced by the symmetric bilinear pairing $\mrm{tr}_{F / \Q_p}(a^{\s} b)$ on $\mc{O}_F$. Indeed, after picking a $\Z_p$-basis $\{1, \a\}$ for $\mc{O}_F$ to identify $\mbf{X}_0 \otimes_{\Z_p} \mc{O}_F \cong \mbf{X}_0^2$, the map in \eqref{equation:moduli_pDiv:RZ:isogeny_from_tensor} is given by the matrix
    \begin{equation}\label{equation:moduli_pDiv:RZ:isogeny_from_tensor:matrix}
    \phi =
    \begin{pmatrix}
    1 & \a \\
    1 & \a^{\s}
    \end{pmatrix} \in M_{2,2}(\mc{O}_F)
    \end{equation}
and $\det \phi$ generates the different ideal of $\mc{O}_F / \Z_p$ (so Smith normal form shows $\deg \phi = |\Delta|_p^{-1}$). Identifying $\mbf{X}_0^{\vee} \otimes_{\Z_p} \mc{O}_F^* \cong \mbf{X}_0^{\vee 2}$ using the basis of $\mc{O}_F^*$ dual to $\{1, \a\}$, the preceding claim about pullback polarizations follows because $({}^t \phi^{\s}) \phi$ (where ${}^t \phi^{\s}$ means conjugate transpose) is the Gram matrix for the basis $\{1,\a\}$ and the trace pairing on $\mc{O}_F$.

If $p \neq 2$, the polarization $\lambda_{\mbf{X}_0} \otimes \lambda_{\mrm{tr}}$ coincides with the polarization on $\mbf{X}_0 \otimes_{\Z_p} \mc{O}_F$ described in \cite[{(6.2)}]{KR11} (inert case, with modification as in \cite[{Footnote 4}]{RSZ17}) and \cite[{(3.4)}]{RSZ17} (ramified case, though we normalize differently).

Suppose $\iota_{\mbf{X}_0} \colon \mc{O}_F \hookrightarrow \End(\mbf{X}_0)$ is an action of signature $(1,0)$. The pair $(\mbf{X}_0, \iota_{\mbf{X}_0})$ exists and is unique up to isomorphism.
In the split case, the element $\iota_{\mbf{X}_0}(e^+)$ is projection to $\pmb{\mu}_{p^{\infty}}$ and $\iota_{\mbf{X}_0}(e^-)$ is projection to $\underline{\Q_p / \Z_p}$. In the ramified case, note that $(\mbf{X}_0, \iota_{\mbf{X}_0})$ is simultaneously of signature $(1,0)$ and $(0,1)$.

Fix $(\mbf{X}_0, \iota_{\mbf{X}_0})$ as above, and form $(\mbf{X}_0^{\s}, \iota_{\mbf{X}_0}^{\s})$. We have
    \begin{equation}\label{equation:moduli_pDiv:RZ:hom_to_sigma}
    \Hom_{\mc{O}_F}(\mbf{X}_0, \mbf{X}_0^{\s}) \cong
    \begin{cases}
    \mc{O}_F & \text{if $F / \Q_p$ is nonsplit} \\
    0 & \text{if $F / \Q_p$ is split}
    \end{cases}
    \end{equation}
as $\mc{O}_F$-modules by precomposition. Using $\End(\mbf{X}_0) \cong \mc{O}_D$ in the nonsplit cases, we find that the $\mc{O}_F$-module $\Hom_{\mc{O}_F}(\mbf{X}_0, \mbf{X}_0^{\s})$ is generated by any isogeny of degree $p$ if $F / \Q_p$ is inert, and is generated by an isomorphism if $F / \Q_p$ is ramified (namely any element $a \in \mc{O}_D^{\times}$ such that conjugation by $a$ induces the nontrivial Galois involution on $F$). We have $\End_{\mc{O}_F}(\mbf{X}_0) = \mc{O}_F$ in all cases.

Suppose $\lambda_{\mbf{X}_0}$ is a principal polarization of $\mbf{X}_0$.
If $F / \Q_p$ is unramified, the triple $(\mbf{X}_0, \iota_{\mbf{X}_0}, \lambda_{\mbf{X}_0})$ is unique (up to isomorphism): given another polarization $\lambda'_{\mbf{X}_0}$, we have $\lambda_{\mbf{X}_0}^{-1} \circ \lambda_{\mbf{X}_0}' \in \Z_p^{\times}$ and know that the norm map $N_{F/\Q_p} \colon \mc{O}_F^{\times} \ra \Z_p^{\times}$ is surjective. 
If $F / \Q_p$ is ramified, the same reasoning shows that there are two choices of $\lambda_{\mbf{X}_0}$ (differing by a $\Z_p^{\times}$ scalar) because $N_{F/\Q_p}(\mc{O}_F^{\times}) \subseteq \Z_p^{\times}$ has index $2$. Fix a choice of $(\mbf{X}_0, \iota_{\mbf{X}_0}, \lambda_{\mbf{X}_0})$.

We now re-impose our running assumption that $p \neq 2$ if $F / \Q_p$ is ramified. Fix any $\circ$-polarized Hermitian $p$-divisible group $(\mbf{X}, \iota_{\mbf{X}}, \lambda_{\mbf{X}})$ of signature $(n-r,r)$ over $\overline{k}$. Such triples exist and are unique up to $F$-linear quasi-isogenies preserving polarizations exactly. 
This uniqueness may be proved via Dieudonn\'e theory: see \cite[{\S 1}]{Vollaard10} (inert case, but we allow $p = 2$ by the same proof) and \cite[{Proposition 3.1}]{RSZ17} \cite[{\S 6}]{RSZ18} (ramified case). In the split case, we have a stronger uniqueness statement.

\begin{lemma}\label{lemma:moduli_pDiv:RZ:unique_framing_object_split}
For $F / \Q_p$ split, the triple $(\mbf{X}, \iota_{\mbf{X}}, \lambda_{\mbf{X}})$ is unique up to isomorphism. This also holds over any algebraically closed field extension $\kappa$ of $\overline{k}$.
\end{lemma}
\begin{proof}
Decompose $\mbf{X} = \mbf{X}^+ \times \mbf{X}^-$ using the idempotents in the $\mc{O}_F = \Z_p \times \Z_p$ action given by $\iota_{\mbf{X}}$. Then $\mbf{X}^+$ and $\mbf{X}^-$ are the unique ordinary $p$-divisible groups over $\kappa$ of height $n$ and the correct dimension ($n - r$ and $r$, respectively).
Uniqueness of $\lambda_{\mbf{X}}$ (up to isomorphism) corresponds to the following fact: there is a unique self-dual Hermitian $\mc{O}_F$-lattice (up to isomorphism) of any given rank.
\end{proof}

For existence, we may construct $(\mbf{X}, \iota_{\mbf{X}}, \lambda_{\mbf{X}})$ as follows. For $F / \Q_p$ unramified, we can take $\mbf{X} = (\mbf{X}_0)^{n - r} \times (\mbf{X}_0^{\s})^r$ (with the product $\mc{O}_F$-action and polarization). 

For $F / \Q_p$ ramified, we can take $\mbf{X} = (\mbf{X}_0)^{n - 2} \times (\mbf{X}_0 \otimes_{\Z_p} \mc{O}_F)$ using the Serre tensor construction. The $\mc{O}_F$-action $\iota_{\mbf{X}}$ is diagonal on $(\mbf{X}_0)^{n - 2}$ and given by the Serre tensor $\mc{O}_F$-action on $\mbf{X}_0 \otimes_{\Z_p} \mc{O}_F$. We can take the product quasi-polarization $\lambda_{\mbf{X}}$ of $\mbf{X}$ given by
    \begin{align}\label{equation:moduli_pDiv:RZ:framing_object_polarization:ramified}
    &
    - \iota_{\mbf{X}_0}(\varpi)^{-2} \circ (\lambda_{\mbf{X}_0} \otimes \lambda_{\mrm{tr}})
    \quad \text{on $\mbf{X}_0 \otimes_{\Z_p} \mc{O}_F$, and} \\
    &
    \underbrace{
    \begin{pmatrix}
    0 & \lambda_{\mbf{X}_0} \circ - \iota_{\mbf{X}_0}(\varpi)^{-1} \\
    \lambda_{\mbf{X}_0} \circ \iota_{\mbf{X}_0}(\varpi)^{-1} & 0
    \end{pmatrix}
    \times \cdots \times
    \begin{pmatrix}
    0 & \lambda_{\mbf{X}_0} \circ - \iota_{\mbf{X}_0}(\varpi)^{-1} \\
    \lambda_{\mbf{X}_0} \circ \iota_{\mbf{X}_0}(\varpi)^{-1} & 0
    \end{pmatrix}
    }_{(n-2)/2 \text{ times}} \quad \text{on $\mbf{X}_0^{n - 2}$.} \notag
    \end{align}
This is the construction of \cite[{\S 3.3}]{RSZ17} (but rescaled).

Given a principally polarized triple $(X_0, \iota_0, \lambda_0)$ of signature $(1,0)$ over some base scheme $S$, a \emph{framing similitude quasi-isogeny $\rho_0$} is an $F$-linear quasi-isogeny $X_{0, \overline{S}} \ra \mbf{X}_{0, \overline{S}}$ such that 
    \begin{equation}
    \rho_0^*(\lambda_{\mbf{X}_0, \overline{S}}) = b \lambda_{0, \overline{S}} \quad \text{for some $b \in \Q_p^{\times}$}
    \end{equation}
where the subscript indicates base-change to $\overline{S} \coloneqq S \times_{\Spec \mc{O}_{\breve{F}}} \Spec \overline{k}$ (and where $b \in \Q_p^{\times}$ really means a section of the constant sheaf). We call $(X_0, \iota_0, \lambda_0, \rho)$ a \emph{framed similitude} tuple. An \emph{isomorphism} of framed similitude tuples $f \colon (X_0, \iota_0, \lambda_0, \rho_0) \ra (X_0', \iota_0', \lambda_0', \rho_0')$ is an $\mc{O}_F$-linear isomorphism of $p$-divisible groups $f \colon X_0 \ra X_0'$ such that $f^*(\lambda_0')$ and $\lambda_0$ agree up to $\Z_p^{\times}$-scalar and also $\rho_0' \circ f_{\overline{S}} = \rho_0$.

Given a $\circ$-polarized triple $(X, \iota, \lambda)$ of signature $(n - r, r)$ over some base scheme $S$, we define a \emph{similitude framing quasi-isogeny} $\rho \colon X_{\overline{S}} \ra \mbf{X}_{\overline{S}}$ in the same way. A \emph{framing quasi-isogeny} $\rho \colon X_{\overline{S}} \ra \mbf{X}_{\overline{S}}$ is given by the stricter requirement $b \in \Z_p^{\times}$. In these two cases, we call the datum $(X, \iota, \lambda, \rho)$ a \emph{framed similitude} tuple and a \emph{framed tuple}, respectively. In both cases, isomorphisms of two such tuples are defined as before: isomorphisms of $p$-divisible groups which are $\mc{O}_F$-linear, preserve polarizations up to $\Z_p^{\times}$, and commute with framings.

\begin{definition}\label{definition:moduli_pDiv:RZ:RZ_spaces}
We consider three \emph{Rapoport--Zink spaces} over $\Spf \mc{O}_{\breve{F}}$, given by the (set-valued) functors
    \begin{align*}
    \mc{N}(1,0)'(S) & \coloneqq \{ \text{isomorphism classes of framed similitude tuples } (X_0, \iota_0, \lambda_0, \rho_0) \text{ over $S$} \}
    \\
    \mc{N}(n-r,r)'(S) & \coloneqq \{ \text{isomorphism classes of framed similitude tuples } (X, \iota, \lambda, \rho) \text{ over $S$} \}
    \\
    \mc{N}(n-r,r)(S) & \coloneqq \{ \text{isomorphism classes of framed tuples } (X, \iota, \lambda, \rho) \text{ over $S$} \}
    \end{align*}
for schemes $S$ over $\Spf \mc{O}_{\breve{F}}$. Here, signature $(1,0)$ and principal polarizations are understood for $\mc{N}(1,0)'$. Signature $(n - r, r)$ and $\circ$-polarizations are understood for $\mc{N}(n - r, r)'$ and $\mc{N}(n - r, r)$.
\end{definition}

These Rapoport--Zink spaces do not depend on the choices of framing objects (up to functorial isomorphism). The functor $\mc{N}(n - r, r)$ is canonically isomorphic to its variant where instead we require framing quasi-isogenies and isomorphisms of framed tuples to preserve polarizations exactly (not just up to $\Z_p^{\times}$ scalar). If $S$ is a formal scheme, we also write e.g. $\mc{N}(n - r, r)(S) \coloneqq \Hom(S, \mc{N}(n - r, r))$.

\begin{lemma}
Each of $\mc{N}(1,0)'$, $\mc{N}(n - r, r)'$, and $\mc{N}(n - r, r)$ is represented by a locally Noetherian formal scheme which is formally locally of finite type and separated over $\Spf \mc{O}_{\breve{F}}$. Each irreducible component of the reduced subschemes is projective over $\overline{k}$.    
\end{lemma}
\begin{proof}
Representability, local Noetherianity, and formally locally of finite type-ness follow via \cite[{Theorem 2.16}]{RZ96}; various closedness statements can be checked via \cite[{Proposition 2.9}]{RZ96}, which holds verbatim with ``isogeny'' replaced by ``homomorphism''. Projectivity of the reduced irreducible components follows from \cite[{Proposition 2.32}]{RZ96}, also using \cite[{Proposition 3.8}]{RSZ17} in the ramified case. Separatedness now follows because this can be checked on underlying reduced subschemes
(then apply the valuative criterion).    
\end{proof}

\begin{lemma}
The formal scheme $\mc{N}(n - r, r)$ is regular and the structure morphism $\mc{N}(n - r, r) \ra \Spf \mc{O}_{\breve{F}}$ is formally smooth of relative dimension $(n - r, r)$.    
\end{lemma}
\begin{proof}
We know the structure map $\mc{N}(n - r, r) \ra \Spf \mc{O}_{\breve{F}}$ is formally smooth of relative dimension $(n - r, r)$ via \cite[{Proposition 1.3}]{Mihatsch22} in the unramified case (also \crefext{III:ssec:ab_var:generic_smoothness}, where we allow $p = 2$) and \cite[{Proposition 3.8}]{RSZ17} for the ramified case. We conclude $\mc{N}(n - r, r)$ is regular because the map to $\Spf \mc{O}_{\breve{F}}$ is formally smooth and formally locally of finite type.
\end{proof}

For $F / \Q_p$ ramified, the Rapoport--Zink space $\mc{N}(n - 1, 1)$ is often called \emph{exotic smooth} in the literature, following the terminology of \cite{RSZ17}.

\begin{lemma}\label{lemma:N(1,0)'_isomorphism}
There is an isomorphism
    \begin{equation}\label{equation:N(1,0)'_isomorphism}
    \begin{tikzcd}[row sep = tiny]
    F^{\times} / \mc{O}_F^{\times} \arrow{r}{\sim} & \mc{N}(1,0)' \\
    a \arrow[mapsto]{r} & (\mbf{X}_0, \iota_{\mbf{X}_0}, \lambda_{\mbf{X}_0}, a \cdot \mrm{id}_{\mbf{X}_0})
    \end{tikzcd}
    \end{equation}
where the left-hand side is viewed as a constant formal scheme over $\Spf \mc{O}_{\breve{F}}$.
\end{lemma}
\begin{proof}
In the nonsplit case, \cite[{Proposition 2.1}]{Howard19} states that $\mc{N}(1,0)'$ is a disjoint union of copies of $\Spf \mc{O}_{\breve{F}}$, so it is enough to check the claim on $\overline{k}$-points. The claim on $\overline{k}$-points follows from uniqueness of the triple $(\mbf{X}_0, \iota_{\mbf{X}_0}, \lambda_{\mbf{X}_0})$ (up to isomorphism preserving polarizations up to $\Z_p^{\times}$ scalar) and the equality $\End_{F}^{0}(\mbf{X}_0) = F$.

The split case holds via the following similar argument. The map $\mc{N}(1,0)'(\kappa) \ra \mc{N}(1,0)'(\kappa')$ is bijective for any extension of algebraically closed fields $\kappa \subseteq \kappa'$ (essentially by uniqueness of the triple $(\mbf{X}_0, \iota_{\mbf{X}_0}, \lambda_{\mbf{X}_0})$, which holds over any algebraically closed field of characteristic $p$). 
So the reduced subscheme $\mc{N}(1,0)'_{\red}$ is isomorphic to a (discrete) disjoint union of copies of $\Spec \overline{k}$.
We also see that the map $F^{\times} / \mc{O}_F^{\times} \ra \mc{N}(1,0)'$ is bijective on $\overline{k}$-points (this follows as in the nonsplit case).
To finish, note that $\mc{N}(1,0)'$ is formally \'etale over $\Spf \mc{O}_{\breve{F}}$ (e.g. by Grothendieck--Messing theory as in \crefext{III:ssec:ab_var:generic_smoothness}).
\end{proof}

\begin{definition}\label{definition:moduli_pDiv:RZ:canonical_lifting}
By the \emph{canonical lifting} of $(\mbf{X}_0, \iota_{\mbf{X}_0}, \lambda_{\mbf{X}_0})$, we mean the tuple $(\mf{X}_0, \iota_{\mf{X}_0}, \lambda_{\mf{X}_0}, \rho_{\mf{X}_0})$ over $\Spf \mc{O}_{\breve{F}}$ corresponding (via Lemma \ref{lemma:N(1,0)'_isomorphism}) to the unique section $\Spf \mc{O}_{\breve{F}} \ra \mc{N}(1,0)'$ associated to the element $1 \in F^{\times} / \mc{O}_F^{\times}$.
\end{definition}

\begin{definition}\label{definition:moduli_pDiv:RZ:N_prime}
We define the open and closed subfunctor $\mc{N}' \subseteq \mc{N}(1,0)' \times \mc{N}(n - r, r)'$ as
    \begin{equation*}
    \mc{N}'(S) \coloneqq \left \{(X_0, \iota_0, \lambda_0, \rho_0, X, \iota, \lambda, \rho) : \begin{array}{l} \rho_0^*(\lambda_{\mbf{X}_0, \overline{S}}) = b_0 \lambda_{0, \overline{S}} \quad \rho^*(\lambda_{\mbf{X}, \overline{S}}) = b \lambda_{\overline{S}} \\ \text{with $b_0 = b$ in $\Q_p^{\times} / \Z_p^{\times}$ } \end{array} \right \}
    \end{equation*}
for schemes $S$ over $\Spf \mc{O}_{\breve{F}}$.
\end{definition}

With $b$ as above and $a \in F^{\times}$ any element with $N_{F / \Q_p}(a) = b$ in $\Q_p^{\times} / \Z_p^{\times}$, there is an isomorphism
    \begin{equation} \label{equation:moduli_pDiv:RZ:unitary_iso}
    \begin{tikzcd}[row sep = tiny]
    \mc{N}' \arrow{r}{\sim} & \mc{N}(1,0)' \times \mc{N}(n-r,r) \\
    (X_0,\iota_0,\lambda_0,\rho_0,X,\iota,\lambda,\rho) \arrow[mapsto]{r} & (X_0,\iota_0,\lambda_0,\rho_0,X,\iota,\lambda, a^{-1} \rho).
    \end{tikzcd}
    \end{equation}
Whenever we write $(X_0,\iota_0,\lambda_0,\rho_0,X,\iota,\lambda,\rho)$ for a (functorial) point of $\mc{N}'$, we mean an object as on the left of \eqref{equation:moduli_pDiv:RZ:unitary_iso} (i.e. $\rho$ preserves polarizations up to $\Q_p^{\times}$ scalar).
    
The functorial assignment $(X, \iota, \lambda, \rho) \mapsto \Lie X$ defines a locally free sheaf $\Lie$ on $\mc{N}(n - r, r)$. In the case of signature $(n - 1, 1)$, there is a unique maximal local direct summand $\mc{F} \subseteq \Lie$ of rank $n - 1$ such that the $\iota$ action on $\mc{F}$ (resp. $\Lie / \mc{F}$) is $\mc{O}_F$-linear (resp. $\s$-linear). 
The ramified case is proved in \cite[{Lemma 2.36}]{LL22II} (and in the unramified case, we have a canonical eigenspace decomposition $\Lie = \mc{F} \oplus (\Lie / \mc{F})$ for the $\mc{O}_F$-action).

Consider the canonical lifting $(\mf{X}_0, \iota_{\mf{X}_0}, \lambda_{\mf{X}_0}, \rho_{\mf{X}_0})$ over $\Spf \mc{O}_{\breve{F}}$ (Definition \ref{definition:moduli_pDiv:RZ:canonical_lifting}). Given any $\Spf \mc{O}_{\breve{F}}$-scheme $S$, we write $\mbb{D}(\mf{X}_{0,S})(S)$ for evaluation of the (covariant) Dieudonn\'e crystal $\mbb{D}(\mf{X}_{0,S})$ at $\mrm{id} \colon S \ra S$, with associated Hodge filtration step $F^0 \mbb{D}(\mf{X}_{0,S})(S)$.
The assignment $S \mapsto F^0 \mbb{D}(\mf{X}_{0,S})(S)$ defines a (trivial) line bundle on $\mc{N}$, which we denote $\Lie_0^{\vee}$. The principal polarization $\lambda_{\mf{X}_0}$ induces an identification $(\Lie \mf{X}_{0,S})^{\vee} \cong (\Lie \mf{X}^{\vee}_{0,S})^{\vee}$ and the latter is $F^0 \mbb{D}(\mf{X}_{0,S})(S)$.

\begin{definition}\label{definition:moduli_pDiv:RZ:tautological}
The \emph{tautological bundle} $\mc{E}$ on $\mc{N}(n - 1, 1)$ is the line bundle whose dual is $\mc{E}^{\vee} \coloneqq \underline{\Hom}(\Lie_0^{\vee}, \Lie / \mc{F})$.
\end{definition}
The definition of $\mc{E}$ is taken from \cite[{Definition 3.4}]{Howard19} (at least in the inert case).
The line bundle $\mc{E}$ on $\mc{N}(n - 1, 1)$ is a local analogue of the global tautological bundle (\cref{definition:ab_var:exotic_smooth:tautological_bundle}).
We are recycling the notation $\mc{E}$ (but the global tautological bundle pulls back to $\mc{E}$ under Rapoport--Zink uniformization, see e.g. \crefext{III:ssec:non-Arch_uniformization:vertical}).

            \subsection{Local special cycles} 
            \label{ssec:moduli_pDiv:special_cycles}
                We define certain local special cycles on Rapoport--Zink spaces, following \cite[{Definition 3.2}]{KR11} (there in the inert case).
Retain notation from Section \ref{ssec:moduli_pDiv:RZ}.

The space of \emph{local special quasi-homomorphisms} means the $F$-module
    \begin{equation}\label{equation:local_special_quasi-hom}
    \mbf{W} \coloneqq \Hom_{F}^0(\mbf{X}_0, \mbf{X}).
    \end{equation}
If $F/\Q_p$ is nonsplit, then $\mbf{W}$ is free of rank $n$ (see also \cite[{Lemma 3.5}]{RSZ17} in the ramified case). If $F / \Q_p$ is split, then $\mbf{W}$ is a free $F$-module of rank $n - r$ (because $\Hom_{\mc{O}_F}(\mbf{X}_0, \mbf{X}_0^{\s}) = 0$ in the split case, in contrast with $\Hom_{\mc{O}_F}(\mbf{X}_0, \mbf{X}_0^{\s}) \cong \mc{O}_F$ in the nonsplit cases).
In the split case only, set $\mbf{W}^{\perp} \coloneqq \Hom^0_F(\mbf{X}_0^{\s}, \mbf{X})$. In the nonsplit cases, set $\mbf{W}^{\perp} = 0$. 

Set
    \begin{equation}
    \mbf{V} = \mbf{W} \oplus \mbf{W}^{\perp} \quad \quad \mbf{V}_0 = \Hom^0_F(\mbf{X}_0, \mbf{X}_0).
    \end{equation}
In all cases, these are free $F$-modules of rank $n$ and $1$, respectively.

We equip $\mbf{W}$, $\mbf{W}^{\perp}$, and $\mbf{V}_0$ with the (non-degenerate) Hermitian pairings $(x,y) = x^{\dagger} y \in \End^0_F(\mbf{X}_0) = F$. We give $\mbf{V}$ the Hermitian form making $\mbf{W}$ and $\mbf{W}^{\perp}$ orthogonal. We have $\varepsilon(\mbf{V}) = (-1)^r$ if $F / \Q_p$ is nonsplit (resp. $\varepsilon(\mbf{V}) = 1$ if $F / \Q_p$ is split). This follows upon inspecting the explicit framing tuples constructed in Section \ref{ssec:moduli_pDiv:RZ} (see \cite[{Lemma 3.5}]{RSZ17} for the ramified case).

\begin{definition}[Kudla--Rapoport local special cycles]
Given any set $L \subseteq \mbf{W}$, there is a associated \emph{local special cycle}
    \begin{equation}
    \mc{Z}(L)' \subseteq \mc{N}' \quad \text{(resp. $\tilde{\mc{Z}}(L) \subseteq \mc{N}(n - r, r)$)}
    \end{equation}
which is the subfunctor consisting of tuples $(X_0, \iota_0, \lambda_0, X, \iota, \lambda, \rho)$ (resp. $(X, \iota, \lambda, \rho)$) over schemes $S$ over $\Spf \mc{O}_{\breve{F}}$ such that, for all $x \in L$, the quasi-homomorphism
    \begin{equation}
    \rho^{-1} \circ x_{\overline{S}} \circ \rho_0 \colon X_{0, \overline{S}} \ra X_{\overline{S}} \quad \text{(resp. $\rho^{-1} \circ x_{\overline{S}} \circ \rho_{\mf{X}_{0},\overline{S}} \colon \mf{X}_{0,\overline{S}} \ra X_{\overline{S}}$)}
    \end{equation}
lifts to a homomorphism $X_0 \ra X$ (resp. $\mbf{X}_{0,S} \ra X$). Here $(\mf{X}_0, \iota_{\mf{X}_0}, \lambda_{\mf{X}_0}, \rho_{\mf{X}_0})$ is the canonical lifting.
\end{definition}

In the preceding definition, such lifts are unique (if they exist) by Drinfeld rigidity for quasi-homomorphisms of $p$-divisible groups.
We know that $\mc{Z}(L)' \subseteq \mc{N}'$ and $\tilde{\mc{Z}}(L) \subseteq \mc{N}(n - r, r)$ are closed subfunctors (hence locally Noetherian formal schemes) by \cite[{Proposition 2.9}]{RZ96} for quasi-homomorphisms.
From the definition, it is clear that $\mc{Z}(L)'$ and $\tilde{\mc{Z}}(L)$ depend only on the $\mc{O}_F$-span of $L$.

The isomorphism $\mc{N}' \xra{\sim} \mc{N}(1,0)' \times \mc{N}(n - r, r)$ of \eqref{equation:moduli_pDiv:RZ:unitary_iso} induces an isomorphism
    \begin{equation}\label{equation:moduli_pDiv:special_cycles:unitary_iso}
    \mc{Z}(L)' \xra{\sim} \mc{N}(1,0)' \times \tilde{\mc{Z}}(L).
    \end{equation}

\begin{lemma}
Let $L \subseteq \mbf{W}$ be any subset. If $\mc{Z}(L)' \neq \emptyset$, then $(x,y) \in \mf{d}^{-1}$ for all $x, y \in L$.
\end{lemma}
\begin{proof}
If $\mc{Z}(L)' \neq \emptyset$, then $\mc{Z}(L)(\overline{k})' \neq \emptyset$ because $\mc{Z}(L)' \ra \Spf \mc{O}_{\breve{F}}$ is formally locally of finite type.
If $\mc{Z}(L)(\overline{k})' \neq \emptyset$ and $x,y \in L$, we find $\mf{d} x^{\dagger} y \in \End_{\mc{O}_F}(\mbf{X}_0) = \mc{O}_{F}$ by the $\circ$-polarization condition defining $\mc{N}(n - r,r)$, where $\mf{d}$ is the different ideal.
\end{proof}

If $F / \Q_p$ is nonsplit, set $\mc{Z}(L) \coloneqq \tilde{\mc{Z}}(L)$ for any subset $L \subseteq \mbf{W}$. If $F / \Q_p$ is split, we will instead define $\mc{Z}(L)$ as a certain open and closed subfunctor (see \eqref{equation:moduli_pDiv:discrete_reduced:split_Z(L)}) for later notational uniformity.

In all cases, we write $\mc{Z}(L)_{\ms{H}} \subseteq \mc{Z}(L)$ (\emph{horizontal special cycle}) for the flat part of $\mc{Z}(L)$, i.e. the largest closed formal subscheme which is flat over $\Spf \mc{O}_{\breve{F}}$.
    
            \subsection{Actions on Rapoport--Zink spaces} 
            \label{ssec:moduli_pDiv:action_RZ}
                Consider the groups
    \begin{align}
    & I_0 \coloneqq \{ \g_0 \in \End^0_F(\mbf{X}_0) : \g_0^{\dagger} \g_0 \in \Q_p^{\times} \} && I \coloneqq \{ \g \in \End^0_F(\mbf{X}) : \g^{\dagger} \g \in \Q_p^{\times} \} \\
    & I' \coloneqq \{ (\g_0, \g) \in I_0 \times I : \g_0^{\dagger} \g_0 = \g^{\dagger} \g \} && I_1 \coloneqq \{ \g \in \End^0_F(\mbf{X}) : \g^{\dagger} \g = 1 \}.
    \end{align}
We have $I_0 = F^\times$ (canonically). Using this identification, there is an isomorphism $I' \ra I_0 \times I_1$ given by $(\g_0, \g) \mapsto (\g_0, \g_0^{-1} \g)$. 
We have actions
    \begin{align}
    & I \acts \mc{N}(n-r,r)' && I_1 \acts \mc{N}(n-r,r) 
    \\
    & (X, \iota, \lambda, \rho) \mapsto (X, \iota, \lambda, \g \circ \rho) && (X, \iota, \lambda, \rho) \mapsto (X, \iota, \lambda, \g \circ \rho) \notag
    \end{align}
    \begin{align}
    & I' \acts \mc{N}'
     \\
    & (X_0, \iota_0, \lambda_0, \rho_0, X, \iota, \lambda, \rho) \mapsto (X_0, \iota_0, \lambda_0, \g_0 \circ \rho_0, X, \iota, \lambda, \g \circ \rho). \notag
    \end{align}
These actions are compatible with the isomorphisms $I' \cong I_0 \times I_1$ and $\mc{N}' \cong \mc{N}(1,0)' \times \mc{N}(n-r,r)$. We have isomorphisms
    \begin{equation}\label{equation:moduli_pDiv:action_RZ:iso_with_unitary_groups}
    \begin{tikzcd}[row sep = tiny, column sep = small]
    I_0 \arrow{r}{\sim} & GU(\mbf{V}_0)
    \end{tikzcd}
    \quad \quad
    \begin{tikzcd}[row sep = tiny, column sep = small]
    I_1 \arrow{r}{\sim} & U(\mbf{W}) \times U(\mbf{W}^{\perp})
    \end{tikzcd}
    \end{equation}
where $\g \in I_1$ acts on $\mbf{V}$ as $x \mapsto \gamma \circ x$, and similarly for $\mbf{V}_0$.

For any subset $L \subseteq \mbf{W}$ with associated local special cycles $\mc{Z}'(L) \subseteq \mc{N}'$ and $\tilde{\mc{Z}}(L) \subseteq \mc{N}(n-r,r)$, the actions of $I'$ and $I_1$ described above satisfy
    \begin{equation}\label{equation:moduli_pDiv:action_RZ:special_cycles}
    (\g_0, \g)(\mc{Z}(L)') = \mc{Z}(\g L \g_0^{-1})' \quad \g(\tilde{\mc{Z}}(L)) = \tilde{\mc{Z}}(\g L).
    \end{equation}
We will also have $\gamma(\mc{Z}(L)) = \mc{Z}(\gamma(L))$ (already checked in the nonsplit cases; in the split case, this will be clear from the definition, see \eqref{equation:moduli_pDiv:discrete_reduced:split_Z(L)}).
    
            \subsection{Discrete reduced subschemes}
            \label{ssec:moduli_pDiv:discrete_reduced}
                In the nonsplit cases (at least if $p \neq 2$), the reduced subscheme $\mc{N}(n - 1, 1)_{\red}$ of $\mc{N}(n - 1, 1)$ admits a stratification by Deligne--Lusztig varieties, described by a certain Bruhat--Tits building \cite{VW11,Wu16}. Later, we will use these results implicitly via citation to \cite{LZ22unitary,LL22II}. 

In this section, we further discuss some cases where the reduced subscheme is discrete (continuing to allow $p = 2$ if $F / \Q_p$ is unramified). 

In the split case, set
    \begin{equation}
    \mbf{L} \coloneqq \Hom_{\mc{O}_F}(\mbf{X}_0, \mbf{X}) \quad \quad \mbf{L}^{\perp} \coloneqq \Hom_{\mc{O}_F}(\mbf{X}_0^{\s}, \mbf{X}).
    \end{equation}
In the nonsplit case, define $\mbf{L}$ in the same way but take $\mbf{L}^{\perp} \coloneqq 0$. Let $K_{1, \mbf{L}} \subseteq U(\mbf{W})$ and $K_{1, \mbf{L}^{\perp}} \subseteq U(\mbf{W}^{\perp})$ be the respective stabilizers.

\begin{lemma}\label{lemma:moduli_pDiv:discrete_reduced}
Consider signature $(n - r, r) = (1,1)$ if $F / \Q_p$ is nonsplit (resp. any signature $(n - r, r)$ if $F / \Q_p$ is split). 
    \begin{enumerate}[(1)]
        \item The framing object $(\mbf{X}, \iota_{\mbf{X}}, \lambda_{\mbf{X}})$ is unique up to isomorphism. This also holds over any algebraically closed field $\kappa$ over $\overline{k}$, at least if $F / \Q_p$ is unramified.
        
        \item The reduced scheme $\mc{N}(n - r, r)_{\red}$ is discrete (i.e. a disjoint union of copies of $\Spec \overline{k}$). If $F / \Q_p$ is inert (resp. ramified), then $\mc{N}(1,1)_{\red}$ is one point (resp. two points).

        \item The lattices $\mbf{L} \subseteq \mbf{W}$ and $\mbf{L}^{\perp} \subseteq \mbf{W}^{\perp}$ are maximal integral lattices. In the nonsplit cases, $\mbf{L} \subseteq \mbf{W} = \mbf{V}$ is the unique maximal integral lattice.
        
        \item The group $I_1 \cong U(\mbf{W}) \times U(\mbf{W}^{\perp})$ acts transitively on $\mc{N}(n - r, r)(\overline{k})$. Consider the resulting surjection
        \begin{equation}\label{equation:moduli_pDiv:discrete_reduced:RZ_identification}
        \begin{tikzcd}[row sep = tiny]
        \mc{N}(n - r, r)(\overline{k}) \arrow{r} & U(\mbf{W}) / K_{1, \mbf{L}} \times U(\mbf{W}^{\perp}) / K_{1, \mbf{L}^{\perp}} \\
        (\mbf{X}, \iota_{\mbf{X}}, \lambda_{\mbf{X}}, (\gamma, \gamma^{\perp})) \arrow[mapsto]{r} & (\g, \g^{\perp}).
        \end{tikzcd}
        \end{equation}
        If $F / \Q_p$ is unramified, this map is a bijection. If $F / \Q_p$ is ramified, this map is $2$-to-$1$. If $F / \Q_p$ is nonsplit, the set $U(\mbf{W}) / K_{1, \mbf{L}} \times U(\mbf{W}^{\perp}) / K_{1, \mbf{L}^{\perp}}$ has size $1$.

        \item Consider the bijective identification
            \begin{equation*}
            \begin{tikzcd}[row sep = tiny]
            U(\mbf{W}) / K_1(\mbf{L}) \times U(\mbf{W}^{\perp}) / K_1(\mbf{L}^{\perp}) \arrow{r} & \arrow{l} \left \{ \begin{array}{l} \text{maximal full-rank integral $\mc{O}_F$-lattices $N \subseteq \mbf{V}$} \\ \text{where $N = M \oplus M^{\perp}$ with} \\ \text{$M \subseteq \mbf{W}$ and $M^{\perp} \subseteq \mbf{W}^{\perp}$} \end{array} \right \}
            \\
            (\gamma, \gamma^{\perp}) \arrow[mapsto]{r} & \gamma(\mbf{L}) \oplus \gamma^{\perp}(\mbf{L}^{\perp}).
            \end{tikzcd}
            \end{equation*}
        Given any subset $L \subseteq \mbf{W}$, the subset $\tilde{\mc{Z}}(L)(\overline{k}) \subseteq \mc{N}(n - r, r)(\overline{k})$ is identified (via \eqref{equation:moduli_pDiv:discrete_reduced:RZ_identification}) with the pre-image of the set of lattices $\{ N : L \subseteq N \}$.
    \end{enumerate}
\end{lemma}
\begin{proof}
\hfill
\begin{enumerate}[(1)]
\item In the inert case, this follows from Dieudonn\'e theory as in \cite[{Proposition 1.10}]{Vollaard10} (but we allow $p = 2$ by the same method), diagonalizability of Hermitian $\mc{O}_F$-lattices, and the following fact: consider the rank $2$ Hermitian $\mc{O}_F$-lattice $\Lambda$ with pairing $(-,-)$ specified by the Gram matrix
    \begin{equation}
    \begin{pmatrix}
    1 & 0 \\
    0 & p
    \end{pmatrix},
    \end{equation}
and also write $(-,-)$ for the induced pairing on $\Lambda \otimes_{\mc{O}_F} W(\kappa)[1/p]$ (which is ``sesquilinear'' for the Frobenius on $W(\kappa)[1/p]$).
If $x \in \Lambda \otimes_{\mc{O}_F} W(\kappa)[1/p]$ is any element with $(x,x) \in W(k)$, then $x \in \Lambda \otimes_{\mc{O}_F} W(\kappa)$ (so $\Lambda \otimes_{\mc{O}_F} W(\kappa)$ satisfies a certain ``unique maximal integral lattice'' property).
This computation shows that $\mc{N}(1,1)(\overline{k})$ is a single point if $F / \Q_p$ is inert.

The ramified case follows from \cite[{Lemma 6.1}]{RSZ17}. The split case was already verified in Lemma \ref{lemma:moduli_pDiv:RZ:unique_framing_object_split}.

\item In the unramified case, part (1) implies $\mc{N}(n - r, r)(\overline{k}) \ra \mc{N}(n - r, r)(\kappa)$ is bijective for every algebraically closed field extension $\kappa$ over $\overline{k}$.
Discreteness then follows because $\mc{N}(n - r, r)_{\red} \ra \Spec \overline{k}$ is locally of finite type. If $F / \Q_p$ is inert, $\mc{N}(1,1)(\overline{k})$ being a single point was already discussed above. The ramified case is \cite[{Lemma 6.1}]{RSZ17}. 

\item By part (1), we may assume $(\mbf{X}, \iota_{\mbf{X}}, \lambda_{\mbf{X}})$ is the explicit tuple constructed in Section \ref{ssec:moduli_pDiv:RZ}. The claim can then be verified explicitly, using \eqref{equation:moduli_pDiv:RZ:hom_to_sigma} and the surrounding discussion. In the split case, $\mbf{L}$ and $\mbf{L}^{\perp}$ will be self-dual. In the inert case, $\mbf{L}$ admits a Gram matrix with basis $\operatorname{diag}(1,p)$. In the ramified case, $\mbf{L}$ admits a Gram matrix with basis $\mrm{diag}(1, -a)$ for some $a \in \Z_p^{\times}$ which is a non-norm, i.e. $a \not \in N_{F / \Q_p}(F^{\times})$. This claim in the ramified case follows from the observation that $\mbf{L}$ is integral and that $\mbf{L} \subseteq \Hom_{\mc{O}_F}(\mbf{X}_0, \mbf{X}_0 \times \mbf{X}_0^{\s}) \subseteq \varpi^{-1} \mbf{L}$ (as the isogeny in \eqref{equation:moduli_pDiv:RZ:isogeny_from_tensor} has kernel contained in the $\varpi$-torsion subgroup).

\item Transitivity of the $I_1$ action is immediate from part (1). Note also $\End(\mbf{X}) \cap (U(\mbf{W}) \times U(\mbf{W}^{\perp})) \subseteq K_{1, \mbf{L}} \times K_{1,\mbf{L}^{\perp}}$ so the displayed map is well-defined. In the nonsplit cases, the assertions follow from parts (2) and (3). Bijectivity in the split case follows because we then have $\End(\mbf{X}) \cap (U(\mbf{W}) \times U(\mbf{W}^{\perp})) = K_{1, \mbf{L}} \times K_{1,\mbf{L}^{\perp}}$.

\item Follows from the previous parts, i.e. $\tilde{\mc{Z}}(L)(\overline{k})$ corresponds to $(\gamma, \gamma^{\perp})$ such that $L \subseteq \gamma(\mbf{L})$. \qedhere
\end{enumerate}
\end{proof}

Suppose $F / \Q_p$ is split and $L \subseteq \mbf{W}$ is any subset (with arbitrary signature $(n - r, r)$). We take $\mc{Z}(L) \subseteq \tilde{\mc{Z}}(L)$
to be the open and closed subfunctor corresponding (via Lemma \ref{lemma:moduli_pDiv:discrete_reduced}(4)) to the locus where $\gamma^{\perp}(\mbf{L}^{\perp}) = \mbf{L}^{\perp}$. By the previous discussion, there is a isomorphism of formal schemes
    \begin{equation}\label{equation:moduli_pDiv:discrete_reduced:split_Z(L)}
    \begin{tikzcd}[row sep = tiny]
    \tilde{\mc{Z}}(L) \arrow{r} & \mc{Z}(L) \times U(\mbf{V}^{\perp})/K_{1,\mbf{L}^{\perp}}
    \\
    (\mbf{X}, \iota_{\mbf{X}}, \lambda_{\mbf{X}}, (\gamma, \gamma^{\perp}))
    \arrow[mapsto]{r} & ((\mbf{X}, \iota_{\mbf{X}}, \lambda_{\mbf{X}}, (\gamma, 1)),\gamma^{\perp}).
    \end{tikzcd}
    \end{equation}
In this case, we have a canonical bijection
    \begin{equation}\label{equation:moduli_pDiv:discrete_reduced:split_Z(L)_points}
    \mc{Z}(L)(\overline{k}) = \{ M \subseteq \mbf{W} : \text{full rank self-dual lattice with $L \subseteq M$} \}
    \end{equation}
via Lemma \ref{lemma:moduli_pDiv:discrete_reduced}.

\begin{lemma}\label{lemma:moduli_pDiv:discrete_reduced:local_special_cycle_red_finite}
Suppose $F / \Q_p$ is split. If $L \subseteq \mbf{W}$ is an $\mc{O}_F$-lattice of full rank (i.e rank $n-r$), then $\mc{Z}(L)(\overline{k})$ is a finite set.
\end{lemma}
\begin{proof}
Our task is to show that the right-hand side of \eqref{equation:moduli_pDiv:discrete_reduced:split_Z(L)_points} is finite. For such $M$, we must have $L \subseteq M \subseteq M^{\vee} \subseteq L^{\vee}$ where $L^{\vee}$ and $M^{\vee}$ denote the dual lattices. If $L \not \subseteq L^{\vee}$ then $\mc{Z}(L)$ is empty. Otherwise, $L^{\vee}/L$ is an $\mc{O}_F$-module of finite length, so there are only finitely many possibilities for $M$.
\end{proof}

            \subsection{Horizontal and vertical decomposition}
            \label{ssec:moduli_pDiv:horizontal_vertical}
                For a locally Noetherian formal scheme $\mc{X}$, viewed as a ringed space with structure sheaf $\mc{O}_{\mc{X}}$,  we write
    \begin{equation}
    K'_0(\mc{X}) \coloneqq K_0(\textit{Coh}(\mc{O}_{\mc{X}})) \quad \quad F_d K'_0(\mc{X}) \subseteq K'_0(\mc{X})
    \end{equation}
for the $K_0$ group of coherent $\mc{O}_{\mc{X}}$-modules and the subgroup generated by coherent sheaves supported in (formal scheme-theoretic) dimension $\leq d$, respectively.
If $\mc{X}$ is moreover formally locally of finite type over $\Spf R$ for a complete discrete valuation ring $R$, we say that $\mc{X}$ is \emph{equidimensional of dimension $n$}
if every open formal subscheme of $\mc{X}$ has dimension $n$.
In this case, if $\mc{Z} \ra \mc{X}$ is an adic finite morphism of locally Noetherian formal schemes, we write
    \begin{equation}
    F^m_{\mc{X}} K'_0(\mc{Z}) \coloneqq F_{n-m} K'_0 (\mc{Z}) \quad \quad \mrm{gr}^m_{\mc{X}} K'_0(\mc{Z}) \coloneqq F^m_{\mc{X}} K'_0(\mc{Z}) / F^{m+1}_{\mc{X}} K'_0(\mc{Z}).
    \end{equation}
We often work with these groups tensor $\Q$, written as $K'_0(\mc{X})_{\Q}$, etc.

The discussion in Section \ref{ssec:moduli_pDiv:RZ} implies that the Rapoport--Zink space $\mc{N}(n-r,r)$ is equidimensional of dimension $(n-r)r+1$. 
For the rest of Section \ref{ssec:moduli_pDiv:horizontal_vertical} we fix signature $(n - 1, 1)$ and use the shorthand $\mc{N} \coloneqq \mc{N}(n - 1, 1)$. The material below is a local analogue of \cref{ssec:part_I:arith_intersections:vertical_classes}.

Assume $F / \Q_p$ is nonsplit for the moment. For any nonzero $x \in \mbf{W}$, the local special cycle $\mc{Z}(x)$ is a Cartier divisor on $\mc{N}$ for any nonzero $x \in \mbf{W}$ (\cite[{Proposition 3.5}]{KR11} (inert) \cite[{Proposition 4.3}]{Howard19} (inert allowing $p = 2$), and also \cite[{Lemma 2.40}]{LL22II} (ramified exotic smooth)). 
For any $x \in \mbf{W}$, set
    \begin{equation}
    {}^{\mbb{L}} \mc{Z}(x) \coloneq
    \begin{cases}
    \mc{O}_{\mc{Z}(x)} & \text{if $x \neq 0$} \\
    (\cdots 0 \ra \mc{E} \xra{0} \mc{O}_{\mc{N}} \ra 0 \cdots) & \text{if $x = 0$}
    \end{cases}
    \end{equation}
in $D^b_{\textit{Coh}(\mc{O}_{\mc{Z}(x)})}(\mc{O}_{\mc{N}})$ (bounded derived category of $\mc{O}_{\mc{N}}$-modules with cohomology sheaves coherent and supported along $\mc{O}_{\mc{Z}(x)}$), where the $\mc{O}_{\mc{N}}$ term is in degree $0$.
For any tuple $\underline{x} \in \mbf{W}^m$, we then consider the \emph{derived local special cycle} 
    \begin{equation}
    {}^\mbb{L} \mc{Z}(\underline{x}) \coloneqq {}^{\mbb{L}} \mc{Z}(x_1) \otimes^{\mbb{L}} \cdots \otimes^{\mbb{L}} {}^{\mbb{L}} \mc{Z}(x_m) \in D^b_{\textit{Coh}(\mc{O}_{\mc{Z}(\underline{x})})}(\mc{O}_{\mc{N}})
    \end{equation}
Its image ${}^{\mbb{L}} \mc{Z}(\underline{x}) \in K'_0(\mc{Z}(\underline{x}))_{\Q}$ lies in $F^m_{\mc{N}} K'_0(\mc{Z}(\underline{x}))_{\Q}$ by multiplicativity of the codimension filtration\footnote{There is a technicality here, as $\mc{N}$ is a formal scheme rather than a scheme. So we instead prove filtration multiplicativity via uniformization (\crefext{III:corollary:non-Arch_uniformization:vertical:local_cycle_filtration}) by reducing to the analogous filtration multiplicativity statement for global special cycles. We will make a few other forward references to \crefext{III:ssec:non-Arch_uniformization:global_to_local} where we verify some properties of local special cycles via uniformization.} and depends only on $\mrm{span}_{\mc{O}_F}(\underline{x})$ (``linear invariance'') by \cite[{Theorem B}]{Howard19} (inert) and \cite[{Proposition 2.33}]{LL22II} (ramified exotic smooth).

Continuing to assume $F / \Q_p$ is nonsplit, assume $\underline{x} \in \mbf{W}^m$ spans a non-degenerate Hermitian $\mc{O}_F$-lattice of rank $m^{\flat}$. We define certain \emph{derived vertical local special cycles} ${}^{\mbb{L}} \mc{Z}(\underline{x})_{\ms{V}} \in \mrm{gr}^m_{\mc{N}} K'_0(\mc{Z}(\underline{x})_{\overline{k}})_{\Q}$ as follows.

For integers $e \gg 0$, we have a scheme-theoretic union decomposition (\crefext{III:lemma:non-Arch_uniformization:global_to_local:scheme_union_decomp})
    \begin{equation}
    \mc{Z}(\underline{x}) = \mc{Z}(\underline{x})_{\ms{H}} \cup \mc{Z}(\underline{x})_{\ms{V}}
    \end{equation}
where $\mc{Z}(\underline{x})_{\ms{H}}$ is the flat part of $\mc{Z}(\underline{x})$, i.e. the largest closed formal subscheme which is flat over $\Spf \mc{O}_{\breve{F}}$, and $\mc{Z}(\underline{x})_{\ms{V}} \coloneqq \mc{Z}(\underline{x})_{\Spf \mc{O}_{\breve{F}} / p^e}$ for $e \gg 0$. Since $\mc{Z}(\underline{x})_{\ms{H}}$ is equidimensional of dimension $n - m^{\flat}$ (\crefext{III:lemma:non-Arch_uniformization:global_to_local:horizontal_description}), and since $\mc{Z}(\underline{x})_{\ms{H}} \cap \mc{Z}(\underline{x})_{\ms{V}}$ has dimension $\leq n - m^{\flat} - 1$, 
there is an induced decomposition
    \begin{equation}
    \mrm{gr}^{m^{\flat}}_{\mc{N}} K'_0(\mc{Z}(\underline{x})) = \mrm{gr}^{m^{\flat}}_{\mc{N}} K'_0(\mc{Z}(\underline{x})_{\ms{H}}) \oplus \mrm{gr}^{m^{\flat}}_{\mc{N}} K'_0(\mc{Z}(\underline{x})_{\overline{k}})
    \end{equation}
independent of $e$ (cf. \cite[{Lemma B.1}]{Zhang21}\footnote{Strictly speaking, our setup for $K'_0$ groups may be different from Zhang's in non quasi-compact settings. The proof of the cited lemma is the same in our setup.}). Here we have used the pushforward d\'evissage isomorphism $K'_0(\mc{Z}(\underline{x})_{\overline{k}}) \xra{\sim} K'_0(\mc{Z}(\underline{x})_{\ms{V}})$ for  for $K'_0$ groups.

If $m = m^{\flat}$, we define ${}^{\mbb{L}} \mc{Z}(\underline{x})_{\ms{V}}$ to be given by the projection
    \begin{equation}\label{equation:moduli_pDiv:horizontal_vertical:derived_vertical}
    \begin{tikzcd}[row sep = tiny]
    \mrm{gr}^{m}_{\mc{N}} K'_0(\mc{Z}(\underline{x}))_{\Q} \arrow{r} & \mrm{gr}^{m}_{\mc{N}} K'_0(\mc{Z}(\underline{x})_{\overline{k}})_{\Q} \\
    {}^{\mbb{L}} \mc{Z}(\underline{x}) \arrow[mapsto]{r} & {}^{\mbb{L}} \mc{Z}(\underline{x})_{\ms{V}}.
    \end{tikzcd}
    \end{equation}
By the linear invariance property for ${}^{\mbb{L}} \mc{Z}(\underline{x})$ discussed above, the class ${}^{\mbb{L}} \mc{Z}(\underline{x})_{\ms{V}}$ depends only on $\mrm{span}_{\mc{O}_F}(\underline{x})$.
    
For possibly $m \neq m^{\flat}$, we say that $\underline{x} = [x_1, \ldots, x_m]$ is in \emph{minimal form} if $\underline{x}^{\flat} \coloneqq [x_{m - m^{\flat} + 1}, \ldots, x_m]$ satisfies $\mrm{span}_{\mc{O}_F}(\underline{x}^{\flat}) = \mrm{span}_{\mc{O}_F}(\underline{x})$. In this case, set $\underline{x}^{\#} \coloneqq [x_1, \ldots, x_{m - m^{\flat}}]$ and define\footnote{One needs to show that the map $\a \mapsto {}^{\mbb{L}} \mc{Z}(\underline{x}^{\#}) \cdot \a$ sends $F^{m^{\flat} + 1}_{\mc{N}} K'_0(\mc{Z}(\underline{x})_{\overline{k}})_{\Q} \ra F^{m + 1}_{\mc{N}} K'_0(\mc{Z}(\underline{x})_{\overline{k}})_{\Q}$. 
This is clear if $m^{\flat} \geq n - 1$, 
but we do not know a proof of this in general as $\mc{N}$ is a formal scheme and not a scheme. Since we are mostly interested in the case $m^{\flat} \geq n - 1$, we do not pursue this point further. Even when $m^{\flat} = n - 1$ and $m = n$, one still needs to check that ${}^{\mbb{L}} \mc{Z}(\underline{x})_{\ms{V}}$ lies in $F^n_{\mc{N}}(\mc{Z}(\underline{x})_{\overline{k}})_{\Q}$ (rather than $F^{n-1}_{\mc{N}}(\mc{Z}(\underline{x})_{\overline{k}})_{\Q}$). This follows e.g. because $\mc{Z}(\underline{x})_{\overline{k}}$ is a Noetherian scheme (\crefext{III:lemma:non-Arch_uniformization:vertical:local_cycle_quasi-compact}) whose reduced irreducible components are projective over $\overline{k}$.
The definition of ${}^{\mbb{L}} \mc{Z}(\underline{x})_{\ms{V}}$ should thus be treated as conditional unless $m = m^{\flat}$ or $m^{\flat} \geq n - 1$.\label{footnote:moduli_pDiv:horizontal_vertical:codimension_multiplicativity}
}
    \begin{equation}
    {}^{\mbb{L}} \mc{Z}(\underline{x})_{\ms{V}} \coloneqq {}^{\mbb{L}} \mc{Z}(\underline{x}^{\#}) \cdot {}^{\mbb{L}} \mc{Z}(\underline{x}^{\flat})_{\ms{V}} \in \mrm{gr}_{\mc{N}}^m K'_0(\mc{Z}(\underline{x})_{\overline{k}})_{\Q}.
    \end{equation}
For $\underline{x}$ possibly not in minimal form, select any $\gamma \in \GL_m(\mc{O}_F)$ such that $\underline{x} \cdot \gamma$ is in minimal form, and set ${}^{\mbb{L}} \mc{Z}(\underline{x})_{\ms{V}} \coloneqq {}^{\mbb{L}} \mc{Z}(\underline{x} \cdot \gamma)_{\ms{V}}$ (note $\mc{Z}(\underline{x}) = \mc{Z}(\underline{x} \cdot \gamma)$).

We claim that ${}^{\mbb{L}} \mc{Z}(\underline{x})_{\ms{V}}$ depends only on $m$ and $\mrm{span}_{\mc{O}_F}(\underline{x})$, and not on the choice of $\underline{x}$ or a minimal form (``linear invariance''). For $\underline{x}$ in minimal form and with notation as above, we already explained that ${}^{\mbb{L}} \mc{Z}(\underline{x}^{\flat})_{\ms{V}}$ depends only on $\mrm{span}_{\mc{O}_F}(\underline{x})$. Recall $\mc{Z}(\underline{x}) = \mc{Z}(\underline{x}^{\flat})$. Consider any element $x_i$ of the tuple $\underline{x}^{\#}$. Then $x_i \in \mrm{span}_{\mc{O}_F}(\underline{x}^{\flat}) = \mrm{span}_{\mc{O}_F}(\underline{x})$. In particular, $\mc{Z}(\underline{x}^{\flat}) \subseteq \mc{Z}(x_i)$. But Grothendieck--Messing theory provides a canonical isomorphism
    \begin{equation}
    \mc{E}|_{\mc{Z}(x_i)} \xra{\sim} \mc{I}(x_i) / \mc{I}(x_i)^2
    \end{equation}
if $x_i \neq 0$, where $\mc{I}(x_i) \subseteq \mc{O}_{\mc{N}}$ is the ideal sheaf of the Cartier divisor $\mc{Z}(x_i) \subseteq \mc{N}$ (follows from \cite[{Definition 4.2}]{Howard19} (inert) and \cite[{Lemma 2.39}]{LL22II}). Hence we have
    \begin{equation}
    {}^{\mbb{L}} \mc{Z}(\underline{x}^{\#})|_{\mc{Z}(\underline{x}^{\flat})} \cong {}^{\mbb{L}} \mc{Z}(\underline{0}_{m - m^{\flat}})|_{\mc{Z}(\underline{x}^{\flat})}  
    \end{equation}
as elements of $D^b_{\textit{Coh}}(\mc{O}_{\mc{Z}(\underline{x}^{\flat})})$, where $\underline{0}_{m - m^{\flat}} \in \mbf{W}^{m - m^{\flat}}$ is the tuple with all entries equal to $0$.
Then we have
    \begin{equation}
    {}^{\mbb{L}} \mc{Z}(\underline{x})_{\ms{V}} = {}^{\mbb{L}} \mc{Z}(\underline{0}_{m - m^{\flat}}) \cdot {}^{\mbb{L}} \mc{Z}(\underline{x}^{\flat})_{\ms{V}} \in \mrm{gr}^m_{\mc{N}} K'_0(\mc{Z}(\underline{x})_{\overline{k}})_{\Q}.
    \end{equation}
We have explained that the right-hand side does not depend on any auxiliary choices.

Next, suppose $F / \Q_p$ is split, and assume $\underline{x} \in \mbf{W}^m$ has $\mc{O}_F$-span which is a lattice of rank $n - 1$ (full rank). 
We have $\mrm{gr}^{n- 1}_{\mc{N}} K'_0(\mc{Z}(\underline{x})_{\overline{k}}) = 0$ for dimension reasons (the reduced subscheme of $\mc{N}$ is dimension $0$, see Section \ref{ssec:moduli_pDiv:discrete_reduced} and \crefext{III:lemma:non-Arch_uniformization:vertical:local_cycle_quasi-compact}). Constructing ${}^{\mbb{L}} \mc{Z}(\underline{x})_{\ms{V}} \in \mrm{gr}^m_{\mc{N}} K'_0(\mc{Z}(\underline{x})_{\overline{k}})$ as above gives the \emph{derived vertical local special cycle} ${}^{\mbb{L}} \mc{Z}(\underline{x})_{\ms{V}} = 0$.

Next, consider $\underline{x} \in \mbf{W}^m$ which is a basis for its $\mc{O}_F$-span $L^{\flat} \coloneqq \mrm{span}_{\mc{O}_F}(\underline{x})$. If $F / \Q_p$ is split, we also assume $m = n - 1$. In this situation, we set ${}^{\mbb{L}} \mc{Z}(L^{\flat})_{\ms{V}} \coloneqq {}^{\mbb{L}} \mc{Z}(\underline{x})_{\ms{V}}$, since the latter depends only on $L^{\flat}$. If $n = 2$ and $m = 1$, we have ${}^{\mbb{L}} \mc{Z}(\underline{x})_{\ms{V}} = 0$ since the reduced subscheme $\mc{N}_{\red}$ has dimension $0$ (Section \ref{ssec:moduli_pDiv:discrete_reduced}) and since $\mc{N}$ has dimension $2$ in this case.
    
            \subsection{Serre tensor and signature \texorpdfstring{$(1,1)$}{(1,1)}}
            \label{ssec:moduli_pDiv:Serre_tensor}
                The case of signature $(1,1)$ plays an important role for describing local special cycles via the Serre tensor construction.

As above, let $\mbf{X}_0$ be the unique supersingular (resp. ordinary) $p$-divisible group over $\overline{k}$ of height $2$ and dimension $1$ if $F / \Q_p$ is nonsplit (resp. split). 
For schemes $S$ over $\Spf \mc{O}_{\breve{F}}$, we consider pairs $(X, \rho)$ where $X$ is a $p$-divisible group over $S$ and $\rho \colon X_{\overline{S}} \ra \mbf{X}_{0,\overline{S}}$ is any quasi-isogeny.

We form the Rapoport--Zink space $\widetilde{\mc{N}}_{2,1}$ over $\Spf \mc{O}_{\breve{F}}$, given by
    \begin{equation}
    \widetilde{\mc{N}}_{2,1}(S) \coloneqq \{ \text{isomorphism classes of framed tuples } (X, \rho) \text{ over $S$} \}.
    \end{equation}
This is a locally Noetherian formal scheme which is formally locally of finite type over $\Spf \mc{O}_{\breve{F}}$ (via the by-now standard representability result \cite[{Theorem 2.16}]{RZ96}). There is an isomorphism of formal schemes $\operatorname{Isog}^0(\mbf{X}_0) \ra \widetilde{\mc{N}}_{2,1}$ given by $\rho \mapsto (\mbf{X}_0, \rho)$, where $\operatorname{Isog}^0(\mbf{X}_0)$ is viewed as a constant formal scheme. Indeed, this follows as in the proof of Lemma \ref{lemma:moduli_pDiv:discrete_reduced} by uniqueness of $\mbf{X}_0$ over any algebraically closed field $\kappa$ (any quasi-endomorphism of $\mbf{X}_0$ descends to $\overline{k}$ as well, as may be checked on isocrystals).

We let $\mc{N}_{2,1} \subseteq \widetilde{\mc{N}}_{2,1}$ be the open and closed locus where the framing $\rho$ is fiberwise an isomorphism.
Then $\mc{N}_{2,1}$ is representable by a formal scheme, and there is a (non-canonical) isomorphism $\mc{N}_{2,1} \cong \Spf \psring{\mc{O}_{\breve{F}}}{t}$ (e.g. by Grothendieck--Messing theory). 

For arbitrary signature $(n -r, r)$ in the split case, we can form $\widetilde{\mc{N}}_{n,r}$ and $\mc{N}_{n,r}$ as above, where we replace $\mbf{X}_0$ with the unique ordinary $p$-divisible group of height $n$ and dimension $r$. The previous asssertions for $\widetilde{\mc{N}}_{2,1}$ and $\mc{N}_{2,1}$ hold in this case as well, except we now have $\mc{N}_{n,r} \cong \Spf \psring{\mc{O}_{\breve{F}}}{t_1, \ldots, t_{(n-r)r}}$ (again by Grothendieck--Messing theory).

\begin{lemma}\label{lemma:moduli_pDiv:Serre_tensor:lift_polarization}
Given any $(X, \rho) \in \mc{N}_{2,1}(S)$, any principal polarization $\lambda_{\mbf{X}_0}$ of $\mbf{X}_0$ lifts uniquely to a principal polarization on $X$.
\end{lemma}
\begin{proof}
Uniqueness follows from Drinfeld rigidity. Any two principal polarizations on $X$ differ by $\Z_p^{\times}$ scalar (since this holds for $\mbf{X}_0$), so it is enough to show existence of a principal polarization on $X$. Since $\mc{N}_{2,1} \cong \Spf \psring{\mc{O}_{\breve{F}}}{t}$, it is enough to check the case where the scheme $S$ is a finite order thickening of $\Spec \overline{k}$. By Serre--Tate, we can view $X$ as the $p$-divisible group of an elliptic curve over $S$ (deforming an elliptic curve over $\Spec \overline{k}$ with $p$-divisible group $\mbf{X}_0$). Any elliptic curve admits a (unique) principal polarization.
\end{proof}

The preceding (possibly standard) argument also appeared in the proof of \cite[{Proposition 6.3}]{RSZ17} (for the same purpose), there in the supersingular case.

Recall the triple $(\mbf{X}_0 \otimes_{\Z_p} \mc{O}_F, \iota, \lambda_{\mbf{X}_0} \otimes \lambda_{\mrm{tr}})$ described in Section \ref{ssec:moduli_pDiv:RZ}, arising from the Serre tensor construction (fixing some choice of $\lambda_{\mbf{X}_0}$).
For any $(X, \rho) \in \mc{N}_{2,1}(S)$, the same construction gives a tuple $(X \otimes_{\Z_p} \mc{O}_F, \iota, \lambda_{\mbf{X}_0} \otimes \lambda_{\mrm{tr}}, \rho \otimes_{\Z_p} \mc{O}_F)$ where $\lambda_{\mbf{X}_0}$ denotes the unique lift to $X$ as in Lemma \ref{lemma:moduli_pDiv:Serre_tensor:lift_polarization} (by abuse of notation), and where $\rho \otimes_{\Z_p} \mc{O}_F \colon X_{\overline{S}} \otimes_{\Z_p} \mc{O}_F \ra \mbf{X}_{0, \overline{S}} \otimes_{\Z_p} \mc{O}_F$.

\begin{lemma}[Serre tensor isomorphism]\label{lemma:moduli_pDiv:Serre_tensor}
For any $F$-linear quasi-isogeny $\phi \colon \mbf{X}_0 \otimes_{\Z_p} \mc{O}_F \ra \mbf{X}$ preserving polarizations exactly, the induced map
    \begin{equation}\label{equation:moduli_pDiv:Serre_tensor:1_1_signature}
    \begin{tikzcd}[row sep = tiny]
    \mc{N}_{2,1} \arrow{r} &  \mc{N}(1,1) \\
    (X, \rho) \arrow[mapsto]{r} & (X \otimes_{\Z_p} \mc{O}_F, \iota, -\mf{d}^2 \cdot (\lambda_{\mbf{X}_0} \otimes \lambda_{\mrm{tr}}), \phi_{\overline{S}} \circ (\rho \otimes_{\Z_p} \mc{O}_F))        
    \end{tikzcd}
    \end{equation}
(defined on $S$-points for schemes $S$ over $\Spf \mc{O}_{\breve{F}}$) is an open and closed immersion whose set-theoretic image is a single point.

If $F / \Q_p$ is split, the inverse is given by restricting $(X, \iota, \lambda, \rho) \mapsto (X^-, (\phi^-_{\overline{S}})^{-1} \circ \rho^-)$ to the appropriate component of $\mc{N}(1,1)$.
\end{lemma}
\begin{proof}
For the ramified case, we refer to \cite[{Proposition 6.3}]{RSZ17}. In the unramified case, the lemma follows by identifying the deformation theory of $\mc{N}_{2,1}$ and $\mc{N}(1,1)$ using Grothendieck--Messing theory, using the eigenspace decomposition for the $\mc{O}_F$-action on the Dieudonn\'e crystals of objects in $\mc{N}(1,1)(S)$ as in \crefext{III:equation:0_eigenspace_deformation_Hodge_filtration} and surrounding discussion (so the deformation problem for $\mc{N}_{2,1}$ identifies with the deformation problem of the ``$-$ eigenspace'' of the Hodge filtration for objects in $\mc{N}(1,1)(S)$ in the notation of loc. cit..).
This is essentially how we verified generic formal smoothness of special cycles in loc. cit..
\end{proof}

When $F / \Q_p$ is split and the signature $(n - r, r)$ is arbitrary, recall that any $(X, \iota, \lambda, \rho) \in \mc{N}(n - r, r)$ admits a decomposition $X = X^+ \times X^-$ and $\rho = \rho^+ \times \rho^-$ where $\rho^{\pm} \colon X^{\pm} \ra \mbf{X}^{\pm}$ using the nontrivial idempotents $e^{\pm} \in \mc{O}_F$.

\begin{lemma}\label{lemma:moduli_pDiv:Serre_tensor:eigenspace_equivalence}
Suppose $F / \Q_p$ is split, and consider arbitrary signature $(n - r, r)$. For any formal scheme $S$ over $\Spf \mc{O}_{\breve{F}}$, the forgetful functor
    \begin{equation}\label{equation:moduli_pDiv:Serre_tensor:eigenspace_equivalence}
    \begin{tikzcd}[row sep = tiny]
    \left \{ \begin{array}{l} \text{groupoid of principally polarized} \\ \text{Hermitian $p$-divisible groups $(X, \iota, \lambda)$} \\ \text{over $S$ of signature $(n-r, r)$} \end{array} \right \}
    \arrow{r} &
    \left \{ \begin{array}{l} \text{groupoid of ordinary $p$-divisible groups} \\ \text{over $S$ of height $n$ and dimension $r$} \end{array} \right \}
    \\
    (X, \iota, \lambda) \arrow[mapsto]{r} & X^-
    \end{tikzcd}
    \end{equation}
is an equivalence of categories. The same holds if we consider the groupoids with morphisms being quasi-isogenies (rather than isomorphisms).
\end{lemma}
\begin{proof}
An explicit quasi-inverse is given by $(X^-) \mapsto (X, \iota, \lambda)$ (over a scheme $S$) with
    \begin{align}\label{equation:split:minus_space_iso}
    & X = (X^-)^{\vee} \times X^- 
    \\
    & \iota(e^+) \colon X \ra (X^-)^{\vee} \quad \quad \iota(e^-) \colon X \ra (X^+)^{\vee} \quad \quad \text{projections} \notag
    \\
    & \lambda = \begin{pmatrix} 0 & 1 \\ -1 & 0 \end{pmatrix} \colon (X^-)^{\vee} \times X^- \ra X^- \times (X^-)^{\vee}. \notag
    \end{align}
This is analogous to the following phenomenon: if $L$ is a free $\mc{O}_F$-module of rank $n$ equipped with a perfect Hermitian pairing, then $U(L) \cong \GL_n(L^-)$ (and similarly with $F$ instead of $\mc{O}_F$).
\end{proof}

\begin{remark}\label{remark:moduli_pDiv:Serre_tensor:split_and_lifts}
Suppose $F / \Q_p$ is split, and suppose $R$ is a complete Noetherian local ring with algebraically closed residue field $\kappa$. If $n \geq 2$ and if $\underline{\Q_p / \Z_p}^{n - 2} \times X^-$ is an ordinary $p$-divisible group of height $n$ and dimension $r$, then $(\mf{X}_0)^{n - 2} \times (X^- \otimes_{\Z_p} \mc{O}_F)$ (with the product $\mc{O}_F$-action and product polarization $(\lambda_{\mf{X}_0})^{n - 2} \times (\lambda_{\mbf{X}_0} \otimes \lambda_{\mrm{tr}})$, for some choice of isomorphism $X^-_{\kappa} \cong \mbf{X}_0$) is a pre-image under the equivalence in \eqref{equation:moduli_pDiv:Serre_tensor:eigenspace_equivalence}. Here $(\mf{X}_0, \iota_{\mf{X}_0}, \lambda_{\mf{X}_0})$ is the canonical lift (over $S$) as in Definition \ref{definition:moduli_pDiv:RZ:canonical_lifting} (forgetting the framing).
\end{remark}

\begin{remark}
If we drop the ordinary hypothesis on both sides of Lemma \ref{lemma:moduli_pDiv:Serre_tensor:eigenspace_equivalence}, the lemma still holds (by the same proof).
\end{remark}

\begin{lemma}\label{lemma:moduli_pdiv:Serre_tensor:split_forget}
Suppose $F / \Q_p$ is split, and form $\widetilde{\mc{N}}_{n,r}$ using the framing object $\mbf{X}^-$. For arbitrary signature $(n - r, r)$, the forgetful map
    \begin{equation}\label{equation:moduli_pdiv:Serre_tensor:split_forget}
    \begin{tikzcd}[row sep = tiny]
    \mc{N}(n - r, r) \arrow{r} & \widetilde{\mc{N}}_{n,r} \\
    (X, \iota, \lambda, \rho) \arrow[mapsto]{r} & (X^-, \rho^-)
    \end{tikzcd}
    \end{equation}
is an isomorphism.
\end{lemma}
\begin{proof}
This is immediate from Lemma \ref{lemma:moduli_pDiv:Serre_tensor:eigenspace_equivalence}. Alternative (less elementary) proof: 
first observe that the forgetful map is an isomorphism on $\kappa$-points for any algebraically closed field $\kappa$ over $\overline{k}$ (see Lemma \ref{lemma:moduli_pDiv:discrete_reduced} and above discussion). As in the proof of Lemma \ref{lemma:moduli_pDiv:Serre_tensor}, the claim now follows from Grothendieck--Messing theory. 
\end{proof}

\begin{remark}\label{remark:moduli_pDiv:Serre_tensor:singleton}
In the situation of Lemma \ref{lemma:moduli_pdiv:Serre_tensor:split_forget}, the open and closed subfunctor $\mc{N}_{n,r} \subseteq \widetilde{\mc{N}}_{n,r} \cong \mc{N}(n - r, r)$ has $\mc{N}_{n,r}(\overline{k})$ being a singleton set, corresponding (via Lemma \ref{lemma:moduli_pDiv:discrete_reduced}) to the lattice $\mbf{L} \oplus \mbf{L}^{\perp} \subseteq \mbf{W} \oplus \mbf{W}^{\perp}$ (i.e. the locus where the framing $\rho$ is a fiberwise isomorphism).
\end{remark}

For $F / \Q_p$ in all cases (inert, ramified, split) and for any $x \in \mbf{W}$, the local special cycle $\mc{Z}(x) \ra \mc{N}(1,1)$ pulls back along the Serre tensor isomorphism (Lemma \ref{lemma:moduli_pDiv:Serre_tensor}) to a certain local special cycle on $\mc{N}_{2,1}$ associated with an element $x' \in \Hom^0(\mbf{X}_0, \mbf{X}_0)$ (arising from adjunction in the Serre tensor construction). The ramified case is explained in \cite[{\S 6.2}]{RSZ17}. The inert and split cases may be formulated in a similar way (we omit a more detailed statement, which we will not need).
This may be viewed as a local version of \cite[{Proposition 14.5}]{KR14} (see also \crefext{IV:ssec:arithmetic_Siegel-Weil:Serre_tensor}).

For $F / \Q_p$ nonsplit (at least if $p \neq 2$), Kudla--Rapoport \cite[{Proposition 8.1}]{KR11} and Rapoport--Smithling--Zhang \cite[{Proposition 7.1}]{RSZ17} use this to describe $\mc{Z}(x)$ in terms of certain \emph{quasi-canonical lifting cycles} on $\mc{N}_{2,1}$, corresponding to closed immersions $\Spf \mc{O}_{\breve{E}_s} \ra \mc{N}_{2,1}$ associated with $(\mf{X}_s, \rho) \in \Spf \mc{O}_{\breve{E}_s}$ where $(\mf{X}_s, \rho)$ arises from a \emph{quasi-canonical lifting} of $(\mbf{X}_0, j)$ for suitable $j \colon \mc{O}_F \ra \End(\mbf{X}_0)$ (in the notation and sense of Section \ref{ssec:can_and_qcan:qcan} below). This was extended by Li--Zhang \cite{LZ22unitary} (inert) and Li--Liu \cite{LL22II} (ramified) to flat parts of $1$-cycles in signature $(n - 1, 1)$, for arbitrary $n$ in the inert case and even $n$ in the ramified case. We will need this result, which we recall in Section \ref{ssec:can_and_qcan:qcan_cycles} below (to the precision we need).

We will need an analogue of the previous paragraph when $F / \Q_p$ is split (allowing $p = 2$). This is accomplished in Section \ref{sec:more_moduli_pDiv_split} below (statement given in Section \ref{ssec:can_and_qcan:qcan_cycles}). Our method in the split case is somewhere different from the proofs cited above.

        \section{More on moduli of \texorpdfstring{$p$}{p}-divisible groups: split}
        \label{sec:more_moduli_pDiv_split}
            Retain notation from Section \ref{sec:moduli_pDiv}. Throughout Section \ref{sec:more_moduli_pDiv_split}, we assume $F / \Q_p$ is split.
    
            \subsection{Lifting theory for ordinary \texorpdfstring{$p$}{p}-divisible groups}
            \label{ssec:more_moduli_pDiv_split:lifting_theory}
                We discuss lifting theory for ordinary $p$-divisible groups over an algebraically closed field $\kappa$ of characteristic $p$. The case of height $2$ dimension $1$ ordinary $p$-divisible groups is discussed in \cite[{Appendix}]{Messing72}. We spell out the case of general height and dimension (which reduces to the results in \cite[{Appendix}]{Messing72}). See also the exposition in \cite{Meusers07} (or the sketch in \cite[{\S 6}]{Gross86}, though we will need some additional material on homomorphisms between liftings.

Take integers $r_1, r_2 \geq 0$. The unique ordinary $p$-divisible group $X$ over $\kappa$ of height $r_1 + r_2$ and dimension $r_1$ is $X = \pmb{\mu}_{p^{\infty}}^{r_1} \times \underline{\Q_p/\Z_p}^{r_2}$.

Next, let $R$ be an adic Noetherian local ring (with maximal ideal being an ideal of definition) with residue field $\kappa$. The $p$-divisible groups $\pmb{\mu}_{p^{\infty}}$ and $\underline{\Q_p/\Z_p}$ lift uniquely to $\Spf R$ (e.g. by Grothendieck--Messing deformation theory), which we still notate as $\pmb{\mu}_{p^{\infty}}$ and $\underline{\Q_p/\Z_p}$.
If $\mf{X}$ is a lift of $X$ over $\Spf R$, its connected-\'etale exact sequence must be $0 \ra \pmb{\mu}_{p^{\infty}}^{r_1} \ra \mf{X} \ra \underline{\Q_p/\Z_p}^{r_2} \ra 0$. Classifying lifts $\mf{X}$ is thus the same as classifying such extensions, which are in canonical bijection with $\Ext^1_{\Spf R}(\underline{\Q_p/\Z_p}^{r_2}, \pmb{\mu}_{p^{\infty}}^{r_1})$ (using also Drinfeld rigidity, as well as the fact $\Hom(\pmb{\mu}_{p^{\infty}}, \underline{\Q_p/\Z_p}) = \Hom(\underline{\Q_p/\Z_p}, \pmb{\mu}_{p^{\infty}}) = 0$). Here, $\Ext^1_{\Spf R}$ is calculated in the abelian category of fppf sheaves of abelian groups over $\Spf R$ (this is also \cite[{Appendix, Corollary (2.3)}]{Messing72}). We typically suppress the $R$-dependence in $\Hom(-,-)$.

Applying $\Hom(-, \pmb{\mu}_{p^{\infty}})$ to the short exact sequence of sheaves (fppf sheaves over $\Spf R$)
    \begin{equation}
    \begin{tikzcd}
    0 \arrow{r} & \underline{\Z} \arrow{r} & \underline{\Z[1/p]} \arrow{r} & \underline{\Q_p/\Z_p} \arrow{r} & 0
    \end{tikzcd}
    \end{equation}
gives a boundary morphism $\delta \colon \Hom(\underline{\Z}, \pmb{\mu}_{p^{\infty}}) \ra \Ext^1_{\Spf R}(\underline{\Q_p/\Z_p}, \pmb{\mu}_{p^{\infty}})$ in the associated long exact sequence. This map $\delta$ is an isomorphism \cite[{Appendix, Proposition (2.5)}]{Messing72}.\footnote{In loc. cit. this is stated for Artinian local rings $R$, but one can pass to the limit and obtain the statement here (compare \cite[{Appendix, Remark (2.2)}]{Messing72}).} By compatibility of $\Ext$ with finite direct sums, it follows that the boundary morphism $\delta \colon \Hom(\underline{\Z}^{r_2}, \pmb{\mu}_{p^{\infty}}^{r_1}) \ra \Ext^1_{\Spf R}(\underline{\Q_p/\Z_p}^{r_2}, \pmb{\mu}_{p^{\infty}}^{r_1})$ is also an isomorphism.

Given an element $\a \in \Hom(\underline{\Z}^{r_2}, \pmb{\mu}^{r_1}_{p^{\infty}})$, we can identify the extension corresponding to $\delta(\a)$ with the bottom row of the diagram
    \begin{equation}
    \begin{tikzcd}
    0 \arrow{r} & \underline{\Z}^{r_2} \arrow{r} \arrow{d}{\a} & \underline{\Z[1/p]}^{r_2} \arrow{r} \arrow{d} & \underline{\Q_p/\Z_p}^{r_2} \arrow{r} \arrow[Equal]{d} & 0
    \\
    0 \arrow{r} & \pmb{\mu}_{p^{\infty}}^{r_1} \arrow{r} & \mf{X} \arrow{r} \arrow[ul, phantom, "\ulcorner", very near start] & \underline{\Q_p/\Z_p}^{r_2} \arrow{r} & 0
    \end{tikzcd}
    \end{equation}
where the rows are exact and the left square is a pushout. This follows from general homological algebra valid in any abelian category (e.g. \cite[\href{https://stacks.math.columbia.edu/tag/010I}{Section 010I}]{stacks-project} and \cite[\href{https://stacks.math.columbia.edu/tag/06XP}{Section 06XP}]{stacks-project}).

Given $r_1, r_1', r_2, r_2' \in \Z_{\geq 0}$, we have
    \begin{align}
    \Hom(\pmb{\mu}_{p^{\infty}}^{r_1} \times \underline{\Q_p/\Z_p}^{r_2}, \pmb{\mu}_{p^{\infty}}^{r_1'} \times \underline{\Q_p/\Z_p}^{r_2'}) & = \Hom(\pmb{\mu}_{p^{\infty}}^{r_1}, \pmb{\mu}_{p^{\infty}}^{r_1'}) \times \Hom(\underline{\Q_p/\Z_p}^{r_2}, \underline{\Q_p/\Z_p}^{r_2'})
    \notag
    \\
    & \cong M_{r_1', r_1}(\Z_p) \times M_{r_2',r_2}(\Z_p),
    \end{align}
since any $p$-divisible group over $\kappa$ of height $1$ has endomorphism ring $\Z_p$. 
Here $M_{s,t}(\Z_p)$ denotes $s \times t$ matrices with entries in $\Z_p$. Given 
    \begin{equation}\label{equation:more_moduli_pDiv_split:lifting_theory:Z_p-linear}
    \a \in \Hom(\underline{\Z}^{r_2}, \pmb{\mu}_{p^{\infty}}^{r_1}) = \Hom_{\Z_p}(\underline{\Z_p}^{r_2}, \pmb{\mu}_{p^{\infty}}^{r_1})
    \quad \quad
    \a' \in \Hom(\underline{\Z}^{r_2'}, \pmb{\mu}_{p^{\infty}}^{r_1'}) = \Hom_{\Z_p}(\underline{\Z_p}^{r_2'}, \pmb{\mu}_{p^{\infty}}^{r_1'})
    \end{equation}
with corresponding lifts $\mf{X}$ and $\mf{X}'$ (of $\pmb{\mu}_{p^{\infty}}^{r_1} \times \underline{\Q_p/\Z_p}^{r_2}$ and $\pmb{\mu}_{p^{\infty}}^{r_1'} \times \underline{\Q_p/\Z_p}^{r_2'}$ respectively) over $\Spf R$, a morphism $(f_1, f_2) \in M_{r'_1,r_1}(\Z_p) \times M_{r'_2, r_2}(\Z_p)$ lifts to a map $f \colon \mf{X} \ra \mf{X}'$ if and only if 
    \begin{equation}\label{equation:split:criterion_hom_lift_ext}
    f_1 \circ \a = \a' \circ f_2,
    \end{equation}
again by general facts about $\Ext$ in abelian categories (compare with the proof of \cite[{Appendix, Proposition (3.3)}]{Messing72}, which discusses the case $r_1 = r_2 = r'_1 = r'_2 = 1$). We will repeatedly use this criterion for lifting to maps $f \colon \mf{X} \ra \mf{X}'$. In \eqref{equation:more_moduli_pDiv_split:lifting_theory:Z_p-linear}, the subscripts $\Z_p$ indicate $\Z_p$-linearity (not the base $\Spf \Z_p$).
    
            \subsection{Quasi-canonical lifting cycles: split}
            \label{ssec:more_moduli_pDiv_split:qcan_cycles_split}
                Throughout Section \ref{ssec:more_moduli_pDiv_split:qcan_cycles_split}, we write $R$ for an adic Noetherian local ring (with maximal ideal being an ideal of definition) equipped with a morphism $\Spf R \ra \Spf \mc{O}_{\breve{F}}$ inducing an isomorphism on residue fields.

Allowing arbitrary signature $(n - r, r)$ for the moment, form the Rapoport--Zink spaces $\mc{N}(n - r, r)$, $\mc{N}_{n,r}$, and $\widetilde{\mc{N}}_{n,r}$ as in Section \ref{ssec:moduli_pDiv:Serre_tensor}. 
With $(\mbf{X}, \iota_{\mbf{X}}, \lambda_{\mbf{X}})$ denoting the framing object for $\mc{N}(n - r, r)$, we take $\mbf{X}^-$ to be the framing object used to define $\mc{N}_{n,r}$ and $\widetilde{\mc{N}}_{n,r}$. There are non-canonical isomorphisms $\mbf{X}^- \cong \pmb{\mu}_{p^{\infty}}^{r} \times \underline{\Q_p / \Z_p}^{n-r}$ and $\mbf{X}_0^- \cong \underline{\Q_p / \Z_p}$.

\begin{definition}\label{definition:more_moduli_pDiv_split:qcan_cycles_split:qcan_cycles}
Given a subset $L^- \subseteq \mbf{W}^- = \Hom^0(\mbf{X}_0^-, \mbf{X}^-)$, consider the associated \emph{local special cycle}
    \begin{equation}
    \mc{Y}(L^-) \subseteq \mc{N}_{n,r}
    \end{equation}
which is the subfunctor consisting of pairs $(X, \rho)$ over schemes $S$ over $\Spf \mc{O}_{\breve{F}}$ such that, for all $x^- \in L^-$, the quasi-homomorphism
    \begin{equation}
    \rho^{-1} \circ x'_{\overline{S}} \circ \rho_{\mf{X}_0, \overline{S}}^- \colon \mf{X}_0^- \ra X_{\overline{S}}
    \end{equation}
lifts to a homomorphism $\mf{X}_0^- \ra X$.
\end{definition}

As in Definition \ref{definition:moduli_pDiv:RZ:canonical_lifting} (also Section \ref{ssec:can_and_qcan:can}), the notation $\mf{X}_0$ refers to the canonical lifting of $\mbf{X}_0$ (and $\mf{X}_0 = \mf{X}_0^+ \times \mf{X}_0^-$ is the decomposition via the nontrivial idempotents $e^{\pm} \in \mc{O}_F$, with $\mbf{X}_0^+ \cong \pmb{\mu}_{p^{\infty}}$ and $\mf{X}_0^- \cong \underline{\Q_p / \Z_p}$).
Again, $\mc{Y}(L^-) \subseteq \mc{N}_{n,r}$ is a closed subfunctor (hence a locally Noetherian formal scheme) by \cite[{Proposition 2.9}]{RZ96} for quasi-homomorphisms.

\begin{lemma}\label{lemma:more_moduli_pDiv_split:qcan_cycles_split:pullback_cycle}
Suppose $L \subseteq \mbf{W} = \Hom_F^0(\mbf{X}_0, \mbf{X})$ is a subset with $L^+ \subseteq \mbf{L}^+ = \Hom(\mbf{X}_0^+, \mbf{X}^+)$. The natural commutative diagram
    \begin{equation}
    \begin{tikzcd}
    \mc{Y}(L^-) \arrow{d} \arrow{rr} \pullbackcorner & & \mc{Z}(L) \arrow{d}
    \\
    \mc{N}_{n,r} \arrow[hook]{r} & \widetilde{\mc{N}}_{n,r} \arrow{r}{\sim} & \mc{N}(n - r, r)
    \end{tikzcd}
    \end{equation}
is Cartesian.
\end{lemma}
\begin{proof}
The lower horizontal arrows are as described in Section \ref{ssec:moduli_pDiv:Serre_tensor} and Lemma \ref{lemma:moduli_pdiv:Serre_tensor:split_forget} (the composite is an open and closed immersion). The lemma amounts to the claim that, for any $(X, \iota, \lambda, \rho) \in X(S)$ (for some scheme $S$ over $\Spf \mc{O}_{\breve{F}}$), if $x = x^+ \times x^- \in \Hom_F^0(\mbf{X}_0, \mbf{X})$ with $x^+ \in \Hom_{\mc{O}_F}(\mbf{X}_0^+, \mbf{X}^+)$, then $x$ lifts to a homomorphism $\mf{X}_0 \ra X$ if and only if $x^-$ lifts to a homomorphism $\mf{X}_0^- \ra X^-$. Stated alternatively, this is the claim that $x^+$ always lifts to a homomorphism $\mf{X}_0^+ \ra X^+$. Since $\mc{N}_{n,r} \cong \Spf \psring{\mc{O}_{\breve{F}}}{t_1, \ldots, t_{(n-r)r}}$,
this is clear because $\mf{X}_0^+ \ra X^+$ automatically factors through a homomorphism to the connected part $(X^+)^0 \cong \pmb{\mu}_{p^{\infty}}^{n -r}$ of $X^+$, over any base $\Spf R$ where $R$ is Noetherian Henselian local ring (alternative proof: apply \eqref{equation:split:criterion_hom_lift_ext}).
\end{proof}

Choose isomorphisms $(\mbf{X}^-)^0 \cong \pmb{\mu}_{p^{\infty}}^r$ and $(\mbf{X}^-)^{\et} \cong \underline{\Q_p / \Z_p}^{n - r}$ for the connected and \'etale parts of $\mbf{X}^-$ respectively. Any element $(X, \rho) \in \mc{N}_{n,r}(\Spf R)$ (i.e. a morphism $\Spf R \ra \mc{N}_{n,r}$) then corresponds to a class $\a \in \Ext^1_{\Spf R}(\underline{\Q_p / \Z_p}^{n - r}, \pmb{\mu}_{p^{\infty}}^r) = \Hom_{\Z_p}(\underline{\Z_p}^{n - r}, \pmb{\mu}_{p^{\infty}}^r)$ via the lifting theory in Section \ref{ssec:more_moduli_pDiv_split:lifting_theory}.

\begin{lemma}\label{lemma:more_moduli_pDiv_split:qcan_cycles_split:hom_factor_through}
Fix any isomorphism $\mbf{X}_0^- \cong \underline{\Q_p/\Z_p}$.
Consider $\varphi \colon \Spf R \ra \mc{N}_{n,r}$, corresponding to $(X, \rho) \in \mc{N}_{n,r}(\Spf R)$ and hence a class $\a' \in \Ext^1_{\Spf R}(\underline{\Q_p / \Z_p}^{n - r}, \pmb{\mu}_{p^{\infty}}^r)$. Given any subset $L^- \subseteq \mbf{L}^-$, the morphism $\varphi$ factors through $\mc{Y}(L^-) \subseteq \mc{N}_{n,r}$ if and only if the map
    \begin{equation}
    x^* \colon \Ext^1_{\Spf R}(\underline{\Q_p / \Z_p}^{n - r}, \pmb{\mu}_{p^{\infty}}^r) \ra \Ext^1_{\Spf R}(\underline{\Q_p / \Z_p}, \pmb{\mu}_{p^{\infty}}^r)
    \end{equation}
satisfies $x^*(\a') = 0$ for all $x \in L^-$.
\end{lemma}
\begin{proof}
In the lemma statement, we have viewed $x \in L^-$ as a morphism $\underline{\Q_p / \Z_p} \ra \underline{\Q_p / \Z_p}^{n - r}$ via the various identifications. The lemma follows from the lifting criterion in \eqref{equation:split:criterion_hom_lift_ext} (in the notation of loc. cit., take $\a = 0$).
\end{proof}

Next, we restrict to the case of signature $(n - 1, 1)$.

\begin{lemma}\label{lemma:split:type_1_and_lifts}
Assume that $R$ is moreover a domain and that $\Spf R \ra \Spf \mc{O}_{\breve{F}}$ is flat.
There is a natural map
    \begin{equation}
    \begin{tikzcd}
    \left \{ \begin{array}{l} \text{cyclic subgroups of order $p^s$} \\ \text{in $\Ext^1_{\Spf R}( \underline{\Q_p/\Z_p}^{n-1}, \pmb{\mu}_{p^{\infty}})$} \end{array} \right \}
    \arrow{r} & 
    \left \{ \begin{array}{l} \text{Full rank integral lattices $M \subseteq \mbf{L}$} \\ \text{such that $t(M) \leq 1$ and $\operatorname{val}(M) = s$} \\ \text{and $M^{+} = \mbf{L}^{+}$} \end{array} \right \}    
    \end{tikzcd}
    \end{equation}
(functorial in $R$ on the left). If $R$ contains a primitive $p^s$-th root of unity, then the map is a bijection. Otherwise, the left-hand side is empty.
\end{lemma}
\begin{proof}
Recall the identification $\Hom_{\Z_p}(\underline{\Z_p}^{n-1}, \pmb{\mu}_{p^{\infty}}) \cong \Ext^1_{\Spf R}(\underline{\Q_p/\Z_p}^{n-1}, \pmb{\mu}_{p^{\infty}})$ from Section \ref{ssec:more_moduli_pDiv_split:lifting_theory}.
Suppose $\a' \in \Hom_{\Z_p}(\underline{\Z_p}^{n-1}, \pmb{\mu}_{p^{\infty}})$ generates a cyclic subgroup of order $p^s$ (possible if and only if $R$ contains a primitive $p^s$-th root of unity). Let $M_{n-1,n-1}(\Z_p)$ act on $\Hom_{\Z_p}(\underline{\Z_p}^{n-1}, \pmb{\mu}_{p^{\infty}})$ by pre-composition. The annihilator of $\a'$ is generated (as a one-sided ideal) by an element $f_2 \in M_{n-1,n-1}(\Z_p)$ which has Smith normal form $\operatorname{diag}(1,\ldots,1,p^s)$.

We have a canonical identification $\mbf{L}^{-} \cong \Hom(\mbf{X}_0^{-}, (\mbf{X}^-)^{\et})$. Via the identification $(\mbf{X}^-)^{\et} \cong \underline{\Q_p/\Z_p}^{n-1}$, we obtain an action of $M_{n-1,n-1}(\Z_p)$ on $\mbf{L}^{-}$ (post-composition).
We then set $M^{-} = f_2 (\mbf{L}^{-})$, and let $M = \mbf{L}^{+} \oplus M^{- }$. Note that $M^{-}$ does not depend on the choice of generator $f_2$.

Conversely, given a lattice $M \subseteq \mbf{L}$ as in the lemma statement, select any $f_2 \in M_{n-1, n-1}(\Z_p)$ satisfying $M^{-} = f_2 (\mbf{L}^{-})$, and note that $f_2$ necessarily has Smith normal form $\operatorname{diag}(1, \ldots, 1, p^s)$. If $R$ contains a primitive $p^s$-th root of unity, then $f_2$ acting on $\Hom_{\Z_p}(\underline{\Z_p}^{n-1}, \pmb{\mu}_{p^{\infty}})$ has kernel which is cyclic of order $p^s$. This gives the inverse map.
\end{proof}

For $s \in \Z_{\geq 0}$, set $\breve{E}_s \coloneqq \breve{F}[\zeta_{p^s}]$ with ring of integers $\mc{O}_{\breve{E}_s} = \mc{O}_{\breve{F}}[\zeta_{p^s}]$, where $\zeta_{p^s}$ is a primitive $p^s$-th root of unity. 
Suppose $M \subseteq \mbf{L}$ is an integral full rank $\mc{O}_F$-lattice satisfying $M^+ = \mbf{L}^+$, with type $t(M) \leq 1$ and $\operatorname{val}(M) = s$.
By Lemma \ref{lemma:split:type_1_and_lifts}, there is an associated cyclic subgroup of $\Ext^1_{\Spf \mc{O}_{\breve{E}_s}}(\underline{\Q_p/\Z_p}^{n-1}, \pmb{\mu}_{p^{\infty}})$. Any generator of this cyclic subgroup defines a morphism $\Spf \mc{O}_{\breve{E}_s} \ra \mc{N}_{n,1}$ (via the lifting theory from Section \ref{ssec:more_moduli_pDiv_split:lifting_theory}). Changing the choice of generator corresponds precisely to the action of $\Gal(\breve{E}_s/\breve{F})$ (by Lubin--Tate theory for $\pmb{\mu}_{p^{\infty}}$). This morphism $\Spf \mc{O}_{\breve{E}_s} \ra \mc{N}_{n, 1}$ must be a closed immersion: if the morphism factors through $\Spf R \ra \mc{N}_{n, 1}$ for some sub $\mc{O}_{\breve{F}}$-algebra $R \subseteq \mc{O}_{\breve{E}_s}$, then Lemma \ref{lemma:split:type_1_and_lifts} implies that $R = \mc{O}_{\breve{E}_s}$.

We write $\Spf \mc{O}_{\breve{E}_s} \cong \mc{Z}(M)^{\qcancirc} \subseteq \mc{N}_{n, 1}$ for the resulting closed subfunctor, and call it a \emph{quasi-canonical lifting cycle}. This closed subfunctor $\mc{Z}(M)^{\qcancirc}$ does not depend on the choices of isomorphisms $(\mbf{X}^-)^0 \cong \pmb{\mu}_{p^{\infty}}$ and $(\mbf{X}^-)^{\et} \cong \underline{\Q_p / \Z_p}^{n - 1}$ appearing in the statement of Lemma \ref{lemma:split:type_1_and_lifts}.

\begin{lemma}\label{lemma:more_moduli_pDiv_split:qcan_cycles_split:qcan_cycle_product_decomp}
With $M$ as above, view $\mc{Z}(M)^{\qcancirc}$ as a morphism $\Spf \mc{O}_{\breve{E}_s} \ra \mc{N}_{n,1}$ corresponding to $(X, \rho) \in \mc{N}_{n,1}(\Spf \mc{O}_{\breve{E}_s})$.

If $n = 1$ then $X \cong \pmb{\mu}_{p^{\infty}}$. If $n \geq 2$ then $X \cong \underline{\Q_p / \Z_p}^{n - 2} \times \mf{X}_s$ (forgetting $\rho$) where $\mf{X}_s$ is a $p$-divisible group of height $2$ and dimension $1$ with $\End(\mf{X}_s) = \Z_p + p^s \mc{O}_F$ (a \emph{quasi-canonical lifting} in the sense of Section \ref{ssec:can_and_qcan:qcan}).
\end{lemma}
\begin{proof}
Let $\a' \in \Hom_{\Z_p}(\underline{\Z_p}^{n - 1}, \pmb{\mu}_{p^{\infty}})$ be the element corresponding to $(X, \rho)$. If $n = 1$ then $\a' = 0$ and $X \cong \pmb{\mu}_{p^{\infty}}$.

If $n \geq 2$, then (after replacing $\rho$ by $\phi \circ \rho$ for some $\phi \in \GL_{n - 1}(\Z_p)$), the lift $(X, \rho)$ of $\underline{\Q_p / \Z_p}^{n - 1} \times \pmb{\mu}_{p^{\infty}}$ is associated with $\a'$ of the form $(0, \ldots, 0, \zeta_{p^s})$ for $\zeta_{p^s} \in \mc{O}_{\breve{E}_s}$ a primitive $p^s$-th root of unity. For some $\mf{X}_s$ as in the lemma statement, we obtain a commutative diagram
    \begin{equation}\label{equation:more_moduli_pDiv_split:qcan_cycles_split:qcan_cycle_product_decomp}
    \begin{tikzcd}
    0 \arrow{r} & \underline{\Z}^{n-1} \arrow{r} \arrow{d}[swap]{(0,\ldots,0,\zeta_{p^s})} & \underline{\Z[1/p]}^{n-1} \arrow{r} \arrow{d} & \underline{\Q_p/\Z_p}^{n-1} \arrow{r} \arrow[Equal]{d} & 0 \\
    0 \arrow{r} & \pmb{\mu}_{p^{\infty}} \arrow{r} & \arrow{r} \underline{\Q_p/\Z_p}^{n-2} \times \mf{X}_s & \underline{\Q_p/\Z_p}^{n-1} \arrow{r} & 0\rlap{~,}
    \end{tikzcd}
    \end{equation}
using the lifting criterion of \eqref{equation:split:criterion_hom_lift_ext} again.
\end{proof}

\begin{lemma}\label{lemma:more_moduli_pDiv_split:qcan_cycles_split:singleton_containment}
Suppose $M \subseteq \mbf{L}$ is an integral full rank $\mc{O}_F$-lattice satisfying $M^+ = \mbf{L}^+$ with type $t(M) \leq 1$ and $\mrm{val}(M) = s$. Let $L^- \subseteq \mbf{W}^-$ be any subset. We have $\mc{Z}(M)^{\qcancirc} \subseteq \mc{Y}(L^-)$ if and only if $L^- \subseteq M$.
\end{lemma}
\begin{proof}
If $n = 1$ then $\mbf{W}^- = 0$ and $\mc{Z}(M)^{\qcancirc} = \mc{Y}(L^-) = \mc{N}_{n,1}$, so the lemma is trivial in this case. We thus assume $n \geq 2$ below.

It is enough to check the case where $L^-$ consists of a single element, i.e. a quasi-homomorphism $x \colon (\mbf{X}_0)^- \ra \mbf{X}^-$ (or equivalently, $x \colon (\mbf{X}_0)^- \ra (\mbf{X}^-)^{\et}$ since $\mbf{X}_0^- \cong \underline{\Q_p / \Z_p}$ is \'etale). If $x \not \in \mbf{L}^-$ then $\mc{Y}(L^-) = \emptyset$ (while $\mc{Z}(M)^{\qcancirc} \neq \emptyset$), so we may assume $x \in \mbf{L}^-$.

Pick any identification $\mbf{X}_0^- \cong \underline{\Q_p / \Z_p}$.  Set $s = \operatorname{val}(M)$. Use the setup and notation in the proof of Lemma \ref{lemma:split:type_1_and_lifts}.

View $\mc{Z}(M)^{\qcancirc}$ as a closed immersion $\Spf \mc{O}_{\breve{E}_s} \ra \mc{N}_{n,1}$, corresponding to an element $\a' \in \Hom_{\Z_p}(\underline{\Z_p}^{n - 1}, \pmb{\mu}_{p^{\infty}})$. By Lemma \ref{lemma:more_moduli_pDiv_split:qcan_cycles_split:hom_factor_through}, our task is to show that $\a' \circ x = 0$ if and only if $\a' \in M^-$. Since $f_2$ generates (in $M_{n - 1, n - 1}(\Z_p)$) the annihilator of $\a'$ (as a one-sided ideal), we see that $\a' \circ x = 0$ if and only if $x \in f_2(\mbf{L}^-) = M^-$ (for example, view $x$ as a column vector and observe that $(x, 0, \ldots, 0) \in M_{n - 1, n - 1}(\Z_p)$ lies in the one-sided ideal generated by $f_2$).
\end{proof}

\begin{definition}\label{definition:more_moduli_pDiv_split:qcan_cycles_split}
Let $M \subseteq \mbf{W}$ be a full rank integral $\mc{O}_F$-lattice, with type $t(M) \leq 1$ and $\operatorname{val}(M) = s$. Select any $\gamma \in U(\mbf{W})$ satisfying $\gamma(\mbf{L}^+) = M^+$ (also write $\gamma$ for $(\gamma,1) \in U(\mbf{W}) \times U(\mbf{W}^{\perp})$, by abuse of notation).

The \emph{quasi-canonical lifting cycle} associated with $M$ is the closed subfunctor
    \begin{equation}
    \Spf \mc{O}_{\breve{E}_s} \cong \mc{Z}(M)^{\qcancirc} \coloneqq \gamma(\mc{Z}(\gamma^{-1}(M))^{\qcancirc}) \subseteq \mc{N}(n - 1, 1)
    \end{equation}
where $\g \in U(\mbf{W}) \times U(\mbf{W}^{\perp})$ acts on $\mc{N}(n - 1, 1)$ as in Section \ref{ssec:moduli_pDiv:action_RZ}.
\end{definition}

In the situation of Definition \ref{definition:more_moduli_pDiv_split:qcan_cycles_split}, the closed subfunctor $\mc{Z}(M)^{\qcancirc}$ does not depend on the choice of $\gamma$. We have also viewed $\mc{N}_{n, 1}$ as an open and closed subfunctor of $\mc{N}(n - 1, 1)$ (as in the lower horizontal arrows in Lemma \ref{lemma:more_moduli_pDiv_split:qcan_cycles_split:pullback_cycle}).

\begin{lemma}\label{lemma:split:qcan_cycle_containments}
If $L \subseteq \mbf{W}$ is any subset and $M \subseteq \mbf{W}$ is any full rank integral lattice with $t(M) \leq 1$, we have $\mc{Z}(M)^{\qcancirc} \subseteq \mc{Z}(L)$ if and only if $L \subseteq M$.
\end{lemma}
\begin{proof}
After acting by $U(\mbf{W})$, it is enough to check the case where $M^+ = \mbf{L}^+$. In this case, we have $\mc{Z}(M)^{\qcancirc} \subseteq \mc{N}_{n,1}$. If $L \not \subseteq \mbf{L}$, then $\mc{Z}(L) \cap \mc{Z}(M)^{\qcancirc} = \emptyset$ by Lemma \ref{lemma:moduli_pDiv:discrete_reduced}(5) and Remark \ref{remark:moduli_pDiv:Serre_tensor:singleton} (and $\mc{Z}(M)^{\qcancirc}$ is nonempty). So assume $L \subseteq \mbf{L}$. Then $\mc{Z}(L) = \mc{Y}(L^-)$ (Lemma \ref{lemma:more_moduli_pDiv_split:qcan_cycles_split:pullback_cycle}). This reduces to the case proved in Lemma \ref{lemma:more_moduli_pDiv_split:qcan_cycles_split:singleton_containment}.
\end{proof}

\begin{corollary}\label{corollary:split:qcan_cycle_inclusion}
Let $L \subseteq \mbf{W}$ be any subset. Form the horizontal (flat) part of the local special cycle $\mc{Z}(L)$, which we denote as $\mc{Z}(L)_{\ms{H}}$. We have an inclusion of closed formal subschemes
    \begin{equation}
    \bigcup_{\substack{L \subseteq M \subseteq M^{*} \\ t(M) \leq 1}} \mc{Z}(M)^{\qcancirc} \subseteq \mc{Z}(L)_{\ms{H}}
    \end{equation}
in $\mc{N}(n - 1, 1)$.
\end{corollary}
\begin{proof}
The union is a scheme-theoretic union (i.e. intersect associated ideal sheaves). The claim follows from Lemma \ref{lemma:split:qcan_cycle_containments} because each $\mc{Z}(M)^{\qcancirc}$ is flat over $\Spf \mc{O}_{\breve{F}}$.
\end{proof}

\begin{lemma}\label{lemma:split:qcan_cycles_distinct}
Let $M \subseteq \mbf{W}$ and $M' \subseteq \mbf{W}$ be integral full-rank $\mc{O}_F$-lattices with $t(M) \leq 1$ and $t(M') \leq 1$. If $M \neq M'$, then $\mc{Z}(M)^{\qcancirc} \neq \mc{Z}(M')^{\qcancirc}$.
\end{lemma}
\begin{proof}
Let $N \subseteq \mbf{W}$ (resp. $N' \subseteq \mbf{W}$) be the unique self-dual full rank lattice such that $N^+ = M^+$ (resp. $N^{\prime +} = M^{\prime +}$). On reduced subschemes, we have $\mc{Z}(M)^{\qcancirc}_{\red} = \mc{Z}(M')^{\qcancirc}_{\red}$ if and only if $N = N'$ by Lemma \ref{lemma:moduli_pDiv:discrete_reduced} (more precisely, Remark \ref{remark:moduli_pDiv:Serre_tensor:singleton}, Definition \ref{definition:more_moduli_pDiv_split:qcan_cycles_split:qcan_cycles}, and the action on special cycles in \eqref{equation:moduli_pDiv:action_RZ:special_cycles}).
So we may assume $N = N'$. Using the $U(\mbf{W})$ action on $\mc{N}(n - 1, 1)$, we also reduce to the case where $N = \mbf{L}$. 

Set $s = \operatorname{val}(M)$ and $s' = \operatorname{val}(M')$, and view $\mc{Z}(M)^{\qcancirc}$ and $\mc{Z}(M')^{\qcancirc}$ as closed immersions $\varphi \colon \Spf \mc{O}_{\breve{E}_s} \ra \mc{N}(n - 1, 1)$ and $\varphi' \colon \Spf \mc{O}_{\breve{E}_{s'}} \ra \mc{N}(n - 1, 1)$. Lemma \ref{lemma:split:type_1_and_lifts} implies that $M = M'$ if and only if both $s = s'$ and the morphisms $\varphi, \varphi'$ are the same up to $\Gal(\breve{E}_s / \breve{F})$-action (this is equivalent to requiring that the corresponding elements of $\Ext^1$ in that lemma generate the same subgroup). This is satisfied if and only if $\mc{Z}(M)^{\qcancirc} = \mc{Z}(M')^{\qcancirc}$.
\end{proof}

\begin{lemma}\label{lemma:split:full_rank_qcan_factor_through}
Let $L \subseteq \mbf{W}$ be a full rank $\mc{O}_F$-lattice. Assume that $R$ is moreover a domain and $\Spf R \ra \Spf \mc{O}_{\breve{F}}$ is flat. Any morphism $\varphi \colon \Spf R \ra \mc{Z}(L)$ factors through some quasi-canonical lifting cycle $\mc{Z}(M)^{\qcancirc}$.
\end{lemma}
\begin{proof}
Again, we may act by $U(\mbf{W})$ on $\mc{N}(n - 1, 1)$ to assume that $\varphi \colon \Spf R \ra \mc{Z}(L)$ factors through the open and closed component $\mc{N}_{n,1} \subseteq \mc{N}(n - 1, 1)$ described in Section \ref{ssec:moduli_pDiv:Serre_tensor} and above. This implies $L \subseteq \mbf{L}$ (as $\mc{Z}(L) \cap \mc{N}_{n,1}$ is otherwise empty, see Lemma \ref{lemma:moduli_pDiv:discrete_reduced}).

Fix isomorphisms as in the statement of Lemma \ref{lemma:split:type_1_and_lifts}. Then $\varphi$ corresponds to some $(X, \rho) \in \mc{N}_{n,1}$, and this lift of $\mbf{X}^-$ corresponds to a class $\a' \in \Ext^1(\underline{\Q_p / \Z_p}^{n - 1}, \pmb{\mu}_{p^{\infty}})$ via the lifting theory in Section \ref{ssec:more_moduli_pDiv_split:lifting_theory}.

By Lemma \ref{lemma:split:type_1_and_lifts}, it is enough to show that $\a'$ is $p^s$-torsion for some $s \in \Z_{\geq 0}$ (then $\varphi$ must factor through $\mc{Z}(M)^{\qcancirc}$ where $M$ is the lattice associated with the cyclic subgroup generated by $\a$). Select $s \geq 0$ such that $p^s \mbf{L} \subseteq L$ (such $s$ exists because $L$ is full rank). Then Lemma \ref{lemma:more_moduli_pDiv_split:qcan_cycles_split:hom_factor_through} implies $p^s a' = 0$, since $\varphi$ factors through $\mc{Z}(L)$ (and hence through $\mc{Y}(L^-)$).
\end{proof}

        \section{Canonical and quasi-canonical liftings}
        \label{sec:can_and_qcan}
            We retain $F / \Q_p$ and accompanying notation as in Section \ref{sec:moduli_pDiv}. In Sections \ref{ssec:can_and_qcan:can} and \ref{ssec:can_and_qcan:qcan}, we allow $p = 2$ even if $F / \Q_p$ is ramified.
We collect some needed facts about canonical and quasi-canonical lifts in all cases (inert, ramified, split). See also \cite{Gross86}, \cite{Wewers07}, \cite{Meusers07}. Our conventions differ slightly from \cite{Wewers07}, due to the phenomenon explained in \cite[{Footnote 7}]{KR11} (there in the inert case, which we also modify to apply in the ramified case).
    
            \subsection{Canonical liftings}
            \label{ssec:can_and_qcan:can}
                As in Section \ref{ssec:moduli_pDiv:RZ}, let $\mbf{X}_0$ be the unique supersingular (resp. ordinary) $p$-divisible group of height $2$ dimension $1$ over $\overline{k}$ if $F / \Q_p$ is nonsplit (resp. split). Let $j \colon \mc{O}_F \hookrightarrow \End(\mbf{X}_0)$ be a ring homomorphism. We reserve the notation $\iota_{\mbf{X}_0}$ to mean a signature $(1,0)$ action, and allow $j$ to have either signature (i.e. $(1,0)$ or $(0,1)$) for its action on $\Lie \mbf{X}_0$.

Let $\breve{E}$ be any finite degree field extension of $\breve{F}$, with ring of integers $\mc{O}_{\breve{E}}$.

The pair $(\mbf{X}_0, j)$ admits a lift $(\mf{X}_0, \iota_{\mf{X}_0}, \rho_{\mf{X}_0})$ over $\Spf \mc{O}_{\breve{E}}$ (i.e. $(\mf{X}_0, \iota_{\mf{X}_0})$ is a $p$-divisible group over $\Spf \mc{O}_{\breve{E}}$ with $\mc{O}_F$-action $\iota_{\mf{X}_0}$, and $\rho_{\mf{X}_0} \colon \mf{X}_{0, \overline{k}} \ra \mbf{X}_0$ is a $\mc{O}_F$-linear isomorphism with respect to $\iota$ and $j$).

In the supersingular case, the pair $(\mf{X}_0, \iota_{\mf{X}_0})$ may be described via Lubin--Tate formal groups. 
In the ordinary case, we have $\mf{X}_0 \cong \pmb{\mu}_{p^{\infty}} \times \underline{\Q_p / \Z_p}$.

By the \emph{signature} of $(\mf{X}_0, \iota_{\mf{X}_0}, \rho_{\mf{X}_0})$ (or $(\mf{X}_0, \iota_{\mf{X}_0})$), we mean the signature of $\iota_{\mf{X}_0}$ acting on $\Lie \mf{X}_0$ (either $(1,0)$ or $(0,1)$). If $F / \Q_p$ is unramified (resp. ramified), then $(\mf{X}_0, \iota_{\mf{X}_0}, \rho_{\mf{X}_0})$ must have the same signature as $(\mbf{X}_0, j)$ (resp. can have either signature).

After fixing a signature, the triple $(\mf{X}_0, \iota_{\mf{X}_0}, \rho_{\mf{X}_0})$ is unique up to unique isomorphism, and we call it the \emph{canonical lifting}\footnote{When $j$ has signature $(1,0)$, what Gross \cite{Gross86} calls a \emph{canonical lifting} is what we call a \emph{canonical lifting of signature $(1,0)$}. This change in terminology allows additional flexibility when discussing quasi-canonical liftings, to account for e.g. \cite[{Footnote 7}]{KR11}.} of $(\mbf{X}_0, j)$. The canonical lifting over $\Spf \mc{O}_{\breve{E}}$ is defined over $\Spf \mc{O}_{\breve{F}}$ (i.e. is the base change of the canonical lift over $\Spf \mc{O}_{\breve{F}}$).

The map $\iota_{\mf{X}_0} \colon \mc{O}_F \ra \End(\mf{X}_0)$ is an isomorphism, since $\End(\mf{X}_0)$ is commutative and $\mc{O}_F$ is self-centralizing in $\End(\mbf{X}_0)$ (in the nonsplit case, note $\End(\mf{X}_0) \hookrightarrow \End(\Lie \mf{X}_0) = \mc{O}_{\breve{E}}$ so $\End(\mf{X}_0)$ must be commutative).

If $(\mf{X}_0^{\s}, \iota_{\mf{X}_0}^{\s})$ is as in \eqref{equation:sigma_twist_triple}, we have $\Hom_{\mc{O}_F}(\mf{X}_0, \mf{X}_0^{\s}) = 0$ because $\End(\mf{X}_0) = \mc{O}_F$.

\begin{example}\label{example:no_lift_quasi_hom}
Assume $F / \Q_p$ is nonsplit, and let $\mf{X}_0$ be the canonical lifting over $\Spf \mc{O}_{\breve{E}}$ (of some fixed signature). Drinfeld rigidity for quasi-homomorphisms implies $\End^0(\mf{X}_{0}) \cong \End^0(\mbf{X}_0) \cong D$, where $D$ is the quaternion division algebra over $\Q_p$. On the other hand, if $\mf{X}'_0$ denotes the $p$-divisible group over $\Spec \mc{O}_{\breve{E}}$ associated with $\mf{X}_0$ (see e.g. \crefext{III:lemma:p_div_completion_equiv}), we have $\End^0(\mf{X}'_0) \cong \End(\mf{X}'_0) \otimes_{\Z_p} \Q_p \cong F$. Thus, by our conventions (explained in \crefext{III:appendix:pDiv_prelim:terminology}), quasi-homomorphisms do not necessarily lift along the equivalence of $p$-divisible groups over $\Spec \mc{O}_{\breve{E}}$ and $\Spf \mc{O}_{\breve{E}}$ from \crefext{III:lemma:p_div_completion_equiv}.
See also \crefext{III:remark:no_lift_quasi_hom}.
\end{example}
            
            \subsection{Quasi-canonical liftings}
            \label{ssec:can_and_qcan:qcan}
                Let $\breve{E}$ and $(\mbf{X}_0, j)$ be as in Section \ref{ssec:can_and_qcan:can}.
For integers $s \geq 0$, let $\mc{O}_{F,s} \coloneqq \Z_p + p^s \mc{O}_F$ be the order of index $p^s$ in $\mc{O}_F$. When $F / \Q_p$ is nonsplit (resp. split) the subgroup $\mc{O}_{F,s}^{\times} \subseteq \mc{O}_F^{\times}$ (resp. $(1 + p^s \Z_p)^{\times} \subseteq \Z_p^{\times}$) has an associated finite totally ramified abelian extension $\breve{E}_s$ of $\breve{F}$ by local class field theory. The index is
    \begin{equation}\label{equation:can_and_qcan:qcan:degree}
    [\breve{E}_s : \breve{F}] = 
    \begin{cases}
    p^s (1 - \eta(p) p^{-1}) & \text{$s \geq 1$} \\
    1 & \text{if $s = 0$}.
    \end{cases}
    \end{equation}
where $\eta(p) \coloneqq -1, 0, 1$ in the inert, ramified, and split cases respectively. 
In the split case, we have $\mc{O}_{\breve{E}_s} = \mc{O}_{\breve{F}}[\zeta_{p^s}]$ where $\zeta_{p^s}$ is a primitive $p^s$-th root of unity.

In all cases, a \emph{quasi-canonical lifting of level $s$} of $(\mbf{X}_0, j)$ is a triple $(\mf{X}_s, \iota_{\mf{X}_s}, \rho_{\mf{X}_s})$ where
    \begin{align*}
    & \mf{X}_s && \text{is a $p$-divisible group over $\Spf \mc{O}_{\breve{E}}$} \\
    & \iota_{\mf{X}_s} \colon \mc{O}_{F,s} \xra{\sim} \End(\mf{X}_s) && \text{is a ring isomorphism} \\
    & \rho_{\mf{X},s} \colon \mf{X}_{s, \overline{k}} \ra \mbf{X}_0 && \text{is a $\mc{O}_{F,s}$-linear isomorphism of $p$-divisible groups over $\overline{k}$.}
    \end{align*}
Note that a quasi-canonical lifting of level $s = 0$ is the same as a canonical lifting. As above, we speak of the \emph{signature} of a quasi-canonical lifting, which means the signature of the action $\iota_{\mf{X}_s} |_{\Lie \mf{X}_s}$.

The signature of $(\mbf{X}_0, j)$ and the signature of a level $s$ quasi-canonical lifting must be
    \begin{equation}
    \begin{cases}
    \text{same} & \text{if $F/\Q_p$ is inert and $s$ is even, or $F / \Q_p$ is split} \\
    \text{opposite} & \text{if $F/\Q_p$ is inert and $s$ is odd} \\
    \text{either signature} & \text{if $F / \Q_p$ is ramified}.
    \end{cases}
    \end{equation}
Quasi-canonical liftings of level $s \geq 0$ exist in all such situations, and are defined over $\Spf \mc{O}_{\breve{E}_s}$. The property of being a level $s$ quasi-canonical lifting is preserved under base change along $\Spf \mc{O}_{\breve{E}'} \ra \Spf \mc{O}_{\breve{E}}$ for any finite degree field extension $\breve{E}'$ over $\breve{E}$.
If $F / \Q_p$ is split, a choice of level $s$ quasi-canonical lifting corresponds to a choice of morphism $\underline{\Z}_p \ra \pmb{\mu}_{p^{\infty}}$ over $\Spf \mc{O}_{\breve{E}}$ of exact order $p^s$ (i.e. a choice of primitive $p^s$-th root of unity in $\breve{E}_s$) via the lifting theory in Section \ref{ssec:more_moduli_pDiv_split:lifting_theory}.

The group $\Gal(\breve{E}_s/\breve{F})$ acts simply transitively on the set of level $s$ quasi-canonical liftings for any fixed signature (if such liftings exist).
By Lubin--Tate theory, this action is compatible with the identification $\Gal(\breve{E}_s/\breve{F}) \cong \mc{O}_F^{\times} / \mc{O}_{F,s}^{\times}$ via local class field theory (normalized to send uniformizers to geometric Frobenius) where $a \in \mc{O}_F^{\times}$ acts on the set of quasi-canonical liftings as $(\mf{X}_s, \iota_{\mf{X}_s}, \rho_{\mf{X}_s}) \mapsto (\mf{X}_s, \iota_{\mf{X}_s}, a \rho_{\mf{X}_s})$. In the split case, we have used the isomorphism
    \begin{equation}
    \begin{tikzcd}[row sep = tiny]
    \mc{O}_F^{\times} / \mc{O}_{F,s}^{\times} \arrow{r} & \Z_p^{\times} / (1 + p^s \Z_p)^{\times} \\
    x \arrow[mapsto]{r} & e^+(x) e^-(x^{-1})
    \end{tikzcd}    
    \end{equation}
if $(\mbf{X}_0, j)$ has signature $(1,0)$ and its reciprocal if $(\mbf{X}_0, j)$ has signature $(0,1)$. In particular, the quasi-canonical liftings of a fixed level $s$ are all isomorphic if the framing $\rho_{\mf{X}_s}$ is forgotten.

Let $(\mf{X}_0, \iota_{\mf{X}_0}, \rho_{\mf{X}_0})$ and $(\mf{X}_s, \iota_{\mf{X}_s}, \rho_{\mf{X}_s})$ be canonical and quasi-canonical lifts over $\Spf \mc{O}_{\breve{E}}$, for some $(\mbf{X}_0, j)$ and $(\mbf{X}_0, j')$ respectively (possibly $j \neq j'$). Then
    \begin{equation}\label{equation:can_and_qcan:qcan:psi_s}
    \Hom(\mf{X}_0, \mf{X}_s) \cong \psi_s \cdot \mc{O}_F
    \end{equation}
(no $\mc{O}_F$-linearity imposed) is a free $\mc{O}_F$-module of rank $1$ (where $\mc{O}_F$ acts by pre-composition), generated by some isogeny $\psi_s$ of degree $p^s$. The isogeny $\psi_s$ is defined over $\Spf \mc{O}_{\breve{E}_s}$. If $\mf{X}_0$ and $\mf{X}_s$ have the same signature, then $\psi_s$ is automatically $\mc{O}_{F,s}$-linear.

When $F / \Q_p$ is split, we may take $\psi_s$ to be the map inducing the map $\mbf{X}_0 \ra \mbf{X}_0$ which is
    \begin{equation}\label{equation:can_and_qcan:qcan:split_psi_s}
    \psi_s|_{\mbf{X}_0^{0}} \colon \mbf{X}_0^0 \xra{\mrm{id}} \mbf{X}_0^0  \quad \quad \psi_s|_{\mbf{X}_0^{\et}} \colon \mbf{X}_0^{\et} \xra{\times p^s} \mbf{X}_0^{\et}.
    \end{equation}
on the connected and \'etale parts, respectively. This follows from the lifting criterion in \eqref{equation:split:criterion_hom_lift_ext}.

For any generator $\psi_s$ of $\Hom(\mf{X}_0, \mf{X}_s)$, we have
    \begin{align}\label{equation:can_and_qcan:qcan:Nakkajima--Taguchi}
    \operatorname{length}_{\mc{O}_{\breve{E}}}(e^* \Omega^1_{\ker \psi_s / \Spec \mc{O}_{\breve{E}}}) & = \frac{1}{2} [\breve{E} : \breve{\Q}_p] \frac{(1 - p^{-s})(1 - \eta_p(p))}{(1 - p^{-1})(p - \eta_p(p))}
    \end{align}
where $\eta(p) \coloneqq -1, 0, 1$ in the inert, ramified, split cases respectively and where $e \colon \Spec \mc{O}_{\breve{E}} \ra \ker \psi_s$ denotes the identity section.
We are passing between $\Spf \mc{O}_{\breve{E}}$ and $\Spec \mc{O}_{\breve{E}}$ as in Appendix \crefext{III:appendix:pDiv_prelim:Spec_v_Spf}.

The nonsplit case of \eqref{equation:can_and_qcan:qcan:Nakkajima--Taguchi} is essentially a computation of Nakkajima and Taguchi \cite{NT91} (see also \cite[{Proposition 10.3}]{KRY04} and its proof). The split case follows from \eqref{equation:can_and_qcan:qcan:split_psi_s}, which implies that $\ker \psi_s$ is \'etale over $\Spec \overline{k}$ (cf. the closely related \cite[{Proposition 10.1}]{KRY04}).

The following constant $\delta_{\mrm{tau}}(s) \in \Q$ (``local change of tautological height'') will be crucial for the formulation of our local main theorems. With notation as above, we define
    \begin{align}\label{equation:can_and_qcan:qcan:local_change_taut}
    \delta_{\mrm{tau}}(s)
    & \coloneqq - \frac{1}{2} v_p(\deg \psi_s) + \frac{1}{[\breve{E} : \breve{\Q}_p]} \operatorname{length}_{\mc{O}_{\breve{E}}}(e^* \Omega^1_{\ker \psi_s / \Spec \mc{O}_{\breve{E}}})
    \\
    & = -\frac{1}{2} \left ( s - \frac{(1 - p^{-s})(1 - \eta_p(p))}{(1 - p^{-1})(p - \eta_p(p))} \right )
    \end{align}
for integers $s \in \Z$, with $\eta_p(p) \coloneqq -1, 0, 1$ in the inert, ramified, split cases respectively. We used \cref{equation:can_and_qcan:qcan:Nakkajima--Taguchi} for the second equality.

The quantity $\delta_{\mrm{tau}}(s)$ depends only on $s$ and $\eta_p(p)$, and does not depend on the choice of $\psi_s$. We also set $\delta_{\mrm{Fal}}(s) \coloneqq - 2 \delta_{\mrm{tau}}(s)$. In \cref{part:local_change_heights} below, we will explain the relation of $\delta_{\mrm{tau}}(s)$ and $\delta_{\mrm{Fal}}(s)$ with local decompositions of ``tautological'' and Faltings heights of special cycles.
    
            \subsection{Quasi-canonical lifting cycles}
            \label{ssec:can_and_qcan:qcan_cycles}
                 We state how certain local special cycles decompose into \emph{quasi-canonical lifting cycles} (Section \ref{ssec:more_moduli_pDiv_split:qcan_cycles_split}). We continue to use the notation in Section \ref{ssec:moduli_pDiv:special_cycles}, now restricting to signature $(n - 1, 1)$. We also assume $p \neq 2$ unless $F / \Q_p$ is split (in the inert case, this is so that we may cite \cite[{Theorem 4.2.1}]{LZ22unitary}).

Suppose $M^{\flat} \subseteq \mbf{W}$ is an integral $\mc{O}_F$-lattice of rank $n - 1$ with $t(M^{\flat}) \leq 1$. Set $s = \lfloor \mrm{val}(M^{\flat}) \rfloor$ (notation as in \cref{ssec:Hermitian_conventions:lattices}).
There is an associated \emph{quasi-canonical lifting cycle} $\mc{Z}(M^{\flat})^{\qcancirc} \subseteq \mc{N}(n - 1, 1)$, which is a certain closed subfunctor such that
    \begin{equation}\label{equation:can_and_qcan:qcan_cycles:underlying_formal_scheme}
    \mc{Z}(M^{\flat})^{\qcancirc} \cong 
    \begin{cases}
    \Spf \mc{O}_{\breve{E}_s} & \text{if $F / \Q_p$ is unramified} \\
    \Spf \mc{O}_{\breve{E}_s} \sqcup \Spf \mc{O}_{\breve{E}_s} & \text{if $F / \Q_p$ is ramified}.
    \end{cases}
    \end{equation}
Suppose $\varphi \colon \Spf \mc{O}_{\breve{E}_s} \ra \mc{N}(n - 1, 1)$ is a morphism representing any component of $\mc{Z}(M^{\flat})^{\qcancirc}$, with corresponding tuple $(X, \iota, \lambda, \rho) \in \mc{N}(n - 1, 1)(\Spf \mc{O}_{\breve{E}_s})$. 
If $n = 1$, then $M^{\flat} = 0$ and $X \cong \mf{X}_0^{\s}$. 
If $n \geq 2$, then there exists a polarization-preserving $\mc{O}_F$-linear isomorphism (forgetting $\rho$)
    \begin{equation}\label{equation:qcan_cycles:pDiv_Serre_tensor}
    X \cong (\mf{X}_0)^{n - 2} \times (\mf{X}_s \otimes_{\Z_p} \mc{O}_F)
    \end{equation}
for some level $s$ quasi-canonical lift $\mf{X}_s$ (and $\mf{X}_0$ being the canonical lift), where $\mf{X}_s \otimes_{\Z_p} \mc{O}_F$ is equipped with the polarization as in \eqref{equation:moduli_pDiv:Serre_tensor:1_1_signature}, where $\mf{X}_0^{n - 2}$ has the diagonal polarization $\lambda_{\mf{X}_0}^{n - 2}$ for some principal polarization $\lambda_{\mf{X}_0}$ on $\mf{X}_0$ if $F / \Q_p$ is unramified, and where $\mf{X}_0^{n - 2}$ has a product polarization as in \eqref{equation:moduli_pDiv:RZ:framing_object_polarization:ramified} (with respect to some principal polarization $\lambda_{\mf{X}_0}$ on $\mf{X}_0$) if $F / \Q_p$ is ramified.

For the inert case of the above assertions, see \cite[{\S 4.2}]{LZ22unitary} (we are using the same notation), and also \cite[{Proposition 8.1}]{KR11} (there for $n = 2$). 

For the ramified case, see \cite[{Proposition 7.1}]{RSZ17} (there for $n = 2$) and also the proof of \cite[{Proposition 2.44}]{LL22II} (also \cite[{Definition 2.45}]{LL22II}; we are using their notation but with $\mc{N}$ replaced by $\mc{Z}$).
In the ramified case, the two components $\mc{Z}(M^{\flat})^{\qcancirc}$ correspond to the two components of $\mc{N}(n - 1, 1)$ (as in in Lemma \ref{lemma:moduli_pDiv:discrete_reduced}, particularly part (5)), i.e. $\mc{Z}(M^{\flat})^{\qcancirc} \ra \mc{N}(n - 1, 1)$ is surjective on underlying topological spaces.

For the split case, $\mc{Z}(M^{\flat})^{\qcancirc}$ was defined in Definition \ref{definition:more_moduli_pDiv_split:qcan_cycles_split}. The assertion $X \cong (\mf{X}_0)^{n - 2} \times (\mf{X}_s \otimes_{\Z_p} \mc{O}_F)$ follows from Lemma \ref{lemma:more_moduli_pDiv_split:qcan_cycles_split:qcan_cycle_product_decomp} (note that $X$ in loc. cit. is $X^-$ in the present notation) and Remark \ref{remark:moduli_pDiv:Serre_tensor:split_and_lifts}.

\begin{proposition}\label{proposition:can_and_qcan:qcan_cycles}
 Let $L^{\flat} \subseteq \mbf{W}$ be an $\mc{O}_F$-lattice of rank $n - 1$. Form the horizontal (flat) part of the local special cycle $\mc{Z}(L^{\flat})$, which we denote as $\mc{Z}(L^{\flat})_{\ms{H}}$. We have an equality of closed formal subschemes
    \begin{equation}
    \mc{Z}(L^{\flat})_{\ms{H}} = \bigcup_{\substack{L^{\flat} \subseteq M^{\flat} \subseteq M^{\flat *} \\ t(M^{\flat}) \leq 1}} \mc{Z}(M^{\flat})^{\qcancirc} 
    \end{equation}
in $\mc{N}(n - 1, 1)$, where the union runs over full rank lattices $M^{\flat} \subseteq L^{\flat}_F$.
\end{proposition}
\begin{proof}
The union is the scheme-theoretic union (i.e. intersect associated ideal sheaves).

The inert case is \cite[{Theorem 4.2.1}]{LZ22unitary}. The ramified case is \cite[{Lemma 2.54}]{LL22II} (if $F / \Q_p$ is ramified, the condition $t(M^{\flat}) \leq 1$ implies $t(M^{\flat}) = 1$ since we have assumed $n$ is even in the ramified case).

For the split case, the inclusion $\subseteq$ is Corollary \ref{corollary:split:qcan_cycle_inclusion}. By Lemma \ref{lemma:split:full_rank_qcan_factor_through}, the inclusion $\supseteq$ will hold if we can verify that $\mc{Z}(L^{\flat})_{\ms{H}} \cong \Spf R$ for some finite flat $\mc{O}_{\breve{F}}$-algebra $R$ with $R \otimes_{\mc{O}_{\breve{F}}} \breve{F}$ reduced (with $R$ not necessarily a domain). We will check this later by passing to global special cycles via uniformization (\crefext{III:lemma:non-Arch_uniformization:global_to_local:horizontal_description}).
\end{proof}

For readers interested in Kr\"amer integral models for $F / \Q_p$ ramified, we mention the analogous \cite[{Theorem 4.2}]{HSY22a}, which we will not need.

    \clearpage


    \part{Local main theorems}
    \label{part:local_main_theorems}

        \section{Sketch of limit argument}
        \label{sec:part_I:sketch}
            For illustration purposes, we further sketch the local limiting strategy described in \cref{ssec:intro:strategy_overview}, particularly \cref{figure:intro:strategy:local_limit_method}. In \cref{sec:part_I:sketch}, we let $F / \Q$ be an imaginary quadratic field, and allow $n$ to be even or odd.

The limit strategy at an odd inert prime is outlined in \cref{ssec:part_I:sketch:non-Archimedean}. For comparison, we also include a sketch of the analogous Archimedean limit strategy in \cref{ssec:part_I:sketch:Archimedean}. We hope that the similarities between the Archimedean and non-Archimedean cases are visible from the sketches below. We treat our non-Archimedean local theorem in detail (via limiting) in \cref{sec:non-Arch_identity} (with inert/ramified/split done essentially simultaneously). The proof of the Archimedean limit will require separate arguments, which are the subject of our companion paper \cite{corank1_ASW_II.pdf}.

            \subsection{Archimedean}
            \label{ssec:part_I:sketch:Archimedean}
                For purposes of exposition in Section \ref{ssec:part_I:sketch:Archimedean}, we consider $T^{\flat} \in \mrm{Herm}_n(\R)$ and $t \in \R_{<0}$. Set $T \coloneqq \mrm{diag}(t, T^{\flat})$. Let $V$ be any signature $(n - 1, 1)$ non-degenerate $F / \Q$ Hermitian space with pairing denoted $(-,-)$ (for any $n \geq 1$).

We prove the limiting identity (left vertical arrow in Figure \ref{figure:intro:strategy:local_limit_method})
    \begin{equation}\label{equation:part_I:sketch:Archimedean:Whittaker_limit}
    \frac{d}{ds} \bigg |_{s = - 1/2} W^*_{T^{\flat}, \infty}(s)^{\circ}_n = \lim_{t \ra 0^{-}} \left ( \frac{d}{ds} \bigg|_{s = 0} W^*_{T,\infty}(s)^{\circ}_n + (\log |t|_{\infty} + \log(4\pi e^{\gamma})) W^*_{T^{\flat},\infty}(-1/2)^{\circ}_n \right )
    \end{equation}
where $\gamma$ is the Euler--Mascheroni constant, and $|-|_{\infty}$ denotes the usual real norm. This formula appears in our companion paper (more generally) as \crefext{II:proposition:local_identities:Archimedean_identity:statement:limiting}. The proof of this limiting formula is the bulk of the work at the Archimedean place. Note the similarity with the non-Archimedean version \cref{equation:part_I:sketch:non-Archimedean:inert_Whittaker_limit} (see also \cref{ssec:part_I:local_Whittaker:limits} for more comparisons).

On the geometric side of Figure \ref{figure:intro:strategy:local_limit_method}, the local $0$-cycles (resp. $1$-cycles) should be interpreted as Green currents of top degree $(n - 1, n - 1)$ (resp. degree $(n - 2, n - 2)$) on the associated Hermitian symmetric domain $\mc{D}$ parameterizing maximal negative definite $\C$-linear subspaces of $V_{\R}$. Consider the signature $(n - 1, 1)$ complex Hermitian space $V_{\R}$, with Hermitian pairing $(-,-)$. Any tuple $\underline{x} \in V_{\R}$ with nonsingular Gram matrix has an associated Kudla Green current $[\xi(\underline{x})]$, studied by Liu \cite{Liu11} in the unitary case. There is a certain local special cycle $\mc{D}(\underline{x}) \subseteq \mc{D}$ (complex submanifold which is the locus of $\C$-lines $z \in \mc{D}$ which are perpendicular to all elements of the tuple $\underline{x}$), arising in the complex uniformization of global special cycles.

Let $\underline{x}^{\flat} = [x_1^{\flat}, \ldots, x_{n - 1}^{\flat}] \in V_{\R}^{n - 1}$ be a tuple with Gram matrix $T^{\flat}$ and consider nonzero $x \in V_{\R} \in \mrm{span}_{\C}(\underline{x}^{\flat})^{\perp}$ in the orthogonal complement. Set
    \begin{equation}
    \underline{x} = [x, x_1^{\flat}, \ldots, x_{n - 1}^{\flat}] \quad \quad t \coloneqq (x,x) \quad \quad T \coloneqq \mrm{diag}(t, T^{\flat}).
    \end{equation}
Liu's Archimedean local theorem \cite[{Theorem 4.1.7}]{Liu11} implies
    \begin{equation}
    \int_{\mc{D}} [\xi(\underline{x})] = \frac{d}{d s} \bigg|_{s = 0} W^*_{T,\infty}(s)^{\circ}_n.
    \end{equation}
We are using the star product construction of $[\xi(\underline{x})]$, which unfolds as
    \begin{equation}\label{equation:intro:strategy:star_product}
    [\xi(\underline{x})] = [\xi(x)] * [\xi(\underline{x}^{\flat})] = \omega(x) \wedge [\xi(\underline{x}^{\flat})] + [\xi(x)] \wedge \delta_{\mc{D}(\underline{x}^{\flat})}
    \end{equation}
where $\omega(x)$ is a $(1,1)$-form associated with $x$ (Kudla--Millson form up to a normalization), $\delta_{\mc{D}(\underline{x}^{\flat})}$ is a Dirac delta current, and $\xi(x)$ is a certain function on $\mc{D}$ with logarithmic singularity along $\mc{D}(x)$. The function $\xi(x)$ is expressed in terms of the exponential integral $\Ei$. We have $\int_{\mc{D}} [\xi(x)] \wedge \delta_{\mc{D}(\underline{x}^{\flat})} = - \Ei(4 \pi t)$ and the limit formulas
    \begin{equation}
    \lim_{x \ra 0} \omega(x) = c_1(\widehat{\mc{E}}^{\vee}) \quad \quad \lim_{u \ra 0^-} (\Ei(u) - \log|u|) = \gamma
    \end{equation}
where $c_1(\widehat{\mc{E}}^{\vee})$ denotes the Chern form of dual tautological bundle on $\mc{D}$ (as in \crefext{II:sec:Hermitian_domain}). Under the assumption that $T^{\flat}$ is positive definite, we have $W^*_{T^{\flat},\infty}(-1/2)^{\circ}_n = \deg \mc{D}(\underline{x}^{\flat}) = 1$ (``local geometric Siegel--Weil'', i.e. $\mc{D}(\underline{x}^{\flat})$ is a single point). We thus have
    \begin{equation}
    \int_{\mc{D}} c_1(\widehat{\mc{E}}^{\vee}) \wedge [\xi(\underline{x}^{\flat})] = \lim_{x \ra 0} \left ( \left ( \int_{\mc{D}} [\xi(\underline{x})] \right ) + \log|t| + \log(4 \pi e^{\gamma}) \right ).
    \end{equation}
Using the functional equation $W^*_{T^{\flat},\infty}(s)^{\circ}_n = W^*_{T^{\flat},\infty}(-s)^{\circ}_n$, the limit formula in \eqref{equation:part_I:sketch:Archimedean:Whittaker_limit} now implies the following theorem.
\begin{theorem*}[Archimedean local version of \cref{theorem:intro:results:main}]
We have
    \begin{equation}\label{equation:intro:strategy:Arch_local_theorem}
    - \frac{d}{ds} \bigg |_{s = 1/2} W^*_{T^{\flat}, \infty}(s)^{\circ}_n = \int_{\mc{D}} c_1(\widehat{\mc{E}}^{\vee}) \wedge [\xi(\underline{x}^{\flat})].
    \end{equation}
\end{theorem*}
This is our main local Archimedean theorem for positive definite $T^{\flat}$ (i.e. the dotted arrow in Figure \ref{figure:intro:strategy:local_limit_method}). Its proof will appear in our companion paper \crefext{II:theorem:local_identities:Archimedean_identity:statement:main_Archimedean}, which also includes a version for non-positive definite $T^{\flat}$. Limiting on the geometric side of Figure \ref{figure:intro:strategy:local_limit_method} was provided by a limiting property of the (normalized) Kudla--Millson form, i.e. $\omega(x) \ra c_1(\widehat{\mc{E}}^{\vee})$ as $x \ra 0$. To compare with the limit of Whittaker function derivatives, we used the special value formula $W^*_{T^{\flat},\infty}(-1/2)^{\circ}_n = 1$ (``local geometric Siegel--Weil'') and the asymptotics of $\Ei$.
            
            \subsection{Non-Archimedean}
            \label{ssec:part_I:sketch:non-Archimedean}
                Let $p$ be an odd prime which is inert in $\mc{O}_F$. For our non-Archimedean main local theorems, we run an argument similar (in spirit) to the Archimedean case (\cref{ssec:part_I:sketch:Archimedean}) where the star product of Green currents is replaced by a derived tensor product of complexes of coherent sheaves on Rapoport--Zink spaces.

Suppose $T^{\flat} \in \mrm{Herm}_{n-1}(\Q_p)$ with $\det T^{\flat} \neq 0$, and consider $t \in \Q_p$ such that $T \coloneqq \mrm{diag}(t, T^{\flat})$ has $\varepsilon(T) = -1$ (defines a nonsplit Hermitian space). For the normalized local Whittaker functions defined in \cref{sec:part_I:local_Whittaker} below, we will prove the limit formula
    \begin{equation}\label{equation:part_I:sketch:non-Archimedean:inert_Whittaker_limit}
        \frac{d}{ds} \bigg|_{s = -1/2} W^*_{T^{\flat},p}(s)^{\circ}_n = \lim_{t \ra 0} \left ( \frac{d}{ds} \bigg|_{s = 0} W^*_{T,p}(s)^{\circ}_n + (\log |t|_{p} - \log p) W^*_{T^{\flat},p}(-1/2)^{\circ}_n \right ).
        \end{equation}
This appears as Proposition \ref{proposition:non-Arch_identity:limits} below. Note the similarity with the Archimedean version \cref{equation:part_I:sketch:Archimedean:Whittaker_limit} (see also \cref{ssec:part_I:local_Whittaker:limits} for comparisons with the ramified and split versions).

In \cref{ssec:part_I:sketch:non-Archimedean}, we set $\mc{N} \coloneqq \mc{N}(n - 1, 1)$ for the Rapoport--Zink space $\mc{N}(n - 1, 1)$ from \cref{sec:moduli_pDiv}. We also use the notation on special cycles from loc. cit., e.g. $\mbf{V}$ is the space of local special quasi-homomorphisms. 

If $\underline{\mbf{x}}$ is a basis for $\mbf{V}$, then $\mc{Z}(\underline{\mbf{x}})$ is a scheme with structure morphism $\mc{Z}(\underline{\mbf{x}}) \ra \Spf \breve{\Z}_p$ which is adic and proper \cite[{Lemma 2.10.1}]{LZ22unitary}. In this case, $\mc{Z}(\underline{\mbf{x}})$ is thus a finite order thickening of its special fiber $\mc{Z}(\underline{\mbf{x}})_{\overline{\F}_p}$, and there is a degree map $\deg_{\overline{\F}_p} \colon \mrm{gr}_0 K'_0(\mc{Z}(\underline{\mbf{x}}))_{\Q} \ra \Q$ given by the composite
    \begin{equation}
    \mrm{gr}_0 K'_0(\mc{Z}(\underline{\mbf{x}}))_{\Q} \xra{\sim} \mrm{gr}_0 K'_0(\mc{Z}(\underline{\mbf{x}})_{\overline{\F}_p})_{\Q} \ra \mrm{gr}_0 K'_0(\Spec \overline{\F}_p)_{\Q} = \Q
    \end{equation}
where the first arrow is induced by the d\'evissage pushforward isomorphism $K'_0(\mc{Z}(\underline{\mbf{x}})_{\overline{\F}_p}) \ra K'_0(\mc{Z}(\underline{\mbf{x}}))$ and the second arrow is pushforward along $\mc{Z}(\underline{\mbf{x}})_{\overline{\F}_p} \ra \Spec \overline{\F}_p$ (e.g. induced by taking Euler characteristics of coherent sheaves on $\mc{Z}(\underline{\mbf{x}})_{\overline{\F}_p}$).

Let $\underline{\mbf{x}}^{\flat} = [\mbf{x}_1^{\flat}, \ldots, \mbf{x}_{n - 1}^{\flat}] \in \mbf{V}^{n - 1}$ be a tuple with Gram matrix $T^{\flat}$ and consider nonzero $\mbf{x} \in \mbf{V} \in \mrm{span}_{F_p}(\underline{\mbf{x}}^{\flat})^{\perp}$ in the orthogonal complement. Set
    \begin{equation}
    \underline{\mbf{x}} = [\mbf{x}, \mbf{x}_1^{\flat}, \ldots, \mbf{x}_{n - 1}^{\flat}] \quad \quad t \coloneqq (\mbf{x}, \mbf{x}) \quad \quad T \coloneqq \mrm{diag}(t, T^{\flat}).
    \end{equation}
Li--Zhang's inert Kudla--Rapoport theorem \cite[{Theorem 1.2.1}]{LZ22unitary} implies
    \begin{equation}
    (\deg_{\overline{\F}_p} {}^{\mbb{L}} \mc{Z}(\underline{\mbf{x}})) \cdot \log p = \frac{1}{2} \frac{d}{ds} \bigg|_{s = 0} W^*_{T,p}(s)^{\circ}_n.
    \end{equation}
As an element of $\mrm{gr}^n_{\mc{N}} K'_0(\mc{Z}(\underline{\mbf{x}}))_{\Q}$, the derived tensor product unfolds as
    \begin{equation}
    {}^{\mbb{L}} \mc{Z}(\underline{\mbf{x}}) = [{}^{\mbb{L}} \mc{Z}(\mbf{x}) \otimes^{\mbb{L}}_{\mc{O}_{\mc{N}}} {}^{\mbb{L}} \mc{Z}(\underline{\mbf{x}}^{\flat})] = [{}^{\mbb{L}} \mc{Z}(\mbf{x}) \otimes^{\mbb{L}}_{\mc{O}_{\mc{N}}} {}^{\mbb{L}} \mc{Z}(\underline{\mbf{x}}^{\flat})_{\ms{V}}] + [{}^{\mbb{L}} \mc{Z}(\mbf{x}) \otimes^{\mbb{L}}_{\mc{O}_{\mc{N}}} {}^{\mbb{L}} \mc{Z}(\underline{\mbf{x}}^{\flat})_{\ms{H}}].
    \end{equation}
Here ${}^{\mbb{L}} \mc{Z}(\underline{\mbf{x}}^{\flat})_{\ms{H}} = [\mc{O}_{\mc{Z}(\underline{\mbf{x}}^{\flat})_{\ms{H}}}] \in \mrm{gr}_{\mc{N}}^{n - 1} K'_0(\mc{Z}(\underline{\mbf{x}}^{\flat})_{\ms{H}})_{\Q}$ is the ``horizontal part'' of ${}^{\mbb{L}} \mc{Z}(\underline{\mbf{x}}^{\flat})$, with $\mc{Z}(\underline{\mbf{x}}^{\flat})_{\ms{H}} \subseteq \mc{Z}(\underline{\mbf{x}}^{\flat})$ denoting the flat part, and ${}^{\mbb{L}} \mc{Z}(\underline{\mbf{x}}^{\flat})_{\ms{V}} \in \mrm{gr}_{\mc{N}}^{n - 1} K'_0(\mc{Z}(\underline{\mbf{x}}^{\flat})_{\overline{\F}_p})_{\Q}$ is the ``vertical part'' of ${}^{\mbb{L}} \mc{Z}(\underline{\mbf{x}}^{\flat})$ (see \cite[{\S 5.2}]{LZ22unitary}; we are using the d\'evissage pushforward isomorphism $K'_0(\mc{Z}(\underline{\mbf{x}}^{\flat})_{\overline{\F}_p}) \ra K'_0(\mc{Z}(\underline{\mbf{x}}^{\flat})_{\ms{V}})$ where $\mc{Z}(\underline{\mbf{x}}^{\flat})_{\ms{V}}$ is the vertical part from loc. cit.).

We show the limit formulas
    \begin{align}
    &\lim_{\mbf{x} \ra 0} (\deg_{\overline{\F}_p} [{}^{\mbb{L}} \mc{Z}(\mbf{x}) \otimes^{\mbb{L}}_{\mc{O}_{\mc{N}}} {}^{\mbb{L}} \mc{Z}(\underline{\mbf{x}}^{\flat})_{\ms{V}}]) = \deg_{\overline{\F}_p}(\mc{E}^{\vee} \cdot {}^{\mbb{L}} \mc{Z}(\underline{\mbf{x}}^{\flat})_{\ms{V}}) \label{equation:intro_strategy:non-Arch_limit:vertical}
    \\
    & \lim_{\mbf{x} \ra 0} \left ( ( \deg_{\overline{\F}_p}[{}^{\mbb{L}} \mc{Z}(\mbf{x}) \otimes^{\mbb{L}}_{\mc{O}_{\mc{N}}} {}^{\mbb{L}} \mc{Z}(\underline{\mbf{x}}^{\flat})_{\ms{H}}]) \cdot \log p- \frac{1}{2} (\log |t|_p - \log p) \cdot  \deg(\mc{Z}(\underline{\mbf{x}}^{\flat})_{\ms{H}}) \right ) \label{equation:intro_strategy:non-Arch_limit:horizontal}
    \\
    & ~ = \sum_{\mc{Z} \hookrightarrow \mc{Z}(\underline{\mbf{x}}^{\flat})_{\ms{H}}} \deg(\mc{Z}) \cdot \delta_{\mrm{tau}}(\mc{Z}) \cdot \log p.\notag
    \end{align}
Here $\mc{E}^{\vee}$ is the certain dual tautological bundle on $\mc{N}$, the sum runs over components $\mc{Z}$ of (the finite scheme associated to) $\mc{Z}(\underline{\mbf{x}}^{\flat})_{\ms{H}}$, the notation $\deg(\mc{Z}(\underline{\mbf{x}}^{\flat})_{\ms{H}})$ means the degree of the adic finite flat morphism $\mc{Z}(\underline{\mbf{x}}^{\flat})_{\ms{H}} \ra \Spf \breve{\Z}_p$ (and similarly for $\deg (\mc{Z})$), and $\delta_{\mrm{tau}}(\mc{Z}) \in \Q$ is the appropriate ``local change of tautological height''. As in \cref{ssec:can_and_qcan:qcan_cycles}, each $\mc{Z}$ is associated with a quasi-canonical lifting of some level $s$, and our notation $\delta_{\mrm{tau}}(\mc{Z})$ here is the $\delta_{\mrm{tau}}(s)$ in \cref{equation:can_and_qcan:qcan:local_change_taut}.

The quantity $\delta_{\mrm{tau}}(\mc{Z})$ (and the closely related ``local change of Faltings height'' $\delta_{\mrm{Fal}}$) arise from our reduction process from mixed characteristic heights to local quantities. This local-to-global reduction is begun in \cref{part:local_change_heights} and completed in our companion paper \cite{corank1_ASW_III.pdf}.

The ``vertical'' limit formula in \eqref{equation:intro_strategy:non-Arch_limit:vertical} follows from a Grothendieck--Messing theory argument (such vertical limiting behavior was observed in the inert case by \cite{LZ22unitary} via computation, and later in the ramified case by \cite{LL22II} via a linear-invariance argument), see Lemma \ref{lemma:non-Arch_identity:limits:vertical_geometric}. We prove the ``horizontal'' limit formula in \eqref{equation:intro_strategy:non-Arch_limit:horizontal} componentwise, i.e. we prove a refined limiting formula for each component $\mc{Z} \hookrightarrow \mc{Z}(\underline{\mbf{x}}^{\flat})_{\ms{H}}$ (Remark \ref{remark:non-Arch_identity:limits:horizontal_geometric}). Each $\mc{Z}$ embeds into a smaller Rapoport--Zink space of dimension $2$, where we make a computation in terms of quasi-canonical liftings.

We have the formula $W^*_{T^{\flat},p}(-1/2)^{\circ}_n = \deg(\mc{Z}(\underline{\mbf{x}}^{\flat})_{\ms{H}} \ra \Spf \breve{\Z}_p)$ (``local geometric Siegel--Weil''; right-hand side denotes degree of the indicated adic finite flat morphism), see Lemma \ref{lemma:non-Arch_identity:statement:special_value} (observed in the inert case by Li--Zhang \cite[{Corollary 4.6.1}]{LZ22unitary}).
Using the functional equation $W^*_{T^{\flat},p}(s)^{\circ}_n = W^*_{T^{\flat},p}(-s)^{\circ}_n$, the limit formula in \eqref{equation:intro:strategy:inert:Whittaker_limit} now implies
    \begin{equation}\label{equation:part_I:sketch:inert_local_theorem}
    - \frac{d}{ds} \bigg|_{s = 1/2} W^*_{T^{\flat},p}(s)^{\circ}_n 
    = \left ( 2 \deg_{\overline{\F}_p}(\mc{E}^{\vee} \cdot {}^{\mbb{L}} \mc{Z}(\underline{\mbf{x}}^{\flat})_{\ms{V}}) 
    + 2 \sum_{\mc{Z} \hookrightarrow \mc{Z}(\underline{\mbf{x}}^{\flat})_{\ms{H}}} \deg(\mc{Z}) \cdot \delta_{\mrm{tau}}(\mc{Z}) \right) \cdot \log p..
    \end{equation}
This is our main non-Archimedean local theorem for odd inert $p$ (i.e. the dotted arrow in Figure \ref{figure:intro:strategy:local_limit_method}), and appeared previously in \cref{ssec:intro:outline}. The analogous statement treating inert/split/ramified simultaneously is \cref{theorem:non-Arch_identity:statement} in this paper (there stated in terms of local densities). Limiting on the geometric side of Figure \ref{figure:intro:strategy:local_limit_method} was provided by the formulas in \eqref{equation:intro_strategy:non-Arch_limit:vertical} and \eqref{equation:intro_strategy:non-Arch_limit:horizontal}.

        \section{Local Whittaker functions}
        \label{sec:part_I:local_Whittaker}
            We will treat Eisenstein series and local Whitaker functions in greater detail in \crefext{IV:part:part_IV:Eisenstein}. Here in \cref{sec:part_I:local_Whittaker}, we fix some notation and state key properties (e.g. functional equations) which we need as ingredients for proving our main local theorems on arithmetic Siegel--Weil in the non-Archimedean cases (\cref{sec:non-Arch_identity}).

Let $F$ be a degree $2$ finite \'etale algebra over a non-Archimedean local field $F_0$. Write $\mc{O}_F$ and $\mc{O}_{F_0}$ for the corresponding rings of integers. Let $p$ (resp. $q$) be the residue characteristic (resp. residue cardinality) of $\mc{O}_{F_0}$. We write $\eta \colon F_0^{\times} \ra \{ \pm 1 \}$ for the quadratic character associated to the extension $F / F_0$. Let $\Delta \subseteq \mc{O}_{F_0}$ be the discriminant ideal associated to $F / F_0$.
    
            \subsection{Additional notation on \texorpdfstring{$U(m,m)$}{U(m,m)}} 
            \label{ssec:part_I:local_Whittaker:group}
                Recall the group $U(m,m)$ defined in \eqref{equation:intro:results:U(m,m)}. There we defined $U(m,m)$ with respect to the degree $2$ finite flat map $\Z \ra \mc{O}_F$, but we can make the same definition with respect to any finite locally free morphism $A \ra B$ of commutative rings such that $B$ is equipped with an involution $b \mapsto \overline{b}$ over $A$. In Section \ref{ssec:part_I:local_Whittaker:group}, we are interested in the case where $A \ra B$ is $\mc{O}_{F_0} \ra \mc{O}_F$, and sometimes use the alternative notation $H \coloneqq U(m,m)$.

We consider
    \begin{align}
    P &\coloneqq  \left \{ h = \begin{pmatrix} * & * \\ 0 & * \end{pmatrix} \in H \right \}  && M \coloneqq  \left \{ m(a) = \begin{pmatrix} a & 0 \\ 0 & {}^t \overline{a}^{-1} \end{pmatrix} : a \in \Res_{B/A} \GL_m \right \} \\
    w & \coloneqq \begin{pmatrix} 0 & 1_m \\ -1_m & 0 \end{pmatrix} && N \coloneqq  \left \{ n(b) = \begin{pmatrix} 1_m & b \\ 0 & 1_m \end{pmatrix} : b \in \mrm{Herm}_m  \right \} 
    \end{align}
where $P$, $M$, and $N$ are subgroups of $H$. We have $P(R) = M(R) N(R)$ for all $A$-algebras $R$.
We consider the standard open compact subgroup
    \begin{equation}
    K^{\circ} \coloneqq H(\mc{O}_{F_0}) \subseteq H(F_0)
    \end{equation}
and have $H(F_0) = P(F_0) K^{\circ}$. 
    
            \subsection{Principal series}
            \label{ssec:part_I:local_Whittaker:principal_series}
                Characters are assumed continuous and unitary unless specified otherwise. Form the group $H = U(m,m)$ as in Section \ref{ssec:part_I:local_Whittaker:group}.
Given $s \in \C$ and a character $\chi \colon F^{\times} \ra \C^{\times}$, we form the local \emph{degenerate principal series}
    \begin{equation}
    I(s,\chi) \coloneqq \Ind_{P(F_0)}^{H(F_0)}(\chi \cdot |-|_{F}^{s+m/2}).
    \end{equation}
This is an unnormalized induction, consisting of smooth and $K^{\circ}$-finite functions $\Phi \colon H(F_0) \ra \C$ satisfying
    \begin{equation}
    \Phi(m(a) n(b) h, s) = \chi(a) |\det a|_F^{s + m/2}
    \end{equation}
for all $m(a) \in M(F_0)$ and $n(b) \in N(F_0)$ and $h \in H(F_0)$. Here we wrote $\chi(a) \coloneqq \chi(\det a)$ for short. A section $\Phi(h,s)$ of $I(s, \chi)$ is \emph{standard} if $\Phi(k,s)$ is independent of $s$ for any fixed $k \in K^{\circ}$. We say that $\Phi$ is \emph{spherical} if $\Phi(hk,s) = \Phi(h,s)$ for any $k \in K^{\circ}$. We write $\Phi^{\circ}$ for the \emph{normalized spherical (standard) section}, i.e. the spherical standard section satisfying $\Phi^{\circ}(1,s) = 1$ for all $s \in \C$.
    
            \subsection{Normalized local Whittaker functions}
            \label{ssec:part_I:local_Whittaker:normalized_Whittaker}
                Now assume $F_0$ has characteristic $0$. Let $\chi \colon F^{\times} \ra \C^{\times}$ and $\psi \colon F_0 \ra \C^{\times}$ be characters with $\psi$ nontrivial. Given a standard section $\Phi \in I(s, \chi)$, and given $T \in \mrm{Herm}_m(F_0)$ with $\det T \neq 0$, there is a \emph{local Whittaker function} defined by the absolutely convergent integral
    \begin{equation}
    W_{T}(h,s,\Phi) \coloneqq \int_{N(F_0)} \Phi(w^{-1} n(b) h, s) \psi(- \mrm{tr}(Tb)) ~dn(b)
    \end{equation}
for $h \in H(F_0)$ and $s \in \C$ with $\Re(s) > m / 2$, where $dn(b)$ is the Haar measure with respect to the pairing $b,b' \mapsto \psi(\mrm{tr}(b b'))$ on $\mrm{Herm}_m(F_0) \cong N(F_0)$. For any fixed $h$, the function $W_T(h,s,\Phi)$ admits holomorphic admits holomorphic continuation to $s \in \C$ \cite[{Corollary 3.6.1}]{Karel79}\cite{KS97}. Extending linearly defines $W_{T}(h, s, \Phi)$ whenever $\Phi$ is a finite meromorphic linear combination of standard sections. 

Next, fix $n \in \Z$ and assume $\chi|_{F_0^{\times}} = \eta^n$. Set $s_0 \coloneqq (n - m)/2$. We also assume $\chi$ is unramified. Note that $n$ must then be even if $F / F_0$ is ramified.

We define a \emph{normalized local Whittaker function}
    \begin{align}
    W_{T}^*(h, s)^{\circ}_n \coloneqq \Lambda_{T}(s)^{\circ}_n W_{T}(h, s, \Phi^{\circ})
    \end{align}
using the a \emph{local normalizing factor}, which we define as
    \begin{align}\label{equation:part_I:Eisenstein:local_Whittaker:non-Arch:normalizing_factor}
    \Lambda_{T}(s)^{\circ}_n & \coloneqq |\Delta|_{F_0}^{- m(m-1)/4} \left ( \prod_{j = 0}^{m - 1} L(2s + m - j, \eta^{j + n})  \right ) |(\det T)\Delta^{\lfloor m/2 \rfloor}|_{F_0}^{-s - s_0}.
    \end{align}
For this particular (shifted) product of $L$-factors, compare \cite[{\S 6}]{HKS96}. The local $L$-factors $L(s,\eta^j)$ are the standard ones as in \cite[{\S 3}]{Tate79}); for the reader's convenience, we recall
    \begin{align*}
    L(s, \xi) & = 
    \begin{cases}    
    (1 - \xi(\varpi_0) q^{-s})^{-1} & \text{if $\xi$ is unramified} \\
    1 & \text{if $\xi$ is ramified}
    \end{cases}
    \end{align*}
for any quasi-character $\xi$ of $F_0$, where $\varpi_0$ denotes a uniformizer of $\mc{O}_{F_0}$.

We have a functional equation
    \begin{align}
    \label{equation:part_I:local_functional_equations:non-Arch:Whittaker_star}
    W_{T}^*(h, s)^{\circ}_n & = \eta((-1)^{m(m-1)/2} \det T)^{n - m - 1} W^*_T(h, s)^{\circ}_n, 
    \end{align}
see \crefext{IV:ssec:Eisenstein:local_functional_equations:non-Arch}.

We use the notation
    \begin{align}\label{equation:part_I:local_Whittaker:normalized_Whittaker:only_s_var}
    W^*_T(s)^{\circ}_n \coloneqq W^*_T(1_{2m},s)^{\circ}_n.
    \end{align}

If $F / F_0$ is ramified, we assume for the rest of \cref{sec:part_I:local_Whittaker} that $\mc{O}_{F_0}$ has odd residue characteristic. There is a \emph{normalized local density} (Laurent) polynomial $\mrm{Den}^*(X,T)_n \in \Z[1/q][X, X^{-1/2}]$ such that
    \begin{equation}\label{equation:part_I:local_Whittaker:local_densities_spherical:Den_star}
    W^*_T(s + s_0)^{\circ}_n = \mrm{Den}^*(q^{-2s},T)_n
    \end{equation}
for all $s \in \C$. We also consider the variant $\mrm{Den}(X,T)_n \in \Z[1/q][X]$ given by
    \begin{equation}
    |(\det T)\Delta^{\lfloor m/2 \rfloor}|_{F_0}^{s + s_0} W^*_T(s + s_0)^{\circ}_n = \mrm{Den}(q^{-2s}, T)_n
    \end{equation}
for all $s \in \C$. We have
    \begin{equation}
    \label{equation:part_I:local_Whittaker:local_densities_spherical:comparison}
    \mrm{Den}^*(X,L)_n
    = 
    (q^{2 s_0} X^{- 1/2})^{\lfloor \mrm{val}(L) \rfloor} \mrm{Den}(X,L)_n
    \end{equation}
as well as
    \begin{align}
    \mrm{Den}(X, T)_{n+1} &= \mrm{Den}(q^{-1} X, T)_n && \text{if $F/F_0$ is split} \notag \\
    \mrm{Den}(X, T)_{n+1} &= \mrm{Den}(-q^{-1} X, T)_n && \text{if $F/F_0$ is inert} \label{equation:part_I:Eisenstein:local_Whittaker:local_densities_spherical:change_of_n} \\
    \mrm{Den}(X, T)_{n+2} &= \mrm{Den}(q^{-2} X, T)_n && \text{if $F/F_0$ is ramified}. \notag
    \end{align}

\begin{remark}\label{remark:part_I:local_Whittaker:local_densities_spherical:square_root_notation}
On the right-hand side of \eqref{equation:part_I:local_Whittaker:local_densities_spherical:Den_star}, we mean evaluating $\mrm{Den}^*(X,T)_n$ at $X^{1/2} = q^{-s}$.
We similarly abuse notation elsewhere. For example, $\mrm{Den}^*(q X, T)_n \in \Z[1 / q^{1/2}][X, X^{-1/2}]$ is obtained from $\mrm{Den}^*(X, T)_n$ by replacing $X^{1/2}$ with $q^{1/2} X^{1/2}$. The notation $\frac{d}{d X} \colon \Q[X, X^{-1/2}] \ra \Q[X, X^{-1/2}]$ means the $\Q$-linear map $X^{j/2} \mapsto (j/2) X^{j/2 - 1}$.
\end{remark}

If $L$ is a non-degenerate Hermitian $\mc{O}_F$-lattice which admits a basis with Gram matrix $T$, we set
    \begin{equation}
    \mrm{Den}^*(X,L)_n \coloneqq \mrm{Den}^*(X,T)_n \quad \quad \mrm{Den}(X,L)_n \coloneqq \mrm{Den}(X,T)_n.
    \end{equation}
The local densities satisfy a certain cancellation property (which we will use): if $L^{\circ}$ is a self-dual Hermitian lattice of rank $n$, then for any non-degenerate Hermitian lattice $L$ and every integer $r \in \Z$ (assume $r$ is even if $F / F_0$ is ramified), we have
    \begin{equation}\label{equation:part_I:local_Whittaker:local_densities_spherical:cancellation}
    \mrm{Den}(X,L \oplus L^{\circ})_{r + n} = \mrm{Den}(X,L)_{r} \quad \quad \mrm{Den}^*(X,L \oplus L^{\circ})_{r + n} = \mrm{Den}^*(X,L)_{r}
    \end{equation}
where $L \oplus L^{\circ}$ is the orthogonal direct sum. This follows from e.g. from ``Cho--Yamauchi'' formulas for local densities (as proved in \cite[{Theorem 3.5.1}]{LZ22unitary} (inert), \cite[{Theorem 2.2}]{FYZ22SW} (split and inert), and \cite[{Lema 2.15}]{LL22II} (ramified)) and the following linear algebra fact: every lattice $M' \subseteq (L \oplus L^{\circ}) \otimes_{\mc{O}_{F}} F$ satisfying $L^{\circ} \subseteq M' \subseteq M^{\prime *}$ admits an orthogonal direct sum decomposition $M' = L^{\circ} \oplus M''$ for some sublattice $M''$.

The functional equation in \eqref{equation:part_I:local_functional_equations:non-Arch:Whittaker_star} implies
    \begin{align}\label{equation:part_I:Eisenstein:local_functional_equations:non-Arch:non-star_density}
    \mrm{Den}(q^{2 s_0} X^{-1},L)_n & = \varepsilon(L)^{n - m - 1} X^{- \lfloor\mrm{val}(L)\rfloor} \mrm{Den}(q^{2 s_0} X,L)_n 
    \\
    \mrm{Den}^*(q^{2 s_0} X^{-1},L)_n & = \varepsilon(L)^{n - m - 1} \mrm{Den}^*(q^{2 s_0} X,L)_n.
    \label{equation:part_I:local_functional_equations:non-Arch:Den_star}
    \end{align}
Both $\varepsilon(L)$ and $\mrm{val}(L)$ were defined in Section \ref{ssec:Hermitian_conventions:lattices}. In the case where $\chi|_{F_0^{\times}}$ is trivial, these non-Archimedean functional equations are essentially \cite[{Corollary 3.2}]{Ikeda08}.

If $m = n$ and $F / F_0$ is unramified, our $\mrm{Den}(X,L)_n$ agrees with the polynomial $\mrm{Den}(X,L)$ from \cite[{\S 3}]{LZ22unitary} and \cite{FYZ22SW}. If $m = n$ and $F / F_0$ is ramified (so $n$ is even, by the assumption that $\chi$ is unramified and $\chi|_{F_0^{\times}} = \eta^n$), our $\mrm{Den}(X^2,L)_n$ is the $\mrm{Den}(X,L)$ of \cite[{\S 2}]{LL22II}.

For any $m$ and $n$, if $T$ defines a self-dual Hermitian lattice when $m$ is even and $F / F_0$ is unramified (resp. an almost self-dual Hermitian lattice when $m$ is odd and $F / F_0$ is ramified), we have
    \begin{equation}
    \label{equation:part_I:Eisenstein:local_Whittaker:local_densities_spherical:often_1}
    W^*_{T}(s)^{\circ}_n = 1 \quad \quad \mrm{Den}^*(X, T)_n = \mrm{Den}(X, T)_n = 1.
    \end{equation}

            \subsection{Limits}
            \label{ssec:part_I:local_Whittaker:limits}
                Take integers $m,n$ with $m \geq 1$, set $s_0 = (n - m)/2$ as above, and set $m^{\flat} = m - 1$.

We consider nonsingular $T \in \mrm{Herm}_m(F_0)$ of the form $T = \mrm{diag}(t, T^{\flat})$ where $T^{\flat} \in \mrm{Herm}_{m^{\flat}}(F_0)$ with $T^{\flat}$ nonsingular, and we study the local Whittaker function $W^*_{T}(s)^{\circ}_n$ as $t \ra 0$. The following limiting identities will be crucial for the proofs of our main local theorems. We collect them here for easier comparison between the inert/ramified/split and Archimedean cases. 

Proofs of the non-Archimedean cases will appear in Section \ref{sec:non-Arch_identity} below. The Archimedean case is (a special case of) one formulation of our main Archimedean local theorem, whose proof will appear in our companion paper \cite{corank1_ASW_II.pdf}.

The Archimedean case involves an analogous normalized Archimedean local Whittaker function $W^*_T(s)^{\circ}_n$. Since the present paper focuses on the non-Archimedean case, we have skipped the definition of $W^*_{T}(s)^{\circ}_n$ (it is defined in a similar way as \cref{ssec:part_I:local_Whittaker:normalized_Whittaker}, but for a certain scalar weight standard section of the principal series, and for a different normalizing factor involving Gamma functions, etc.). A definition may be found in \crefext{II:ssec:part_II:local_Whittaker_Archimedean:group} (also \crefext{IV:ssec:Eisenstein:local_Whittaker:Archimedean}).

Although the following limit formulas below look very similar in all cases, our proofs of these identities are fairly different in the non-Archimedean vs Archimedean cases.

If $F_0$ is non-Archimedean and $F/F_0$ is inert, Proposition \ref{proposition:non-Arch_identity:limits} implies
    \begin{equation}
    \frac{d}{ds} \bigg|_{s = -1/2} W^*_{T^{\flat}}(s)^{\circ}_n = \lim_{t \ra 0} \left ( \frac{d}{ds} \bigg|_{s = 0} W^*_{T}(s)^{\circ}_n + (\log |t|_{F_0} - \log q_v) W^*_{T^{\flat}}(-1/2)^{\circ}_n \right )
    \end{equation}
if the limit is taken over nonzero $t \in F_0$ with $\varepsilon(\mrm{diag}(t,T^{\flat})) = - 1$.

If $F_0$ is non-Archimedean and $F/F_0$ is split, Proposition \ref{proposition:non-Arch_identity:limits} implies
    \begin{equation}
    \frac{d}{ds} \bigg|_{s = -1/2} W^*_{T^{\flat}}(s)^{\circ}_n = \lim_{t \ra 0} \left ( \log q_v \cdot W^*_{T}(0)^{\circ}_n + (\log |t|_{F_0} - \log q_v) \cdot W^*_{T^{\flat}}(-1/2)^{\circ}_n \right )
    \end{equation}
if the limit is taken over nonzero $t \in F_0$.

If $F_0$ is non-Archimedean and $F/F_0$ is ramified, Proposition \ref{proposition:non-Arch_identity:limits} implies
    \begin{equation}
    2 \frac{d}{ds} \bigg|_{s = -1/2} W^*_{T^{\flat}}(s)^{\circ}_n = \lim_{t \ra 0} \left ( \frac{d}{ds} \bigg|_{s = 0} W^*_{T}(s)^{\circ}_n + (\log |t|_{F_0} - \log q_v) \cdot W^*_{T^{\flat}}(-1/2)^{\circ}_n \right )
    \end{equation}
if the limit is taken over nonzero $t \in F_0$ with $\varepsilon(\mrm{diag}(t,T^{\flat})) = -1$.

If $F/F_0$ is $\C / \R$, \crefext{II:proposition:local_identities:Archimedean_identity:statement:limiting} gives
    \begin{equation}
    \frac{d}{ds} \bigg |_{s = - 1/2} W^*_{T^{\flat}}(s)^{\circ}_n = \lim_{t \ra 0^{\pm}} \left ( \frac{d}{ds} \bigg|_{s = 0} W^*_{T}(s)^{\circ}_n + (\log |t|_{F_0} + \log(4\pi) - \Gamma'(1)) W^*_{T^{\flat}}(-1/2)^{\circ}_n \right )
    \end{equation}
where the sign on $0^{\pm}$ is $-$ (resp. $+$) if $T^{\flat}$ is positive definite (resp. not positive definite). If $T^{\flat} \in \mrm{Herm}_{m^{\flat}}(\R)$ is not positive definite, \crefext{II:proposition:local_identities:Archimedean_identity:statement:limiting} also proves a similar limiting statement for arbitrary $m^{\flat}$ (i.e. not necessarily $m^{\flat} = n - 1$).
    
        \section{Non-Archimedean local main theorems}
        \label{sec:non-Arch_identity}
            Let $F_0$ be a non-Archimedean local field of characteristic $0$, residue cardinality $q$, and residue characteristic $p$. Let $F$ be a finite \'etale $F_0$-algebra of degree $2$. We use notation $\breve{F}$ and $\breve{F}_0$ as in Sections \ref{sec:moduli_pDiv}--\ref{sec:can_and_qcan} (there with $F_0 = \Q_p$), so that $[\breve{F} : \breve{F}_0] = 1$ (resp. $[\breve{F} : \breve{F}_0] = 2$) if $F / F_0$ is unramified (resp. ramified).

Notation on Hermitian lattices from Section \ref{ssec:Hermitian_conventions:lattices} will be used freely.
For a non-degenerate Hermitian $\mc{O}_F$-lattice $L$, we use the shorthand $\mrm{val}'(L) \coloneqq \lfloor \mrm{val}(L) \rfloor \in \Z_{\geq 0}$, as well as $\mrm{val}'(x) \coloneqq \lfloor \mrm{val}(x) \rfloor$ for any $x \in L$ (i.e. $\mrm{val}'(L) = \mrm{val}(L) - 1/2$ if $F / F_0$ is ramified and $L$ has odd rank, and $\mrm{val}'(L) = \mrm{val}(L)$ otherwise). 
Fix an integer $n \geq 1$, and assume $n$ is even if $F / F_0$ is ramified.

If $F_0 = \Q_p$, we form the associated Rapoport--Zink space $\mc{N} \coloneqq \mc{N}(n - 1, 1)$ (Section \ref{ssec:moduli_pDiv:RZ}). Recall the space of local special quasi-homomorphisms $\mbf{W} \subseteq \mbf{V}$ (Section \ref{ssec:moduli_pDiv:special_cycles}). 
Recall that $\mbf{W}$ and $\mbf{V}$ are non-degenerate Hermitian $F$-modules of rank $n$ if $F / \Q_p$ is nonsplit (resp. rank $n - 1$ and rank $n$ is $F / \Q_p$ is split). Recall $\varepsilon(\mbf{V}) = -1$ if $F / \Q_p$ is nonsplit (resp. $\varepsilon(\mbf{V}) = 1$ if $F / \Q_p$ is split).

            \subsection{Statement of local main theorem}
            \label{ssec:non-Arch_identity:statement}
                We first define the geometric side of our main local identity, taking $F_0 = \Q_p$. We also assume $p \neq 2$ if $F / \Q_p$ is nonsplit.
Let $L^{\flat} \subseteq \mbf{W}$ be any non-degenerate Hermitian $\mc{O}_F$-lattice of rank $n - 1$. Form the associated local special cycle $\mc{Z}(L^{\flat}) \subseteq \mc{N}$. Recall that the flat part $\mc{Z}(L^{\flat})_{\ms{H}} \subseteq \mc{Z}(L^{\flat})$ decomposes into quasi-canonical lifting cycles $\mc{Z}(M^{\flat})^{\circ}$ for certain lattices $M^{\flat}$ (Proposition \ref{proposition:can_and_qcan:qcan_cycles}). Recall also the derived vertical local special cycle ${}^{\mbb{L}} \mc{Z}(L^{\flat})_{\ms{V}} \in \mrm{gr}^{n- 1}_{\mc{N}} K'_0(\mc{Z}(L^{\flat}))_{\Q}$ (Section \ref{ssec:moduli_pDiv:horizontal_vertical}).

\begin{definition}\label{definition:non-Arch_identity:statement:local_intersection_number}
Given a non-degenerate Hermitian $\mc{O}_F$-lattice $L^{\flat} \subseteq \mbf{W}$ of rank $n - 1$, the associated \emph{local intersection number} is
    \begin{equation}\label{equation:non-Arch_identity:statement:local_intersection_number:horizontal_vertical}
    \mrm{Int}(L^{\flat})_n \coloneqq \mrm{Int}_{\ms{H}}(L^{\flat})_n + \mrm{Int}_{\ms{V}}(L^{\flat})_n
    \end{equation}
where
    \begin{equation}\label{equation:non-Arch_identity:statement:local_intersection_number:horizontal_decomp}
    \mrm{Int}_{\ms{H}}(L^{\flat})_n \coloneqq \sum_{\substack{L^{\flat} \subseteq M^{\flat} \subseteq M^{\flat *} \\ t(M^{\flat}) \leq 1}} \mrm{Int}_{\ms{H}}(M^{\flat})^{\circ}_n
    \end{equation}
with the sum running over full rank lattices $M^{\flat} \subseteq L^{\flat}_F$, where 
    \begin{equation}
    \mrm{Int}_{\ms{H}}(M^{\flat})^{\circ}_n \coloneqq 2 \deg \mc{Z}(M^{\flat})^{\circ} \cdot \delta_{\mrm{tau}}(\mrm{val}'(M^{\flat}))
    \end{equation}
for any non-degenerate integral lattice $M^{\flat} \subseteq \mbf{W}$ with $t(M^{\flat}) \leq 1$, and where
    \begin{equation} \label{definition:non-Arch_identity:statement:local_intersection_number:vertical}
    \mrm{Int}_{\ms{V}}(L^{\flat})_n \coloneqq 2 [\breve{F} : \breve{\Q}_p]^{-1} \deg_{\overline{k}}({}^{\mbb{L}} \mc{Z}(L^{\flat})_{\ms{V}} \cdot \mc{E}^{\vee}).
    \end{equation}
\end{definition}

In our companion paper, we relate these local intersection numbers with global intersection numbers (end of \crefext{III:ssec:non-Arch_uniformization:vertical,III:ssec:non-Arch_uniformization:horizontal}).

The quantity
    \begin{equation}
    \deg \mc{Z}(M^{\flat})^{\circ}
    =
    \begin{cases}
    [\breve{F} : \breve{\Q}_p] p^{\mrm{val}'(M^{\flat})} (1 - \eta(p) p^{-1}) & \text{if $\mrm{val}'(M^{\flat}) \geq 1$} \\
    [\breve{F} : \breve{\Q}_p] & \text{if $\mrm{val}'(M^{\flat}) = 0$}
    \end{cases}
    \end{equation}
is the degree of the adic finite flat morphism $\mc{Z}(M^{\flat})^{\circ} \ra \Spf \mc{O}_{\breve{F}}$, where $\eta(p) \coloneqq -1, 0, 1$ in the inert, ramified, and split cases respectively (see \eqref{equation:can_and_qcan:qcan:degree}; the extra factor of $[\breve{F} : \breve{F}_0]$ accounts for the two components of $\mc{Z}(M^{\flat})^{\circ}$ when $F / F_0$ is ramified, see \eqref{equation:can_and_qcan:qcan_cycles:underlying_formal_scheme}). Recall that $\delta_{\mrm{tau}}(s)$ is the ``local change of tautological height'' defined in \eqref{equation:can_and_qcan:qcan:local_change_taut}, and recall that $\mc{E}^{\vee}$ is the dual tautological bundle on $\mc{N}$ (Definition \ref{definition:moduli_pDiv:RZ:tautological}). In \eqref{definition:non-Arch_identity:statement:local_intersection_number:vertical}, we understand $\mc{E}^{\vee} = [\mc{O}_{\mc{N}}] - [\mc{E}] \in K'_0(\mc{N})$ so that ${}^{\mbb{L}} \mc{Z}(L^{\flat})_{\ms{V}} \cdot \mc{E}^{\vee} \in F^n_{\mc{N}} K'_0 (\mc{Z}(L^{\flat})_{\overline{k}})_{\Q}$. For $L^{\flat}$ as above, recall that $\mc{Z}(L^{\flat})_{\overline{k}}$ is a scheme proper over $\Spec \overline{k}$ \crefext{III:lemma:non-Arch_uniformization:vertical:local_cycle_quasi-compact}, so there is a degree map $\deg_{\overline{k}} \colon F^n_{\mc{N}} K'_0(\mc{Z}(L^{\flat})_{\overline{k}}) \ra \Z$.

We refer to $\mrm{Int}_{\ms{H}}(L^{\flat})_n$ as the ``horizontal part'' of the local intersection number (coming from the flat part $\mc{Z}(L^{\flat})_{\ms{H}}$) and we refer to $\mrm{Int}_{\ms{V}}(L^{\flat})_n$ as the ``vertical part'' of the local intersection number (coming from $\mc{Z}(L^{\flat})_{\overline{k}}$, supported in positive characteristic).

We next define the automorphic side of our main local identity. For this, we allow $F_0$ to be an arbitrary finite extension of $\Q_p$ (allowing $p = 2$ if $F / F_0$ is unramified).
If $L^{\flat}$ is a non-degenerate Hermitian $\mc{O}_F$-lattice of rank $n - 1$, we set
    \begin{equation}
    \partial \mrm{Den}^*(L^{\flat})_n \coloneqq -2 [\breve{F} : \breve{F}_0] \frac{d}{d X}\bigg|_{X = 1} \mrm{Den}^*(q^2 X, L^{\flat})_n
    \end{equation}
where $\mrm{Den}^*(X,L^{\flat})_n \in \Z[1/q][X, X^{-1/2}]$ is a normalized local density \eqref{equation:part_I:local_Whittaker:local_densities_spherical:Den_star}. We are abusing notation as in Remark \ref{remark:part_I:local_Whittaker:local_densities_spherical:square_root_notation}, i.e. $\mrm{Den}^*(q^2 X, L^{\flat})_n$ means to evaluate $\mrm{Den}^*(X, L^{\flat})_n$ at $X^{1/2}$ being $q X^{1/2}$. We also set 
    \begin{equation}
    \mrm{Den}^*(L^{\flat})_n \coloneqq [\breve{F} : \breve{F}_0] \cdot \mrm{Den}^*(q^2, L^{\flat})_n.    
    \end{equation}

Suppose $M^{\flat}$ is a non-degenerate integral Hermitian $\mc{O}_F$-lattice of rank $n - 1$ with $t(M^{\flat}) \leq 1$. If $M^{\flat}$ is maximal integral,\footnote{The symbol $\circ$ indicates ``primitive'' here (for quasi-canonical lifting cycles), while $\circ$ indices ``spherical'' in Section \ref{ssec:part_I:local_Whittaker:normalized_Whittaker}. There is no notation clash as written, but we hope this remark helps to avoid confusion.} we set $\partial \mrm{Den}^*_{\ms{H}}(M^{\flat})^{\circ}_n \coloneqq \partial \mrm{Den}^*_{\ms{H}}(M^{\flat})_n$. Otherwise, we define $\partial \mrm{Den}^*_{\ms{H}}(M^{\flat})^{\circ}_n$ inductively so that the relation
    \begin{equation}\label{equation:non-Arch_identity:statement:horizontal_density_decomp}
    \partial \mrm{Den}^*(M^{\flat})_n = \sum_{M^{\flat} \subseteq N^{\flat} \subseteq N^{\flat * }} \partial \mrm{Den}^*_{\ms{H}}(N^{\flat})_n^{\circ}
    \end{equation}
is satisfied (induct on $\mrm{val}(M^{\flat})$), where the sum runs over lattices $N^{\flat} \subseteq M^{\flat}_F$. Given any non-degenerate integral Hermitian $\mc{O}_F$-lattice $L^{\flat}$ of rank $n - 1$, we then define $\partial \mrm{Den}^*_{\ms{V}}(L^{\flat})_n$ so that the relation
    \begin{equation}\label{equation:non-Arch_identity:statement:density_decomp}
    \partial \mrm{Den}^*(L^{\flat}) = \bigg( \sum_{\substack{L^{\flat} \subseteq M^{\flat} \subseteq M^{\flat *} \\ t(M^{\flat}) \leq 1}} \mrm{\partial} \mrm{Den}^*_{\ms{H}}(M^{\flat})_n^{\circ} \bigg ) + \partial \mrm{Den}^*_{\ms{V}}(L^{\flat})_n
    \end{equation}
is satisfied, where the sum runs over lattices $M^{\flat} \subseteq L^{\flat}_F$.

\begin{theorem}\label{theorem:non-Arch_identity:statement}
Suppose $F_0 = \Q_p$ and that $p \neq 2$ unless $F / \Q_p$ is split. For any non-degenerate Hermitian $\mc{O}_F$-lattice $L^{\flat} \subseteq \mbf{W}$ of rank $n - 1$, we have
    \begin{equation}\label{equation:theorem:non-Arch_identity:statement}
    \mrm{Int}(L^{\flat})_n = \partial \mrm{Den}^*(L^{\flat})_n.
    \end{equation}
Moreover, we have
    \begin{equation}\label{equation:theorem:non-Arch_identity:statement:refined}
    \mrm{Int}_{\ms{H}}(M^{\flat})^{\circ}_n = \partial \mrm{Den}^*_{\ms{H}}(M^{\flat})^{\circ}_n 
    \quad \quad 
    \mrm{Int}_{\ms{V}}(L^{\flat})_n = \partial \mrm{Den}^*_{\ms{V}}(L^{\flat})_n.
    \end{equation}
where $M^{\flat} \subseteq \mbf{W}$ is any non-degenerate integral Hermitian $\mc{O}_F$-lattice of rank $n - 1$ with $t(M^{\flat}) \leq 1$.
\end{theorem}

On account of the decompositions in \eqref{equation:non-Arch_identity:statement:local_intersection_number:horizontal_vertical}, \eqref{equation:non-Arch_identity:statement:local_intersection_number:horizontal_decomp}, and \eqref{equation:non-Arch_identity:statement:density_decomp}, it is clearly enough to prove the refined identities in \eqref{equation:theorem:non-Arch_identity:statement:refined}. 
The theorem is also clear if $L^{\flat}$ is not integral, since both sides of \eqref{equation:theorem:non-Arch_identity:statement} are zero in this case (the special cycle $\mc{Z}(L^{\flat})$ will be empty, and $\mrm{Den}(X,L^{\flat})_n$ will be identically zero as discussed in \crefext{IV:ssec:Eisenstein:local_Whittaker:local_densities}).

We also record a special value formula (as observed in the inert case by Li and Zhang \cite[{Corollary 4.6.1}]{LZ22unitary}) for later use. Its proof will appear in Section \ref{ssec:non-Arch_identity:horizontal_identity}.

\begin{lemma}\label{lemma:non-Arch_identity:statement:special_value}
Suppose $F_0 = \Q_p$ and that $p \neq 2$ unless $F /\Q_p$ is split. For any non-degenerate Hermitian $\mc{O}_F$-lattice $L^{\flat} \subseteq \mbf{W}$ of rank $n - 1$, we have
    \begin{equation}
    \deg \mc{Z}(L^{\flat})_{\ms{H}} = \mrm{Den}^*(L^{\flat})_n.
    \end{equation}
\end{lemma}

In the preceding lemma statement, $\deg \mc{Z}(L^{\flat})_{\ms{H}}$ means the degree of the adic finite flat morphism $\mc{Z}(L^{\flat})_{\ms{H}} \ra \Spf \mc{O}_{\breve{F}}$ of formal schemes.
    
            \subsection{Horizontal identity}
            \label{ssec:non-Arch_identity:horizontal_identity}
                We will need Cho--Yamauchi formulas for local densities (unitary version, as proved in \cite[{Theorem 3.5.1}]{LZ22unitary} (inert) \cite[{Theorem 2.2(3)}]{FYZ22SW} (split) \cite[{Lemma 2.15}]{LL22II} (ramified)). For this, we allow $F_0$ to be an arbitrary finite extension of $\Q_p$ (allowing $p = 2$ if $F / F_0$ is unramified). Then, if $L$ is any non-degenerate Hermitian $\mc{O}_F$-lattice of rank $n$ (still assuming $n$ even if $F / F_0$ is ramified), we have
    \begin{align}\label{equation:non-Arch_identity:horizontal_identity:Cho--Yamauchi}
    \mrm{Den}(X,L)_n & = \sum_{L \subseteq M \subseteq M^*} X^{\ell(M / L)} \mrm{Den}(X,M)_n^{\circ}
    \\
    \mrm{Den}(X,M)_n^{\circ} & \coloneqq \prod_{i = 0}^{t(M) - 1} (1 - \eta^i(\varpi_0) q^i X) \label{equation:non-Arch_identity:horizontal_identity:Cho--Yamauchi:"primitive"}
    \end{align}
where $\eta(\varpi_0) \coloneqq \eta^i(\varpi_0) \coloneqq -1, 0, 1$ if $i$ is odd (resp. $\eta^i(\varpi_0) \coloneqq 1$ if $i$ is even) in the inert, ramified, split cases respectively, and $\mrm{Den}(X,L)_n \in \Z[X]$ is the local density polynomial normalized as in Section \ref{ssec:part_I:local_Whittaker:normalized_Whittaker}. The displayed sum runs over lattices $M \subseteq L_F$.

Suppose $L^{\flat}$ is a Hermitian $\mc{O}_F$-lattice of rank $n - 1$ (still assuming $n$ even if $F / F_0$ is ramified). If $F /F_0$ is unramified, we have $\mrm{Den}(X, L^{\flat})_n = \mrm{Den}(\eta(\varpi_0) q^{-1} X, L^{\flat})_{n - 1}$ \eqref{equation:part_I:Eisenstein:local_Whittaker:local_densities_spherical:change_of_n} and we set $\mrm{Den}(X, L^{\flat})^{\circ}_n \coloneqq \mrm{Den}(\eta(\varpi_0) q^{-1}, L^{\flat})^{\circ}_{n - 1}$ if $L^{\flat}$ is also integral.

If $F / F_0$ is ramified, we have
    \begin{align}
    \mrm{Den}(X, L^{\flat})_n & = \sum_{L^{\flat} \subseteq M^{\flat} \subseteq M^{\flat *}} (q^{-1} X)^{\ell(M^{\flat}/L^{\flat})} \mrm{Den}(X, M^{\flat})^{\circ}_n
    \\
    \mrm{Den}(X, M^{\flat})^{\circ}_n & \coloneqq  \prod_{i = 0}^{\frac{t(M^{\flat})-3}{2}} (1 - q^{2i} X)
    \end{align}
where the sum runs over lattices $M^{\flat} \subseteq L^{\flat}_F$ (may be verified using \cite[{Lemma 2.15}]{LL22II}).

If $M^{\flat}$ is a non-degenerate integral Hermitian $\mc{O}_F$-lattice of rank $n - 1$ with $t(M^{\flat}) \leq 1$, set
    \begin{equation}
    \mrm{Den}^*(M^{\flat})^{\circ}_n \coloneqq
    [\breve{F} : \breve{F}_0] q^{\mrm{val}'(M^{\flat})} \mrm{Den}(1, M^{\flat})^{\circ}_n
    \end{equation}
We have
    \begin{align} \label{equation:non-Arch_identity:horizontal_identity:decomp}
    [\breve{F} : \breve{F}_0] \mrm{Den}^*(q^2, L^{\flat})_n = [\breve{F} : \breve{F}_0] \mrm{Den}^*(1, L^{\flat})_n = \sum_{L^{\flat} \subseteq M^{\flat} \subseteq M^{\flat *}} \mrm{Den}^*(M^{\flat})^{\circ}_n
    \end{align}
where the sum runs over lattices $M^{\flat} \subseteq L^{\flat}_F$. The first equality follows from the functional equation \eqref{equation:part_I:local_functional_equations:non-Arch:Den_star}, and the second equality follows from the Cho--Yamauchi formulas (and \eqref{equation:part_I:local_Whittaker:local_densities_spherical:comparison}). Note $\mrm{Den}^*(M^{\flat})_n = \mrm{Den}^*(M^{\flat})^{\circ}_n$ if $M^{\flat}$ is maximal integral.

\begin{proof}[Proof of Lemma \ref{lemma:non-Arch_identity:statement:special_value}]
Follows from \eqref{equation:non-Arch_identity:horizontal_identity:decomp}.
Note $\mrm{Den}^*(M^{\flat})^{\circ}_n = \deg \mc{Z}(M^{\flat})^{\circ}$ if $t(M^{\flat}) \leq 1$, and $\mrm{Den}^*(M^{\flat})^{\circ}_n = 0$ if $t(M^{\flat}) \geq 2$.
\end{proof}

\begin{proposition}\label{proposition:non-Arch_identity:horizontal_identity}
Assume $F_0 = \Q_p$ and that $p \neq 2$ unless $F / \Q_p$ is split. For any rank $n - 1$ non-degenerate Hermitian $\mc{O}_F$-lattice $M^{\flat} \subseteq \mbf{W}$ with $t(M^{\flat}) \leq 1$, we have
    \begin{equation}
    \mrm{Int}_{\ms{H}}(M^{\flat})^{\circ}_n = \partial \mrm{Den}^*_{\ms{H}}(M^{\flat})^{\circ}_n 
    \end{equation}
\end{proposition}
\begin{proof}
By definition, the quantity $\mrm{Int}_{\ms{H}}(M^{\flat})^{\circ}_n$ depends only on $\mrm{val}(M^{\flat})$. Since $t(M^{\flat}) \leq 1$, we may write $M^{\flat} = L^{\flat \prime} \oplus L^{\flat \prime \prime}$ (orthogonal direct sum) where $L^{\flat \prime}$ is self-dual of rank $n - 2$ and $\mrm{val}(L^{\flat \prime \prime}) = \mrm{val}(M^{\flat})$ (in the unramified case, this follows upon diagonalizing $M^{\flat}$; in the ramified case, this follows from picking a ``standard basis'' as in \cite[{Lemma 2.12}]{LL22II}). Using the cancellation property of local densities explained in \eqref{equation:part_I:local_Whittaker:local_densities_spherical:cancellation}, we thus reduce to the case $n = 2$ (which we now assume).

By the inductive decompositions in \eqref{equation:non-Arch_identity:statement:local_intersection_number:horizontal_decomp} and \eqref{equation:non-Arch_identity:statement:horizontal_density_decomp}, it is enough to show $\mrm{Int}_{\ms{H}}(M^{\flat})_n = \partial \mrm{Den}^*_{\ms{H}}(M^{\flat})_n$ (induct on $\mrm{val}(M^{\flat})$). We have $\partial \mrm{Den}^*(M^{\flat})_n = \partial \mrm{Den}^*_{\ms{H}}(M^{\flat})_n$ by construction, since $t(M^{\flat}) \leq 1$ (i.e. compare \eqref{equation:non-Arch_identity:statement:horizontal_density_decomp} and \eqref{equation:non-Arch_identity:statement:density_decomp}).

Set $b = \mrm{val}'(M^{\flat})$. Using the Cho--Yamauchi formulas, we find
    \begin{equation}\label{equation:non-Arch_idnentity:horizontal_identity:rank_one_Den}
    \mrm{Den}^*(q^2 X, M^{\flat})_n = X^{-b / 2} \sum_{j = 0}^b (q X)^j \quad \quad \partial \mrm{Den}^*(M^{\flat})_n = [\breve{F} : \breve{F}_0] \sum_{j = 0}^b (b - 2j) q^j
    \end{equation}
in all cases. The preceding formulas are valid even if $F_0 \neq \Q_p$ (and also valid if $p = 2$ whenever $F / F_0$ is unramified), hence why we wrote $q$ instead of $p$.

We have
    \begin{align}
    \mrm{Int}_{\ms{H}}(M^{\flat})_n & = 2 [\breve{F} : \breve{\Q}_p] \sum_{M^{\flat} \subseteq N^{\flat} \subseteq N^{\flat *}} p^s(1 - \eta(p) p^{-1}) \delta_{\mrm{tau}}(s)
    \\
    & = - [\breve{F} : \breve{\Q}_p] \sum_{M^{\flat} \subseteq N^{\flat} \subseteq N^{\flat *}} p^s(1 - \eta(p) p^{-1}) \left ( s - \frac{(1 - p^{-s})(1 - \eta(p))}{(1 - p^{-1})(p - \eta(p))} \right ) \notag
    \end{align}
where the sum runs over lattices $N^{\flat} \subseteq M^{\flat}_F$, where $s \coloneqq \mrm{val}'(N^{\flat})$, and where $\eta(p) \coloneqq -1, 0, 1$ in the inert, ramified, split cases respectively.

We prove the identity $\mrm{Int}_{\ms{H}}(M^{\flat})_n = \partial \mrm{Den}^*(M^{\flat})_n$ by induction on $b$. The case $b = 0$ is clear, as both quantities are $0$. 
    
Next suppose $b \geq 1$ and that $M^{\flat \prime}$ (resp. $M^{\flat \prime \prime}$) is a rank one non-degenerate lattice with $\mrm{val}'(M^{\flat \prime}) = b - 1$ (resp. $\mrm{val}'(M^{\flat \prime \prime}) = b - 2$). If $b - 2 \leq 0$, set $\mrm{Int}_{\ms{H}}(M^{\flat \prime \prime})_n \coloneqq 0$ (in which case $\partial \mrm{Den}^*(M^{\flat \prime \prime})_n = 0$ as well).

We have
    \begin{align}
    \partial \mrm{Den}^*(M^{\flat})_n - \partial \mrm{Den}^*(M^{\flat \prime})_n & = [\breve{F} : \breve{F}_0] (- b q^{b} + \sum_{j = 0}^{b - 1} q^j)
    \\
    \partial \mrm{Den}^*(M^{\flat})_n - \partial \mrm{Den}^*(M^{\flat \prime \prime})_n & = [\breve{F} : \breve{F}_0] (- b q^{b} - b q^{b - 1} + 2 \sum_{j = 0}^{b - 1} q^j).
    \end{align}
If $F / \Q_p$ is inert, we find
    \begin{align}
    \mrm{Int}_{\ms{H}}(M^{\flat})_n - \mrm{Int}_{\ms{H}}(M^{\flat \prime \prime})_n & = - [\breve{F} : \breve{\Q}_p] p^b (1 + p^{-1}) \left ( b - 2 \frac{(1 - p^{-b})}{(1 - p^{-1})(p + 1)} \right ) \notag
    \\
    & = [\breve{F} : \breve{\Q}_p] (- b p^{b} - b p^{b-1} + 2 \sum_{j = 0}^{b - 1} p^j).
    \end{align}
If $F / \Q_p$ is ramified, we find
    \begin{align}
    \mrm{Int}_{\ms{H}}(M^{\flat})_n - \mrm{Int}_{\ms{H}}(M^{\flat \prime})_n & = - [\breve{F} : \breve{\Q}_p] p^b \left ( b - \frac{(1 - p^{-b})}{(1 - p^{-1})p} \right )
    \\
    & = [\breve{F} : \breve{\Q}_p] (- b p^{b} + \sum_{j = 0}^{b - 1} p^j).
    \end{align}
If $F / \Q_p$ is split, we find
    \begin{align}
    \mrm{Int}_{\ms{H}}(M^{\flat})_n - \mrm{Int}_{\ms{H}}(M^{\flat \prime})_n & = - [\breve{F} : \breve{\Q}_p] \sum_{j = 0}^b p^j(1 - p^{-1}) j
    \\
    & = [\breve{F} : \breve{\Q}_p] (- b p^{b} + \sum_{j = 0}^{b - 1} p^j).
    \end{align}
This proves the lemma in all cases, by induction on $b$.
\end{proof}

\begin{corollary}
Theorem \ref{theorem:non-Arch_identity:statement} holds when $n = 2$.
\end{corollary}
\begin{proof}
If $n = 2$, Proposition \ref{proposition:non-Arch_identity:horizontal_identity} shows $\mrm{Int}_{\ms{H}}(L^{\flat})_n = \partial \mrm{Den}^*_{\ms{H}}(L^{\flat})_n = \partial \mrm{Den}^*(L^{\flat})_n$. We have $\mrm{gr}^1_{\mc{N}} K'_0(\mc{Z}(L^{\flat})_{\overline{k}})_{\Q} = 0$ because $\mc{Z}(L^{\flat})_{\overline{k}}$ is a scheme and because the reduced subscheme $\mc{N}_{\red} \subseteq \mc{N}$ is $0$-dimensional (a disjoint union of copies of $\Spec \overline{k}$), see Lemma \ref{lemma:moduli_pDiv:discrete_reduced}. Hence $\mrm{Int}_{\ms{V}}(L^{\flat})_n = 0$ since ${}^{\mbb{L}} \mc{Z}(L^{\flat})_{\ms{V}} \in \mrm{gr}^1_{\mc{N}} K'_0(\mc{Z}(L^{\flat})_{\overline{k}})_{\Q}$.
\end{proof}
                
            \subsection{Induction formula}
            \label{ssec:non-Arch_identity:induction_formula}
                Throughout Sections \ref{ssec:non-Arch_identity:induction_formula} and \ref{ssec:non-Arch_identity:limits}, we allow $F_0$ to be an arbitrary finite extension of $\Q_p$ (allowing $p = 2$ if $F / F_0$ is unramified). We take the following setup for the rest of of Section \ref{sec:non-Arch_identity} (i.e. the notations $n$, $V$, $L$, $L'$, $L''$, $L^{\flat}$, $x$, $x'$, and $x''$ are all reserved unless otherwise indicated).

\begin{setup}\label{setup:non-Arch_identity:induction_formula:lattices}
Let $V$ be a non-degenerate Hermitian $F$-module of rank $n$, with pairing $(-,-)$. Assume $n$ is even if $F / F_0$ is ramified. Let $L^{\flat} \subseteq V$ be a non-degenerate Hermitian $\mc{O}_F$-lattice of rank $n - 1$. Let $x, x', x'' \in V$ be nonzero and orthogonal to $L^{\flat}$ with $\langle x \rangle \subseteq \langle x' \rangle \subseteq \langle x'' \rangle$ and $\mrm{length}_{\mc{O}_F}( \langle x' \rangle / \langle x \rangle ) = \mrm{length}_{\mc{O}_F}( \langle x'' \rangle / \langle x' \rangle ) = 1$. Set
    \begin{equation}
    L \coloneqq L^{\flat} \oplus \langle x \rangle \quad \quad L' \coloneqq L^{\flat} \oplus \langle x' \rangle \quad \quad L'' \coloneqq L^{\flat} \oplus \langle x'' \rangle
    \end{equation}
\end{setup}

The notation $L''$ and $x''$ will only appear in our proof of the induction formula (Proposition \ref{proposition:non-Arch_identity:induction_formula}) for $F / F_0$ split.

\begin{proposition}[Induction formula]\label{proposition:non-Arch_identity:induction_formula}
If $\mrm{val}(x) > a_{\mrm{max}}(L^{\flat})$ in the nonsplit cases (resp. if $\mrm{val}(x) > 2 a_{\mrm{max}}(L^{\flat})$ in the case $F / F_0$ is split), we have
    \begin{equation}
    \mrm{Den}(X, L)_n =
    \begin{cases}
    X^2 \mrm{Den}(X, L')_n + (1 - X) \mrm{Den}(q^2 X, L^{\flat})_n & \text{if $F / F_0$ is inert} \\
    X \mrm{Den}(X, L')_n + (1 - X) \mrm{Den}(q^2 X, L^{\flat})_n & \text{if $F / F_0$ is ramified} \\
    X \mrm{Den}(X, L')_n + \mrm{Den}(q^2 X, L^{\flat})_n & \text{if $F / F_0$ is split.}
    \end{cases}
    \end{equation}
\end{proposition}

In the inert case, this is \cite[{Theorem 5.1}]{Terstiege13} (strictly speaking, there is a blanket $p \neq 2$ assumption there), which is a unitary analogue of \cite[{Theorem 2.6(1)}]{Katsurada99} (orthogonal groups); see also \cite[{Proposition 3.7.1}]{LZ22unitary} (there stated allowing $p = 2$) for a statement closer to ours.

Using the Cho--Yamauchi formulas, we give a uniform proof of the inert and ramified cases (Lemma \ref{lemma:non-Arch_identity:induction_formula:weak}). Our lower bounds on $\mrm{val}(x)$ are possibly nonsharp (e.g. in the inert case, we only show the induction formula when $\mrm{val}(x) > 2 a_{\mrm{max}}(L^{\flat})$) but this makes no difference for the proof of Theorem \ref{theorem:non-Arch_identity:statement}, where we will take $\mrm{val}(x) \ra \infty$ (Proposition \ref{proposition:non-Arch_identity:limits}).

The case when $F / F_0$ is split is more difficult for us, and the same proof only shows a weaker version of the induction formula (stated in Lemma \ref{lemma:non-Arch_identity:induction_formula:weak}), which is insufficient for our purposes. Extracting the induction formula from this weak version is the subject of Section \ref{ssec:non-Arch_identity:more_induction_formula_split}.

For the proof of Theorem \ref{theorem:non-Arch_identity:statement}, only the statement of the induction formula and the definitions in \eqref{equation:non-Arch_identity:induction_formula:Den_"prim"_definitions:1} and \eqref{equation:non-Arch_identity:induction_formula:Den_"prim"_definitions:2} will be needed.

We first record a few preparatory lemmas. As in Section \ref{ssec:Hermitian_conventions:lattices}, we fix a uniformizer $\varpi \in \mc{O}_F$ and a generator $u \in \mc{O}_F$ of the different ideal such that $\varpi^{\s} = - \varpi$ and $u^{\s} = -u$. 
We say a quantity \emph{stabilizes}, e.g. for $\mrm{val}(x) > C$ (for some constant $C$) if that quantity does not depend on $x$ if $\mrm{val}(x) > C$. When $F / F_0$ is nonsplit, given an $\mc{O}_F$-module $M$ and $m \in M$, we say e.g. that $m$ is \emph{exact $\varpi^e$-torsion} for $e \geq 0$ if $\varpi^e m = 0$ but $\varpi^{e-1} m \neq 0$ (and if $e = 0$, the only exact $\varpi^e$-torsion element is $0$). We use similar terminology for $\mc{O}_{F_0}$-modules and exact $\varpi_0^e$-torsion elements, etc..

\begin{lemma}\label{lemma:non-Arch_identity:induction_formula:type_and_rank}
Let $M$ be a non-degenerate integral Hermitian $\mc{O}_F$-lattice of rank $m$. Suppose elements $w_1, \ldots, w_r \in M$ have $\mc{O}_F$-span $M$. Write $T$ for the associated Gram matrix. Then $t(M) + \rank((uT) \otimes \mc{O}_F / \varpi) = m$.
\end{lemma}
\begin{proof}
If $F / F_0$ is split, the rank continues to make sense because $T$ is Hermitian (e.g. diagonalize the Hermitian form). In the unramified cases, the lemma follows by diagonalizing the Hermitian form. In the ramified case, the lemma follows by putting $M$ into ``standard form'' (i.e. an orthogonal direct sum of rank one lattices and rank two hyperbolic lattices) as in \cite[{Definition 2.11}]{LL22II}. We are allowing $m$ even or odd.
\end{proof}

\begin{lemma}\label{lemma:non-Arch-identity:limits:type-change_intersection}
Let $M$ be a non-degenerate integral Hermitian $\mc{O}_F$-lattice of rank $m$. Suppose $M_F = W' \oplus W''$ is an orthogonal decomposition with $W''$ of rank $1$.
    \begin{enumerate}[(1)]
        \item We have $t(M) - 1 \leq t(M \cap W') \leq t(M) + 1$.
        \item Let $M' \subseteq W'$ and $M'' \subseteq W''$ be the images of $M$ under the projections $M_F \ra W'$ and $M_F \ra W''$. Assume that $M'$ and $M''$ are integral and that $\mrm{val}(M'') > 0$. Then we have $t(M) = t(M') + 1$.
    \end{enumerate}
\end{lemma}
\begin{proof}
\hfill
\begin{enumerate}[(1)]
    \item The ramified case follows from \cite[{Lemma 2.23(2)}]{LL22II}. The inert case when $t(M) = 0$ is \cite[{Lemma 4.5.1}]{LZ22unitary}. The same proof works in general for arbitrary $F / F_0$ in arbitrary characteristic: select any basis $(w_1, \ldots, w_{m - 1})$ of $N \cap W'$, extend to a basis $(w_1, \ldots, w_m)$ of $M$ with Gram matrix $T$, then use the formulas $t(M) + \rank((u T) \otimes \mc{O}_F / \varpi) = m$ and $t(M \cap W') + \rank((u T^{\flat}) \otimes \mc{O}_F / \varpi) = m - 1$.
    \item Let $\underline{w} = [w_1, \ldots, w_m]$ be any basis of $M$, and let $T = (\underline{w}, \underline{w})$ be the corresponding Gram matrix. Let $\underline{w}' = [w_1', \ldots, w'_m]$ be the projection of $\underline{w}$ to $W'$, with Gram matrix $T' = (\underline{w}', \underline{w}')$. Since $\mrm{val}(M'') > 0$, we see $(u T) \otimes \mc{O}_F / \varpi = (u T') \otimes \mc{O}_F / \varpi$. Applying Lemma \ref{lemma:non-Arch_identity:induction_formula:type_and_rank} twice (once for $M$ and $r = m$ and once for $M'$ and $r = m$) proves the claim. \qedhere
\end{enumerate}
\end{proof}
Set
    \begin{align}\label{equation:non-Arch_identity:induction_formula:Den_"prim"_definitions:1}
    \mrm{Den}_{L^{\flat}, x}(X)^{\circ}_n & \coloneqq \sum_{\substack{L \subseteq M \subseteq M^* \\ M \cap L^{\flat}_F = L^{\flat}}} X^{\ell(M / L)} \mrm{Den}(X, M)^{\circ}_n 
    \\
    \mrm{Den}_{L^{\flat}, x}(X)_n & \coloneqq \mrm{Den}(X,L)_n = \sum_{L^{\flat} \subseteq M^{\flat} \subseteq M^{\flat *}} \mrm{Den}_{M^{\flat}, x}(X)^{\circ}_n
    \label{equation:non-Arch_identity:induction_formula:Den_"prim"_definitions:2}
    \end{align}
where the first sum runs over lattices $M \subseteq V$ and the second sum runs over lattices $M^{\flat} \subseteq L^{\flat}_F$. Note that the only dependence on $x$ is on $\mrm{val}(x)$ (since $\mrm{Den}_{L^{\flat}, x}(X)^{\circ}_n$ only depends on the isomorphism class of the Hermitian lattice $L$).

\begin{lemma}\label{lemma:non-Arch_identity:induction_formula:stabilizes}
The polynomial
    \begin{equation}\label{equation:non-Arch_identity:induction_formula:stabilizes:0}
    f_x(X) \coloneqq \begin{cases}
    \mrm{Den}_{L^{\flat},x}(X)^{\circ}_n - X^2 \mrm{Den}_{L^{\flat}, x'}(X)^{\circ}_n & \text{if $F / F_0$ is inert} \\
    \mrm{Den}_{L^{\flat}, x}(X)^{\circ}_n - X \mrm{Den}_{L^{\flat}, x'}(X)^{\circ}_n & \text{if $F / F_0$ is ramified} \\
    \mrm{Den}_{L^{\flat}, x}(X)^{\circ}_n - 2 X \mrm{Den}_{L^{\flat}, x'}(X)^{\circ}_n + X^2 \mrm{Den}_{L^{\flat}, x''}(X)^{\circ}_n & \text{if $F / F_0$ is split} \\
    \end{cases}
    \end{equation}
(an element of $\Z[X]$) stabilizes for $\mrm{val}(x) > 2 [\breve{F} : \breve{F}_0]^{-1} a_{\mrm{max}}(L^{\flat})$.
\end{lemma}
\begin{proof}
The notation $f_x(X)$ is temporary, used only for this lemma. If $F / F_0$ is split, write $\varpi = \varpi_1 \varpi_2$ for $\varpi_i \in \mc{O}_F$ with $\mrm{val}(\varpi_1 x) = \mrm{val}(\varpi_2 x) = \mrm{val}(x) + 1$. We may assume $x'' = \varpi^{-1} x$ in this case. Note $\mrm{Den}_{L^{\flat}, x'}(X)^{\circ}_n = \mrm{Den}_{L^{\flat}, \varpi_1^{-1} x}(X)^{\circ}_n = \mrm{Den}_{L^{\flat}, \varpi_2^{-1} x}(X)^{\circ}_n$.

Inspecting \eqref{equation:non-Arch_identity:induction_formula:Den_"prim"_definitions:1} shows
    \begin{align}
    f_x(X) & =
    \sum_{\substack{L \subseteq M \subseteq M^* \\ M \cap L^{\flat}_F = L^{\flat} \\ \varpi^{-1} x \not \in M }} X^{\ell(M / L)} \mrm{Den}(X, M)^{\circ}_n
    && \text{if $F / F_0$ is nonsplit} \label{equation:non-Arch_identity:induction_formula:stabilizes:difference:1}
    \\
    f_x(X) & = \sum_{\substack{L \subseteq M \subseteq M^* \\ M \cap L^{\flat}_F = L^{\flat} \\ \varpi_1^{-1} x \not \in M \\ \varpi_2^{-1} x \not \in M}} X^{\ell(M / L)} \mrm{Den}(X, M)^{\circ}_n
    && \text{if $F / F_0$ is split} \label{equation:non-Arch_identity:induction_formula:stabilizes:difference:2}
    \end{align}
where the sums run over lattices $M \subseteq V$. For each such $M$, we know $L^{\flat} \subseteq M$ is a saturated sublattice, hence $M = L^{\flat} \oplus \langle \xi \rangle$ (not necessarily orthogonal direct sum) for some $\xi \in V$.

If $L^{\flat}$ is not integral, then the lemma is trivial as the polynomials of the lemma statement are $0$. We assume $L^{\flat}$ is integral for the rest of the proof.

If $F / F_0$ is nonsplit, each lattice $M$ appearing in \eqref{equation:non-Arch_identity:induction_formula:stabilizes:difference:1} is of the form $M = L^{\flat} \oplus \langle y + \varpi^{-e} x \rangle$ for a uniquely determined element $y \in L^{\flat *} / L^{\flat}$, where $e \in \Z_{\geq 0}$ is such that $y \in L^{\flat *} / L^{\flat}$ is of exact $\varpi^e$-torsion. Conversely, an element $y \in L^{\flat *} / L^{\flat}$ gives rise to an $M$ appearing in \eqref{equation:non-Arch_identity:induction_formula:stabilizes:difference:1} if and only if $\mrm{val}(y + \varpi^{-e} x) \geq 0$. If $\mrm{val}(x) > 2[\breve{F} : \breve{F}_0]^{-1} a_{\mrm{max}}(L^{\flat})$, then $\mrm{val}(\varpi^{-e} x) > 0$, so $\mrm{val}(y + \varpi^{-e} x) \geq 0$ holds if and only if $\mrm{val}(y) \geq 0$.  

If $F / F_0$ is split, the preceding paragraph holds upon replacing $\varpi^{-e}$ with $\varpi_1^{-e_1} \varpi_2^{-e_2}$ for $e_1, e_2 \in \Z_{\geq 0}$ such that $y \in L^{\flat *} / L^{\flat}$ is of exact $\varpi_1^{e_1} \varpi_2^{e_2}$-torsion (i.e. $\varpi_1^{e_1} \varpi_2^{e_2} y \in L^{\flat}$ but $\varpi_1^{e_1 - 1} \varpi_2^{e_2} \not \in L^{\flat}$ and $\varpi_1^{e_1} \varpi_2^{e_2 - 1} \not \in L^{\flat}$).

In the previous notation, we thus have
    \begin{equation}\label{equation:non-Arch_identity:induction_formula:stabilizes:difference:y_version}
    f_x(X) = \sum_{\substack{y \in L^{\flat *} / L^{\flat} \\ \mrm{val}(y) \geq 0}} X^{\ell((L^{\flat} + \langle y \rangle)/L^{\flat})} \mrm{Den}(X,M)^{\circ}_n
    \end{equation}
where the sum runs over $y$, and $M = L^{\flat} \oplus \langle y + \varpi^{-e} x \rangle$ in the nonsplit case (resp. $M = L^{\flat} \oplus \langle y + \varpi_1^{- e_1} \varpi_2^{-e_2} x \rangle$ in the split case).

In the notation of \eqref{equation:non-Arch_identity:induction_formula:stabilizes:difference:y_version}, we have $t(M) = t(L^{\flat} + \langle y \rangle)$ + 1 by Lemma \ref{lemma:non-Arch-identity:limits:type-change_intersection}(2) (using $\mrm{val}(\varpi^{-e} x) > 0$ in the nonsplit case and $\mrm{val}(\varpi_1^{-e_2} \varpi_2^{-e_2} x) > 0$ in the split case). Hence we have
    \begin{equation}\label{equation:non-Arch_identity:induction_formula:stabilizes:difference:Den}
    \mrm{Den}(X,M)^{\circ}_n = \prod_{i = 0}^{t(M^{\flat} + \langle y \rangle)}(1 - \eta^i(\varpi_0) q^i X)
    \end{equation}
(see definition in \eqref{equation:non-Arch_identity:horizontal_identity:Cho--Yamauchi:"primitive"}), and now the right-hand side of \eqref{equation:non-Arch_identity:induction_formula:stabilizes:difference:y_version} clearly depends only on $L^{\flat}$ (and not on $x$).
\end{proof}

\begin{lemma}\label{lemma:non-Arch_identity:induction_formula:weak}
With notation as above, assume $\mrm{val}(x) > 2 [\breve{F} : \breve{F}_0]^{-1} a_{\mrm{max}}(L^{\flat})$. We have
    \begin{equation}\label{equation:non-Arch_identity:induction_formula:weak}
    (1 - X) \mrm{Den}(q^2 X, L^{\flat})_n = 
    \begin{cases}
    \mrm{Den}(X,L) - X^2 \mrm{Den}(X,L') & \text{if $F / F_0$ is inert} \\
    \mrm{Den}(X,L) - X \mrm{Den}(X,L') & \text{if $F / F_0$ is ramified} \\
    \mrm{Den}(X,L) - 2 X \mrm{Den}(X,L') + X^2 \mrm{Den}(X, L'') & \text{if $F / F_0$ is split.} \\
    \end{cases}
    \end{equation}
\end{lemma}
\begin{proof}
Combining \eqref{equation:non-Arch_identity:induction_formula:stabilizes:difference:y_version} and \eqref{equation:non-Arch_identity:induction_formula:stabilizes:difference:Den}, we find that the right-hand side of \eqref{equation:non-Arch_identity:induction_formula:weak} is given by
    \begin{equation}\label{equation:non-Arch_identity:induction_formula:weak:1}
    \sum_{L^{\flat} \subseteq M^{\flat} \subseteq M^{\flat *}} X^{\ell(M^{\flat} / L^{\flat})} \sum_{\substack{y \in M^{\flat *} / M^{\flat} \\ \mrm{val}(y) \geq 0}} X^{\ell((M^{\flat} + \langle y \rangle)/M^{\flat})} \prod_{i = 0}^{t(M^{\flat} + \langle y \rangle)} (1 - \eta^i(\varpi_0) q^i X)
    \end{equation}
in all cases, where the outer sum runs over lattices $M^{\flat} \subseteq L^{\flat}_F$.

Collecting the terms with $M^{\flat} + \langle y \rangle = N^{\flat}$ for fixed integral lattices $N^{\flat} \subseteq L^{\flat}_F$, we find that \eqref{equation:non-Arch_identity:induction_formula:weak:1} is equal to
    \begin{equation*}
    \sum_{L^{\flat} \subseteq N^{\flat} \subseteq N^{\flat *}} X^{\ell(N^{\flat} / L^{\flat})} \prod_{i = 0}^{t(N^{\flat})} (1 - \eta^i(\varpi_0) q^i X) \sum_{\substack{L^{\flat} \subseteq M^{\flat} \subseteq N^{\flat} \\ N^{\flat} / M^{\flat} \text{cyclic}}} (\text{number of generators of $N^{\flat} / M^{\flat}$}).
    \end{equation*}
where the outer sum runs over lattices $N^{\flat} \subseteq L^{\flat}_F$ and the inner sum runs over lattices $M^{\flat}$. We have
    \begin{align*}
    \sum_{\substack{L^{\flat} \subseteq M^{\flat} \subseteq N^{\flat} \\ N^{\flat} / M^{\flat} \text{cyclic}}} (\text{number of generators of $N^{\flat} / M^{\flat}$}) &= \sum_{\substack{N^{\flat *} \subseteq M^{\flat *} \subseteq L^{\flat *} \\ M^{\flat *} / N^{\flat *} \text{cyclic}}} (\text{number of generators of $M^{\flat *} / N^{\flat *}$})
    \\
    & = |L^{\flat *} / N^{\flat *}| = |N^{\flat} / L^{\flat}| = q^{\ell(N^{\flat} / L^{\flat})}
    \end{align*}
so \eqref{equation:non-Arch_identity:induction_formula:weak:1} is equal to
    \begin{equation}
    \sum_{L^{\flat} \subseteq N^{\flat} \subseteq N^{\flat *}} (q X)^{\ell(N^{\flat} / L^{\flat})} \prod_{i = 0}^{t(N^{\flat})} (1 - \eta^i(\varpi_0) q^i X).
    \end{equation}
Inspecting the Cho--Yamauchi formulas (and surrounding discussion) at the beginning of Section \ref{ssec:non-Arch_identity:horizontal_identity} shows that the displayed expression is equal to $(1 - X) \mrm{Den}(q^2 X, L^{\flat})_{n}$ in all cases (if $F / F_0$ is ramified, note that $t(N^{\flat})$ is always odd because $N^{\flat}$ has rank $n - 1$, which we have assumed is odd in the ramified case).
\end{proof}
    
            \subsection{More on induction formula: split}
            \label{ssec:non-Arch_identity:more_induction_formula_split}
                Suppose $F / F_0$ is split. To prove the induction formula (Proposition \ref{proposition:non-Arch_identity:induction_formula}), it remains only to show that $\mrm{Den}(X,L) - X \mrm{Den}(X,L')$ stabilizes for $\mrm{val}(x) > 2 a_{\mrm{max}}(L^{\flat})$, as Lemma \ref{lemma:non-Arch_identity:induction_formula:weak} then shows $(1 - X) (\mrm{Den}(X,L) - X \mrm{Den}(X, L')) = (1 - X) \mrm{Den}(q^2 X, L^{\flat})_n$.

We define some more notation (only used in Section \ref{ssec:non-Arch_identity:more_induction_formula_split}). Fix a uniformizer $\varpi_0$ of $\mc{O}_{F_0}$, and consider the elements
    \begin{equation}
    \varpi_1 = (\varpi_0, 1) \quad \quad \varpi_2 = (1, - \varpi_0) \quad \quad e_1 = (1,0) \quad \quad e_2 = (0,1)
    \end{equation}
in $\mc{O}_F = \mc{O}_{F_0} \times \mc{O}_{F_0}$. Given an $\mc{O}_F$-module $M$, we set $M_1 \coloneqq e_1 M$ and $M_2  \coloneqq e_2 M$ (so $M = M_1 \oplus M_2$). We similarly write $y_1 \coloneqq e_1 y$ and $y_2 \coloneqq e_2 y$ for $y \in M$. If $M$ is a non-degenerate Hermitian $\mc{O}_F$-lattice, we set $M_1^* \coloneqq e_2 M^*$ and $M_2^* \coloneqq e_1 M^*$. If $M$ is moreover integral, the Hermitian pairing induces an identification $M_2^* / M_1 \cong \Hom_{\mc{O}_{F_0}}(M_1^* / M_2, F_0 / \mc{O}_{F_0})$.

For integers $t \geq 0$, we set
    \begin{equation}
    \mf{m}(t, X) \coloneqq \prod_{i = 0}^{t - 1} (1 - q^i X)
    \end{equation}
so that $\mrm{Den}(X,M)_n^{\circ} = \mf{m}(t(M), x)$ for any integral non-degenerate Hermitian $\mc{O}_F$-lattice $M$ of rank $n$. If $\mc{T}$ is a finite length $\mc{O}_{F_0}$-module, we set 
    \begin{equation}
    t_0(\mc{T}) \coloneqq \dim_{\F_q}(\mc{T} \otimes_{\mc{O}_{F_0}} \F_q)
    \quad \quad \ell_0(\mc{T}) \coloneqq \mrm{length}_{\mc{O}_{F_0}}(\mc{T}).
    \end{equation}

\begin{lemma}\label{lemma:more_induction_formula_split:diff_stabilizes}
Consider the polynomial
    \begin{equation}\label{equation:more_induction_formula_split:diff_stabilizes}
    h_{\mrm{diff},x}(X) \coloneqq \sum_{\substack{L \subseteq M \subseteq M^* \\ M_1 \cap L^{\flat}_F  = L_1^{\flat} \\ M_2 \cap L^{\flat}_F \neq L_2^{\flat} \\ M / L \text{ is cyclic} \\ \varpi_1^{-1} x \not \in M}} X^{\ell(M / L)} \mf{m}(t(M), X)
    \end{equation}
where the sum runs over lattices $M \subseteq V$ (satisfying the displayed conditions). This sum stabilizes for $\mrm{val}(x) > 2 a_{\mrm{max}}(L^{\flat})$.
\end{lemma}
\begin{proof}
Each lattice $M$ in the sum is of the form $M = L + \langle \xi \rangle$ for a unique element $\xi = y + \varpi_1^{-e_1} \varpi_2^{-e_2}  x \in L^* / L$ with $y \in L^{\flat *}$, such that $\mrm{val}(\xi) \geq 0$, and with $e_1, e_2 \in \Z_{\geq 0}$.

Assume $\mrm{val}(x) > 2 a_{\mrm{max}}(L^{\flat})$. We claim that $\mrm{val}(y) \geq 0$ (in the notation above). The additional conditions on $M$ imply that $y_1 \in L_2^{\flat *} / L_1^{\flat}$ is of exact $\varpi_1^{e_1}$-torsion and that $\varpi_2^{e_2} y_2 \not \in L^{\flat}$. We thus have $e_1 \leq a_{\mrm{max}}(L^{\flat})$ and $e_2 < a_{\mrm{max}}(L^{\flat})$, so $\mrm{val}(\varpi_1^{-e_1} \varpi_2^{-e_2}  x) > 0$ when $\mrm{val}(x) > 2 a_{\mrm{max}}(L^{\flat})$. This implies that $\mrm{val}(y) \geq 0$ as well.

Consider the $F$-linear (non-unitary) automorphism $\phi \colon V \ra V$ which is the identity on $L^{\flat}_F$ and sends $x \mapsto \varpi_2 x$. Then $M \mapsto \phi(M)$ is a bijection from the set of lattices appearing in the sum for $h_{\mrm{diff},x}(X)$ to the set of lattices appearing in the sum for $h_{\mrm{diff},\varpi_2 x}(X)$ (we remind the reader that $L$ depends on $x$ as well).

In the above setup, an application of Lemma \ref{lemma:non-Arch-identity:limits:type-change_intersection}(2) shows $t(M) = t(\phi(M)) = t(L^{\flat} + \langle y \rangle) + 1$. We also find $\ell(M / L) = \ell(\phi(M) / (L^{\flat} \oplus \langle \varpi_2 x \rangle )) = \ell((L^{\flat} + \langle y \rangle) / L^{\flat})$. This shows $h_{\mrm{diff},x}(X) = h_{\mrm{diff}, \varpi_2 x}(X)$ (compare the $M$ term and the $\phi(M)$ term). This proves the lemma, as the $x$-dependence of $h_{\mrm{diff},x}(X)$ is only on $\mrm{val}(x)$.
\end{proof}

\begin{lemma}\label{lemma:more_induction_formula_split:quotient}
Let $\mc{T}$ be a finite length $\mc{O}_{F_0}$-module, and suppose $\mc{T}$ is $\varpi_0^e$-torsion. For any integer $b > e$, form the $\mc{O}_{F_0}$-module $A = \mc{T} \oplus (\varpi_0^{-b} \mc{O}_{F_0} /\mc{O}_{F_0})$. Consider $u = t + w \in A$ with $t \in \mc{T}$ and $w \in \varpi_0^{-b} \mc{O}_{F_0} / \mc{O}_{F_0}$ both of exact $\varpi_0^r$-torsion. There is a (non-canonical) isomorphism
    \begin{equation}
    A/(u) \cong (\mc{T}/(t)) \oplus (\varpi_0^{-b} \mc{O}_{F_0} /\mc{O}_{F_0}).
    \end{equation}
\end{lemma}
\begin{proof}
This follows from the structure theorem for finitely generated modules over the discrete valuation ring $\mc{O}_{F_0}$. For example, we can select elements $e_1, \ldots, e_m \in \mc{T}$ such that $\mc{T} = \langle e_1 \rangle \oplus \cdots \oplus \langle e_n \rangle$ and such that $t = \varpi_0^s e_1$ for some $s \geq 0$. The case $r = 0$ is trivial, so take $r \geq 1$. Then $r + s \leq e$. If $w' \in \varpi_0^{-b} \mc{O}_{F_0} / \mc{O}_{F_0}$ is such that $\varpi_0^s w' = w$ there is an isomorphism
    \begin{equation}
    \mc{T} \oplus (\varpi_0^{-b} \mc{O}_{F_0} / \mc{O}_{F_0}) \ra A
    \end{equation}
sending $e_1 \mapsto e_1 + w'$, $e_i \mapsto e_i$ for $i \geq 2$, and $z \mapsto z$ (for any generator $z$ of $\varpi_0^{-b} \mc{O}_{F_0} / \mc{O}_{F_0}$). This isomorphism takes $t$ to $t + w$.
\end{proof}

Given a finite torsion cyclic $\mc{O}_{F_0}$-module $N \cong \mc{O}_{F_0} / \varpi_0^ a \mc{O}_{F_0}$, we set $\mrm{ord}(N) \coloneqq a$.

\begin{lemma}\label{lemma:non-Arch_identity:more_induction_formula_split:torsion_counting}
Let $\mc{T}$ be a finite length $\mc{O}_{F_0}$-module, and assume $\mc{T}$ is $\varpi^e$-torsion for some $e \geq 0$. For any integer $b \geq 0$, form the $\mc{O}_{F_0}$-module $A_b \coloneqq \mc{T} \oplus (\varpi_0^{-b} \mc{O}_{F_0} / \mc{O}_{F_0})$.
The polynomial
    \begin{equation}
    \a_b \coloneqq \sum_{\substack{\text{cyclic submodules} \\ N \subseteq A_b}} X^{\mrm{ord}(N)} \mf{m}(t_0(A_b / N), X)
    \end{equation}
stabilizes for $b > e$.
\end{lemma}
\begin{proof}
Applying $- \otimes_{\mc{O}_{F_0}} \F_q$ to the exact sequence
    \begin{equation}
    0 \ra N \ra A_b \ra A_b / N \ra 0
    \end{equation}
shows
    \begin{equation}
    t_0(A_b / N) =
    \begin{cases}
    t_0(A_b)  & \text{if $N \subseteq \varpi_0 A_b$} \\
    t_0(A_b) - 1 & \text{if $N \not \subseteq \varpi_0 A_b$}
    \end{cases}
    \end{equation}
for any cyclic submodule $N \subseteq A_b$. We also have $t_0(A_b) = t_0(\mc{T}) + 1$ if $b > 1$. 

There is a natural inclusion $A_b \ra A_{b+1}$. For any cyclic submodule $N \subseteq A_b$, we have
    \begin{equation}
    t_0(A_b / N) =
    \begin{cases}
    t_0(A_{b+1}/N) - 1 & \text{if $N = \langle t + \varpi_0^{-b} \rangle$ with $t \in \varpi_0 \mc{T}$} \\
    t_0(A_{b+1}/N) & \text{otherwise}
    \end{cases}
    \end{equation}
where $\varpi_0^{-b} \in \varpi_0^{-b} \mc{O}_{F_0} / \mc{O}_{F_0}$. Assume $b > e$. Then, in the first case above, the element $t \in \mc{T}$ is uniquely determined by $N$ (using $b > e$). The cyclic submodules $N \subseteq A_{b+1}$ with $N \not \subseteq A_b$ are of the form $N \langle t + \varpi_0^{-b-1} \rangle$ for a unique $t \in \mc{T}$.

We thus have
    \begin{align}\label{equation:non-Arch_identity:more_induction_formula_split:alpha_difference}
    & \a_{b+1} - \a_b
    \\
    & = \sum_{\substack{t \in \mc{T} \\ N = \langle t + \varpi_0^{-b - 1} \rangle}} X^{\mrm{ord}(N)} \mf{m}(t_0(A_{b+1} / N), X) + \sum_{\substack{t \in \varpi_0 \mc{T} \\ N = \langle t + \varpi_0^{- b} \rangle}} X^{\mrm{ord}(N)} \mf{m}(t_0(A_{b+1}/N), X)
    \\
    & \mathrel{\hphantom{=}} - \sum_{\substack{t \in \varpi_0 \mc{T} \\ N = \langle t + \varpi_0^{- b} \rangle }} X^{\mrm{ord}(N)} \mf{m}(t_0(A_b / N), X)
    \end{align}
where the sums run over $t \in \mc{T}$ or $t \in \varpi \mc{T}$, as indicated. We compute
    \begin{equation}
    \sum_{\substack{t \in \mc{T} \\ N = \langle t + \varpi_0^{-b - 1} \rangle}} X^{\mrm{ord}(N)} \mf{m}(t_0(A_{b+1} / N), X) = |\mc{T}| X^{b+1} \mf{m}(t_0(\mc{T}), X)
    \end{equation}
where $|\mc{T}|$ is the cardinality of $\mc{T}$. For any integer $a \geq 0$, we have the identity $\mf{m}(a+1, X) - \mf{m}(a, x) = -q^a X \mf{m}(a, X)$, so we compute
    \begin{align}
    & \sum_{\substack{t \in \varpi_0 \mc{T} \\ N = \langle t + \varpi_0^{- b} \rangle}} X^{\mrm{ord}(N)} \mf{m}(t_0(A_{b+1}/N), X) - \sum_{\substack{t \in \varpi_0 \mc{T} \\ N = \langle t + \varpi_0^{- b} \rangle }} X^{\mrm{ord}(N)} \mf{m}(t_0(A_b / N), X)
    \\
    & = - |\varpi_0 \mc{T}| q^{t_0(\mc{T})} X^{b+1} \mf{m}(t_0(\mc{T}), X).
    \end{align}
But the exact sequence
    \begin{equation}
    0 \ra \varpi_0 \mc{T} \ra \mc{T} \ra \mc{T} / \varpi_0 \mc{T} \ra 0
    \end{equation}
shows that $|\mc{T}| = |\varpi_0 \mc{T}| q^{t_0(\mc{T})}$ since $t_0(\mc{T}) = \dim_{\F_q} \mc{T} / \varpi_0 \mc{T}$ by definition. Substituting into \eqref{equation:non-Arch_identity:more_induction_formula_split:alpha_difference} shows $\a_{b+1} - \a_b = 0$.
\end{proof}

\begin{lemma}\label{lemma:non-Arch_identity:more_induction_formula_split:stabilizes}
The polynomial $\mrm{Den}_{L^{\flat}, x}(X)^{\circ}_n - X \mrm{Den}_{L^{\flat}, x'}(X)^{\circ}_n$ stabilizes for $\mrm{val}(x) > 2 a_{\mrm{max}}(L^{\flat})$.
\end{lemma}
\begin{proof}
As the $x'$ dependence of $\mrm{Den}_{L^{\flat}, x'}(X)^{\circ}_n$ is only on $\mrm{val}(x')$, we may assume $x' = \varpi_1^{-1} x$ without loss of generality. Assume $\mrm{val}(x) > 2 a_{\mrm{max}}(L^{\flat})$. The lemma is trivial if $L^{\flat}$ is not integral (the polynomial is $0$), so assume $L^{\flat}$ is integral.

Inspecting \eqref{equation:non-Arch_identity:induction_formula:Den_"prim"_definitions:1} shows that $\mrm{Den}_{L^{\flat}, x}(X)^{\circ}_n - X \mrm{Den}_{L^{\flat}, x'}(X)^{\circ}_n$ is equal to 
    \begin{equation}\label{equation:non-Arch_identity:more_induction_formula_split:stabilizes:1}
    \sum_{\substack{L \subseteq M \subseteq M^* \\ M \cap L^{\flat}_F = L^{\flat} \\ \varpi_1^{-1} x \not \in M }} X^{\ell(M / L)} \mf{m}(t(M), X).
    \end{equation}
where the sum runs over lattices $M \subseteq V$ (similar reasoning was used at the beginning of the proof of Lemma \ref{lemma:non-Arch_identity:induction_formula:stabilizes}). For each $M$ in the above sum, note that $M / L$ is cyclic (again, $L^{\flat} \subseteq M$ is a saturated sublattice, so there is a direct sum decomposition $M = L^{\flat} \oplus \langle \xi \rangle$ (not necessarily orthogonal) for some $\xi \in V$). By Lemma \ref{lemma:more_induction_formula_split:diff_stabilizes}, it is enough to show that
    \begin{equation}\label{equation:non-Arch_identity:more_induction_formula_split:stabilizes:with_error}
    \sum_{\substack{L \subseteq M \subseteq M^* \\ M_1 \cap L^{\flat}_F = L_1^{\flat} \\ M / L \text{ is cyclic} \\ \varpi_1^{-1} x \not \in M }} X^{\ell(M / L)} \mf{m}(t(M), X).
    \end{equation}
stabilizes for $\mrm{val}(x) > 2 a_{\mrm{max}}(L^{\flat})$, where the sum runs over lattices $M \subseteq V$ (because the difference between \eqref{equation:non-Arch_identity:more_induction_formula_split:stabilizes:with_error} and \eqref{equation:non-Arch_identity:more_induction_formula_split:stabilizes:1} is \eqref{equation:more_induction_formula_split:diff_stabilizes}).

We find that \eqref{equation:non-Arch_identity:more_induction_formula_split:stabilizes:with_error} equals
    \begin{equation}\label{equation:non-Arch_identity:more_induction_formula_split:stabilizes:with_error:2}
    \sum_{\substack{L_1 \subseteq M_1 \subseteq L_2^* \\ M_1 \cap L^{\flat}_F = L^{\flat}_1 \\ \varpi_1^{-1} x_1 \not \in M_1}} \sum_{\substack{L_2 \subseteq M_2 \subseteq M_1^* \\ M_2 / L_2 \text{ is cyclic}}} X^{\ell(M / L)} \mf{m}(t(M), X)
    \end{equation}
where the outer sum runs over lattices $M_1 \subseteq V_1$, the right-most sum runs over lattices $M_2 \subseteq V_2$, and $M = M_1 \oplus M_2$. Note that the lattices $M_1$ always satisfy $M_1 / L_1$ being cyclic, because $M_1 \cap L^{\flat}_F = L_1^{\flat}$ implies $M_1 = L^{\flat}_1 \oplus \langle y_1 + \varpi_1^{-e_1} x \rangle$ where $y_1 \in L_2^{\flat *} / L_1^{\flat}$ is of exact $\varpi_1^{e_1}$-torsion.

To prove the lemma, it is enough to check that \eqref{equation:non-Arch_identity:more_induction_formula_split:stabilizes:with_error:2} does not change if $x$ is replaced with $\varpi_2 x$. The set of lattices $M_1 \subseteq V_1$ appearing in the outer sum is indexed elements $y_1 \in L_2^{\flat *} / L_1^{\flat}$ (since $e_1$ is determined by $y_1$, in the above notation), and hence does not change if $x$ is replaced by $\varpi_2 x$ (here using $\mrm{val}(x) > a_{\mrm{max}}(L^{\flat})$ to ensure $M_1 \subset L_2^*$ for any choice of $y_1$). Note also that $\ell_0(M_1 / L_1) = e_1$ and hence does not change when $x$ is replaced by $\varpi_2 x$.

For the rest of the proof, fix an $M_1$ as in the outer sum of \eqref{equation:non-Arch_identity:more_induction_formula_split:stabilizes:with_error:2}. We will show that the inner sum of \eqref{equation:non-Arch_identity:more_induction_formula_split:stabilizes:with_error:2} does not change if $x$ is replaced by $\varpi_2 x$.

Set $A = M_1^* / L_2$. The inner sum is
    \begin{equation}\label{equation:non-Arch_identity:more_induction_formula_split:inner_sum}
    X^{\ell_0(M_1 / L_1)} \sum_{\substack{\text{cyclic submodules} \\ N \subseteq A}} X^{\mrm{ord}(N)} \mf{m}(t_0(A / N), X).
    \end{equation}
We already discussed that the factor $X^{\ell_0(M_1 / L_1)}$ does not change when $x$ is replaced by $\varpi_2 x$. On the other hand, we have $A \cong \Hom_{\mc{O}_{F_0}}(L_2^* / M_1, F_0 / \mc{O}_{F_0})$ so $A \cong L_2^* / M_1$ (non-canonically). If $b \coloneqq \mrm{val}(x)$ and $\mc{T} \coloneqq L_2^{\flat *} / L_1^{\flat}$, then Lemma \ref{lemma:more_induction_formula_split:quotient} shows $A \cong (\mc{T} / \langle y_1 \rangle) \oplus (\varpi_0^{-b} \mc{O}_{F_0} / \mc{O}_{F_0})$, where $y_1$ is associated to $M_1$ as above (since the submodule $(M_1 / L_1) \subseteq L_2^* / L_1$ is cyclic and generated by $y_1 + \varpi^{-e_1} x_1$ where $y_1$ is of exact $\varpi_1^{e_1}$-torsion).

Now Lemma \ref{lemma:non-Arch_identity:more_induction_formula_split:torsion_counting} implies that the sum in \eqref{equation:non-Arch_identity:more_induction_formula_split:inner_sum} does not change if $x$ is replaced by $\varpi_2 x$.
\end{proof}

\begin{proof}[Proof of Proposition \ref{proposition:non-Arch_identity:induction_formula} in split case]
Assume $F / F_0$ is split. As remarked at the beinning of Section \ref{ssec:non-Arch_identity:more_induction_formula_split}, it is enough to show that $\mrm{Den}(X,L) - X \mrm{Den}(X,L')$ stabilizes for $\mrm{val}(x) > 2 a_{\mrm{max}}(L^{\flat})$. We have
    \begin{equation}
    \mrm{Den}(X,L) - X \mrm{Den}(X,L') = \sum_{L^{\flat} \subseteq M^{\flat} \subseteq M^{\flat *}} \mrm{Den}_{L^{\flat},x}(X)^{\circ}_n - X \mrm{Den}_{L^{\flat}, x'}(X)^{\circ}_n 
    \end{equation}
by definition (see \eqref{equation:non-Arch_identity:induction_formula:Den_"prim"_definitions:2}), so Lemma \ref{lemma:non-Arch_identity:more_induction_formula_split:stabilizes} proves the claimed stabilization.
\end{proof}
    
            \subsection{Limits}
            \label{ssec:non-Arch_identity:limits}
                We continue in the setup of Section \ref{ssec:non-Arch_identity:induction_formula} but now assume $\varepsilon(V) = - 1$ if $F / F_0$ is nonsplit. Recall also the definitions in \eqref{equation:non-Arch_identity:induction_formula:Den_"prim"_definitions:1} and \eqref{equation:non-Arch_identity:induction_formula:Den_"prim"_definitions:2}.

Let $M^{\flat} \subseteq L^{\flat}_F$ be any non-degenerate Hermitian $\mc{O}_F$-lattice of rank $n - 1$ with $t(M^{\flat}) \leq 1$. If $F / F_0$ is nonsplit, set
    \begin{align}
    & \partial \mrm{Den}_{L^{\flat}}(x)_n \coloneqq - [\breve{F} : \breve{F}_0] \frac{d}{d X} \bigg|_{X = 1} \mrm{Den}(X,L)_n
    & & \partial \mrm{Den}_{M^{\flat}, \ms{H}}(x)^{\circ}_n \coloneqq - [\breve{F} : \breve{F}_0] \frac{d}{d X} \bigg|_{X = 1} \mrm{Den}_{M^{\flat}, x}(X)^{\circ}_n \notag
    \\
    & \partial \mrm{Den}_{L^{\flat}, \ms{H}}(x)_n \coloneqq \sum_{\substack{L \subseteq N \subseteq N^* \\ N^{\flat} = N \cap L^{\flat}_F \\ t(N^{\flat}) \leq 1}} \partial \mrm{Den}_{N^{\flat}, \ms{H}}(x)^{\circ}_n
    && \partial \mrm{Den}_{L^{\flat}, \ms{V}}(x)_n \coloneqq \partial \mrm{Den}_{L^{\flat}}(x)_n - \partial \mrm{Den}_{L^{\flat}, \ms{H}}(x)_n. \notag
    \end{align}
If $F / F_0$ is split, set
    \begin{align}
    & \mrm{Den}_{L^{\flat}}(x)_n \coloneqq \mrm{Den}(X,L)_n \bigg|_{X = 1}
    && \mrm{Den}_{M^{\flat}, \ms{H}}(x)^{\circ}_n \coloneqq \mrm{Den}_{M^{\flat}, x}(X)^{\circ}_n \bigg|_{X = 1} \notag
    \\
    & \mrm{Den}_{L^{\flat}, \ms{H}}(x)_n \coloneqq \sum_{\substack{L \subseteq N \subseteq N^* \\ N^{\flat} = N \cap L^{\flat}_F \\ t(N^{\flat}) \leq 1}} \mrm{Den}_{N^{\flat}, \ms{H}}(x)^{\circ}_n
    && \mrm{Den}_{L^{\flat}, \ms{V}}(x)_n \coloneqq \mrm{Den}_{L^{\flat}}(x)_n - \mrm{Den}_{L^{\flat}, \ms{H}}(x)_n. \notag
    \end{align}
The above sums run over lattices $N \subseteq V$ (so $N^{\flat}$ varies). These definitions also apply for any $x \not \in L^{\flat}_F$ (not necessarily perpendicular to $L^{\flat}_F$), as long as we take $L = L^{\flat} + \langle x \rangle$.

\begin{lemma}\label{lemma:non-Arch_identity:limits:split_vertical_den_vanishes}
If $F / F_0$ is split, then $\mrm{Den}_{L^{\flat},\ms{V}}(x)_n = 0$ for all $x$.
\end{lemma}
\begin{proof}
Inspecting \eqref{equation:non-Arch_identity:horizontal_identity:Cho--Yamauchi:"primitive"} shows that $\mrm{Den}(X,M)^{\circ}_n = 0$ unless $M = M^*$. Lemma \ref{lemma:non-Arch-identity:limits:type-change_intersection} implies $\mrm{Den}_{N^{\flat}, x}(X)^{\circ}_n |_{X = 1} = 0$ unless $t(N^{\flat}) \leq 1$, i.e. $\mrm{Den}_{L^{\flat}}(x)_n = \mrm{Den}_{L^{\flat}, \ms{H}}(x)$.
\end{proof}

Given $x \in V$ with $(x,x) \neq 0$, we set $\mrm{val}''(x) \coloneqq \mrm{val}'(x)$ if $F / F_0$ is not inert (resp. $\mrm{val}''(x) \coloneqq (\mrm{val}(x)-1)/2$ if $F / F_0$ is inert) to save space. We say a limit \emph{stabilizes} if the argument of the limit becomes constant.

\begin{proposition}\label{proposition:non-Arch_identity:limits}
If $F / F_0$ is nonsplit, we have
    \begin{equation}
    \partial \mrm{Den}^*(L^{\flat})_n = 2 [\breve{F} : \breve{F}_0]^{-1} \lim_{x \ra 0} \left ( \partial \mrm{Den}_{L^{\flat}}(x')_n - \mrm{val}''(x) \mrm{Den}^*(L^{\flat})_n \right ).
    \end{equation}
If $F / F_0$ is split, we have
    \begin{equation}
    \partial \mrm{Den}^*(L^{\flat})_n = \lim_{x \ra 0} \left ( \mrm{Den}_{L^{\flat}}(x')_n - \mrm{val}(x) \mrm{Den}^*(L^{\flat})_n \right ).
    \end{equation}
The expressions are $0$ if $L^{\flat}$ is not integral, and all limits stabilize for $\mrm{val}(x) \gg 0$. If $L^{\flat}$ is integral and $F / F_0$ is nonsplit (resp. split), then the limits stabilize when $\mrm{val}(x) > a_{\mrm{max}}(L^{\flat})$ (resp. $\mrm{val}(x) > 2 a_{\mrm{max}}(L^{\flat})$).
\end{proposition}
\begin{proof}
We emphasize that we are following Setup \ref{setup:non-Arch_identity:induction_formula:lattices}; in particular, we have $x' \ra 0$ as $x \ra 0$. Assume $L^{\flat}$ is integral (as the lemma is otherwise clear) and assume $\mrm{val}(x) > a_{\mrm{max}}(L^{\flat})$. The key input is the induction formula from Proposition \ref{proposition:non-Arch_identity:induction_formula}.

\emph{Case $F / F_0$ is nonsplit:} Multiply the induction formula from Proposition \ref{proposition:non-Arch_identity:induction_formula} by $X^{-\mrm{val}'(L^{\flat})/2}$, and call the resulting expression ($*$) (temporary notation). Taking one derivative of ($*$) at $X = 1$ yields
    \begin{equation}
    \mrm{Den}^*(L^{\flat})_n = \partial \mrm{Den}_{L^{\flat}}(x)_n - \partial \mrm{Den}_{L^{\flat}}(x')_n.
    \end{equation}
Here we used $\mrm{Den}(1,L)_n = \mrm{Den}(1,L')_n = 0$ because $\varepsilon(V) = -1$ causes a sign in the functional equation \eqref{equation:part_I:Eisenstein:local_functional_equations:non-Arch:non-star_density}. Taking two derivatives of ($*$) at $X = 1$ yields the identity
    \begin{align}\label{equation:proposition:non-Arch_identity:limits:nonsplit}
    & \mrm{val}'(L^{\flat}) \partial \mrm{Den}_{L^{\flat}}(x)_n + [\breve{F} : \breve{F}_0] \frac{d^2}{d X^2} \bigg|_{X = 1} \mrm{Den}(X,L)_n 
    \\
    & = (\mrm{val}'(L^{\flat}) - 4 [\breve{F} : \breve{F}_0]^{-1}) \partial \mrm{Den}_{L^{\flat}}(x')_n + [\breve{F} : \breve{F}_0] \frac{d^2}{d X^2} \bigg|_{X = 1} \mrm{Den}(X,L')_n + \partial \mrm{Den}^*(L^{\flat})_n. \notag
    \end{align}
Again using $\varepsilon(V) = -1$, we apply the functional equation for $\mrm{Den}(X,L)$ \eqref{equation:part_I:Eisenstein:local_functional_equations:non-Arch:non-star_density} to find
    \begin{align}
    \frac{d^2}{d X^2} \bigg|_{X = 1} \mrm{Den}(X,L)_n & = (\mrm{val}(L) - 1) \frac{d}{d X}\bigg|_{X = 1} \mrm{Den}(X,L)_n 
    \\
    & = - (\mrm{val}(L) - 1) [\breve{F} : \breve{F}_0]^{-1} \partial \mrm{Den}_{L^{\flat}}(x)_n
    \end{align}
(the second equality is by definition) and similarly for $L'$. We also have
    \begin{equation}
    \mrm{val}'(L^{\flat}) = \mrm{val}(L) - 2 [\breve{F} : \breve{F}_0]^{-1} \mrm{val}''(x) - 1
    \quad \quad
    \mrm{val}(L) = \mrm{val}(L') + 2 [\breve{F} : \breve{F}_0]^{-1}.
    \end{equation}
Substituting all displayed equations into \eqref{equation:proposition:non-Arch_identity:limits:nonsplit} proves the claim.

\emph{Case $F / F_0$ is split:} Evaluating the induction formula from Proposition \ref{proposition:non-Arch_identity:induction_formula} at $X = 1$ yields
    \begin{equation}
    \mrm{Den}^*(L^{\flat})_n = \mrm{Den}(1,L)_n - \mrm{Den}(1, L')_n.
    \end{equation}
Multiplying both sides of the induction formula by $X^{-\mrm{val}(L^{\flat}) / 2}$ and taking one derivative at $X = 1$, we find
    \begin{align}\label{equation:proposition:non-Arch_identity:limits:split}
    & \mrm{val}(L^{\flat}) \mrm{Den}(1,L)_n - 2 \frac{d}{d X} \bigg |_{X = 1} \mrm{Den}(X,L)_n
    \\
    & = (\mrm{val}(L^{\flat}) - 2) \mrm{Den}(1, L')_n - 2 \frac{d}{d X} \bigg |_{X = 1} \mrm{Den}(X, L')_n + \partial \mrm{Den}^*(L^{\flat})_n. \notag
    \end{align}
    
The functional equation \eqref{equation:part_I:Eisenstein:local_functional_equations:non-Arch:non-star_density} implies
    \begin{equation}
    2 \frac{d}{d X} \bigg|_{X = 1} \mrm{Den}(X,L)_n = \mrm{val}(L) \mrm{Den}(1,L)_n
    \end{equation}
and similarly for $L'$. We also have
    \begin{equation}
    \mrm{val}(L^{\flat}) = \mrm{val}(L) - \mrm{val}(x) \quad \quad \mrm{val}(L) = \mrm{val}(L') + 1.
    \end{equation}
Substituting all displayed equations into \eqref{equation:proposition:non-Arch_identity:limits:split} proves the claim.
\end{proof}

\begin{corollary}\label{corollary:non-Arch_identity:limits}
Let $M^{\flat} \subseteq L^{\flat}_F$ be any full rank integral lattice with $t(M^{\flat}) \leq 1$. If $F / F_0$ is nonsplit, the following formulas hold.
    \begin{enumerate}[(1)]
        \item $\partial \mrm{Den}^*_{\ms{V}}(L^{\flat})_n = 2 [ \breve{F} : \breve{F}_0]^{-1} \lim_{x \ra 0} \partial \mrm{Den}_{L^{\flat}, \ms{V}}(x)_n$
        \item $\partial \mrm{Den}^*_{\ms{H}}(L^{\flat})_n = 2 [ \breve{F} : \breve{F}_0]^{-1} \lim_{x \ra 0} \left ( \partial \mrm{Den}_{L^{\flat}, \ms{H}}(x')_n - \mrm{val}''(x) \mrm{Den}^*(L^{\flat})_n \right )$
        \item $\partial \mrm{Den}^*_{\ms{H}}(M^{\flat})^{\circ}_n = 2 [ \breve{F} : \breve{F}_0]^{-1} \lim_{x \ra 0} \left ( \partial \mrm{Den}_{M^{\flat}, \ms{H}}(x')^{\circ}_n - \mrm{val}''(x) \mrm{Den}^*(M^{\flat})^{\circ}_n \right )$.
    \end{enumerate}
If $F / F_0$ is split, the following formulas hold.
    \begin{enumerate}[(1)]
    \item $\partial \mrm{Den}^*_{\ms{V}}(L^{\flat})_n = \lim_{x \ra 0} \mrm{Den}_{L^{\flat}, \ms{V}}(x)_n$
    \item $\partial \mrm{Den}^*_{\ms{H}}(L^{\flat})_n = \lim_{x \ra 0} \left ( \mrm{Den}_{L^{\flat}, \ms{H}}(x')_n - \mrm{val}(x) \mrm{Den}^*(L^{\flat})_n \right )$
    \item $\partial \mrm{Den}^*_{\ms{H}}(M^{\flat})^{\circ}_n = \lim_{x \ra 0} \left ( \mrm{Den}_{M^{\flat}, \ms{H}}(x')^{\circ}_n - \mrm{val}(x) \mrm{Den}^*(M^{\flat})^{\circ}_n \right )$.
    \end{enumerate}
All limits stabilize for $\mrm{val}(x) \gg 0$. The expressions (1) and (2) are $0$ if $L^{\flat}$ is not integral. If $L^{\flat}$ is integral and $F / F_0$ is nonsplit (resp. split), then the limits in (1) and (2) stabilize when $\mrm{val}(x) > a_{\mrm{max}}(L^{\flat})$ (resp. $\mrm{val}(x) > 2 a_{\mrm{max}}(M^{\flat})$). If $F / F_0$ is nonsplit (resp. split) , the limits in (3) stabilizes when $\mrm{val}(x) > a_{\mrm{max}}(M^{\flat})$ (resp. $\mrm{val}(x) > 2 a_{\mrm{max}}(M^{\flat})$).
\end{corollary}
\begin{proof}
Denote the result of Proposition \ref{proposition:non-Arch_identity:limits} as (0).
We have (3) $\implies$ (2) (in all cases, nonsplit or split), by summing over $M^{\flat}$ containing $L^{\flat}$. We have (0) $\implies$ (3) by taking $L^{\flat} = M^{\flat}$ and inducting on $\mrm{val}(M^{\flat})$ (starting with the base cases of $M^{\flat}$ being maximal integral (still with $t(M^{\flat}) \leq 1$), in which case $\partial \mrm{Den}^*(M^{\flat})_n = \partial \mrm{Den}^*_{\ms{H}}(M^{\flat})_n = \partial \mrm{Den}^*_{\ms{H}}(M^{\flat})^{\circ}_n$, and similarly for $\mrm{Den}^*(M^{\flat})_n$, as well as $\partial \mrm{Den}_{M^{\flat}}(x')_n$ (nonsplit) and $\mrm{Den}_{M^{\flat}}(x')_n$ (split)).
Since (0) = (1) $+$ (2), we conclude that (0) $\implies$ (1) as well.
\end{proof}

The following lemma is the geometric counterpart of Corollary \ref{corollary:non-Arch_identity:limits}(1) (in the special case when $\a = {}^{\mbb{L}} \mc{Z}(L^{\flat})_{\ms{V}}$ for a non-degenerate Hermitian $\mc{O}_F$-lattice $L^{\flat} \subseteq \mbf{W}$).

\begin{lemma}\label{lemma:non-Arch_identity:limits:vertical_geometric}
Take $F_0 = \Q$, and assume $p \neq 2$ if $F / \Q_p$ is ramified. Let $\mc{Z} \ra \Spec \overline{k}$ be a proper scheme equipped with a closed immersion $\mc{Z} \hookrightarrow \mc{N}$. Given any $\a \in \mrm{gr}_1 K'_0(\mc{Z})$, we have
    \begin{equation}
    \deg_{\overline{k}}( \a \cdot \mc{E}^{\vee}) = \lim_{w \ra 0} \deg_{\overline{k}}( \a \cdot {}^{\mbb{L}} \mc{Z}(w))
    \end{equation}
where the limit runs over $w \in \mbf{W}$. The limit stabilizes for $w$ satisfying $\mc{Z} \subseteq \mc{Z}(w)$.
\end{lemma}
\begin{proof}
We may assume $F / \Q_p$ is nonsplit, as otherwise $\mrm{gr}_1 K'_0(\mc{Z}) = 0$ (Section \ref{ssec:moduli_pDiv:discrete_reduced}) for dimension reasons so the lemma is trivial.

For any fixed $w \in \mbf{W}$, there exists $e \gg 0$ such that $\mc{Z} \subseteq \mc{Z}(p^e w)$ (over a quasi-compact base scheme, $p^e$ times any quasi-homomorphism is a homomorphism for $e \gg 0$). Hence $\mc{Z} \subseteq \mc{Z}(w)$ for all $w \in \mbf{W}$ lying in a sufficiently small neighborhood of $0$.

Assume $w \in \mbf{W}$ is such that $\mc{Z} \subseteq \mc{Z}(w)$. Write $\mc{I}(w) \subseteq \mc{O}_{\mc{N}}$ for the ideal sheaf of $\mc{Z}(w)$ (recall that $\mc{Z}(w)$ is a Cartier divisor, see Section \ref{ssec:moduli_pDiv:horizontal_vertical}). 
The lemma now follows from the ``linear invariance'' argument in the proof of \cite[{Lemma 2.55(3)}]{LL22II} (valid in the inert case as well, using \cite{Howard19}). Alternatively, the proof of linear invariance (particularly \cite[{Definition 4.2}]{Howard19} (inert) \cite[{Lemma 2.39}]{LL22II} (ramified)) exhibits a canonical isomorphism $\mc{E} \otimes \mc{O}_{Z(w)} \cong \mc{I}(w) \otimes \mc{O}_{Z(w)}$ via Grothendieck--Messing theory.
\end{proof}

\begin{proof}[Proof of Theorem \ref{theorem:non-Arch_identity:statement}]
The horizontal part of the theorem was already verified in Proposition \ref{proposition:non-Arch_identity:horizontal_identity}, so it remains to show $\mrm{Int}_{\ms{V}}(L^{\flat})_n = \partial \mrm{Den}^*_{\ms{V}}(L^{\flat})_n$.

If $F / \Q_p$ is split, then ${}^{\mbb{L}} \mc{Z}(L^{\flat})_{\ms{V}} = 0$ and so $\mrm{Int}_{\ms{V}}(L^{\flat})_n = 0$. Applying Corollary \ref{corollary:non-Arch_identity:limits}(1) with $V = \mbf{V}$, we find $\partial \mrm{Den}^*_{\ms{V}}(L^{\flat})_n = 0$ since $\mrm{Den}_{L^{\flat}, \ms{V}}(x)_n = 0$ for all $x$ (Lemma \ref{lemma:non-Arch_identity:limits:split_vertical_den_vanishes}).

Next assume $F / \Q_p$ is nonsplit. 
For any $w \in \mbf{W}$ not in $L^{\flat}_F$, we have $\deg_{\overline{k}}({}^{\mbb{L}} \mc{Z}(L^{\flat})_{\ms{V}} \cdot {}^{\mbb{L}} \mc{Z}(w)) = \partial \mrm{Den}_{L^{\flat}, \ms{V}}(w)_n$ by \cite[{Theorem 8.2.1}]{LZ22unitary} (inert) and \cite[{Theorem 2.7}]{LL22II} (the ``vertical'' parts of the main results of loc. cit..).
Lemma \ref{lemma:non-Arch_identity:limits:vertical_geometric} implies $\mrm{Int}_{\ms{V}}(L^{\flat})_n = 2 [\breve{F} : \breve{F}_0]^{-1} \lim_{w \ra 0} \partial \mrm{Den}_{L^{\flat}, \ms{V}}(w)_n$. Restricting to $w$ perpendicular to $L^{\flat}_F$, the limiting formula in Corollary \ref{corollary:non-Arch_identity:limits}(1) now implies $\mrm{Int}_{\ms{V}}(L^{\flat})_n = \partial \mrm{Den}^*_{\ms{V}}(L^{\flat})_n$.
\end{proof}

\begin{remark}\label{remark:non-Arch_identity:limits:horizontal_geometric}
Suppose $F_0 = \Q_p$, suppose $F / \Q_p$ is nonsplit, and assume $p \neq 2$. Let $M^{\flat} \subseteq \mbf{V}$ be a non-degenerate integral $\mc{O}_F$-lattice of rank $n - 1$ with $t(M^{\flat}) \leq 1$. As above, let $\mc{Z}(M^{\flat})^{\circ} \subseteq \mc{N}$ be the associated quasi-canonical lifting cycle. 

For any nonzero $w \in \mbf{V}$ not in $M^{\flat}_F$, we have $\deg_{\overline{k}}(\mc{Z}(M^{\flat})^{\circ} \cap \mc{Z}(w)) = \partial \mrm{Den}_{M^{\flat}, \ms{H}}(w)^{\circ}_n$ by \cite[{Proposition 8.4}]{KR11} (inert, see also \cite[{Corollary 5.4.6, Theorem 6.1.3}]{LZ22unitary}) and \cite[{Corollary 2.46}]{LL22II} (ramified), i.e. the ``horizontal'' parts of the main results of loc. cit..

The ``horizontal part'' of our main theorem showed $\mrm{Int}_{\ms{H}}(M^{\flat})^{\circ}_n = \partial \mrm{Den}^*_{\ms{H}}(M^{\flat})^{\circ}_n$ (Proposition \ref{proposition:non-Arch_identity:horizontal_identity}). Using also the special value formula in Lemma \ref{lemma:non-Arch_identity:statement:special_value}, our limiting result Corollary \ref{corollary:non-Arch_identity:limits}(3) is equivalent to the geometric statement
    \begin{align}
    & 2 \deg \mc{Z}(M^{\flat})^{\circ} \cdot \delta_{\mrm{tau}}(\mrm{val}'(M^{\flat})) \notag
    \\
    & = 2 [\breve{F} : \breve{\Q}_p]^{-1} \lim_{x \ra 0} \left ( \deg_{\overline{k}}(\mc{Z}(M^{\flat})^{\circ}_n \cap \mc{Z}(x)) - \mrm{val}''(\varpi x) \deg \mc{Z}(M^{\flat})^{\circ}\right ) \label{equation:non-Arch_identity:limits:horizontal_geometric}
    \end{align}
(limiting over nonzero $x$ perpendicular to $M^{\flat}$) where $\delta_{\mrm{tau}}(\mrm{val}'(M^{\flat}))$ is the ``local change of tautological height'', as in \eqref{equation:can_and_qcan:qcan:local_change_taut} (which is $-1/2$ times the ``local change of Faltings height'' $\delta_{\mrm{Fal}}(\mrm{val}'(M^{\flat}))$).

To prove our main theorem, we verified \eqref{equation:non-Arch_identity:limits:horizontal_geometric} indirectly by the computations in Section \ref{ssec:non-Arch_identity:horizontal_identity}. Direct computations are also possible.
\end{remark}

    \clearpage


    \part{Local change of heights}
    \label{part:local_change_heights}
        Fix an imaginary quadratic field $F / \Q$. Write $\Delta$ for the discriminant and $\s$ for the nontrivial involution. We allow $2 \mid \Delta$ in Sections \ref{sec:Faltings_and_taut} and \ref{sec:qcan_heights} unless otherwise specified. We set $F_p \coloneqq F \otimes_{\Q} \Q_p$ and $\mc{O}_{F_p} \coloneqq \mc{O}_F \otimes_{\Z} \Z_p$.

Throughout Sections \ref{sec:Faltings_and_taut} and \ref{sec:qcan_heights}, we write $E$ for a number field, with ring of integers $\mc{O}_E$. Given a prime $p$, we set $\mc{O}_{E, (p)} \coloneqq \mc{O}_E \otimes_{\Z} \Z_{(p)}$. We use $\breve{E}$ to denote a finite degree field extension of $\breve{\Q}_p$, with ring of integers $\mc{O}_{\breve{E}}$. We write $\mf{d}_p \subseteq \mc{O}_{F_p}$ for the different ideal of $\mc{O}_{F_p/\Q_p}$. We abuse notation and also mean $\mf{d}_p \coloneqq \sqrt{\Delta}$, which is a generator of the different ideal.

By a \emph{place} $\breve{w}$ of $E \otimes_{\Q} \breve{\Q}_p$, we mean a prime ideal of $E \otimes_{\Q} \breve{\Q}_p$ (equivalently, an element of $\Hom_{\breve{\Q}_p}(E \otimes_{\Q} \breve{\Q}_p, \C_p)$ up to automorphisms of $\C_p$). We write $\breve{E}_{\breve{w}}$ for the residue field of $\breve{w}$, with ring of integers $\mc{O}_{\breve{E}_{\breve{w}}}$. We use the shorthand $\breve{w} \mid p$ to indicate a place of $E \otimes_{\Q} \breve{\Q}_p$, and may use subscripts (e.g. $X'_{\breve{w}}$ and $\phi_{\breve{w}}$ in Section \ref{ssec:qcan_heights:minimal_isogenies}) to indicate base-change from $\Spec \mc{O}_{E,(p)}$ to $\Spf \mc{O}_{\breve{E}_{\breve{w}}}$.

Whenever an $\mc{O}_F$-action or $F$-action is mentioned (e.g. on a sheaf of modules on $\Spec \mc{O}_E$), we assume that $\mc{O}_E$ (resp. $\mc{O}_{\breve{E}}$) is equipped with morphism $\mc{O}_F \ra \mc{O}_E$ (resp. $\mc{O}_F \ra \mc{O}_{\breve{E}}$).

We write $\mf{X}_s$ for a level $s \geq 0$ quasi-canonical lifting of signature $(1,0)$ over $\Spec \mc{O}_{\breve{E}}$ with its $\mc{O}_{F_p}$-action $\iota_{\mf{X}_s}$, as explained in Section \ref{sec:can_and_qcan}. The framing $\rho_{\mf{X}_s}$ of loc. cit. is unimportant in Sections \ref{sec:Faltings_and_taut} and \ref{sec:qcan_heights} (and will be omitted). As before, the notation $\mf{X}_s^{\s}$ means $\mf{X}_s$ but with $\mc{O}_{F_p}$-action pre-composed by $\s$.

Given a group scheme $G$ over a base $S$, we typically write $e \colon S \ra G$ for the identity section. We abuse notation and use ``$e$'' simultaneously for different group schemes.

        \section{Faltings and ``tautological'' heights}
        \label{sec:Faltings_and_taut}
        
            \subsection{Heights}
            \label{ssec:Faltings_and_taut:heights}
                Suppose $A \ra \Spec \mc{O}_E$ is a semi-abelian N\'eron model of an abelian variety over $E$.
The \emph{Faltings height} of $A$ (or its generic fiber $A_E$) is
    \begin{equation}
    h_{\mrm{Fal}}(A_E) \coloneqq h_{\mrm{Fal}}(A) \coloneqq \frac{1}{[E : \Q]} \widehat{\deg}(\widehat{\omega}_{A})
    \end{equation}
where $\widehat{\omega}_{A} = (\omega_{A}, \norm{-}) = (e^* \bigwedge^n \Omega^1_{A/\mc{O}_E}, \norm{-})$ is the Hermitian line bundle with norm $\norm{-}$ normalized as in \eqref{equation:arith_cycle_classes:Hodge_bundles:Faltings_metric}. The usual arithmetic degree $\widehat{\deg}$ was recalled in \cref{ssec:part_I:arith_intersections:Hermitian_bundles}. Any abelian variety over a number field has everywhere potentially semi-abelian reduction, and the \emph{Faltings height} of any abelian variety $B$ over $\Spec E$ is defined so that $h_{\mrm{Fal}}(B)$ is remains constant under finite field extensions $E \ra E'$. (We only consider stable Faltings height, as defined above.)

We also consider certain ``tautological heights'' to describe the arithmetic intersections numbers appearing in Equation \eqref{equation:intro:results:if_proper} (middle term). The terminology we introduce for this (e.g. ``Kr\"amer datum") is likely nonstandard.

\begin{definition}\label{definition:Falting_and_taut:heights:Kramer_datum}
\hfill
\begin{enumerate}[(1)]
    \item Given a scheme $S$ over $\Spec \mc{O}_F$, a \emph{Kr\"amer datum} (of signature $(n - 1, 1)$) is a tuple $(A, \iota, \mc{F})$ where $A \ra S$ is an abelian scheme, where $\iota \colon \mc{O}_F \ra \End(A)$ an action of signature $(n - 1, 1)$, and where $\mc{F} \subseteq \Lie A$ is a $\iota$-stable local direct summand of rank $n - 1$ such that the $\mc{O}_F$ action via $\iota$ on $\mc{F}$ (resp. $(\Lie A) / \mc{F}$) is $\mc{O}_F$-linear (resp. $\s$-linear). We say that $\mc{F}$ is the associated \emph{Kr\"amer hyperplane}.
    
    \item Given a formal scheme $S$ over $\Spf \mc{O}_{F_p}$, a \emph{local Kr\"amer datum} (of signature $(n - 1,1)$) is a tuple $(X, \iota, \mc{F})$ where $X$ is a $p$-divisible group over $S$ of height $2n$ and dimension $n$, where $\iota \colon \mc{O}_{F_p} \ra \End(X)$ is an action of signature $(n - 1, 1)$, and where $\mc{F} \subseteq \Lie X$ is a $\iota$-stable local direct summand of rank $n - 1$ such that the $\mc{O}_F$ action via $\iota$ on $\mc{F}$ (resp. $(\Lie X) / \mc{F}$) is $\mc{O}_F$-linear (resp. $\s$-linear). We say that $\mc{F}$ is the associated \emph{Kr\"amer hyperplane}.
    
    \item A \emph{quasi-polarized} Kr\"amer datum (resp. \emph{quasi-polarized} local Kr\"amer datum) is a tuple $(A, \iota, \lambda, \mc{F})$ (resp. $(X, \iota, \lambda, \mc{F})$) where $(A, \iota, \lambda)$ is a Hermitian abelian scheme (\crefext{definition:ab_var:integral_models:Hermitian_abelian_scheme}) (resp. $(X, \iota, \lambda)$ is a Hermitian $p$-divisible group (Definition \ref{definition:moduli_pDiv:RZ:Hermitian_p-divisible}, but we allow $p = 2$ even if $F / \Q_p$ is ramified)) and $(A, \iota, \mc{F})$ is a Kr\"amer datum (resp. $(X, \iota, \mc{F})$ is a Kr\"amer datum).
\end{enumerate}
\end{definition}

The name ``Kr\"amer datum'' refers to the Kr\"amer model, as studied (locally) in \cite{Kramer03}.
For an understood Kr\"amer datum $(A, \iota, \lambda, \mc{F})$, we will use the shorthand $\ms{E}^{\vee} \coloneqq (\Lie A)/\mc{F}$ (cf. the ``tautological bundles'' of \cref{definition:ab_var:exotic_smooth:tautological_bundle}). We use the same notation $\ms{E}^{\vee} \coloneqq (\Lie X) / \mc{F}$ given an understood local Kr\"amer datum $(X, \iota, \mc{F})$. In both cases, the sheaf $\ms{E}^{\vee}$ is locally free of rank $1$, and we call it the associated \emph{Kr\"amer hyperplane quotient}.

\begin{definition}
A \emph{morphism} (resp. \emph{isogeny}) of Kr\"amer data $(A_1, \iota_1, \mc{F}_1) \ra (A_2, \iota_2, \mc{F}_2)$ is an $\mc{O}_F$-homomorphism (resp. isogeny) $A_1 \ra A_2$ such that $\mrm{im}(\mc{F}_1) \subseteq \mc{F}_2$, where $\mrm{im}(\mc{F}_1)$ is the image of $\mc{F}_1$ under $\Lie A_1 \ra \Lie A_2$. A \emph{morphism} (resp. \emph{isogeny}) of local Kr\"amer data is defined in the same way.
\end{definition}

\begin{lemma}\label{lemma:Faltings_and_taut:height:canonical_Kramer_data}
Let $S$ be a scheme over $\Spec \mc{O}_F$. Assume either that $S$ is a scheme over $\Spec \mc{O}_F[1/\Delta]$ or that $S = \Spec R$ where $R$ is a Dedekind domain with fraction field of characteristic $0$.
\begin{enumerate}[(1)]
    \item Suppose $A \ra S$ is an abelian scheme with an action $\iota \colon \mc{O}_F \ra \End(X)$ of signature $(n - 1, 1)$. Then the pair $(A, \iota)$ extends uniquely to a Kr\"amer datum $(A, \iota, \mc{F})$ over $S$.
    \item Given pairs $(A_1, \iota_1)$ and $(A_2, \iota_2)$ as above, any $\mc{O}_F$-linear homomorphism (resp. isogeny) $A_1 \ra A_2$ induces a morphism (resp. isogeny) of Kr\"amer data.
    \item If $S$ is a scheme over $\Spec \mc{O}_F[1 / \Delta]$, the exact sequence
    \begin{equation}
    0 \ra \mc{F} \ra \Lie A \ra \ms{E}^{\vee} \ra 0
    \end{equation}
has a unique $\mc{O}_F$-linear splitting.
\end{enumerate}
\end{lemma}
\begin{proof}
If $S$ is a scheme over $\Spec \mc{O}_F[1/\Delta]$, the claims hold because there is a unique decomposition $\Lie A = (\Lie A)^+ \oplus (\Lie A)^-$ characterized by $\iota$ acting $\mc{O}_F$-linearly on the rank $n - 1$ subbundle $(\Lie A)^+$ (resp. $\s$-linearly on the rank $1$ bundle $(\Lie A)^-$).

Suppose instead that $S = \Spec R$ is a Dedekind domain with fraction field $K$ of characteristic $0$. By localizing, it is enough to verify the lemma when $R$ is a discrete valuation ring. Then part (1) amounts to the following fact: given a finite free $R$-module $M$ and any $K$-subspace $W \subseteq M \otimes K$, there is a unique summand $M' \subseteq M$ such that $W = M' \otimes K$ (namely $M' = M \cap W$; note that $M' \subseteq M$ is a saturated sublattice). We are applying this when $M = \Lie A$ and $W = (\Lie A \otimes K)^+$, in the notation above (and taking $\mc{F} = M'$). The signature $(n - 1, 1)$ condition forces the $\mc{O}_F$-action on $(\Lie A) / \mc{F}$ to be $\s$-linear. 
These considerations also verify the claim in part (2) (since it holds in the generic fiber).
\end{proof}

\begin{lemma}\label{lemma:Faltings_and_taut:height:canonical_local_Kramer_data}
Let $S$ be a formal scheme over $\Spf \mc{O}_{F_p}$. Assume either that $p$ is unramified in $\mc{O}_F$ or that $S = \Spf R$ for an adic ring which is a Dedekind domain with fraction field of characteristic $0$. Then the following conclusions hold.
    \begin{enumerate}[(1)]
    \item Suppose $X$ is a $p$-divisible group of height $2n$ over $S$ with an action $\iota \colon \mc{O}_{F_p} \ra \End(X)$ of signature $(n - 1, 1)$. Then the pair $(X, \iota)$ extends uniquely to a local Kr\"amer datum $(X, \iota, \mc{F})$ over $S$.
    \item Given pairs $(X_1, \iota_1)$ and $(X_2, \iota_2)$ as above, any $\mc{O}_{F_p}$-linear homomorphism (resp. isogeny) $X_1 \ra X_2$ induces a morphism (resp. isogeny) of local Kr\"amer data.
    \item If $p$ is unramified, the exact sequence
        \begin{equation}
        0 \ra \mc{F} \ra \Lie A \ra \ms{E}^{\vee} \ra 0
        \end{equation}
    has a unique $\mc{O}_{F_p}$-linear splitting.
    \end{enumerate}
\end{lemma}
\begin{proof}
This may be proved in the same way as Lemma \ref{lemma:Faltings_and_taut:height:canonical_Kramer_data}.
If $S = \Spf R$ for $R$ a Dedekind domain with fraction field of characteristic $0$, note that $R$ must be a complete discrete valuation ring.
\end{proof}

In the situations of Lemma \ref{lemma:Faltings_and_taut:height:canonical_Kramer_data} and \ref{lemma:Faltings_and_taut:height:canonical_Kramer_data}, we also use the alternative terminology \emph{dual tautological bundle} for the Kr\"amer hyperplane quotient $\ms{E}^{\vee}$.

If $(A, \iota, \lambda)$ is a Hermitian abelian scheme of signature $(n - 1, 1)$ over $\Spec \mc{O}_E$ with associated quasi-polarized Kr\"amer datum $(A, \iota, \lambda, \mc{F})$, we thus obtain a Hermitian line bundle $\widehat{\ms{E}^{\vee}} = (\ms{E}^{\vee}, \norm{-})$ on $\Spec \mc{O}_E$ as follows: the metric $\norm{-}$ is given by restricting the metric on $\Lie A$ induced by $\lambda$ (which we take to be normalized as in \eqref{equation:arith_cycle_classes:Hodge_bundles:4_pi_gamma_normalization}) along the $\mc{O}_F$-linear splitting $\ms{E}^{\vee}[1/\Delta] \hookrightarrow (\Lie A)[1/\Delta]$ (where $(-)[1/\Delta]$ means restriction to $\Spec \mc{O}_E[1/\Delta]$). We say $\widehat{\ms{E}}^{\vee}$ is the associated \emph{metrized dual tautological bundle}. We also make the same construction over $\Spec \mc{O}_{E,(p)}$ and $\Spec E$.

\begin{definition}
Let $(A, \iota, \lambda)$ be a Hermitian abelian scheme of signature $(n - 1, 1)$ over $\Spec \mc{O}_E$. The associated \emph{tautological height} is
    \begin{equation}
    h_{\mrm{tau}}(A_E) \coloneqq h_{\mrm{tau}}(A) \coloneqq \frac{1}{[E : \Q]} \widehat{\deg}(\widehat{\ms{E}}^{\vee}).
    \end{equation}
\end{definition}

The tautological height depends on the auxiliary data in the definition (not just $A_E$ or $A$), which we have suppressed from notation. If $(A, \iota, \lambda)$ is a Hermitian abelian scheme of signature $(n - 1, 1)$ over $\Spec E$ such that $A$ has everywhere potentially good reduction, we define the \emph{tautological height} $h_{\mrm{tau}}(A)$ so that it is invariant under finite degree field extension $E \ra E'$.

\begin{remark}\label{remark:Faltings_and_taut:heights:punctured_heights}
If we instead work over $\mc{O}_E[1/N]$ for some integer $N \geq 1$, we may define \emph{Faltings height} and \emph{tautological height} as above, but where $\widehat{\deg}$ now takes values in $\R_N \coloneqq \R / \sum_{p \mid N} \Q \cdot \log p$ (as explained in Section \ref{ssec:part_I:arith_intersections:Hermitian_bundles}).
\end{remark}
    
            \subsection{Change along global isogenies}
            \label{ssec:Faltings_and_taut:isogeny_change_global}
                Let $A_1 \ra \Spec \mc{O}_E$ and $A_2 \ra \Spec \mc{O}_E$ be semi-abelian N\'eron models of abelian varieties over $E$. We have
    \begin{equation}
    [E : \Q] (h_{\mrm{Fal}}(A_2) - h_{\mrm{Fal}}(A_1)) = \widehat{\deg}(\widehat{\omega}_{A_2}) - \widehat{\deg}(\widehat{\omega}_{A_1}) = -\widehat{\deg}(\underline{\Hom}(\widehat{\omega}_{A_2}, \widehat{\omega}_{A_1})).
    \end{equation}
Any isogeny $\phi \colon A_1 \ra A_2$ defines a section $\phi$ of the Hermitian line bundle $\underline{\Hom}(\widehat{\omega}_{A_2}, \widehat{\omega}_{A_1})$, which gives
    \begin{align}\label{equation:Faltings_and_taut:isogeny_change_global:Faltings}
    h_{\mrm{Fal}}(A_2) - h_{\mrm{Fal}}(A_1) & = \frac{1}{[E : \Q]} \left ( \log \norm{\phi}_{\infty} + \sum_{v < \infty} \log \norm{\phi}_v \right )
    \\
    & = \frac{1}{2} \log(\deg \phi) - \frac{1}{[E : \Q]} \log|e^* \Omega^1_{\ker \phi / \mc{O}_E}|
    \end{align}
(sum is over places $v$ of $E$) as in \cite[{Lemma 5}]{Faltings86}, where $|e^* \Omega^1_{\ker \phi / \mc{O}_E}|$ denotes the cardinality of the finite length $\mc{O}_E$-module $e^* \Omega^1_{\ker \phi / \mc{O}_E}$. Note $|e^*\Omega^1_{\ker \phi / \mc{O}_E}| = |\coker(\phi^* \colon \omega_{A_2} \ra \omega_{A_1})| = |\coker(\phi_* \colon \Lie A_1 \ra \Lie A_2)|$. Also note
    \begin{equation}\label{equation:Faltings_and_taut:isogeny_change_global:Faltings:log_p}
    h_{\mrm{Fal}}(A_2) - h_{\mrm{Fal}}(A_1) = \sum_p a_p \log p = \sum_{p \mid \deg \phi} a_p \log p
    \end{equation}
for some $a_p \in \Q$ independent of $\phi$.

Given Hermitian abelian schemes $(A_1, \iota_1, \lambda_1)$ and $(A_2, \iota_2, \lambda_2)$ of signature $(n - 1, 1)$ over $\Spec \mc{O}_E$ with associated Hermitian line bundles $\widehat{\ms{E}}^{\vee}_1$ and $\widehat{\ms{E}}^{\vee}_2$, we similarly have
    \begin{equation}
    h_{\mrm{tau}}(A_2) - h_{\mrm{tau}}(A_1) = \frac{1}{[E : \Q]} \widehat{\deg}(\underline{\Hom}(\widehat{\ms{E}}_{1}^{\vee}, \widehat{\ms{E}}_{2}^{\vee})).
    \end{equation}
Any $\mc{O}_F$-linear isogeny $\phi \colon A_1 \ra A_2$ defines a section $\phi$ of the Hermitian line bundle $\underline{\Hom}(\widehat{\ms{E}}_{1}^{\vee}, \widehat{\ms{E}}_{2}^{\vee})$, and we have
    \begin{align}\label{equation:Faltings_and_taut:isogeny_change_global:taut}
    h_{\mrm{tau}}(A_2)- h_{\mrm{tau}}(A_1) & = \frac{1}{[E : \Q]} \left ( - \log \norm{\phi}_{\infty} - \sum_{v < \infty} \log \norm{\phi}_v \right ) 
    \\
    & = \frac{1}{[E : \Q]} \left ( - \log \norm{\phi}_{\infty} +  \log |\coker(\phi_* \colon \ms{E}_1^{\vee} \ra \ms{E}_2^{\vee})| \right ).
    \end{align}
    
            \subsection{Change along local isogenies: Faltings}
            \label{ssec:Faltings_and_taut:isogeny_change_local_Faltings}
                Given an isogeny $\phi \colon A_1 \ra A_2$ of abelian schemes over $\Spec \mc{O}_{E,(p)}$, we define the \emph{semi-global change of Faltings height}
    \begin{equation}
    \delta_{\mrm{Fal},(p)}(\phi) \coloneqq - \frac{1}{2}\log|\deg \phi|_p - \frac{1}{[E : \Q]} \log|e^* \Omega^1_{\ker \phi / \mc{O}_{E,(p)}}|
    \end{equation}
where $|-|_p$ is the usual $p$-adic norm. We have $\delta_{\mrm{Fal},(p)}(\phi) \in \Q \cdot \log p$. The formula for change of Faltings height \eqref{equation:Faltings_and_taut:isogeny_change_global:Faltings} shows that $\delta_{\mrm{Fal},(p)}(\phi) = a_p \log p$, in the notation of \eqref{equation:Faltings_and_taut:isogeny_change_global:Faltings:log_p}.
In particular, $\delta_{\mrm{Fal},(p)}(\phi)$ does not depend on the choice of isogeny $\phi$ (and depends only on $A_1$ and $A_2$).
If $A_1$ and $A_2$ have everywhere potentially good reduction, we have
    \begin{equation}\label{equation:Faltings_and_taut:isogeny_change_local_Faltings:decomp_to_semi-global}
    h_{\mrm{Fal}}(A_{2,E}) - h_{\mrm{Fal}}(A_{1,E}) = \sum_{\ell} \delta_{\mrm{Fal},(\ell)}(\phi) = \sum_{\ell \mid \deg \phi} \delta_{\mrm{Fal},(\ell)}(\phi)
    \end{equation}
where $\phi$ also denotes the induced isogeny on N\'eron models over $\Spec \mc{O}_{E,(\ell)}$ for each prime $\ell$ (after enlarging $E$ if necessary).

Given any isogeny $\phi \colon X_1 \ra X_2$ of $p$-divisible groups over $\Spf \mc{O}_{\breve{E}}$, we have $(\Lie X_i)^{\vee} \cong e^* \Omega^1_{X_i[p^N]/\Spec \mc{O}_{\breve{E}}}$ (canonically) for $N \gg 0$ by \cite[{Corollary II.3.3.17}]{Messing72} (passing to the limit over $\mc{O}_{\breve{E}} / p^k \mc{O}_{\breve{E}}$ as $k \ra \infty$), so there is a canonical exact sequence
    \begin{equation}\label{equation:Faltings_and_taut:isogeny_change_local:coker_and_differentials}
    0 \ra (\Lie X_2)^{\vee} \xra{\phi^*} (\Lie X_1)^{\vee} \ra e^* \Omega^1_{\ker \phi / \Spec \mc{O}_{\breve{E}}} \ra 0
    \end{equation}
of finite free $\mc{O}_{\breve{E}}$-modules (note that $\Lie X_1 \ra \Lie X_2$ is injective, e.g. by \crefext{III:lemma:isogeny_criterion}). If $X_1$ and $X_2$ are moreover height $2n$ and dimension $n$, we define the \emph{local change of Faltings height}
    \begin{equation}
    \breve{\delta}_{\mrm{Fal}}(\phi) \coloneqq \frac{1}{2} \log(\deg \phi) - \frac{1}{[\breve{E} : \breve{\Q}_p]} \mrm{length}_{\mc{O}_{\breve{E}}}(e^* \Omega^1_{\ker \phi / \Spec \mc{O}_{\breve{E}}}) \cdot \log p.
    \end{equation}
We have $\breve{\delta}_{\mrm{Fal}}(\phi) = \Q \cdot \log p$, as well as
    \begin{align}
    & \breve{\delta}_{\mrm{Fal}}(\phi' \circ \phi) = \breve{\delta}_{\mrm{Fal}}(\phi') + \breve{\delta}_{\mrm{Fal}}(\phi) && \breve{\delta}_{\mrm{Fal}}([N]) = 0
    \end{align}
where $\phi' \colon X_2 \ra X_3$ is any isogeny of $p$-divisible groups and $[N] \colon X_1 \ra X_1$ is the multiplication-by-$N$ isogeny (follows from \eqref{equation:Faltings_and_taut:isogeny_change_local:coker_and_differentials}). 
Unlike $\delta_{\mrm{Fal},(p)}(-)$ from above, the quantity $\breve{\delta}_{\mrm{Fal}}(\phi)$ may depend on the isogeny $\phi$.

Given isogenous abelian schemes over $\Spec \mc{O}_{E,(p)}$ and an isogeny $\phi \colon X_1 \ra X_2$ of the associated $p$-divisible groups, set
    \begin{equation}\label{equation:Faltings_and_taut:isogeny_change_local_Faltings:pDiv_semi-global_def}
    \delta_{\mrm{Fal},(p)}(\phi) \coloneqq \frac{1}{[E : \Q]} \sum_{\breve{w} \mid p} [\breve{E}_{\breve{w}} : \breve{\Q}_p] \breve{\delta}_{\mrm{Fal}}(\phi_{\breve{w}}).
    \end{equation}
where $\phi_{\breve{w}}$ denotes the base-change of $\phi$ to $\Spf \mc{O}_{\breve{E}_{\breve{w}}}$.

\begin{lemma}\label{lemma:Faltings_and_taut:isogeny_change_local_Faltings:global-to-local}
Let $A_1, A_2$ be isogenous abelian schemes over $\Spec \mc{O}_{E,(p)}$. Let $X_i$ be the associated $p$-divisible groups. Given any isogenies $\tilde{\phi} \colon A_1 \ra A_2$ and $\phi \colon X_1 \ra X_2$, we have
    \begin{equation}\label{equation:Faltings_and_taut:isogeny_change_local:global-to-local}
    \delta_{\mrm{Fal},(p)}(\tilde{\phi}) = \delta_{\mrm{Fal},(p)}(\phi).
    \end{equation}
\end{lemma}
\begin{proof}
The lemma is clear if $\phi$ is the isogeny associated with $\tilde{\phi}$. If $\phi' \colon X_1 \ra X_2$ is another isogeny, we have $[p^N] \circ \phi = \phi' \circ \phi''$ for some isogeny $\phi'' \colon X_1 \ra X_1$ \crefext{III:lemma:isogeny_criterion}). By additivity of $\breve{\delta}_{\mrm{Fal}}$ and since $\breve{\delta}_{\mrm{Fal}}([p^N]) = 0$, it is enough to show $\delta_{\mrm{Fal},(p)}(\phi) = 0$ if $X_1 = X_2$. For this purpose, we may also assume $A_1 = A_2$.

Write $A \coloneqq A_1$ and $X \coloneq X_1$ to lighten notation. As usual, $A_E$ and $X_E$ denote the respective generic fibers (over $\Spec E$). We write $\mrm{Isog}(A)$ and $\mrm{Isog}(X)$ for the set of self-isogenies of $A$ and $X$.

We have canonical identifications
    \begin{equation}
    \End(X) = \End(X_E) = \End(T_p(X_E)).
    \end{equation}
The first equality holds by a theorem of Tate \cite[{Theorem 4}]{Tate67} (base-change along $\Spec E \ra \Spec \mc{O}_{E,(p)}$ is fully faithful on $p$-divisible groups) and the second equality holds because $X_E$ is an \'etale $p$-divisible group. Here, the notation $\End(T_p(X_E))$ means endomorphisms of $T_p(X_E)$ as a Galois module.

Equip the finite $\Z_p$-module $\End(T_p(X_E))$ with the $p$-adic topology, and give $\mrm{Isog}(X)$ the subspace topology.
We have $\delta_{\mrm{Fal},(p)}(\phi \circ \phi') = \delta_{\mrm{Fal},(p)}(\phi)$ for any $\phi \in \mrm{Isog}(X)$ and $\phi' \in \End(X)$ with $\phi' \equiv 1 \pmod{p}$, since any such $\phi'$ is an automorphism of $X$. The map $\mrm{Isog}(X) \ra \R$ given by $\phi \mapsto \delta_{\mrm{Fal},(p)}(\phi)$ is thus locally constant.

We also have canonical identifications
    \begin{equation}
    \End(A) \otimes_{\Z} \Z_p = \End(A_E) \otimes_{\Z} \Z_p = \End(T_p(A_E)) = \End(T_p(X_E)).
    \end{equation}
The first equality holds by the N\'eron mapping property, and the second equality holds by Faltings's theorem \cite[{\S 5 Corollary 1}]{Faltings86}. Hence $\mrm{Isog}(A) = \mrm{Isog}(X) \cap \End(A)$ is a dense subset of $\mrm{Isog}(X)$. Since $\delta_{\mrm{Fal},(p)}(\phi) = 0$ for any $\phi \in \mrm{Isog}(A)$, this proves the lemma.
\end{proof}

\begin{corollary}
In the situation of Lemma \ref{lemma:Faltings_and_taut:isogeny_change_local_Faltings:global-to-local}, the quantity $\delta_{\mrm{Fal},(p)}(\phi)$ does not depend on the choice of isogeny $\phi \colon X_1 \ra X_2$.
\end{corollary}
\begin{proof}
In the notation of the lemma, this follows immediately from $\tilde{\phi}$ independence of $\delta_{\mrm{Fal},(p)}(\tilde{\phi})$ (discussed above).
\end{proof}

We will use Lemma \ref{lemma:Faltings_and_taut:isogeny_change_local_Faltings:global-to-local} to compute Faltings heights without producing isogenies on abelian varieties, only isogenies on underlying $p$-divisible groups over $\Spec \mc{O}_{E,(p)}$. The analogous lemma for tautological height (Lemma \ref{lemma:Faltings_and_taut:isogeny_change_local_taut:global-to-local}) serves a similar purpose.
    
            \subsection{Change along local isogenies: tautological}
            \label{ssec:Faltings_and_taut:isogeny_change_local_taut}
                To locally decompose the change of tautological height along an isogeny, we impose an additional condition.

\begin{definition}\label{definition:Faltings_and_taut:isogeny_change_local_taut:special}
\hfill
\begin{enumerate}[(1)]
    \item A Hermitian abelian scheme $(A, \iota, \lambda)$ of signature $(n - 1, 1)$ over $E$ is \emph{special} if $A$ is $\mc{O}_F$-linearly isogenous to a product of elliptic curves, each with $\mc{O}_F$-action. A Hermitian abelian scheme of signature $(n - 1, 1)$ over $\Spec \mc{O}_E$ or $\Spec \mc{O}_{E,(p)}$ is \emph{special} if its generic fiber is special.
    \item A Hermitian $p$-divisible group $(X, \iota, \lambda)$ of signature $(n - 1, 1)$ over $\Spf \mc{O}_{\breve{E}}$ is \emph{special} if $X$ is $\mc{O}_{F_p}$-linearly isogenous to $\mf{X}_0^{n - 1} \times \mf{X}_0^{\s}$.
\end{enumerate}
\end{definition}

We only use the term ``special" this way in Sections \ref{sec:Faltings_and_taut} and \ref{sec:qcan_heights} (but we have global special cycles in mind, cf. \crefext{III:lemma:corank_1_isogeny_decomposition}). The norm $\norm{-}_{\infty}$ below is as in \eqref{equation:Faltings_and_taut:isogeny_change_global:taut}.

\begin{lemma}\label{lemma:Faltings_and_taut:isogeny_change_local_taut:infty_norm}
Let $(A_1, \iota_1, \lambda_1)$ and $(A_2, \iota_2, \lambda_2)$ be special Hermitian abelian schemes of signature $(n - 1, 1)$ over $\Spec E$. For any $\mc{O}_F$-linear isogeny $\phi \colon A_1 \ra A_2$, we have $\norm{\phi}_{\infty}^2 \in \Q_{>0}$.
\end{lemma}
\begin{proof}
Given such $\phi$, form a diagram
    \begin{equation}
    B_1 \times B_1^{\perp} \xra{\phi_1} A_1 \xra{\phi} A_2 \xra{\phi_2} B_2 \times B_2^{\perp}
    \end{equation}
where each $\phi_i$ is an $\mc{O}_F$-linear isogeny, each $B_i$ is a product of $(n
- 1)$ elliptic curves each with $\mc{O}_F$-action of signature $(1,0)$, and each $B_i^{\perp}$ is an elliptic curve with $\mc{O}_F$-action of signature $(0,1)$. Signature incompatibility implies that $\lambda_1$ pulls back to a diagonal quasi-polarization $\lambda_{B_1} \times \lambda_{B_1^{\perp}}$ on $B_1 \times B_1^{\perp}$ (e.g. $\Hom_{\mc{O}_F}(B_1^{\s}, B_1^{\perp \vee}) = \Hom_{\mc{O}_F}(B_1^{\perp \s}, B_1^{\vee}) = 0$). Similarly, $\lambda_2$ pulls back along the quasi-isogeny $\phi_2^{-1}$ to a diagonal quasi-polarization $\lambda_{B_2} \times \lambda_{B_2^{\perp}}$.

With these quasi-polarizations, we have $\norm{\phi_2 \circ \phi \circ \phi_1}_{\infty} = \norm{\phi_2}_{\infty} \norm{\phi}_{\infty} \norm{\phi_1}_{\infty} = \norm{\phi}_{\infty}$ since $\norm{\phi_1}_{\infty} = \norm{\phi_2}_{\infty} = 1$ (because $\phi_1$ and $\phi_2$ preserve quasi-polarizations, by construction).
On the other hand, if $\phi' \colon B_1^{\perp} \ra B_2^{\perp}$ is the induced isogeny (signature incompatibility again implies $\Hom_{\mc{O}_F}(B_1, B_2^{\perp}) = \Hom_{\mc{O}_F}(B_1^{\perp}, B_2) = 0$), we must have $\norm{\phi_2 \circ \phi \circ \phi_1}_{\infty} = \norm{\phi'}_{\infty}$ (the latter norm is taken with respect to $\lambda_{B_1^{\perp}}$ and $\lambda_{B_2^{\perp}}$). For each embedding $\t \colon E \ra \C$, the quantity $\norm{\phi'}^2_{\t}$ must be the element of $\Q_{>0}$ satisfying
    \begin{equation}
    \phi^{\prime *} \lambda_{B_2^{\perp}} = \norm{\phi'}_{\t}^2 \lambda_{B_1^{\perp}},
    \end{equation}
(quasi-polarizations on elliptic curves are unique up to $\Q_{>0}$ scalar), so we have $\norm{\phi'}_{\infty}^2 = \prod_{\t \colon E \ra \C} \norm{\phi'}_{\t}^2 \in \Q_{>0}$.
\end{proof}

For the rest of Section \ref{ssec:Faltings_and_taut:isogeny_change_local_taut}, we let $(A_i, \iota_i, \lambda_i)$ for $i = 1, 2$ be special Hermitian abelian schemes of signature $(n - 1, 1)$ over $\Spec \mc{O}_{E,(p)}$, with associated Kr\"amer hyerplanes $\mc{F}_i$ and dual tautological bundles $\ms{E}^{\vee}_i$. We also let $(X_i, \iota_i, \lambda_i)$ for $i = 1, 2, 3$ be special Hermitian $p$-divisible groups of signature $(n - 1, 1)$ over $\Spf \mc{O}_{\breve{E}}$, and reuse the notation $\mc{F}_i$ and $\ms{E}^{\vee}_i$ for the respective Kr\"amer hyerplanes and dual tautological bundles.

Given $(X_1, \iota_1, \lambda_1)$ and an $\mc{O}_{F_p}$-linear isogeny $Y_1 \times Y_1^{\perp} \ra X_1$ with $Y_1$ being a product of $n - 1$ canonical liftings of signature $(1,0)$ and $Y_1^{\perp}$ being a canonical lifting of signature $(0,1)$, there is an induced decomposition
    \begin{equation}\label{equation:Faltings_and_taug:isogeny_change_local_taut:T_p-decomp}
    T_p(X_1)^0 = T_p(Y_1)^0 \oplus T_p(Y_1^{\perp})^0
    \end{equation}
on rational Tate modules (of the generic fibers). Equip $Y_1 \times Y_1^{\perp}$ with the pullback of $\lambda_1$. This gives a product quasi-polarization $\lambda_{Y_1} \times \lambda_{Y_1^{\perp}}$ on $Y_1 \times Y_1^{\perp}$ (by signature incompatibility as in the abelian scheme case, i.e. $\Hom_{\mc{O}_{F_p}}(Y_1^{\s}, Y_1^{\perp \vee}) = \Hom_{\mc{O}_{F_p}}(Y_1^{\perp \s}, Y_1^{\vee}) = 0$).
Hence the decomposition in \eqref{equation:Faltings_and_taug:isogeny_change_local_taut:T_p-decomp} is orthogonal for the Hermitian pairing on $T_p(X_1)^0$.

Consider $(X_2, \iota_2, \lambda_2)$ with $\mc{O}_{F_p}$-linear isogeny $Y_2 \times Y_2^{\perp} \ra X_2$ as above and, and suppose $\phi \colon X_1 \ra X_2$ is an $\mc{O}_{F_p}$-linear isogeny. Then the induced map $\phi_* \colon T_p(X_1)^0 \ra T_p(X_2)^0$ sends $T_p(Y_1)^0$ to $T_p(Y_2)^0$ and similarly for $T_p(Y_i^{\perp})^0$ (again by signature incompatibility, i.e. $\Hom_{\mc{O}_{F_p}}(Y_1, Y_2^{\perp}) = \Hom_{\mc{O}_{F_p}}(Y_1^{\perp}, Y_2) = 0$). In particular, the decomposition in \eqref{equation:Faltings_and_taug:isogeny_change_local_taut:T_p-decomp} does not depend on the choice of $Y_1 \times Y_1^{\perp} \ra X_1$.

Any $\mc{O}_{F_p}$-linear isogeny $\phi \colon X_1 \ra X_2$ thus gives a nonzero element $\phi \in \Hom_{F_p}(T_p(Y_1^{\perp})^0, T_p(Y_2^{\perp})^0)$. We then set
    \begin{equation}\label{equation:Faltings_and_taut:isogeny_change_local_taut:infinity_p_norm}
    \norm{\phi}_{\infty,p} \coloneqq \norm{\phi}
    \end{equation}
where $\norm{-}$ on the right means the norm for the (one-dimensional and non-degenerate) $F_p$-Hermitian space $\Hom_{F_p}(T_p(Y_1^{\perp})^0, T_p(Y_2^{\perp})^0)$.

We may now proceed as in the Faltings height case.
Given an $\mc{O}_F$-linear isogeny $\phi \colon A_1 \ra A_2$, we define the \emph{semi-global change of tautological height}
    \begin{equation}
    \delta_{\mrm{tau},(p)}(\phi) \coloneqq \frac{1}{[E : \Q]} \left (\log | \norm{\phi}_{\infty} |_p +  \log|\coker(\phi_* \colon \ms{E}_1^{\vee} \ra \ms{E}_2^{\vee})| \right )
    \end{equation}
where $|-|_p$ is the usual $p$-adic norm (well-defined by Lemma \ref{lemma:Faltings_and_taut:isogeny_change_local_taut:infty_norm}).
We have $\delta_{\mrm{tau},(p)}(\phi) \in \Q \cdot \log p$. Since $A_1$ and $A_2$ have everywhere potentially good reduction (implied by the special hypothesis: elliptic curves with $\mc{O}_F$-action over number fields have everywhere potentially good reduction)
the formula for change of tautological height \eqref{equation:Faltings_and_taut:isogeny_change_global:taut} implies
    \begin{equation}\label{equation:Faltings_and_taut:isogeny_change_local_taut:decomp_to_semi-global}
    h_{\mrm{tau}}(A_{2,E}) - h_{\mrm{tau}}(A_{1,E}) = \sum_{\ell} \delta_{\mrm{tau},(\ell)}(\phi) = \sum_{\ell \mid \deg \phi} \delta_{\mrm{tau},(\ell)}(\phi)
    \end{equation}
where $\phi$ also denotes the induced isogeny on N\'eron models over $\Spec \mc{O}_{E,(\ell)}$ for each prime $\ell$ (after enlarging $E$ if necessary). In particular, $\delta_{\mrm{tau},(p)}(\phi)$ does not depend on the choice of isogeny $\phi$.

Given any $\mc{O}_{F_p}$-linear isogeny $\phi \colon X_1 \ra X_2$, we define the \emph{local change of tautological height}
    \begin{equation}
    \breve{\delta}_{\mrm{tau}}(\phi) \coloneqq \log \norm{\phi}_{\infty, p} + \frac{1}{[\breve{E} : \breve{\Q}_p]} \mrm{length}_{\mc{O}_{\breve{E}}}(\coker(\phi_* \colon \ms{E}_1^{\vee} \ra \ms{E}_2^{\vee})) \cdot \log p.
    \end{equation}
We have $\breve{\delta}_{\mrm{tau}}(\phi) \in \Q \cdot \log p$, as well as
    \begin{align}\label{equation:isogeny_change_local:taut:linear_relations}
    & \breve{\delta}_{\mrm{tau}}(\phi' \circ \phi) = \breve{\delta}_{\mrm{tau}}(\phi') + \breve{\delta}_{\mrm{tau}}(\phi)
    && \breve{\delta}_{\mrm{tau}}([N]) = 0
    \end{align}
where $\phi' \colon X_2 \ra X_2$ is any $\mc{O}_{F_p}$-linear isogeny and $[N] \colon X_1 \ra X_1$ is the multiplication-by-$N$ isogeny. 
For use in later calculations, we note the identity
    \begin{align}\label{equation:Faltings_and_taut:isogeny_change_local_taut:length_additivity}
    & \operatorname{length}_{\mc{O}_{\breve{E}}}(\coker(\phi_* \colon \Lie(X_1) \ra \Lie(X_2)))
    \\
    & = \operatorname{length}_{\mc{O}_{\breve{E}}}(\coker(\phi_* \colon \mc{F}_1 \ra \mc{F}_2)) + \operatorname{length}_{\mc{O}_{\breve{E}}}(\coker(\phi_* \colon \ms{E}^{\vee}_1 \ra \ms{E}^{\vee}_2)) \notag
    \end{align}
(by the snake lemma).

\begin{lemma}\label{lemma:isogeny_change_local_taut:nonsplit_independence}
If $F_p / \Q_p$ is nonsplit, we have $\breve{\delta}_{\mrm{tau}}(\phi) = \breve{\delta}_{\mrm{tau}}(\phi')$ for any two $\mc{O}_F$-linear isogenies $\phi, \phi' \colon X_1 \ra X_2$.
\end{lemma}
\begin{proof}
Set $X = \mf{X}_0^{n - 1} \times \mf{X}_0^{\s}$, and equip $X$ with any $\mc{O}_F$-action-compatible quasi-polarization. Select any $\mc{O}_F$-linear isogeny $\phi'' \colon X \ra X_1$. Using the additivity property $\breve{\delta}_{\mrm{tau}}(\phi \circ \phi'') = \breve{\delta}_{\mrm{tau}}(\phi) + \breve{\delta}_{\mrm{tau}}(\phi'')$ and similarly for $\phi'$, this reduces us to the case where $X = X_1$.

As in the proof of Lemma \ref{lemma:Faltings_and_taut:isogeny_change_local_Faltings:global-to-local}, there exists an isogeny $\phi'' \colon X \ra X$ such that $[p^N] \circ \phi = \phi' \circ \phi''$ for some $N \geq 0$, so the additivity properties of $\breve{\delta}_{\mrm{tau}}$ reduce us to showing $\breve{\delta}_{\mrm{tau}}(\phi) = 0$ when $(X_1, \iota_1, \lambda_1) = (X_2, \iota_2, \lambda_2)$.

Since $\Hom_{\mc{O}_F}(\mf{X}_0, \mf{X}_0^{\s}) = \Hom_{\mc{O}_F}(\mf{X}_0^{\s}, \mf{X}_0) = 0$, we must have $\phi = f \times f^{\perp}$ where $f \colon \mf{X}_0^{n - 1} \ra \mf{X}_0$ and $f^{\perp} \colon \mf{X}_0^{\s} \ra \mf{X}_0^{\s}$. We find $\breve{\delta}_{\mrm{tau}}(\phi) = \breve{\delta}_{\mrm{tau}}(f^{\perp}) = 0$ since $f^{\perp} \colon \mf{X}_0^{\s} \ra \mf{X}_0^{\s}$ is an automorphism times $[p^N]$ for some $N \geq 0$.
\end{proof}

\begin{remark}\label{remark:Faltings_and_taut:isogeny_change_local_taut:split_complication}
If $F_p / \Q_p$ is split, then Lemma \ref{lemma:isogeny_change_local_taut:nonsplit_independence} fails (consider multiplication by $(1,p)$ and $(p,1)$ in $\mc{O}_{F_p} \cong \Z_p \times \Z_p$). This is the reason for Lemma \ref{lemma:Faltings_and_taut:isogeny_change_local_taut:global-to-local} below, which allows us to uniformly treat all cases of $F_p / \Q_p$.
\end{remark}

Continuing to allow $F_p/\Q_p$ inert/ramified/split, now suppose that $(X_i, \iota_i, \lambda_i)$ is the Hermitian $p$-divisible group associated with $(A_i, \iota_i, \lambda_i)$, for $i = 1, 2$. Since each $(A_i, \iota_i, \lambda_i)$ is special, there automatically exists an $\mc{O}_F$-linear isogeny $A_1 \ra A_2$ after possibly replacing $E$ by a finite extension (by the theory of complex multiplication for elliptic curves). 
Given any $\mc{O}_{F_p}$-linear isogeny $\phi \colon X_1 \ra X_2$, set
    \begin{equation}\label{equation:Faltings_and_taut:isogeny_change_local_taut:pDiv_semi-global_def}
    \delta_{\mrm{tau},(p)}(\phi) \coloneqq \frac{1}{[E : \Q]} \sum_{\breve{w} \mid p} [\breve{E}_{\breve{w}} : \breve{\Q}_p] \breve{\delta}_{\mrm{tau}}(\phi_{\breve{w}})
    \end{equation}
where $\phi_{\breve{w}}$ denotes the base-change of $\phi$ to $\Spf \mc{O}_{\breve{E}_{\breve{w}}}$.

\begin{lemma}\label{lemma:Faltings_and_taut:isogeny_change_local_taut:global-to-local}
Suppose that $(X_i, \iota_i, \lambda_i)$ is the Hermitian $p$-divisible group associated with $(A_i, \iota_i, \lambda_i)$, for $i = 1, 2$.
For any $\mc{O}_F$-linear isogenies $\tilde{\phi} \colon A_1 \ra A_2$ and $\phi \colon X_1 \ra X_2$, we have
    \begin{equation}
    \delta_{\mrm{tau},(p)}(\tilde{\phi}) = \delta_{\mrm{tau}, (p)}(\phi).
    \end{equation}
\end{lemma}
\begin{proof}
This may be proved exactly as in Lemma \ref{lemma:Faltings_and_taut:isogeny_change_local_Faltings:global-to-local}, now requiring isogenies and endomorphisms to be $\mc{O}_F$-linear.
\end{proof}

\begin{corollary}
In the situation of Lemma \ref{lemma:Faltings_and_taut:isogeny_change_local_taut:global-to-local}, the quantity $\delta_{\mrm{tau},(p)}(\phi)$ does not depend on the choice of isogeny $\phi \colon X_1 \ra X_2$.
\end{corollary}
\begin{proof}
In the notation of the lemma, this follows immediately from $\tilde{\phi}$ independence of $\delta_{\mrm{tau},(p)}(\tilde{\phi})$ (discussed above).
\end{proof}
    
            \subsection{Serre tensor}
            \label{ssec:Faltings_and_taut:Serre_tensor}
                We compute local changes of Faltings and tautological heights for isogenies involving the Serre tensor $p$-divisible groups $\mf{X}_s \otimes_{\Z_p} \mc{O}_{F_p}$. These results will later be used to compute heights of arithmetic special $1$-cycles.

Given $s \in \Z_{\geq 0}$ and a quasi-canonical lifting $\mf{X}_s$ over $\Spf \mc{O}_{\breve{E}}$, we write $\lambda_{\mf{X}_s}$ for an understood principal polarization of $\mf{X}_s$. Recall that $\lambda_{\mf{X}_s}$ exists and is unique up to $\Z_p^{\times}$ scalar (Lemma \ref{lemma:moduli_pDiv:Serre_tensor:lift_polarization} and its proof). As in Section \ref{ssec:moduli_pDiv:RZ}, we consider the map $\lambda_{\mrm{tr}} \colon \mc{O}_{F_p} \ra \mc{O}_{F_p}^{*}$ determined by the $\Z_p$-bilinear pairing $\mrm{tr}_{F_p / \Q_p}(x^{\s} y)$ on $\mc{O}_{F_p}$, where $\mc{O}_{F_p}^{*} \coloneqq \Hom_{\Z_p}(\mc{O}_{F_p}, \Z_p)$.

We equip $\mf{X}_s \otimes_{\Z_p} \mc{O}_{F_p}$ with its Serre tensor $\mc{O}_{F_p}$-action $\iota$ and the polarization $-\iota(\mf{d}_p^2)^{-1} \circ (\lambda_{\mf{X}_s} \otimes \lambda_{\mrm{tr}})$. We equip $\mf{X}_0 \times \mf{X}_0^{\s}$ with its diagonal $\mc{O}_{F_p}$ action $\iota_{\mf{X}_0} \times \iota_{\mf{X}_0}^{\s}$ (of signature $(1,1)$) and the diagonal quasi-polarization $-\iota(\mf{d}_p^2)^{-1} \circ (\lambda_{\mf{X}_0} \times \lambda_{\mf{X}_0})$.

\begin{lemma}\label{lemma:Faltings_and_taut:Serre_tensor:product}
For the $\mc{O}_{F_p}$-linear isogeny
    \begin{equation}
    \begin{tikzcd}[row sep = tiny]
    \mf{X}_0 \otimes_{\Z_p} \mc{O}_{F_p} \arrow{r}{\phi} & \mf{X}_0 \times \mf{X}_0^{\s}
    \\
    x \otimes a \arrow[mapsto]{r} & (\iota_{\mf{X}_0}(a) x, \iota_{\mf{X}_0}(a^{\s}) x)
    \end{tikzcd}
    \end{equation}
we have $\breve{\delta}_{\mrm{Fal}}(\phi) = 0$. Assuming $p \neq 2$ if $F_p / \Q_p$ is ramified, we also have $\breve{\delta}_{\mrm{tau}}(\phi) = 0$. 
\end{lemma}
\begin{proof}
We already know $\deg \phi = |\Delta|_p^{-1}$ (see \eqref{equation:moduli_pDiv:RZ:isogeny_from_tensor:polarizations} and surrounding discussion). 

Pick any $\mc{O}_{F_p}$-linear isomorphism $\Lie \mf{X}_0 \cong \mc{O}_{\breve{E}}$. Then the map $\phi_* \colon \Lie(\mf{X}_0 \otimes_{\Z_p} \mc{O}_{F_p}) \ra \Lie(\mf{X}_0 \times \mf{X}_0^{\s})$ may be identified with the map of $\mc{O}_{\breve{E}}$-modules $f \colon \mc{O}_{\breve{E}} \otimes_{\Z_p} \mc{O}_{F_p} \ra \mc{O}_{\breve{E}} \oplus \mc{O}_{\breve{E}}$ given by $f(x \otimes a) = (a x, a^{\s} x)$. Thus $\phi_*$ is given by the matrix in \eqref{equation:moduli_pDiv:RZ:isogeny_from_tensor:matrix} (the same matrix describing $\phi$ after identifying $\mf{X}_0 \otimes_{\Z_p} \mc{O}_{F_p} \cong \mf{X}_0^2$ using a $\Z_p$-basis of $\mc{O}_{F_p}$). That matrix has determinant which generates the different ideal $\mf{d}_p$, hence 
    \begin{equation}
    \operatorname{length}_{\mc{O}_{\breve{E}}}(\coker(\phi_* \Lie(\mf{X}_0 \otimes_{\Z_p} \mc{O}_{F_p}) \ra \Lie(\mf{X}_0 \times \mf{X}_0^{\s}))) = \frac{1}{2} [\breve{E} : \breve{\Q}_p] v_p(\Delta). 
    \end{equation}
This gives $2 \breve{\delta}_{\mrm{Fal}}(\phi) = \log \deg \phi - v_p(\Delta) \log p = 0$.

We also know that $\phi^*(\lambda_{\mf{X}_0} \times \lambda_{\mf{X}_0}) = \lambda_{\mf{X}_0} \otimes \lambda_{\mrm{tr}}$ (see discussion surrounding \eqref{equation:moduli_pDiv:RZ:isogeny_from_tensor:polarizations} again). Thus $\norm{\phi}_{\infty, p} = 1$, in the notation of \eqref{equation:Faltings_and_taut:isogeny_change_local_taut:infinity_p_norm}.

Let $\mc{F}_1 \subseteq \Lie (\mf{X}_0 \otimes_{\Z_p} \mc{O}_{F_p})$ and $\mc{F}_2 \subseteq \Lie(\mf{X}_0 \times \mf{X}_0^{\s})$ be the (unique) associated Kr\"amer hyperplanes, with associated Kr\"amer hyperplane quotients $\ms{E}^{\vee}_1 $ and $\ms{E}^{\vee}_2$. 
If $F_p / \Q_p$ is unramified, then $\phi$ is an isomorphism, hence $\coker(\phi_* \colon \ms{E}^{\vee}_1 \ra \ms{E}^{\vee}_2) = 0$. 
If $F_p / \Q_p$ is ramified, assume $p \neq 2$ and select a uniformizer $\varpi \in \mc{O}_{F_p}$ satisfying $\varpi^{\s} = -\varpi$. 
Then $(\varpi \otimes 1 + 1 \otimes \varpi) \in \mc{O}_{\breve{E}} \otimes_{\Z_p} \mc{O}_{F_p}$ is a generator of $\mc{F}_1$. 
We thus find $\coker(\phi_* \colon \mc{F}_1 \ra \mc{F}_2) \cong \mc{O}_{\breve{E}} / \varpi \mc{O}_{\breve{E}}$. 

By \eqref{equation:Faltings_and_taut:isogeny_change_local_taut:length_additivity}, the previous computations imply $\coker(\phi_* \colon \ms{E}^{\vee}_1 \ra \ms{E}^{\vee}_2) = 0$, and hence
    \begin{equation}
    \breve{\delta}_{\mrm{tau}}(\phi) = \log \norm{\phi}_{\infty, p} + \frac{1}{[\breve{E} : \breve{\Q}_p]} \mrm{length}_{\mc{O}_{\breve{E}}}(\coker(\phi_* \colon \ms{E}^{\vee}_1 \ra \ms{E}^{\vee}_2)) \cdot \log p = 0.
    \end{equation}
\end{proof}

For any given integer $s \in \Z_{\geq 0}$, recall the constants $\delta_{\mrm{tau}}(s), \delta_{\mrm{Fal}}(s) \in \Q$ (``local change of `tautological' and Faltings heights'') as defined in \cref{equation:can_and_qcan:qcan:local_change_taut} and surrounding text.

\begin{lemma}\label{lemma:Faltings_and_taut:Serre_tensor:qcan_isog}
Let $\psi_s \colon \mf{X}_0 \ra \mf{X}_s$ be any isogeny of degree $p^s$. For the $\mc{O}_{F_p}$-linear isogeny
    \begin{equation}
    \begin{tikzcd}[row sep = tiny]
    \mf{X}_0 \otimes_{\Z_p} \mc{O}_{F_p} \arrow{r}{\phi} & \mf{X}_s \otimes_{\Z_p} \mc{O}_{F_p}
    \\
    x \otimes a \arrow[mapsto]{r} & \psi_s(x) \otimes a
    \end{tikzcd}
    \end{equation}
we have 
    \begin{align}
    & \breve{\delta}_{\mrm{Fal}}(\phi) = -2 \breve{\delta}_{\mrm{tau}}(\phi) = -2 \delta_{\mrm{tau}}(s) \cdot \log p.
    \end{align}
\end{lemma}
\begin{proof}
Recall that $\psi_s$ is unique up to pre-composition by elements of $\mc{O}_{F_p}^{\times}$ \eqref{equation:can_and_qcan:qcan:psi_s}. Write $\mc{F}_1$ and $\mc{F}_2$ (resp. $\ms{E}^{\vee}_1$ and $\ms{E}^{\vee}_2$) for the associated Kr\"amer hyerplanes (resp. dual tautological bundles) of $\mf{X}_0 \otimes_{\Z_p} \mc{O}_{F_p}$ and $\mf{X}_s \otimes_{\Z_p} \mc{O}_{F_p}$ respectively.

We have $\deg \phi = (\deg \psi_s)^2 = p^{2s}$. Since quasi-polarizations on $\mf{X}_0$ are unique up to $\Q_p^{\times}$ scalar (follows from Drinfeld rigidity and the corresponding statement for $\mbf{X}_0$ in Section \ref{ssec:moduli_pDiv:RZ}), we have $\psi_s^* \lambda_{\mf{X}_s} = b \lambda_{\mf{X}_0}$ for some $b \in p^s \Z_p^{\times}$. Hence we have $\phi^*(\lambda_{\mf{X}_s} \otimes \lambda_{\mrm{tr}}) = b (\lambda_{\mf{X}_0} \otimes \lambda_{\mrm{tr}})$, so $\norm{\phi}_{\infty,p} = p^{-s/2}$.

Pick any identifications $\Lie \mf{X}_0 \cong \Lie \mf{X}_s \cong \mc{O}_{\breve{E}}$ of $\mc{O}_{\breve{E}}$-modules. With these identifications, the map $\psi_{s,*} \colon \Lie \mf{X}_0 \ra \Lie \mf{X}_s$ is multiplication by some $c \in \mc{O}_{\breve{E}}$ satisfying $[\breve{E} : \breve{\Q}_p] v_p(c) = \operatorname{length}_{\mc{O}_{\breve{E}}}(\coker (\psi_{s,*} \colon \Lie \mf{X}_0 \ra \Lie \mf{X}_s))$.

We also obtain identifications $\Lie (\mf{X}_0 \otimes_{\Z_p} \mc{O}_{F_p}) \cong \Lie (\mf{X}_s \otimes_{\Z_p} \mc{O}_{F_p})$ of $\mc{O}_{\breve{E}} \otimes_{\Z_p} \mc{O}_{F_p}$-modules, with induced identifications $\mc{F}_1 \cong \mc{F}_2$ and $\ms{E}^{\vee}_1 \cong \ms{E}_2^{\vee}$. Then $\phi_* \colon \Lie (\mf{X}_0 \otimes_{\Z_p} \mc{O}_{F_p}) \ra \Lie (\mf{X}_s \otimes_{\Z_p} \mc{O}_{F_p})$ is identified with multiplication by $c$, and hence $\phi_* \colon \ms{E}^{\vee}_1 \ra \ms{E}_2^{\vee}$ must also be multiplication by $c$.
Hence
    \begin{align}
    \operatorname{length}_{\mc{O}_{\breve{E}}}(\coker (\phi_* \colon \Lie(\mf{X}_0 \otimes_{\Z_p} \mc{O}_{F_p}) \ra \Lie(\mf{X}_s \otimes_{\Z_p} \mc{O}_{F_p}))) & = 2 [\breve{E} : \breve{\Q}_p] v_p(c)
    \\
    \operatorname{length}_{\mc{O}_{\breve{E}}}(\coker (\phi_* \colon \ms{E}^{\vee}_1 \ra \ms{E}^{\vee}_2)) & = [\breve{E} : \breve{\Q}_p] v_p(c).
    \end{align}
The lemma now follows from the formula for $\mrm{length}_{\mc{O}_{\breve{E}}}(\coker (\psi_{s,*} \colon \Lie \mf{X}_0 \ra \Lie \mf{X}_s))$ in \eqref{equation:can_and_qcan:qcan:Nakkajima--Taguchi}.
\end{proof}

        \section{Heights and quasi-canonical liftings}
        \label{sec:qcan_heights}
    
            \subsection{A descent lemma}
            \label{ssec:qcan_heights:descent_lemma}
                To compute Faltings and tautological heights, we will produce isogenies of $p$-divisible groups over $\Spec \mc{O}_{E,(p)}$ from isogenies over $\Spf \mc{O}_{\breve{E}_{\breve{w}}}$ for any choice of $\breve{w} \mid p$. We now explain this descent procedure, in a more general setup.

\begin{lemma}\label{lemma:qcan_heights:descent_lemma}
Let $S' \ra S$ be a morphism of schemes whose scheme-theoretic image is all of $S$. Suppose $X$ is a $p$-divisible group over $S$ which satisfies $\End^{0}(X) = \End^{0}(X_{S'})$. Let $Y$ and $Z$ be $p$-divisible groups over $S$ which are isogenous to $X$. The base-change maps
    \begin{align*}
    & \! \Hom^0(Y,Z) \ra \Hom^0(Y_{S'}, Z_{S'}) & & \! \Hom(Y,Z) \ra \Hom(Y_{S'}, Z_{S'}) \\
    & \mrm{Isog}^0(Y,Z) \ra \mrm{Isog}^0(Y_{S'}, Z_{S'}) & & \mrm{Isog}(Y,Z) \ra \mrm{Isog}(Y_{S'}, Z_{S'})
    \end{align*}
are bijections. 
\end{lemma}
\begin{proof}
Choose isogenies $\phi_Y \colon X \ra Y$ and $\phi_Z \colon X \ra Z$. There is a commutative diagram
    \[
    \begin{tikzcd}
    \a \arrow[mapsto]{d} & \End^0(X) \arrow{r}{\sim} \arrow{d}[vertLabel]{\sim} & \End^0(X_{S'}) \arrow{d}[vertLabelSwap]{\sim} & \a' \arrow[mapsto]{d} \\
    \phi_Z \circ \a \circ \phi_Y^{-1} & \Hom^0(Y,Z) \arrow{r} & \Hom^0(Y_{S'}, Z_{S'}) & \phi_{Z,S'} \circ \a' \circ \phi_{Y,S'}^{-1}
    \end{tikzcd}
    \]
where horizontal arrows are base-change. The vertical arrows are isomorphisms, and the upper horizontal arrow is an isomorphism by hypothesis. Hence the bottom arrow is an isomorphism. Suppose $\b \in \Hom^0(Y,Z)$ is any quasi-homomorphism. The functor $T \mapsto \{ \phi \in \Hom(T,S) : \phi^* \b \text{ is a homomorphism} \}$ is represented by a closed subscheme of $T$, see \cite[{Proposition 2.9}]{RZ96}. If $\b|_{S'}$ is a homomorphism, then $\b$ must also be a homomorphism, since the smallest closed subscheme of $S$ through which $S'$ factors is all of $S$ (by hypothesis). Hence $\Hom(Y,Z) \ra \Hom(Y_{S'},Z_{S'})$ is an isomorphism. The statements about (quasi-)isogenies follow from an essentially identical argument, replacing $\End$ and $\Hom$ with $\mrm{Isog}$, and noting $\mrm{Isog}^0(X) = (\End^0(X))^\times$ (e.g. by \crefext{III:lemma:qisog_invertible}).
\end{proof}

\begin{remark}
We will be interested in the case where $S = \Spec \mc{O}_{E,(p)}$ and $S' = \Spec \mc{O}_{\breve{E}}$ for some finite extension $\breve{E}$ of $\breve{E}_{\breve{w}}$ for some $\breve{w} \mid p$. In this case, Lemma \ref{lemma:qcan_heights:descent_lemma} admits an alternative proof: a quasi-homomorphism of $p$-divisible groups over $\Spec E$ is a homomorphism if and only if the map on rational Tate modules preserves (integral) Tate modules, and this can be checked after base-change to $\Spec \breve{E}$. Then apply the theorem of Tate \cite[{Theorem 4}]{Tate67} which states that the generic fiber functor for $p$-divisible groups over $\Spec \mc{O}_{E,(p)}$ (similarly, for $\Spec \mc{O}_{\breve{E}}$) is fully faithful.
\end{remark}

\begin{lemma}\label{lemma:qcan_heights:descent_lemma:product_decomp_pdiv}
Let $X$ be a $p$-divisible group over a formal scheme $S$. Suppose there is a decomposition $X = X_1 \times X_2$ as fppf sheaves of abelian groups (on $(\Sch/S)_{fppf}$). Then $X_1$ and $X_2$ are both $p$-divisible groups.
\end{lemma}
\begin{proof}
Write $e_1, e_2 \in \End(X)$ for the projections to $X_1$ and $X_2$ respectively. As being a $p$-divisible group can be checked locally on $(\Sch/S)_{fppf}$, assume $S$ is a usual scheme.

It is clear that the multiplication by $p$ map $[p] \colon X \ra X$ is a surjection if and only if $[p] \colon X_1 \ra X_1$ and $[p] \colon X_2 \ra X_2$ are surjections. We also have $X[p^n] = X_1[p^n] \times X_2[p^n]$ for all $n \geq 1$. Thus the natural map $\varinjlim X[p^n] \ra X$ is an isomorphism if and only if $\varinjlim X_1[p^n] \ra X_1$ and $\varinjlim X_2[p^n] \ra X_2$ are isomorphisms.

Next, note $X_1[p] = \ker(e_2 \colon X[p] \ra X[p])$ and similarly $X_2[p] = \ker(e_1 \colon X[p] \ra X[p])$. Since $X[p]$ is representable by a finite locally free scheme over $S$, we conclude that $X_1[p]$ and $X_2[p]$ are represented by schemes which are finite and finitely presented over $S$. We also have short exact sequences
    \begin{align*}
    & 0 \ra X_1[p] \ra X[p] \ra X_2[p] \ra 0 \\
    & 0 \ra X_2[p] \ra X[p] \ra X_1[p] \ra 0
    \end{align*}
so \crefext{III:lemma:surj_finite_locally_free_fppf} implies that $X_1[p]$ and $X_2[p]$ are finite locally free over $S$.
\end{proof}

\begin{corollary}\label{corollary:qcan_heights:descent_lemma:hom_lemma_product_decomp}
Let $S' \ra S$ and $X$ be as in Lemma \ref{lemma:qcan_heights:descent_lemma}. Suppose $Y$ and $Z$ are $p$-divisible groups over $S$ isogenous to $X$. 

If $Y_{S'} = Y'_1 \times \cdots \times Y'_r$ and $Z_{S'} = Z'_1 \times \cdots \times Z'_r$ for $p$-divisible groups $Y'_i$ and $Z'_i$ over $S'$, then there are unique decompositions $Y = Y_1 \times \cdots \times Y_r$ and $Z = Z_1 \times \cdots \times Z_r$ such that $Y_i|_{S'} = Y'_i$ and $Z_i|_{S'} = Z'_i$ for all $i$. For any $i$, the base-change maps
    \begin{align*}
    & \! \Hom^0(Y_i, Z_i) \ra \Hom^0(Y_{i, S'}, Z_{i, S'}) & & \! \Hom(Y_i,Z_i) \ra \Hom(Y_{i, S'}, Z_{i, S'}) \\
    & \mrm{Isog}^0(Y_i, Z_i) \ra \mrm{Isog}^0(Y_{i, S'}, Z_{i, S'}) & & \mrm{Isog}(Y_i ,Z_i) \ra \mrm{Isog}(Y_{i, S'}, Z_{i, S'})
    \end{align*}
are bijective.
\end{corollary}
\begin{proof}
The decomposition $Y_{S'} = Y'_1 \times \cdots \times Y'_r$ corresponds to a system of orthogonal idempotents $d'_1, \ldots, d'_r \in \End(Y_{S'})$, i.e. $d_i^{\prime 2} = d'_i$ for all $i$ and $d'_i d'_j = 0$ for all $i \neq j$. Lifting to a decomposition $Y = Y_1 \times \cdots \times Y_r$ is the same as lifting $\{d'_i\}_{i}$ to a system of orthogonal idempotents $\{d_i\}_{i}$ in $\End(Y)$. Such a lift exists and is unique by Lemma \ref{lemma:qcan_heights:descent_lemma}. The same applies for $Z$, and we write $\{e'_i\}_i$ and $\{e_i\}_i$ for the corresponding systems of idempotents. 
Using Lemma \ref{lemma:qcan_heights:descent_lemma}, we have 
    \begin{align*}
    & \Hom^{0}(Y_i,Z_i) = d_i \Hom^{0}(Y,Z) e_i = d'_i \Hom^{0} (Y_{S'}, Z_{S'}) e'_i = \Hom^{0}(Y_{i, S'},Z_{i, S'})
    \\
    & \Hom(Y_i,Z_i) = d_i \Hom(Y,Z) e_i = d'_i \Hom (Y_{S'}, Z_{S'}) e'_i = \Hom(Y_{i, S'}, Z_{i, S'}).
    \end{align*}
The statement about $\mrm{Isog}^{0}$ then follows from \crefext{III:lemma:qisog_invertible}, and the statement about $\mrm{Isog}$ follows from the relation $\mrm{Isog}(-,-) = \mrm{Isog}^{0}(-,-) \cap \Hom(-,-)$.
\end{proof}
            
            \subsection{Minimal isogenies}
            \label{ssec:qcan_heights:minimal_isogenies}
                Given any abelian scheme $A \ra S$ over some base $S$, we can form the \emph{Serre tensor} abelian scheme $A \otimes_{\Z} \mc{O}_F$ given by $(A \otimes_{\Z} \mc{O}_F)(T) \coloneqq A(T) \otimes_{\Z} \mc{O}_F$ for $S$-schemes $T$. There is a natural action of $\mc{O}_F$ on $A \otimes_{\Z} \mc{O}_F$, as we have discussed for $p$-divisible groups (``Serre tensor''). If $\lambda \colon A \ra A^{\vee}$ is a quasi-polarization, then $\lambda \otimes \lambda_{\mrm{tr}} \colon A \otimes_{\Z} \mc{O}_F \ra A^{\vee} \otimes_{\Z} \mc{O}_F^{\vee} \cong (A \otimes_{\Z} \mc{O}_F)^{\vee}$ is a polarization, where $\lambda_{\mrm{tr}} \colon \mc{O}_F \ra \mc{O}_F^{\vee}$ is induced by the trace pairing, as above.

Let $A_0 \ra \Spec \mc{O}_{E,(p)}$ be any (relative) elliptic curve with $\mc{O}_F$-action $\iota_0$ of signature $(1,0)$, and let $\lambda_0$ be the unique principal polarization of $A_0$. For $n \geq 2$, set
    \begin{equation}
    A \coloneqq A_0^{n - 2} \times (A_0 \otimes_{\Z} \mc{O}_F)
    \end{equation}
with $\mc{O}_F$ action $\iota$ which is diagonal on $A_0^{n - 2}$ and the Serre tensor action on $A_0 \otimes_{\Z} \mc{O}_F$, and polarization $\lambda_0^{n - 2} \times (|\Delta|^{-1} (\lambda_0 \otimes \lambda_{\mrm{tr}}))$. Then $(A, \iota, \lambda)$ is a special Hermitian abelian scheme of signature $(n - 1, 1)$. We write $(X, \iota, \lambda)$ for the associated special Hermitian $p$-divisible group of signature $(n - 1, 1)$, with
    \begin{equation}
    X = X_0^{n - 2} \times (X_0 \otimes_{\Z_p} \mc{O}_{F_p})
    \end{equation}
where $X_0$ is the $p$-divisible group of $A_0$. For any $\breve{w} \mid p$, the base-change $X_{0, \breve{w}}$ is a canonical lifting. The preceding notation (e.g. for $A_0$ and $X_0$) will be fixed for all of Section \ref{ssec:qcan_heights:minimal_isogenies}.

In Proposition \ref{proposition:qcan_heights:minimal_isogenies} and Corollary \ref{corollary:qcan_heights:minimal_isogenies:local_decomp} below, we equip $\mf{X}_0^{n - 2} \times (\mf{X}_s \otimes_{\Z_p} \mc{O}_{F_p})$ with the diagonal $\mc{O}_{F_p}$ action (which is the Serre tensor action on $\mf{X}_s \otimes_{\Z_p} \mc{O}_{F_p}$) and a product quasi-polarization, for some quasi-polarization of $\mf{X}_0^{n - 2}$ and the quasi-polarization $-\iota(\mf{d}_p^2)^{-1} \circ (\lambda_{\mf{X}_s} \otimes \lambda_{\mrm{tr}})$ on $(\mf{X}_s \otimes_{\Z_p} \mc{O}_{F_p})$.

\begin{proposition}\label{proposition:qcan_heights:minimal_isogenies}
Let $(A', \iota', \lambda')$ be a special Hermitian abelian scheme of signature $(n - 1, 1)$ over $\Spec \mc{O}_{E,(p)}$, with associated Hermitian $p$-divisible group $(X', \iota', \lambda')$. Replace $E$ with a finite extension if necessary, so that $A$ and $A'$ are $\mc{O}_F$-linearly isogenous.

Suppose there exists a $\mc{O}_{F_p}$-linear quasi-polarization preserving isomorphism
    
    \begin{equation}\label{equation:qcan_heights:minimal_isogenies:local_decomp}
    X'_{\Spf \mc{O}_{\breve{E}}} \cong \mf{X}_0^{n - 2} \times (\mf{X}_s \otimes_{\Z_p} \mc{O}_{F_p})
    \end{equation}
over $\Spf \mc{O}_{\breve{E}}$, where $\breve{E}$ is a finite extension of $\mc{O}_{\breve{E}_{\breve{w}'}}$, for some $\breve{w}' \mid p$ and $s \geq 0$.
Fix an isomorphism $X_{0, \Spf \mc{O}_{\breve{E}}} \cong \mf{X}_0$.

\begin{enumerate}[(1)]
    \item Then there exists an $\mc{O}_{F_p}$-linear quasi-polarization preserving isomorphism
        \begin{equation}\label{equation:qcan_heights:minimal_isogenies:semi-global_decomp}
        X' \cong X_0^{n - 2} \times (X_s \otimes_{\Z_p} \mc{O}_{F_p})
        \end{equation}
    over $\Spec \mc{O}_{E,(p)}$, for some $p$-divisible group $X_s$ of height $2$ and dimension $1$ with fixed identification $X_{s, \Spf \mc{O}_{\breve{E}}} \cong \mf{X}_s$, such that \eqref{equation:qcan_heights:minimal_isogenies:semi-global_decomp} recovers \eqref{equation:qcan_heights:minimal_isogenies:local_decomp} upon base-change to $\Spf \mc{O}_{\breve{E}}$. 
    On the right-hand side of \eqref{equation:qcan_heights:minimal_isogenies:semi-global_decomp}, the polarization is the product of a polarization on $X_0^{n - 2}$ and a quasi-polarization $- (\mf{d}^2_p)^{-1} \cdot (\lambda_s \otimes_{\Z_p} \lambda_{\mrm{tr}})$ on $X_s \otimes_{\Z_p} \mc{O}_{F_p}$, where $\lambda_s$ is a principal polarization on $X_s$.
    
    \item For any $\breve{w} \mid p$, the base-change $X_{s, \breve{w}}$ is a quasi-canonical lifting of level $s$, and hence there is an identification as in \eqref{equation:qcan_heights:minimal_isogenies:semi-global_decomp} for all $\breve{w} \mid p$.

    \item There exists an isogeny $\psi_s \colon X_0 \ra X_s$ of degree $p^s$. The $\mc{O}_{F_p}$-linear product isogeny $\phi \colon X \ra X'$ given by
        \begin{equation}
        \phi \coloneqq \mrm{id}_{X_0^{n - 2}} \times (\psi_s \otimes 1) \colon X_0^{n - 2} \times (X_0 \otimes_{\Z_p} \mc{O}_{F_p}) \ra X_0^{n - 2} \times (X_s \otimes_{\Z_p} \mc{O}_{F_p})
        \end{equation}
    over $\Spec \mc{O}_{E,(p)}$ satisfies
        \begin{align}
        & \breve{\delta}_{\mrm{Fal}}(\phi_{\breve{w}}) = -2 \breve{\delta}_{\mrm{tau}}(\phi_{\breve{w}}) = -2 \delta_{\mrm{tau}}(s) \cdot \log p && \text{for all $\breve{w} \mid p$} \label{equation:qcan_heights:minimal_isogenies:local_change}
        \\
        & \delta_{\mrm{Fal}, (p)}(\phi) = -2 \delta_{\mrm{tau}, (p)}(\phi) = -2 \delta_{\mrm{tau}}(s) \cdot \log p. \label{equation:qcan_heights:minimal_isogenies:semi-global_change}
        \end{align}
\end{enumerate}
\end{proposition}
\begin{proof}
Note that $X$ satisfies the hypotheses of Lemma \ref{lemma:qcan_heights:descent_lemma} with $S = \Spec \mc{O}_{E,{(p)}}$ and $S' = \Spec \mc{O}_{\breve{E}_{\breve{w}}}$ for any $\breve{w} \mid p$, as $\End(X) \cong M_{n,n}(\mc{O}_{F_p})$ over both $\Spec \mc{O}_{E,(p)}$ and $\Spf \mc{O}_{\breve{w}}$ for any $\breve{w}$. The same holds for $S' = \Spec \mc{O}_{\breve{E}}$. Again, we pass between $\Spf \mc{O}_{\breve{E}_{\breve{w}}}$ and $\Spec \mc{O}_{\breve{E}_{\breve{w}}}$ as in \crefext{III:appendix:pDiv_prelim:Spec_v_Spf}.

The proposition then follows from repeated applications of Lemma \ref{lemma:qcan_heights:descent_lemma} and Corollary \ref{corollary:qcan_heights:descent_lemma:hom_lemma_product_decomp}, as we now explain.

\begin{enumerate}[(1)]
    \item and (2) Corollary \ref{corollary:qcan_heights:descent_lemma:hom_lemma_product_decomp} implies that \eqref{equation:qcan_heights:minimal_isogenies:local_decomp} descends to a $\mc{O}_{F_p}$-linear product decomposition $X' \cong X_0^{\prime n - 2} \times (X_s \otimes_{\Z_p} \mc{O}_{F_p})$ over $\Spec \mc{O}_{E, (p)}$ for some $X'_0$ descending $\mf{X}_0$ (first pick any identification of $p$-divisible groups $\mf{X}_s \otimes_{\Z_p} \mc{O}_{F_p} \cong \mf{X}_s^2$, then descend the $\mc{O}_{F_p}$-action), and the fully-faithfulness in Corollary \ref{corollary:qcan_heights:descent_lemma:hom_lemma_product_decomp} implies $\End(X_s) = \mc{O}_{F_p,s}$ (with $\mc{O}_{F_p,s} = \Z_p + p^s \mc{O}_{F_p}$ as in Section \ref{ssec:can_and_qcan:qcan}) over $\Spec \mc{O}_{E,(p)}$ and also over $\Spf \mc{O}_{\breve{w}}$ for any $\breve{w} \mid p$. The fully-faithfulness in Corollary \ref{corollary:qcan_heights:descent_lemma:hom_lemma_product_decomp} also implies that the fixed $\mc{O}_{F_p}$-linear isomorphism $X_{0, \breve{w}} \ra \mf{X}_0$ lifts to an isomorphism $X_0 \ra X'_0$. The polarization on $\mf{X}_0^{n - 2} \times (\mf{X}_s \otimes_{\Z_p} \mc{O}_{F_p})$ descends to $X_0^{n - 2} \times (X_s \otimes_{\Z_p} \mc{O}_{F_p})$ by Corollary \ref{corollary:qcan_heights:descent_lemma:hom_lemma_product_decomp} again (applied to $X'$ and $X^{\prime \vee}$; note that the property of being a polarization is represented by a closed subfunctor of $\Spec \mc{O}_{E,(p)}$, hence can be checked in the generic fiber or over $\Spec \breve{E}$)
    \addtocounter{enumi}{1}
    
    \item If $\psi_s \colon \mf{X}_0 \ra \mf{X}_s$ is any isogeny of degree $p^s$ (exists and is unique up to precomposition by $\mc{O}_{F_p}^{\times}$, as discussed in Section \ref{ssec:can_and_qcan:qcan}), we apply Corollary \ref{corollary:qcan_heights:descent_lemma:hom_lemma_product_decomp} to descend to an isogeny $\psi_s \colon X_0 \ra X_s$ of degree $p^s$. Equation \eqref{equation:qcan_heights:minimal_isogenies:local_change} now follows from Lemma \ref{lemma:Faltings_and_taut:Serre_tensor:qcan_isog}. Equation \eqref{equation:qcan_heights:minimal_isogenies:semi-global_change} follows from this (by the definitions in \eqref{equation:Faltings_and_taut:isogeny_change_local_Faltings:pDiv_semi-global_def} and \eqref{equation:Faltings_and_taut:isogeny_change_local_taut:pDiv_semi-global_def}). \qedhere
\end{enumerate}
\end{proof}

We will use the following reformulation (tailored to our intended application for global heights via local special cycles). In the corollary statement and proof, $A_0^{\s}$ and $A_0^{n - 1} \times A_0^{\s}$ are equipped with the product quasi-polarizations $- \Delta^{-1} \lambda_0$ and $-\Delta^{-1} (\lambda_0 \times \cdots \times \lambda_0)$ (where $A_0^{\s} = A_0$ but with $\mc{O}_F$-action $\iota \circ \s$, as above).

\begin{corollary}\label{corollary:qcan_heights:minimal_isogenies:local_decomp}
Let $S$ be a reduced scheme which is finite flat over $\Spec \mc{O}_F$. Let $(A', \iota', \lambda', \mc{F}')$ be a quasi-polarized Kr\"amer datum over $S$ (of signature $(n - 1, 1)$) for $n \geq 2$, with associated metrized line bundles $\widehat{\omega}$ and $\widehat{\ms{E}}^{\vee}$ on $S$. Assume that $(A', \iota', \lambda')$ is special at all generic points of $S$.
Let $(X', \iota', \lambda')$ be the associated Hermitian $p$-divisible group. 

Suppose we are given a finite \'etale surjection
    \begin{equation}
    \coprod_j \mc{Z}_j \ra S \times_{\Spec \Z} \Spec \breve{\Z}_p
    \end{equation}
such that each restricted map $\Theta_j \colon \mc{Z}_j \ra S \times_{\Spec \Z} \Spec \breve{\Z}_p$ has constant degree $\deg(j)$ onto its image. Assume that $\Theta_j$ and $\Theta_{j'}$ have disjoint images for $j \neq j'$.

For each irreducible component $\mc{Z} \hookrightarrow \coprod_j \mc{Z}_j$, write $\breve{E}_{\mc{Z}}$ for the residue field of its generic point. Assume there exists an isomorphism of Hermitian $p$-divisible groups
    \begin{equation}
    X'|_{\Spf \mc{O}_{\breve{E}_{\mc{Z}}}} \cong \mf{X}_0^{n - 2} \times (\mf{X}_{s_{\mc{Z}}} \otimes_{\Z_p} \mc{O}_{F_p})
    \end{equation}
for all $\mc{Z}$, where $s_{\mc{Z}} \in \Z_{\geq 0}$ is an integer depending on $\mc{Z}$.

We then have
    \begin{align}\label{equation:qcan_heights:minimal_isogenies:degree_decomp}
    \widehat{\deg}(\widehat{\ms{E}}^{\vee}) - (\deg_{\Z} S) \cdot h_{\mrm{tau}}(A_0^{\s}) & =  \sum_{j \in J} \frac{1}{\deg(j)} \sum_{\mc{Z} \hookrightarrow \mc{Z}_j} (\deg_{\breve{\Z}_p} \mc{Z}) \delta_{\mrm{tau}}(s_{\mc{Z}}) \log p
    \\
    \widehat{\deg}(\widehat{\omega}) - (\deg_{\Z} S) \cdot h_{\mrm{Fal}}(A_0^{n - 1} \times A_0^{\s}) & = \sum_{j \in J} \frac{1}{\deg(j)} \sum_{\mc{Z} \hookrightarrow \mc{Z}_j}  (\deg_{\breve{\Z}_p} \mc{Z}) \delta_{\mrm{Fal}}(s_{\mc{Z}}) \log p \notag
    \end{align}
modulo $\sum_{\ell \neq p} \Q \cdot \log \ell$, where the inner sums run over all irreducible components $\mc{Z} \hookrightarrow \mc{Z}_j$.
\end{corollary}
\begin{proof}
In the corollary statement, the expression ``modulo $\sum_{\ell \neq p} \Q \cdot \log \ell$'' means an equality of elements in the additive quotient $\R / (\sum_{\ell \neq p} \Q \cdot \log \ell)$. The notation $\deg_{\Z} S$ (resp. $\deg_{\breve{\Z}_p} \mc{Z}$) denotes the degree of $S \ra \Spec \Q$ (resp. $\mc{Z} \ra \Spec \breve{\Z}_p$) in the generic fiber.

By additivity, we immediately reduce to the case where $S$ is irreducible. Then $J$ consists of a single element $j$. By normalization, we may assume $S = \Spec \mc{O}_E$ for a number field $E$. We may also enlarge $E$ as necessary so that $(A_0, \iota, \lambda_0)$ also extends to $\Spec \mc{O}_E$, and such that there exists an $\mc{O}_F$-linear isogeny $\phi \colon A_0^{n - 2} \times (A_0 \otimes_{\Z} \mc{O}_F) \ra A$. We also consider the $\mc{O}_F$-linear isogeny $\phi' \colon A_0^{n - 2} \times (A_0 \otimes_{\Z} \mc{O}_F) \ra A_0^{n - 2} \times (A_0 \times A_0^{\s})$ which is the identity on $A_0^{n - 2}$ and given by $(x \otimes a) \mapsto (ax, a^{\s} x)$ for $(A_0 \otimes_{\Z} \mc{O}_F) \ra A_0 \times A_0^{\s}$.

Since $\delta_{\mrm{tau}}(\phi') = 0$ (Lemma \ref{lemma:Faltings_and_taut:Serre_tensor:product}, along with the local decomposition \eqref{equation:Faltings_and_taut:isogeny_change_local_taut:pDiv_semi-global_def}, also Lemma \ref{lemma:Faltings_and_taut:isogeny_change_local_taut:global-to-local}), the decomposition in \eqref{equation:Faltings_and_taut:isogeny_change_local_taut:decomp_to_semi-global} shows
    \begin{align}
    \widehat{\deg}(\widehat{\ms{E}}^{\vee}) - [E : \Q] \cdot h_{\mrm{tau}}(A_0^{n - 1} \times A_0^{\s}) = [E : \Q] \delta_{\mrm{tau},(p)}(\phi) \mod \sum_{\ell \neq p} \Q \cdot \log \ell.
    \end{align}
We have $\sum_{\mc{Z} \hookrightarrow \mc{Z}_j} (\deg_{\breve{\Z}_p} \mc{Z}) = \deg(j) \cdot [E : \Q]$. Applying Proposition \ref{proposition:qcan_heights:minimal_isogenies}(3) (combined with the ``isogeny independence'' result of Lemma \ref{lemma:Faltings_and_taut:isogeny_change_local_taut:global-to-local}) now shows $\delta_{\mrm{tau},(p)}(\phi) = \delta_{\mrm{tau}}(s_{\mc{Z}}) \log p $ for any $\mc{Z} \hookrightarrow \mc{Z}_j$. This also shows that all $s_{\mc{Z}}$ are equal (when $S$ is irreducible): the quantity $\delta_{\mrm{tau}}(s)$ takes distinct values for distinct $s \in \Z_{\geq 0}$ (in the nonsplit cases, note $\delta_{\mrm{tau}}(s)$ has strictly decreasing $p$-adic valuation as $s$ increases, for $s > 2$). We also have $h_{\mrm{tau}}(A_0^{n - 1} \times A_0^{\s}) = h_{\mrm{tau}}(A_0^{\s})$ (straightforward from the definition). This verifies \eqref{equation:qcan_heights:minimal_isogenies:degree_decomp} for $\widehat{\deg}(\widehat{\ms{E}}^{\vee})$ and the tautological height.

Since $\delta_{\mrm{Fal}}(\phi') = 0$ (Lemma \ref{lemma:Faltings_and_taut:Serre_tensor:product}, along with the local decomposition \eqref{equation:Faltings_and_taut:isogeny_change_local_Faltings:pDiv_semi-global_def}, also Lemma \ref{lemma:Faltings_and_taut:isogeny_change_local_Faltings:global-to-local}), \eqref{equation:Faltings_and_taut:isogeny_change_local_taut:decomp_to_semi-global} similarly shows
    \begin{align}
    \widehat{\deg}(\widehat{\omega}) - [E : \Q] \cdot h_{\mrm{Fal}}(A_0^{n - 1} \times A_0^{\s}) = [E : \Q] \delta_{\mrm{Fal},(p)}(\phi) \mod \sum_{\ell \neq p} \Q \cdot \log \ell.
    \end{align}
Applying Proposition \ref{proposition:qcan_heights:minimal_isogenies}(3) (combined with the ``isogeny independence'' result of Lemma \ref{lemma:Faltings_and_taut:isogeny_change_local_Faltings:global-to-local}) verifies \eqref{equation:qcan_heights:minimal_isogenies:degree_decomp} for $\widehat{\deg}(\widehat{\omega})$ and the Faltings height, just as for tautological height above.
\end{proof}

In the situation above, we have 
    \begin{align}\label{equation:qcan_heights:minimal_isogenies:CM_heights}
    h_{\mrm{tau}}(A_0^{\s}) = h_{\mrm{tau}}^{\mrm{CM}} \quad \quad
    h_{\mrm{Fal}}(A_0^{n - 1} \times A_0^{\s}) = n \cdot h_{\mrm{Fal}}^{\mrm{CM}}
    \end{align}
in the notation of \eqref{equation:part_I:arith_cycle_classes:Hodge_bundles:taut_height_constants} and \eqref{equation:part_I:arith_cycle_classes:Hodge_bundles:CM_Faltings_height}.

\begin{remark}
In Proposition \ref{proposition:qcan_heights:minimal_isogenies}(3), it was important that $\psi_s$ was an isogeny of minimal degree $p^s$. If $\psi_s$ were replaced by an arbitrary isogeny $f \colon \mf{X}_0 \ra \mf{X}_s$, we would not be able to determine $\breve{\delta}_{\mrm{Fal}}(f)$ or $\breve{\delta}_{\mrm{tau}}(f)$ using only $\deg f$ in the case when $F_p / \Q_p$ is split (due to Remark \ref{remark:Faltings_and_taut:isogeny_change_local_taut:split_complication}).
\end{remark}

    \clearpage
    

    \appendix

    \part*{Appendices}

        \section{\texorpdfstring{$K_0$}{K0} groups} \label{appendix:K0}
    
            \subsection{\texorpdfstring{$K_0$}{K0} groups for Deligne--Mumford stacks} 
            \label{appendix:K0:K0_stacky}
                Suppose $\mc{X}$ is a Noetherian Deligne--Mumford stack. There are at least two different ways one might define $K_0$ groups for $\mc{X}$. One way is to define a $K$-theory spectrum for $\mc{X}$ using the $K$-theory spectra of schemes in the small \'etale site of $\mc{X}$, as in \cite[{\S 2}]{Gillet09}. This is the approach used in \cite{HM22}. Another way is to simply mimic a definition of $K_0$ for schemes and consider perfect complexes on the small \'etale site of $\mc{X}$. These two approaches will in general result in different $K_0$ groups \cite[{Remark A.2.4}]{HM22}. At least if $\mc{X}$ is regular (and, say, with the additional running hypotheses of \cite[{Appendix A}]{HM22}), there is a map from the latter $K_0$ group to the former $K_0$ group \cite[{(A.7),(A.8)}]{HM22}.

In this paper, we take the latter approach and mimic constructions for schemes to define $K_0(\mc{X})$. Our definitions and notation will be analogous to those for schemes in \cite[\href{https://stacks.math.columbia.edu/tag/0FDE}{Section 0FDE}]{stacks-project}. When defining dimension/codimension filtrations on $K_0'(\mc{X})$ (with notation and hypotheses as below), we will require existence of a finite flat cover by a scheme (enough for our intended application). A similar approach appears in \cite[{Appendix A}]{YZ17} (at least for $K'_0$), but there the stacks are over a base field. We need a slightly more general setup which allows base schemes such as $\Spec R$ for Dedekind domains $R$.

Suppose $\mc{X}$ is a Deligne--Mumford stack. By an \emph{$\mc{O}_{\mc{X}}$-module}\footnote{This is one of the only places where our conventions differ from the Stacks project \cite[\href{https://stacks.math.columbia.edu/tag/06TF}{Chapter 06TF}]{stacks-project}, which mostly works with sheaves on big sites (say, fppf and \'etale) for general algebraic stacks. Restriction from these two big sites to the small \'etale site (for Deligne--Mumford stacks) induces equivalences on categories of quasi-coherent sheaves. But the equivalences are not compatible with pushforward, and are also not compatible with exactness for $\mc{O}_{\mc{X}}$-modules on big sites versus small sites.}
, we mean a sheaf of modules on the small \'etale site\footnote{The small \'etale site is as defined in \cite[{Definition 4.10}]{DM69}: the underlying category has objects which are pairs $(U,f)$ for $f \colon U \ra \mc{X}$ an \'etale morphism (definition as in \cite[\href{https://stacks.math.columbia.edu/tag/0CIL}{Definition 0CIL}]{stacks-project}) from a scheme $U$, and morphisms are pairs $(g, \xi) \colon (U,f) \ra (U',f')$ where $g \colon U \ra U'$ is a $1$-morphism and $\xi \colon f \ra f' \circ g$ is a $2$-isomorphism.} of $\mc{X}$. 
Similarly, \emph{quasi-coherent $\mc{O}_{\mc{X}}$-modules} will mean quasi-coherent sheaves of modules on the small \'etale site. When $\mc{X}$ is locally Noetherian, we will also speak of \emph{coherent $\mc{O}_{\mc{X}}$-modules} on the small \'etale site, which are the same as finitely presented quasi-coherent $\mc{O}_{\mc{X}}$-modules in this situation.

Suppose $\mc{X}$ is a locally Noetherian Deligne--Mumford stack. The category $\textit{Coh}(\mc{O}_{\mc{X}})$ of coherent $\mc{O}_{\mc{X}}$-modules forms a weak Serre subcategory of the abelian category $\textit{Mod}(\mc{O}_{\mc{X}})$ of $\mc{O}_{\mc{X}}$-modules (reduce to the case of small \'etale sites of schemes and apply \cite[\href{https://stacks.math.columbia.edu/tag/05VG}{Lemma 05VG}, \href{https://stacks.math.columbia.edu/tag/0GNB}{Lemma 0GNB}]{stacks-project}). We may form derived categories such as
    \begin{equation}
    D(\mc{O}_{\mc{X}}) \quad D^b(\mc{O}_{\mc{X}}) \quad D_{\textit{perf}}(\mc{O}_{\mc{X}}) \quad D^b_{\textit{Coh}}(\mc{O}_{\mc{X}}) \quad D^b(\textit{Coh}(\mc{O}_{\mc{X}}))
    \end{equation}
which denote the derived category of $\mc{O}_{\mc{X}}$-modules, bounded derived category of $\mc{O}_{\mc{X}}$-modules, derived category of perfect objects (definition as in \cite[\href{https://stacks.math.columbia.edu/tag/08G4}{Section 08G4}]{stacks-project}) in $D(\mc{O}_{\mc{X}})$, bounded derived category of $\mc{O}_{\mc{X}}$-modules with coherent cohomology, and the bounded derived category of coherent $\mc{O}_{\mc{X}}$-modules, respectively.

If $\mc{X}$ happens to be a scheme, then $D_{\textit{perf}}(\mc{O}_{\mc{X}})$ and $D^b_{\textit{Coh}}(\mc{O}_{\mc{X}})$ and $D^b(\textit{Coh}(\mc{O}_{\mc{X}}))$ will agree with the usual constructions using the small Zariski site instead of the small \'etale site, via comparison results such as \cite[\href{https://stacks.math.columbia.edu/tag/08HG}{Lemma 08HG}, \href{https://stacks.math.columbia.edu/tag/071Q}{Lemma 071Q}, \href{https://stacks.math.columbia.edu/tag/05VG}{Lemma 05VG}]{stacks-project}.

\begin{definition}
Let $\mc{X}$ be a locally Noetherian Deligne--Mumford stack. We set
    \begin{equation}
    K_0(\mc{X}) \coloneqq K_0(D_{\textit{perf}}(\mc{O}_{\mc{X}})) \quad \quad K'_0(\mc{X}) \coloneqq K_0(\textit{Coh}(\mc{O}_{\mc{X}})).
    \end{equation}
\end{definition}
Above, the left expression means $K_0$ of a triangulated category and the right expression means $K_0$ of an abelian category. If $\mc{X}$ is a locally Noetherian Deligne--Mumford stack, we have canonical identifications
    \begin{equation}
    K_0(\textit{Coh}(\mc{O}_{\mc{X}})) = K_0(D^b(\textit{Coh}(\mc{O}_{\mc{X}}))) = K_0(D^b_{\textit{Coh}}(\mc{O}_{\mc{X}}))
    \end{equation}
as in \cite[\href{https://stacks.math.columbia.edu/tag/0FDF}{Lemma 0FDF}]{stacks-project} (the case of schemes) by general facts about derived categories (see also \cite[\href{https://stacks.math.columbia.edu/tag/0FCS}{Lemma 0FCS}]{stacks-project}).

If $\mc{X}$ is a quasi-compact locally Noetherian Deligne--Mumford stack, there is an inclusion $D_{\textit{perf}}(\mc{O}_{\mc{X}}) \ra D^b_{\textit{Coh}}(\mc{O}_{\mc{X}})$ and a corresponding group homomorphism $K_0(\mc{X}) \ra K_0'(\mc{X})$. If $\mc{X}$ is a regular locally Noetherian Deligne--Mumford stack (not necessarily quasi-compact), there is an inclusion $D^b_{\textit{Coh}}(\mc{O}_{\mc{X}}) \ra D_{\textit{perf}}(\mc{O}_{\mc{X}})$ and a corresponding group homomorphism $K_0'(\mc{X}) \ra K_0(\mc{X})$. If $\mc{X}$ is a locally Noetherian Deligne--Mumford stack which is both quasi-compact and regular, we have $D_{\textit{perf}}(\mc{O}_{\mc{X}}) = D^b_{\textit{Coh}}(\mc{O}_{\mc{X}})$ and a corresponding isomorphism 
    \begin{equation}
    K_0(\mc{X}) \xra{\sim} K_0'(\mc{X}).    
    \end{equation}
These claims follow from the corresponding facts for schemes \cite[\href{https://stacks.math.columbia.edu/tag/0FDC}{Lemma 0FDC}]{stacks-project} and comparison results mentioned previously.

The derived tensor product $\otimes^{\mbb{L}}$ on $D(\mc{O}_{\mc{X}})$ gives $K_0(\mc{X})$ the structure of a commutative ring. Compatibility of $\otimes^{\mbb{L}}$ with the case when $\mc{X}$ is also a scheme follows from the displayed equation in the proof of \cite[\href{https://stacks.math.columbia.edu/tag/08HF}{Lemma 08HF}]{stacks-project} (comparison between the small Zariski and small \'etale sites).

We next describe dimension and codimension filtrations. Our setup for dimension theory is as in \cite[\href{https://stacks.math.columbia.edu/tag/02QK}{Section 02QK}]{stacks-project}. That is, we work over a locally Noetherian and universally catenary base scheme $S$ with a dimension function $\delta \colon |S| \ra \Z$ (which we typically suppress). Typical setups will be $S = \Spec R$ for $R$ a field or Dedekind domain, where $\delta$ is the dimension function sending closed points to $0$. Any Deligne--Mumford stack $\mc{X}$ which is quasi-separated and locally of finite type over $S$ inherits a dimension function $\delta_{\mc{X}} \colon |\mc{X}| \ra \Z$ (work \'etale locally to pass to the case of schemes; the case of algebraic spaces is \cite[\href{https://stacks.math.columbia.edu/tag/0EDS}{Section 0EDS}]{stacks-project}). If $\mc{X}$ is equidimensional of dimension $n$, then $n - \delta_{\mc{X}}$ is also the codimension function (given by dimensions of local rings on \'etale covers by schemes).
For a scheme $X$ which is locally of finite type over $S$, consider the full subcategory $\textit{Coh}_{\leq d}(\mc{O}_X) \subseteq \textit{Coh}(\mc{O}_X)$ consisting of coherent $\mc{O}_X$-modules $\mc{F}$ with $\dim(\mrm{Supp}(\mc{F})) \leq d$. Then there is an increasing \emph{dimension filtration} on $K'_0(X) = K_0(\textit{Coh}(\mc{O}_X))$ given by the image
    \begin{equation}
    F_d K'_0(X) \coloneqq \mrm{im} (K_0(\textit{Coh}_{\leq d}(\mc{O}_X)) \ra K_0(\textit{Coh}(\mc{O}_X)))
    \end{equation}
as in \cite[\href{https://stacks.math.columbia.edu/tag/0FEV}{Section 0FEV}]{stacks-project}. We similarly consider the full subcategory $\textit{Coh}^{\geq m}(\mc{O}_X) \subseteq \textit{Coh}(\mc{O}_X)$ of coherent sheaves supported in codimension $\geq m$, and form the decreasing \emph{codimension filtration}
    \begin{equation}
    F^m K'_0(X) \coloneqq \mrm{im} (K_0(\textit{Coh}^{\geq m}(\mc{O}_X)) \ra K_0(\textit{Coh}(\mc{O}_X))).
    \end{equation}
When $X$ is equidimensional of dimension $n$, we have $F^m K'_0(X) = F_{n-m} K'_0(X)$.

For the case of Deligne--Mumford stacks, one could consider naive dimension/codimension filtrations on $K_0'(\mc{X})$ by mimicking the definition for schemes. This may not be a well-behaved notation, and we instead take the filtration defined in \cite[{A.2.3}]{YZ17} (with $\Q$-cofficients).

\begin{definition}\label{definition:stack_dimension_filtration}
For $S$ as above, let $\mc{X}$ be a Deligne--Mumford stack which is quasi-separated and locally of finite type over $S$. Suppose there exists a finite flat surjection $\pi \colon U \ra \mc{X}$ from a scheme $U$. Pick such a morphism $\pi$. 

The \emph{dimension filtration} on $K'_0(\mc{X})_{\Q}$ is the increasing filtration given by
    \begin{equation}
    F_d K'_0(\mc{X})_{\Q} \coloneqq \{ \b \in K'_0(\mc{X})_{\Q} : \pi^* \b \in F_d K'_0(U)_{\Q} \} \subseteq K'_0(\mc{X})_{\Q}
    \end{equation}
for $d \in \Z$. If $\mc{X}$ is equidimensional, we also consider the decreasing \emph{codimension filtration} on $K'_0(\mc{X})_{\Q}$ given by
    \begin{equation}
    F^m K'_0(\mc{X})_{\Q} \coloneqq \{ \b \in K'_0(\mc{X})_{\Q} : \pi^* \b \in F^m K'_0(U)_{\Q} \} \subseteq K'_0(\mc{X})_{\Q}.
    \end{equation}
for $m \in \Z$.
\end{definition}

For $\mc{X}$ as in the preceding definition, the filtrations just defined give rise to graded pieces $\mrm{gr}_d K'_0(\mc{X})_{\Q} \coloneqq F_d K'_0(\mc{X})_{\Q} / F_{d-1} K'_0(\mc{X})_{\Q}$ and $\mrm{gr}^m K'_0(\mc{X})_{\Q} \coloneqq F^m K'_0(\mc{X})_{\Q} / F^{m+1} K'_0(\mc{X})_{\Q}$. If $\mc{X}$ as above is equidimensional of dimension $n$, we have $F^m K'_0(\mc{X})_{\Q} = F_{n-m} K'_0(\mc{X})_{\Q}$ for all $m \in \Z$.

\begin{lemma}
With notation as in Definition \ref{definition:stack_dimension_filtration}, the filtrations $F_d K'_0(\mc{X})_{\Q}$ and $F^m K'_0(\mc{X})_{\Q}$ do not depend on the choice of finite flat surjection $\pi \colon U \ra \mc{X}$. If $\mc{X}$ is a scheme, these filtrations recovers the usual filtrations.
\end{lemma}
\begin{proof}
Suppose $X$ is a scheme which is locally of finite type over $S$. If $Z_d(X)$ is the group of $d$-cycles on $X$, recall that there is an identification 
    \begin{equation}
    K_0(\textit{Coh}_{\leq d}(\mc{O}_X) / \textit{Coh}_{\leq {d-1}}(\mc{O}_X)) \xra{\sim} Z_d(X)
    \end{equation}
which is compatible with flat pullback of constant relative dimension and finite pushforward \cite[\href{https://stacks.math.columbia.edu/tag/02S9}{Lemma 02S9}, \href{https://stacks.math.columbia.edu/tag/0FDR}{Lemma 0FDR}]{stacks-project} (see also \cite[\href{https://stacks.math.columbia.edu/tag/02MX}{Lemma 02MX}]{stacks-project}). For any finite flat surjection $\pi \colon U \ra X$ of constant degree $a$, the map $\pi_* \pi^* \colon Z_d(X) \ra Z_d(X)$ is multiplication by $a$.
It follows that $\pi_* \pi^* \colon F_d K'_0(X) / F_{d-1} K'_0(X) \ra F_d K'_0(X) / F_{d-1} K'_0(X)$ is multiplication by $a$. This is an isomorphism after tensoring by $\Q$. When $X$ is equidimensional, this gives the corresponding statement for the codimension filtration as well. This verifies the lemma when $\mc{X}$ is a scheme.

Let $\mc{X}$ is a Deligne--Mumford stack as in the lemma statement. Let $\pi \colon U \ra \mc{X}$ and $\pi' \colon U' \ra \mc{X}$ be two finite flat surjections, for schemes $U$ and $U'$. Consider the fiber product $U \times_{\mc{X}} U'$ with its finite flat projections to $U$ and $U'$. We then apply the preceding discussion to see that the filtrations do not depend on the choice of finite flat surjection.

These arguments are essentially the same as in \cite[{A.2.3}]{YZ17} (the arguments of loc. cit. are over a base field, so we have used different references).
\end{proof}

\begin{remark}
As in \cite[{A.2.3}]{YZ17}, it is possible to have $F_d K'_0(\mc{X})_{\Q} \neq 0$ for $d < 0$ in the situation of Definition \ref{definition:stack_dimension_filtration}. 
\end{remark}

We are mainly interested in $K_0$ groups for the purpose of intersection theory, so we next discuss degree theory over a field. Suppose $S = \Spec k$ for a field $k$, and suppose $\mc{X}$ is a Deligne--Mumford stack which is proper over $S$. Again assuming that $\mc{X}$ admits a finite flat surjection from a scheme, there is a graded group homomorphism $\mrm{gr}_* K'_0(\mc{X})_{\Q} \ra \mrm{Ch}_*(\mc{X})_{\Q}$ as defined in \cite[{A.2.6}]{YZ17} (pass to a finite flat surjection to reduce to the case of schemes).

There is a degree map $\deg \colon \mrm{Ch}_0(\mc{X})_{\Q} \ra \Q$ on $0$-cycles which may be described as follows. Suppose $\mc{Z}$ is a quasi-separated finite type Deligne--Mumford stack over $\Spec k$ with separated diagonal, and assume the underlying topological space $|\mc{Z}|$ is a single point. If $V \ra \mc{Z}$ is any finite flat surjection from a scheme $V$ (which exists as in Remark \ref{remark:existence_finite_covers_stack}), one can check that $V$ is finite over $\Spec k$ and we take
    \begin{equation}\label{equation:zero_cycle_degree_over_field}
    \deg(\mc{Z}) \coloneqq \deg_k(\mc{Z}) \coloneqq \frac{\deg_k (V)}{\deg_{\mc{Z}} (V)}
    \end{equation}
where $\deg_k(V)$ (resp. $\deg_{\mc{Z}}(V)$) is the degree of the finite flat morphism $V \ra \Spec k$ (resp. $V \ra \mc{Z}$).
It is straightforward to see that $\deg_k(\mc{Z})$ does not depend on the choice of $V \ra \mc{Z}$ (compare \cite[{Definition 1.15}]{Vistoli89}). This generalizes immediately to the case where $|\mc{Z}|$ is instead a discrete finite set (add the degrees of its components). When $\mc{Z} = \emptyset$, we take $\deg(\mc{Z}) \coloneqq 0$.

There is an induced degree map 
    \begin{equation}\label{equation:stacky_degree_K0}
    \deg \colon \mrm{gr}_0 K'_0(\mc{X})_{\Q} \ra \Q.
    \end{equation}
Consider a class $\b = \sum_{i} b_i [\mc{F}_i] \in F_0 K'_0(\mc{X})_{\Q}$ where each $\mc{F}_i$ is a coherent sheaf on $\mc{X}$ (we do not assume $[\mc{F}_i] \in F_0 K'_0(\mc{X})_{\Q}$ for any given $i$). Select any finite flat surjection $\pi \colon U \ra \mc{X}$ with $U$ a scheme. If $\pi$ has constant degree $a$, we have
    \begin{equation}\label{equation:euler_char_degree}
    \deg(\b) = \frac{1}{a} \deg(\pi^* \b) = \frac{1}{a} \sum_i b_i \cdot \chi(\pi^* \mc{F}_i)
    \end{equation}
where $\chi$ denotes Euler characteristic. We can give a similar description for general finite flat surjections $\pi$ by decomposing $\mc{X}$ into its connected components. On account of \eqref{equation:euler_char_degree}, we may write $\chi(\b) \coloneqq \deg(\b)$ and think of $\chi \colon \mrm{gr}_0 K'_0(\mc{X})_{\Q} \ra \Q$ as a ``stacky Euler characteristic'' (compare usage in \cite[{Definition 11.4}]{KR14}). We caution, however, that we have only defined $\chi$ on $\mrm{gr}_0 K'_0(\mc{X})_{\Q}$ and have not defined $\chi(\mc{F})$ for a general coherent sheaf $\mc{F}$ on $\mc{X}$.

We conclude this subsection with a lemma which we will use to decompose $K'_0(\mc{X})$ in terms of irreducible components of $\mc{X}$. A similar lemma for formal schemes is \cite[{Lemma B.1}]{Zhang21}.

\begin{lemma}\label{lemma:K'0_component_decomp}
Let $\mc{X}$ be a locally Noetherian Deligne--Mumford stack. Let $\pi_1 \colon \mc{Z}_1 \ra \mc{X}$ and $\pi_2 \colon \mc{Z}_2 \ra \mc{X}$ be closed immersions of Deligne--Mumford stacks with corresponding ideal sheaves $\mc{I}_1$ and $\mc{I}_2$. Assume that the diagonals of $\mc{X}$, $\mc{Z}_1$, and $\mc{Z}_2$ are representable by schemes (e.g. if $\mc{X}$ is separated).

Assume that $\mc{X} = \mc{Z}_1 \cup \mc{Z}_2$ scheme-theoretically (meaning $\mc{I}_1 \cap \mc{I}_2 = 0$). There are mutually inverse isomorphisms
    \begin{equation}
    \begin{tikzcd}[row sep = tiny]
    \frac{K'_0(\mc{X})}{K'_0(\mc{Z}_1 \cap \mc{Z}_2)} \arrow{r} & \arrow{l} \frac{K'_0(\mc{Z}_1)}{K'_0(\mc{Z}_1 \cap \mc{Z}_2)} \oplus \frac{K'_0(\mc{Z}_2)}{K'_0(\mc{Z}_1 \cap \mc{Z}_2)}
    \\
    {[\mc{F}]} \arrow[mapsto]{r} & ({[\pi_1^* \mc{F}_1]}, {[\pi_2^* \mc{F}_2]}) \\
    {[\pi_{1,*} \mc{F}_1]} + {[\pi_{2,*} \mc{F}_2]} & \arrow[mapsto]{l} ({[\mc{F}_1]}, {[\mc{F}_2]}) \rlap{~.}
    \end{tikzcd}
    \end{equation}
Here, $\mc{F}$, $\mc{F}_1$, and $\mc{F}_2$ stand for coherent sheaves on $\mc{X}$, $\mc{Z}_1$, and $\mc{Z}_2$ respectively.
\end{lemma}
\begin{proof}
The condition about diagonals is included for technical convenience. Some additional explanation on notation in the lemma statement: the symbol $\mc{Z}_1 \cap \mc{Z}_2$ denotes the closed substack $\mc{Z}_1 \times_{\mc{X}} \mc{Z}_2$ of $\mc{X}$, with associated ideal sheaf $\mc{I}_1 + \mc{I}_2$, and we have also written $K'_0(\mc{X})/K'_0(\mc{Z}_1 \cap \mc{Z}_2) \coloneqq \coker(K'_0(\mc{Z}_1 \cap \mc{Z}_2) \ra K'_0(\mc{X}))$ etc. (the latter map may not be injective).

Consider the short exact sequence
    \begin{equation}
    0 \ra \mc{O}_{\mc{X}}/(\mc{I}_1 \cap \mc{I}_2) \ra \mc{O}_{\mc{X}}/\mc{I}_1 \oplus \mc{O}_{\mc{X}}/\mc{I}_2 \ra \mc{O}_{\mc{X}}/(\mc{I}_1 + \mc{I}_2) \ra 0.
    \end{equation}
Tensoring by any coherent sheaf $\mc{F}$ on $\mc{X}$, we find that $\mrm{Tor}_1^{\mc{O}_{\mc{X}}}(\mc{F}, \mc{O}_{\mc{X}}/\mc{I}_1)$ is an $\mc{O}_{\mc{X}}/(\mc{I}_1 + \mc{I}_2)$-module, and similarly with $\mc{I}_2$ instead of $\mc{I}_1$. This shows that the displayed projection maps $\mc{F} \mapsto \pi_1^* \mc{F}$ and $\mc{F} \mapsto \pi_2^* \mc{F}$ are well-defined (i.e. that they are additive in short exact sequences and hence descend to the given quotients of $K'_0$-groups). Since $\mrm{Tor}_1^{\mc{O}_{\mc{X}}}(\mc{F}, \mc{O}_{\mc{X}}/(\mc{I}_1 + \mc{I}_2))$ is an $\mc{O}_{\mc{X}}/(\mc{I}_1 + \mc{I}_2)$-module, the Tor long exact sequence of the displayed short exact sequence also shows that $[\mc{F}] = [\mc{F} \otimes_{\mc{O}_{\mc{X}}} \mc{O}_{\mc{X}}/\mc{I}_1] + [\mc{F} \otimes_{\mc{O}_{\mc{X}}} \mc{O}_{\mc{X}}/\mc{I}_2]$ in $K'_0(\mc{X}) / K'_0(\mc{Z}_1 \cap \mc{Z}_2)$.
\end{proof}
                
            \subsection{\texorpdfstring{$K_0$}{K0} groups with supports along finite morphisms} \label{appendix:K0:K0_relative_supports}
                Suppose $X$ is a separated regular Noetherian scheme.
There is an established intersection theory for $K_0$ groups with supports along closed subsets of $X$, and the intersection pairing is multiplicative with respect to codimension filtrations (after tensoring by $\Q$) \cite{GS87}. However, we will need a slightly more general setup which allows for ``supports along finite morphisms''. This is needed because the special cycles $\mc{Z}(T) \ra \mc{M}$ (\cref{ssec:part_I:arith_intersections:special_cycles}) are not literally cycles but are instead finite unramified morphisms.

Intersection theory with supports along finite morphisms is also discussed in \cite[{Appendix A.4}]{HM22} for a similar purpose. They are not able to show the codimension filtration is multiplicative in general \cite[{Remark A.4.2}]{HM22}, but can show multiplicativity for intersections against classes of codimension $1$ (in the case of supports along finite unramified morphisms) \cite[{Proposition A.4.4}]{HM22}.

We have two main objectives in this section (besides fixing notation). Our first objective is to comment on another situation where the codimension filtration is multiplicative (namely, when the finite supports become disjoint unions of closed immersions after finite flat base change to a regular scheme). The (short) proof reduces to the case of supports along closed immersions. This is relevant for us because of \crefext{III:lemma:special_cycles_level_structure}, which says that each special cycle $\mc{Z}(T) \ra \mc{M}$ becomes a disjoint union of closed immersions after finite \'etale base change, at least after inverting the prime $p$ in the cited lemma.  
For $\mc{M}$ associated to a Hermitian lattice of signature $(n - r, r)$, intersecting special cycles over $\mc{M}$ involves multiplicativity for classes of codimension $r$ (not covered by \cite[{Proposition A.4.4}]{HM22} when $r > 1$).

Our second objective is to explain intersection theory with supports along finite morphisms of Deligne--Mumford stacks in terms of the $K_0$ groups of \cref{appendix:K0:K0_stacky}. A stacky theory is also considered in \cite[{Appendix A.4}]{HM22}, but the $K_0$ groups we use are slightly different (as discussed at the beginning of \cref{appendix:K0:K0_stacky}). The setup we consider agrees with \cite[{Appendix A.4}]{HM22} for schemes.

\begin{lemma}\label{lemma:K0_relative}
Consider a $2$-commutative diagram of algebraic stacks
    \begin{equation}
    \begin{tikzcd}
    \mc{Z} \times_{\mc{X}} \mc{W} \arrow{r} \arrow{d} \arrow{dr}{h} & \mc{W} \arrow{d}{g} \\
    \mc{Z} \arrow{r}{f} & \mc{X}
    \end{tikzcd}
    \end{equation}
with outer square $2$-Cartesian, where $\mc{X}$ is a separated regular Noetherian Deligne--Mumford stack and the morphisms $f$ and $g$ (and hence $h$) are finite.

There is a bilinear pairing
    \begin{equation}
    \begin{tikzcd}[row sep = tiny]
    K'_0(\mc{Z}) \times K'_0(\mc{W}) \arrow{r} & K'_0(\mc{Z} \times_{\mc{X}} \mc{W}) \\
    (\mc{F}, \mc{G}) \arrow[mapsto]{r} & \sum_{i} (-1)^i \mrm{Tor}_i^{\mc{O}_{\mc{X}}}(f_*\mc{F}, g_*\mc{G})
    \end{tikzcd}
    \end{equation}
where $\mc{F}$ and $\mc{G}$ stand for coherent $\mc{O}_{\mc{Z}}$-modules and coherent $\mc{O}_{\mc{W}}$-modules, respectively. We have a commutative diagram
    \begin{equation}
    \begin{tikzcd}
    K'_0(\mc{Z}) \times K'_0(\mc{W}) \arrow{r} \arrow{d}{f_* \times g_*} & K'_0(\mc{Z} \times_{\mc{X}} \mc{W}) \arrow{d}{h_*} \\
    K'_0(\mc{X}) \times K'_0(\mc{X}) \arrow{r} & K'_0(\mc{X})
    \end{tikzcd}
    \end{equation}
where vertical arrows are pushforward and the lower horizontal arrow is the bilinear pairing from the ring structure on $K'_0(\mc{X}) \cong K_0(\mc{X})$.
\end{lemma}
\begin{proof}
If $\mc{F}$ is a coherent $\mc{O}_{\mc{Z}}$-module and $\mc{G}$ is a coherent $\mc{O}_{\mc{W}}$-module, we may form the object $(f_* \mc{F} \otimes^{\mbb{L}} g_* \mc{G})$ in $D_{\textit{perf}}(\mc{O}_{\mc{X}})$. For each object $U \ra \mc{X}$ in the small \'etale site of $\mc{X}$ (i.e. $U$ is a scheme with an \'etale morphism to $\mc{X}$), the restriction $(f_* \mc{F} \otimes^{\mbb{L}} g_* \mc{G})|_U \in D_{\textit{perf}}(\mc{O}_{U})$ carries natural $\mc{O}_{U}$-linear actions of $(f_* \mc{O}_{\mc{Z}})|_U $ and $(f_* \mc{O}_{\mc{W}})|_U$.
The resulting cohomology sheaves $\mrm{Tor}_i^{\mc{O}_{\mc{X}}}(f_*\mc{F}, g_* \mc{G}) = H^{-i}(f_* \mc{F} \otimes^{\mbb{L}} g_* \mc{G})$ (a priori coherent $\mc{O}_{\mc{X}}$-modules) are thus sheaves of $(f_* \mc{O}_{\mc{Z}}) \otimes_{\mc{O}_{\mc{X}}} (g_* \mc{O}_{\mc{W}})$-algebras. There is a canonical isomorphism $(f_* \mc{O}_{\mc{Z}}) \otimes_{\mc{O}_{\mc{X}}} (g_* \mc{O}_{\mc{W}}) \ra h_* \mc{O}_{\mc{Z} \times_{\mc{X}} \mc{W}}$ of $\mc{O}_{\mc{X}}$-algebras. Since $h$ is affine, we obtain a lift (up to canonical isomorphism) of each $\mrm{Tor}_i^{\mc{O}_{\mc{X}}}(f_* \mc{F}, g_* \mc{G})$ to a coherent sheaf of $\mc{O}_{\mc{Z} \times_{\mc{X}} \mc{W}}$ modules (to pass between quasi-coherent $h_* \mc{O}_{\mc{Z} \times_{\mc{X}} \mc{W}}$-modules and quasi-coherent $\mc{O}_{\mc{Z} \times_{\mc{X}} \mc{W}}$-modules, we may take an \'etale surjection of $\mc{X}$ from a scheme, use the corresponding result for the small \'etale site of schemes which is \cite[\href{https://stacks.math.columbia.edu/tag/08AI}{Lemma 08AI}]{stacks-project}, and reduce to a statement about glueing data on the small \'etale sites of $\mc{Z} \times_{\mc{X}} \mc{W}$ and $\mc{X}$).

The procedure just described descends to $K'_0$ groups and gives the pairing in the lemma statement.
\end{proof}

We think of the map $K'_0(\mc{Z}) \times K'_0(\mc{W}) \ra K'_0(\mc{Z} \times_{\mc{X}} \mc{W})$ from the preceding lemma as an \emph{intersection pairing} ``with supports along finite morphisms''.

Next, fix a base scheme $S$ with dimension function $\delta$ as in \cref{appendix:K0:K0_stacky}. Suppose $\mc{X}$ is a Deligne--Mumford stack which is quasi-separated and locally of finite type over $S$. We assume that $\mc{X}$ is equidimensional of dimension $n$, and we also assume that $\mc{X}$ admits a finite flat surjection from a scheme in order to define dimension and codimension filtrations as in Definition \ref{definition:stack_dimension_filtration}.

Consider a finite morphism $f \colon \mc{Z} \ra \mc{X}$ from a Deligne--Mumford stack $\mc{Z}$. We define a ``relative codimension'' filtration on $K'_0(\mc{Z})_{\Q}$ by setting
    \begin{equation}
    F^m_{\mc{X}} K'_0(\mc{Z})_{\Q} \coloneqq F_{n-m} K'_0(\mc{Z})_{\Q}.
    \end{equation}
We similarly set $\mrm{gr}^m_{\mc{X}} K'_0(\mc{Z})_{\Q} \coloneqq F^m_{\mc{X}} K'_0(\mc{Z})_{\Q} / F^{m+1}_{\mc{X}} K'_0(\mc{Z})_{\Q}$. The subscript $\mc{X}$ is meant to remind of the dependence on $\mc{X}$.

\begin{lemma}\label{lemma:K0_relative_with_filtrations}
Let $\mc{X}$ be a regular Noetherian Deligne--Mumford stack which is separated and finite type over $S$. Assume that $\mc{X}$ is equidimensional. Let $f \colon \mc{Z} \ra \mc{X}$ and $g \colon \mc{W} \ra \mc{X}$ be finite morphisms from Deligne--Mumford stacks $\mc{Z}$ and $\mc{W}$.

Assume that there exists a finite flat surjection $\pi \colon U \ra \mc{X}$ with $U$ a regular Noetherian scheme, such that $\mc{Z} \times_{\mc{X}} U \ra U$ and $\mc{W} \times_{\mc{X}} U \ra U$ are both disjoint unions of closed immersions. Then the intersection pairing of Lemma \ref{lemma:K0_relative} restricts to a pairing
    \begin{equation}
    F^s_{\mc{X}} K'_0(\mc{Z})_{\Q} \times F^t_{\mc{X}} K'_0(\mc{W})_{\Q} \ra F^{s+t}_{\mc{X}} K'_0(\mc{Z} \times_{\mc{X}} \mc{W})_{\Q}
    \end{equation}
for any $s,t \in \Z$.
\end{lemma}
\begin{proof}
We use the shorthand $\mc{Z}_U \coloneqq \mc{Z} \times_{\mc{X}} U$ and $\mc{W}_U \coloneqq \mc{W} \times_{\mc{X}} U$. If we abuse notation and also write $\pi$ for the natural projections $\mc{Z}_U \ra \mc{Z}$ and $\mc{W}_U \ra \mc{W}$ and $\mc{Z}_U \times_U \mc{W}_U \ra \mc{Z} \times_{\mc{X}} \mc{W}$, we have $(\pi^* \a) \cdot (\pi^* \b) = \pi^* (a \cdot \b)$ for any $\a \in K'_0(\mc{Z})_{\Q}$ and $\b \in K'_0(\mc{W})_{\Q}$.
By definition of the dimension filtration (Definition \ref{definition:stack_dimension_filtration}), it is enough to check that the intersection pairing over $U$ restricts to
    \begin{equation}
    F^s_{U} K'_0(\mc{Z}_U)_{\Q} \times F^t_{U} K'_0(\mc{W}_U)_{\Q} \ra F_{U}^{s+t} K'_0(\mc{Z}_U \times_{U} \mc{W}_U)_{\Q}
    \end{equation}
(i.e. respects filtrations). We have thus reduced to the case where $\mc{X}$ is a scheme and $\mc{Z} \ra \mc{X}$ and $\mc{W} \ra \mc{X}$ are disjoint unions of closed immersions, and we assume these conditions hold for the rest of the proof. Write $\mc{Z} = \coprod_i \mc{Z}_i$ where each $\mc{Z}_i \ra \mc{X}$ is a closed immersion of schemes, and similarly write $\mc{W} = \coprod_j \mc{W}_j$. By a result of Gillet--Soul\'e \cite[{Proposition 5.5}]{GS87}, the pairing $F^s_{\mc{X}} K'_0(\mc{Z}_i)_{\Q} \times F^t_{\mc{X}} K'_0(\mc{W}_j)_{\Q} \ra K'_0(\mc{Z}_i \times_{\mc{X}} \mc{W}_j)_{\Q}$ factors through $F^{s+t}_{\mc{X}} K'_0(\mc{Z}_i \times_{\mc{X}} \mc{W}_j)_{\Q}$. We may decompose $F_{\mc{X}}^s K'_0(\mc{Z})_{\Q} = \bigoplus_i F_{\mc{X}}^s K'_0(\mc{Z}_i)_{\Q}$ and $F^t_{\mc{X}} K'_0(\mc{W})_{\Q} = \bigoplus_j F^t_{\mc{X}} K'_0(\mc{W}_j)_{\Q}$. Commutativity of the diagram
    \begin{equation}
    \begin{tikzcd}
    F^s_{\mc{X}} K'_0(\mc{Z}_i)_{\Q} \times F^t_{\mc{X}} K'_0(\mc{W}_j)_{\Q} \arrow{r} \arrow{d} & F^{s+t}_{\mc{X}} K'_0(\mc{Z}_i \times_{\mc{X}} \mc{W}_j)_{\Q} \arrow{d} \\
    ( \bigoplus_i F^s_{\mc{X}} K'_0(\mc{Z}_i)_{\Q} ) \times ( \bigoplus_j F^t_{\mc{X}} K'_0(\mc{W}_j)_{\Q} ) \arrow{r} & K'_0(\mc{Z} \times_{\mc{X}} \mc{W})_{\Q}
    \end{tikzcd}
    \end{equation}
for each $i,j$ gives the claim.
\end{proof}

    \clearpage


    \phantomsection
    \addcontentsline{toc}{part}{References}
    \renewcommand{\addcontentsline}[3]{}
    \printbibliography

\end{document}